\documentclass[12pt,oneside]{book}
\newif\ifprint
\printfalse 

\usepackage{amsmath,amssymb,mathtools}
\usepackage{amsthm}
\usepackage{amsrefs}
\usepackage{libertinus}
\usepackage{geometry}
\usepackage{enumerate,enumitem,stmaryrd}
\usepackage{hyperref,framed}
\usepackage{caption,subcaption}
\usepackage[dvipsnames]{xcolor}
\usepackage{tikz}
\usepackage{tikz-cd}
\usetikzlibrary{positioning}
\usetikzlibrary{calc}
\usepackage{pgfplots}
\pgfplotsset{compat=1.18}
\usetikzlibrary{arrows.meta}
\usepackage{tabularx}
\usepackage{booktabs}  
\usepackage{makeidx}
\makeindex

\theoremstyle{definition}
\newtheorem{thm}{Theorem}[chapter]
\newtheorem{ex}[thm]{Examples}
\newtheorem{1ex}[thm]{Example}
\theoremstyle{remark}
\newtheorem*{note}{Note}

\newcounter{example}[chapter]

\newenvironment{examples}{
  \begin{enumerate}[label=\textbf{(\alph*)},leftmargin=*]
}{
  \end{enumerate}
}

\newenvironment{exercises}{
  \begin{enumerate}[leftmargin=*]
}{
  \end{enumerate}
}

\setlist[itemize]{leftmargin=*, itemsep=0.3em}
\setlist[enumerate]{itemsep=0.3em}

\newcommand{\N}{\mathbb{N}}
\renewcommand{\H}{\mathbb{H}}
\newcommand{\Z}{\mathbb{Z}}
\newcommand{\F}{\mathbb{F}}
\newcommand{\Q}{\mathbb{Q}}
\newcommand{\R}{\mathbb{R}}
\newcommand{\C}{\mathbb{C}}
\newcommand{\im}{\mathrm{im}}
\newcommand{\lcm}{\mathrm{lcm}}
\newcommand{\calC}{\mathcal{C}}
\newcommand{\GL}{\operatorname{GL}}
\newcommand{\Aut}{\operatorname{Aut}}
\newcommand{\Frac}{\operatorname{Frac}}
\newcommand{\Hom}{\operatorname{Hom}}

\newcommand{\Spec}{\operatorname{Spec}}
\newcommand{\Sub}{\operatorname{Sub}}
\newcommand{\Id}{\operatorname{Id}}
\newcommand{\adj}{\operatorname{adj}}
\newcommand{\ord}{\operatorname{ord}}
\newcommand{\ch}{\operatorname{char}}

\DeclarePairedDelimiter\pow{\llbracket}{\rrbracket}
\DeclarePairedDelimiter\conv{\{}{\}}
\DeclarePairedDelimiter\lau{(\!(}{)\!)}

\ifprint
    \newcommand{\col}[1]{\textbf{#1}}
    \newcommand{\sel}{\textbf{$\star$}}
    \newcommand{\house}{%
      \tikz[scale=0.09, baseline=-0.5ex]{
        \draw[thick, black] (0,0) -- (1,1.2) -- (2,0) -- cycle;
        \draw[thick, black] (0.4,0) -- (0.4,-1) -- (1.6,-1) -- (1.6,0);}}
    \hypersetup{
    colorlinks=false,
    pdfborder={0 0 0}
  }

\else
    \definecolor{DarkMagenta}{RGB}{150,0,120}
    \newcommand{\col}[1]{\textbf{\textcolor{DarkMagenta}{#1}}}
    \newcommand{\sel}{\textcolor{DarkMagenta}{$\star$}}
    \newcommand{\house}{%
      \tikz[scale=0.09, baseline=-0.5ex]{
        \draw[thick, DarkMagenta] (0,0) -- (1,1.2) -- (2,0) -- cycle;
        \draw[thick, DarkMagenta] (0.4,0) -- (0.4,-1) -- (1.6,-1) -- (1.6,0);}}
    \hypersetup{
      colorlinks=true,
      linkcolor=RoyalBlue,
      citecolor=RoyalBlue,
      urlcolor=RoyalBlue
    }
\fi

\ifprint
  \tikzset{
    fone/.style  = {ultra thick, black},                 
    ftwo/.style  = {ultra thick, black, dashed},         
    funity/.style = {thick, black, dash pattern=on 6pt off 3pt},
    fzero/.style  = {thick, black, dotted}
  }
\else
  \tikzset{
    fone/.style  = {ultra thick, teal},
    ftwo/.style  = {ultra thick, orange},
    funity/.style = {thick, black, dash pattern=on 6pt off 3pt},
    fzero/.style  = {thick, black, dotted}
  }
\fi

\title{A Course in Ring Theory}
\author{David Krumm}
\date{}

\begin{document}

\frontmatter
\maketitle

\ifprint
\thispagestyle{empty}

\begin{center}
{\small
\textcopyright\ \the\year\ David Krumm\\[1em]
\vspace{1em}

}
\end{center}
\clearpage
\fi

\section*{Preface}
\addcontentsline{toc}{section}{Preface}

This book offers a rigorous and conceptually unified development of ring theory. Beyond an exposition of the underlying theory, the book includes numerous examples and exercises, sample homework problems, and thematic quizzes. These materials are meant to serve both as a resource for a self-guided course and as an aid for students and instructors in a formal academic setting. Importantly, the book does not assume prior knowledge of group theory.

The early chapters introduce fundamental notions such as rings, subrings, and polynomial rings, with an emphasis on both structure and examples. The later chapters build toward more sophisticated concepts such as prime and maximal ideals, universal properties, and unique factorization domains.

Each chapter includes exercises that reinforce core concepts while encouraging deeper thought. Exercises marked with a star (\sel) are solved in a later review chapter, while those marked with a house (\house) are suggested as sample homework problems.

The style of writing is intended to be rigorous and succinct, omitting details where deemed appropriate. Substantial engagement is expected from readers who wish to verify all unproven claims.

\medskip

\textbf{Prerequisites.} A solid foundation in proof-based mathematics is assumed. Readers should be comfortable with logical reasoning, set theory, functions, equivalence relations, partially ordered sets, and elementary concepts in number theory.

More specifically, the book assumes familiarity with the following topics:
\begin{itemize}
  \item The number systems $\Z$, $\Q$, $\R$, and $\C$, as well as modular arithmetic.
  \item Proofs by mathematical induction.
  \item Fundamental set-theoretic concepts: union, intersection, difference, functions, and equivalence relations.
  \item Basic linear algebra: especially matrices and their addition and multiplication.
  \item Elementary calculus: continuous and differentiable functions; convergence of sequences and series.
  \item Partially ordered sets (posets); joins (suprema) and meets (infima).
  \item The extended real number system, including $\infty$ and $-\infty$.
  \item Basic number theory: greatest common divisors and prime factorization of integers. The division algorithm: every integer division yields a quotient and a remainder of smaller absolute value than the divisor.
  \item A few non-essential exercises make use of concepts from topology or complex analysis.
\end{itemize}

\medskip

\textbf{Notation.} Throughout this book, $\N$ denotes the set of natural numbers, including zero:
\[
\N = \{0, 1, 2, \ldots\}.
\]

\medskip

\textbf{Further Reading.} The bibliography contains several references covering both ring theory and group theory. All were consulted in the preparation of this book, and are recommended as sources of additional exercises and different perspectives on the subject. In particular, Gallian \cite{Gallian} includes numerous historical notes on abstract algebra.

\medskip

\textbf{Comments.} Please email any comments regarding the book to

\medskip
\centerline{\href{mailto:david.krumm@gmail.com}{\texttt{david.krumm@gmail.com}}.}
\medskip

Notices of typos or any other errors are especially appreciated.

\medskip
\textbf{Acknowledgments.} Many thanks to Dino Lorenzini for helpful comments on an early draft of the book.

\section*{List of Symbols}
\addcontentsline{toc}{section}{List of Symbols}

\begin{tabularx}{\textwidth}{@{}lX@{}}
\textbf{Symbol} & \textbf{Meaning} \\
\toprule
\( \N \) & The set of natural numbers \( \{0,1,2,\dots\} \) \\
\( \Z \) & The ring of integers \\
\( \Q \) & The field of rational numbers \\
\( \R \) & The field of real numbers \\
\( \C \) & The field of complex numbers \\
\( \F_p \) & The finite field with \( p \) elements \\
\( A[x] \) & Polynomial ring over ring \( A \) \\
\( A\pow{x} \) & Power series ring over ring \( A \) \\
\( A\lau{x} \) & Laurent series ring over ring \( A \) \\
\( A^{\times} \) & The unit group of ring $A$. \\
\( a \mid b \) & \( a \) divides \( b \), i.e., \( b = ac \) for some $c$ \\
\( a \nmid b \) & \( a \) does not divide \( b \) \\
\( a \sim b \) & \( a \) is associate to $b$, i.e., $a=ub$ for some unit $u$ \\
\(a \equiv b \pmod I\) & \(a - b \in I\); congruence modulo an ideal \(I\) \\
\( (a) \) & The principal ideal generated by \( a \) \\
\( A/I \) & The quotient ring of \( A \) by the ideal \( I \) \\
$\Frac(D)$ & The field of fractions of an integral domain \(D\) \\
\( \cong \) & Ring isomorphism \\
\( \ker(f) \) & Kernel of a ring homomorphism \( f \) \\
\( \mathrm{im}(f) \) & Image of the map \( f \) \\
\( \gcd(a,b) \) & Any greatest common divisor of \( a \) and \( b \) \\
\( \mathrm{lcm}(a,b) \) & Any least common multiple of \( a \) and \( b \) \\
\( \mathrm{id} \) & Identity map \\
\( 0_A \) & The additive identity of a ring $A$ \\
\( 1_A \) & The multiplicative identity of a ring $A$ \\
\bottomrule
\end{tabularx}

\tableofcontents

\mainmatter
\chapter{Course Introduction}\label{intro_chap}

Ring theory studies certain types of sets equipped with operations of addition and multiplication. This includes $\Z$, the set of integers, as well as $M_2(\R)$, the set of $2\times 2$ matrices with real entries, whose operations are the sum and product of matrices defined in linear algebra. 

The central concept of ring theory is, naturally, that of a \emph{ring}, defined formally later in this chapter. Rings are ubiquitous throughout mathematics---they appear in algebraic geometry, number theory, algebraic topology, representation theory, mathematical physics, and other areas. They also play a key role in modern applications such as \href{https://en.wikipedia.org/wiki/Cryptography}{cryptography} and \href{https://en.wikipedia.org/wiki/Coding_theory}{coding theory}, contributing to technologies including blockchains and quantum computing.

This book focuses on describing the theoretical foundations of ring theory. Readers who become comfortable with the material presented here will be well prepared to explore the many applications of the subject.

A topic closely related to ring theory is \emph{group theory}, which studies a certain type of set equipped with a single operation. While there is a beautiful story to be told about the interplay between ring theory and group theory, that story is not the subject of this book. Our aim here is to develop ring theory on its own terms, assuming no prior knowledge of group theory. The concept of a group will be defined and used repeatedly, but no theorems about groups will be needed.

Before proceeding to develop ring theory, we remind the reader of the assumed background and other assumptions detailed in the book's preface. In particular, some statements are left unproven when the proofs are deemed to be relatively straightforward.

\section*{Binary operations}

We begin course material by defining the concept of a binary operation on a set. If $S$ is an arbitrary set, a \col{binary operation}\index{binary operation} on $S$ is a function $S\times S\to S$. Informally, we can think of the operation as a rule for combining two elements of $S$ into one element of $S$. Note, however, that the order in which the two elements are selected may affect the outcome of the operation. The prototypical examples of binary operations are addition and multiplication on the set of natural numbers.
 
Throughout the course we will use \col{infix notation}\index{infix notation}, such as $2+3$, rather than \col{function notation}, such as $+(2,3)$, when working with binary operations. Thus, if\, $\square$\, is a binary operation on a set $S$ (where the square represents any symbol used to denote a binary operation), we'll write $a\,\square\, b$ instead of $\square(a,b)$ for elements $a,b\in S$.

The objects we're aiming to define in this chapter, namely groups and rings, involve binary operations with specific properties defined below. If $*$ is a binary operation on a set $S$, we say that $*$ is \col{associative}\index{associative operation} if it satisfies \[a*(b*c)=(a*b)*c\quad\text{for all}\;\;a,b,c\in S.\] Associativity allows us to unambiguously interpret expressions such as $a*b*c$.

We say that $\ast$ is \col{commutative}\index{commutative operation} if \[a\ast b=b\ast a\quad\text{for all}\;\;a,b\in S.\] Commutativity allows us to ignore the order in which elements are listed when performing operations. Elements $x,y\in S$ are said to \col{commute} if $x\ast y=y\ast x$. Consequently, the operation $\ast$ is commutative if and only if every pair of elements of $S$ commutes.

An element $e\in S$ is called an \col{identity element}\index{identity element} of $\ast$, or simply an \col{identity} of $*$, if it satisfies \[e\ast x=x\ast e=x\quad\text{for every}\;\; x\in S.\] Such an element need not exist, but if $*$ does have an identity element, then it has exactly one identity; thus, we may unambiguously refer to \emph{the} identity of the operation. The proof of uniqueness of identities is left to the reader.
\newpage
\begin{ex}[Binary operations]\;
\begin{examples}
\item Addition and multiplication on $\N$ are both associative and commutative. The identity element of addition is 0, and that of multiplication is 1.
\item Subtraction is a binary operation on $\Z$. It has no identity element, and is neither associative nor commutative.
\item Let $X$ be a set and let $P(X)$ denote the power set of $X$, i.e., the set of all subsets of $X$. We may consider intersection as a binary operation on $P(X)$, given that the intersection of any two subsets of $X$ is again a subset of $X$. This operation is commutative, associative, and has an identity element, namely $X$. Similarly, union is a commutative and associative binary operation on $P(X)$, and its identity element is the empty set, $\emptyset$.
\item Let $X$ be a set and let $S=X^X$ be the set of all functions $X\to X$. Composition of functions is then a binary operation on $S$, for if $f\in S$ and $g\in S$ then $f\circ g$ is a function $X\to X$, so $f\circ g\in S$. This operation is associative but will generally not be commutative. It does have an identity element, namely the identity function on $X$ given by $e(x)=x$ for all $x\in X$.
\end{examples}
\end{ex}

Suppose that $\ast$ is a binary operation on a set $S$, and that $e\in S$ is an identity of the operation (and therefore the only identity). An element $x\in S$ is called \col{invertible}\index{invertible element} if there exists an element $y\in S$ such that \[x\ast y=y\ast x=e.\] In that case, $y$ is called an \col{inverse}\index{inverse element} of $x$. For example, $3$ and $-3$ are inverses with respect to addition on $\Z$. However, the same is not true in $\N$ because $3$ has no inverse with respect to addition in $\N$.

In general, it is possible for an element of $S$ to have more than one inverse, but if $\ast$ is associative then inverses are unique; see Exercise \ref{inverse_unique}. The bulk of the operations we'll be considering in this course are in fact associative, so it's safe to assume that invertible elements will always have unique inverses.

\section*{Groups and rings}

At this point we have all the concepts needed to define groups and rings. A \col{group}\index{group} is formally defined as an ordered pair $(G,\ast)$ consisting of a set $G$ and a binary operation $*$ on $G$ satisfying the \col{group axioms}\index{group axioms}, namely:
\begin{itemize}
\item The operation $\ast$ is associative and has an identity element.
\item Every element of $G$ is invertible (and thus has a unique inverse).
\end{itemize}

If the operation $*$ is commutative, then $G$ is called a \col{commutative group}\index{commutative group} or an \col{Abelian group}\index{Abelian group}, the latter in honor of \href{https://en.wikipedia.org/wiki/Niels_Henrik_Abel}{Niels Henrik Abel} (1802-1829).

\begin{ex}[Groups]\;
\begin{examples}
\item A standard example of a group is $(\Z,+)$, which is an Abelian group since addition is a commutative operation. In this group the identity element is $0$, and the inverse of an integer $n$ is its negative, $-n$.
\item The set $\Z$ is not a group under multiplication. Even though multiplication is associative and has an identity, not every integer has a multiplicative inverse (for instance, 2 is not invertible since there is no integer $n$ such that $2n=1$.)
\item For an example of a group that is not commutative, let $M$ be the set of all $2\times 2$ matrices with real entries and nonzero determinant, and let $*$ denote matrix multiplication. Then $(M,*)$ is a group in which the rule $ab=ba$ does not apply to all elements $a,b$.
\end{examples}
\end{ex}

A \col{ring}\index{ring} is a triple $(R,+,\cdot)$ consisting of a set $R$ and two binary operations on $R$, called addition and multiplication, satisfying the \col{ring axioms}\index{ring axioms}, namely:
\begin{itemize}
\item The pair $(R,+)$ is an Abelian group.
\item Multiplication is associative and has an identity element.
\item The following distributive properties hold: for all $a,b,c\in R$,
\[a\cdot (b+c)=(a\cdot b)+(a\cdot c)\quad\text{and}\quad (a+b)\cdot c=(a\cdot c)+(b\cdot c).\]
\end{itemize}

The identity elements for addition and multiplication in $R$ are called, respectively, the \col{additive identity} and the \col{multiplicative identity} of $R$.

Since $R$ is a group under addition, every element $x\in R$ has a unique inverse with respect to addition; this element is denoted by $-x$ and called the \col{additive inverse} of $x$. If $x$ has an inverse with respect to multiplication (which is not always the case), this inverse is called the \col{multiplicative inverse} of $x$ and denoted by $x^{-1}$.

Though less common, a broader concept of ring exists in which multiplication is not assumed to have an identity element. From that point of view, the structure defined by the axioms above is called a \col{ring with unity}, the term \emph{unity} referring to the identity of multiplication. As optional reading, see the article of Poonen \cite{Poonen_unity}, which provides several arguments in support of assuming that all rings have a unity.

Note that addition is commutative by definition, but the same is not true for multiplication. A ring in which multiplication is commutative is called a \col{commutative ring}\index{commutative ring}. 
The identity element of addition in a ring will be denoted by 0, and the multiplicative identity by 1. If the particular ring, say $R$, needs to be specified for clarity, we may also write $0_R$ and $1_R$. It is important to keep in mind, however, that these elements may not be the integers 0 and 1. For instance, in a ring of matrices, the symbol 1 may denote the identity matrix (the diagonal matrix with ones along the diagonal).

\begin{ex}[Rings]\;
\begin{examples}
\item Standard examples of rings include the sets $\Z$, $\Q$, $\R$, and $\C$, each with their usual operations of addition and multiplication. All of these rings are commutative. 
\item For an example of a non-commutative ring, consider the set of all $2\times 2$ matrices with real entries, the operations being the usual addition and multiplication of matrices. In this ring, the element 1 is the identity matrix and the element 0 is the zero matrix.
\item Here is a rather different example of a ring. Let $X$ be a set and $S=P(X)$ its power set. We will form a ring with $S$ as its underlying set. The operation of multiplication on $S$ is set intersection, and the operation of addition is the symmetric difference of sets, defined by \[a\,\Delta\, b=(a\cup b)\setminus(a\cap b)\] for subsets $a,b\subseteq X$. Though it takes some work to check all the ring axioms, basic set theory can be used to prove that $(S,\Delta,\cap)$ is a commutative ring in which the additive identity is $\emptyset$ and the multiplicative identity is $X$. Moreover, the additive inverse of a set $a\in S$ is its complement: $-a=X\setminus a$.
\end{examples}
\end{ex}

We will be studying rings in much greater detail throughout the course, but it is sufficient for now to understand the definition of this concept, as well as that of a group, and to be able to use the group and ring axioms to determine whether a particular structure is a group or a ring.

\section*{Exercises}

\begin{note}\;
\begin{itemize}
\item Solutions to starred exercises (\sel) are included in Chapter~\ref{review1chap}.
\item A house (\house) marks exercises suggested as sample homework.
\end{itemize}
\end{note}

\medskip

\begin{exercises}
\item Consider the following two binary operations on $\Z$:
\[
a\circ b=a\cdot b+a-b,\quad\text{and}\quad a\ast b=b.
\]
For each operation, determine whether it is commutative, associative, has an identity element, and which elements (if any) are invertible.
\item\label{inverse_unique} Let $\ast$ be an associative binary operation on a set $S$, and let $x\in S$ be an invertible element. Show that $x$ has a unique inverse in $S$.
\item(\house) Consider the set
\[
X=\{(a,b)\in\R^2:a\ne 0\}
\]
and the binary operation on $X$ given by
\[
(a,b)\ast(c,d)=(ac,ad+b).
\]
Determine whether $(X,\ast)$ is a group.
\item(\sel)\label{chap1_sel} Let $(R,+,\cdot)$ be a ring and $x,y \in R$. Use the ring axioms to prove the following:
\begin{enumerate}
\item $0 \cdot x = x \cdot 0 = 0$,
\item $(-x)\cdot y = x \cdot (-y) = -(x \cdot y)$,
\item $(-x)\cdot (-y) = x \cdot y$.
\end{enumerate}
\item(\house) Define operations on $\R^2$ by $(a,b)+(c,d)=(a+c,b+d)$ and
\[
(a,b)\odot(c,d)=(ac-bd,ad+bc).
\]
Determine whether $(\R^2,+,\odot)$ is a ring.
\item Consider a set $(R,+,\cdot)$ endowed with two binary operations satisfying all the ring axioms, except that addition is not assumed to be commutative. Show that $+$ is necessarily commutative.
\end{exercises}

\chapter{Introduction to Rings}\label{rings_intro_chap}

As the formal definition of a ring in terms of the ring axioms was given in Chapter 1, we begin our exploration of rings with an informal description: a ring is a set whose elements can be added and multiplied in such a way that the rules of elementary algebra are satisfied. Examples of such rules are
\[x+0=x,\;\; x\cdot 0=0,\;\;x\cdot 1=x,\;\;x+y=y+x,\;\;x\cdot (y+z)=x\cdot y+x\cdot z.\]

In order to ease notation, we will often write $xy$ instead of $x\cdot y$. The additive identity in a ring may also be called the \col{zero element}\index{zero element} of the ring, and the multiplicative identity the \col{unity}\index{unity of a ring} of the ring. Though a ring was formally defined as a triple $(R,+,\cdot)$, we typically refer to a ring only by its underlying set, $R$, with the understanding that this set comes equipped with two operations.

In any ring we may define an operation of \col{subtraction}\index{subtraction in a ring} by \[x-y:=x+(-y)\quad\text{for all elements}\; x,y.\] Recall that $-y$ denotes the additive inverse of $y$, which exists as a consequence of the ring axioms. With this definition, multiplication distributes with respect to subtraction, so we can append the following rules to the above list:
\[x(y-z)=xy-xz\quad\text{and}\quad (y-z)x=yx-zx.\]

\begin{ex}[Rings] All of the following examples are instances of more general constructions to be formally discussed later in the course.
\begin{examples}
\item For a positive integer $n$, we denote by $\Z_n$ the ring of integers modulo $n$. A simple way to describe this ring is to define \[\Z_n:=\{0,1,2,\ldots, n-1\}\] and define addition and multiplication on this set as the usual operations in $\Z$ followed by reduction modulo $n$. For example, to evaluate $x=2\cdot 3+7$ in $\Z_8$, we compute $2\cdot 3+7=6+7=13$ in $\Z$ and then reduce modulo 8 to obtain $x=5$. Every ring $\Z_n$ is commutative because $\Z$ is commutative.
\item For a positive integer $n$, the ring of $n\times n$ matrices with real entries will be denoted by $M_n(\R)$. Addition and multiplication in this ring are the usual operations on matrices from linear algebra. The zero element in this ring is the matrix of all zeros, and the unity is the identity matrix, denoted by $I_n$. The rings $M_n(\R)$ are non-commutative for $n>1$. Note that if $m\ne n$, the set of $m\times n$ matrices does not form a ring because multiplication of two such matrices is undefined. 
\item The above observations can also be used to define rings $M_n(\C)$ for $n\ge 1$, whose elements are matrices with complex entries. Once again, these rings are non-commutative if $n>1$.
\item The ring of polynomials with real coefficients is denoted by $\R[x]$. Elements of this ring have the form \[a_0+a_1x+a_2x^2+\cdots+a_nx^n,\] where the $a_i$'s are real numbers called the \col{coefficients} of the polynomial. Having $n=0$ is allowed, so every real number can be viewed as a polynomial called a \col{constant polynomial}. Operations in $\R[x]$ are carried out by distributing and combining like terms; for instance, \[x(2-x^2)+3x=2x-x^3+3x=5x-x^3.\] The zero element and unity in $\R[x]$ are the constant polynomials 0 and 1. The ring $\R[x]$ is commutative.
\item We can similarly define the ring $\C[x]$, whose elements are polynomials with complex coefficients.
\end{examples}
\end{ex}

Theorem \ref{identi} below establishes several basic algebraic identities that hold in any ring, and which can be deduced from the ring axioms. These identities involve iterated, or repeated, addition and multiplication, defined as follows. Let $R$ be a ring. Given $a\in R$ and a natural number $n$, we define $a^n$ recursively by setting $a^0=1_R$ and $a^n=a\cdot a^{n-1}$ for $n\ge 1$. More concretely, \[a^n=a\cdot a\cdot a\cdots a\] with $n$ factors equal to $a$. Similarly, we define $na$ by setting $0\cdot a=0_R$ and $na=a+(n-1)a$, so that \[na=a+a+\cdots+a.\] With these definitions, all the identities in the following theorem can be proved by induction.

\begin{framed}
\begin{thm}\label{identi} Let $R$ be a ring. The following identities hold for all elements $a,b\in R$ and all $m,n\in\N$:
\begin{enumerate}
\item[(a)] \;$a^ma^n=a^{m+n}$ \;and\; $(a^m)^n=a^{mn}$.
\item[(b)] \;$n(a+b)=na+nb$,\; $(n+m)a=na+ma$,\; $(mn)a=m(na)$.
\item[(c)] (The binomial theorem) If $ab=ba$, then \[(a+b)^n=\sum_{i=0}^n{n\choose i}a^ib^{n-i}.\]
\end{enumerate}
\end{thm}
\end{framed}

We arrive now at two fundamental concepts in ring theory: units and zero divisors. If $R$ is a ring, a \col{unit}\index{unit in a ring} in $R$ is an element having a multiplicative inverse. Thus, $a\in R$ is a unit if there exists an element $b\in R$ such that $ab=ba=1$. Such an element $b$ is unique and will be denoted by $a^{-1}$. The set of all units in $R$ is denoted by $R^{\times}$ and called the \col{unit group}\index{unit group} of $R$. As the name suggests, $R^{\times}$ is a group under the operation of multiplication in $R$.
\newpage
\begin{ex}[Units]\;
\begin{examples} 
\item The unit group of $\Z$ is $\Z^{\times}=\{\pm 1\}$. 
\item In $\R$, every nonzero element is a unit.
\item In $\Z_9$, 2 is a unit and its inverse is $2^{-1}=5$.
\item In $M_n(\R)$, the units are the matrices with nonzero determinant.
\item The unit group of $\R[x]$ is $\R[x]^{\times}=\R^{\times}=\R\setminus\{0\}$.
\end{examples}
\end{ex}

A \col{zero divisor}\index{zero divisor} in a ring $R$ is a \underline{nonzero} element $a\in R$ such that there exists a \underline{nonzero} $b\in R$ with $ab=0$ or $ba=0$. This may seem like an odd concept at first if your intuition is based on rings such as the integers or the real numbers, in which the only way for a product of elements to be 0 is for one of the elements to be 0. However, for arbitrary rings this is not necessarily the case.

\begin{ex}[Zero divisors]\;
\begin{examples} 
\item In the rings $\Z$ and in $\R[x]$ there are no zero divisors.
\item In $\Z_6$, 2 is a zero divisor since $2\cdot 3=0$. For the same reason, 3 is a zero divisor.
\item In $M_2(\R)$, the matrix $\begin{pmatrix}1&1\\2&2\end{pmatrix}$ is a zero divisor since
\[\begin{pmatrix}1&1\\2&2\end{pmatrix}\begin{pmatrix}1&1\\-1&-1\end{pmatrix}=\begin{pmatrix}0&0\\0&0\end{pmatrix}.\]
\end{examples}
\end{ex}

\begin{note} In the context of non-commutative rings, it is appropriate to distinguish between \emph{left} zero divisors and \emph{right} zero divisors, depending on the side on which an element $x$ can be multiplied by some $y$ to obtain 0. In this course, the above distinction will not be needed, so a zero divisor will refer to an element that is either a left zero divisor or a right zero divisor (or both).
\end{note}

Units and zero divisors will be relevant to our discussions throughout the duration of the course, so there will be ample opportunity to develop familiarity with these types of elements. For now, we just need one observation about them, which is that units are never zero divisors:

\begin{framed}
\begin{thm}Let $R$ be a ring. The sets of units and zero divisors in $R$ are disjoint.\end{thm}
\end{framed}
\begin{proof}
Let $u\in R$ be a unit. Assuming that $u$ is a zero divisor, we seek a contradiction. By definition of zero divisor, there is a nonzero element $x\in R$ such that $ux=0$ or $xu=0$. We assume the latter case, with the former being handled similarly. Multiplying by $u^{-1}$ and using associativity of multiplication, we have $(xu)u^{-1}=0\cdot u^{-1}=0$, so $x(uu^{-1})=0$. Since $uu^{-1}=1$, this implies $x=0$, yielding the desired contradiction.
\end{proof}

We now define integral domains and fields, two very commonly referenced types of rings. An \col{integral domain}\index{integral domain}, often also called a \col{domain}\index{domain}, is a commutative ring with no zero divisors, and in which $0\ne 1$. (If the latter condition seems strange, we will soon see that it can occur that $0=1$ in a ring.) Equivalently, a domain is a commutative ring $R$ in which $0\ne 1$ and where the following implication holds true: for all $a,b\in R$, \[ab=0\implies a=0\;\;\text{or}\;\; b=0.\] From this fact we deduce the \col{cancellation property}\index{cancellation property} of integral domains: if $a,b,c$ belong to a domain and $ab=ac$ with $a\ne 0$, then $b=c$. The proof of this property is left to the reader. The term \emph{integral domain} applied to a ring is meant to indicate that the ring is in some ways similar to the ring of integers, $\Z$. We will see later that integral domains provide a good framework in which to develop concepts such as \emph{prime element}, generalizing the familiar concept of a prime number in $\Z$. 

\begin{ex}[Integral domains]\;
\begin{examples}
\item The rings $\Z$ and $\R[x]$ are domains.
\item The rings $M_n(\R)$ are not domains for $n>1$, as domains are by definition commutative.
\item The ring $\Z_p$ is a domain if and only if $p$ is prime. For one direction of this claim, note that if $p$ factors as $p=ab$ with $1<a,b<p$, then $a$ and $b$ are nonzero elements of $\Z_p$ such that $ab=0$ in $\Z_p$. Hence, $\Z_p$ is not a domain.
\end{examples}
\end{ex}

\begin{note} When $p$ is prime, the notation $\Z_p$ conflicts with standard notation used for the ring of \emph{p-adic integers}, a ring commonly used in number theory. For this reason, we introduce the notation $\F_p$ to denote the ring $\Z_n$ when $n=p$ is prime. This choice of notation is due to the fact that $\F_p$ is a \emph{field}, as explained below.
\end{note}

A \col{field}\index{field} is a commutative ring in which $0\ne 1$ and in which every nonzero element is a unit. Equivalently, a field is an integral domain in which every nonzero element is a unit. 

On any field $F$ there is an operation of \col{division}\index{division in a field} defined by setting \[\frac{x}{y}:=xy^{-1}\;\;\text{for all}\;\;x,y\in F\;\;\text{with}\;\; y\ne 0.\]
The notation $x/y$ may also be used instead of $\frac{x}{y}$. The following elementary identities are then satisfied:
\[\frac{a}{b}\cdot\frac{c}{d}=\frac{ac}{bd},\quad\frac{a}{b}+\frac{c}{d}=\frac{ad+bc}{bd}.\]

We may therefore informally describe fields as sets in which the four basic arithmetic operations ($+,-,\cdot,/$) can be performed and satisfy the rules of elementary algebra. 

\begin{ex}[Fields]\;
\begin{examples}
\item The rings $\Q,\R,$ and $\C$ are fields.
\item The rings $\Z$ and $\R[x]$ are domains but not fields.
\end{examples}
\end{ex}

The above examples in (b) show that, although every field is a domain, the converse statement does not hold. Interestingly, though, the converse \emph{is} true in the finite case:

\begin{framed}
\begin{thm}\label{finite_domain} Every finite integral domain is a field.\end{thm}
\end{framed}

\begin{proof} We use the following basic fact from set theory: if $X$ is a finite set and $f:X\to X$ is an injective map, then $f$ is also surjective. Suppose that $F$ is a finite integral domain and let $a\in F$ be nonzero. Since $F$ is a domain, the cancellation property implies that the map $F\to F$ given by $x\mapsto ax$ is injective. Moreover, since $F$ is finite, the map is also surjective. In particular, there exists $b\in F$ such that $ab=1$, proving that $a$ is a unit.
\end{proof}

One important consequence of Theorem \ref{finite_domain} is that the ring $\F_p$ is a field, as claimed earlier. Once the proper concepts are in place, we will be able to prove that (in some precise sense) $\F_p$ is the \emph{only} field with exactly $p$ elements.

We end with a simple example of a computation in the field $\F_7$:\[\frac{2}{3}+\frac{1}{2}\cdot\frac{3}{4}=\frac{2}{3}+\frac{3}{1}=2\cdot 3^{-1}+3.\]

Since $3^{-1}=5$ in $\F_7$, we compute $2\cdot 5+3=3+3=6$.

\section*{Exercises}

\begin{note}\;
\begin{itemize}
\item Solutions to starred exercises (\sel) are included in Chapter~\ref{review1chap}.
\item A house (\house) marks exercises suggested as sample homework.
\end{itemize}
\end{note}

\medskip

\begin{exercises}
\item Find $a^{-1}$ for every nonzero $a\in\F_{11}$. Compute $2/5-7/9$ in $\F_{11}$.
\item Find an example of a ring and an element in the ring that is neither a unit nor a zero divisor.
\item(\sel)\label{chap2_sel} If $a$ is an element of an integral domain such that $a^2=1$, show that $a=1$ or $a=-1$. Find a counterexample to this statement in some ring that is not a domain.
\item Let $x,y$ be units in a ring. Show that $(xy)^{-1}=y^{-1}x^{-1}$.
\item(\house)\label{Zpdef} Let $p$ be a prime integer and let $\Z_{(p)}$ denote the set of all rational numbers of the form $a/b$ with $\gcd(a,b)=1$ and $b$ not divisible by $p$. Show that $\Z_{(p)}$ is a ring with addition and multiplication inherited from $\Q$, and determine its unit group. Is $\Z_{(p)}$ a domain? Is it a field?
\item\label{idempotent_def} An element $x$ of a ring is called \col{idempotent}\index{idempotent} if $x^2=x$. Show that if $x\ne 1$ is idempotent, then $x$ is a zero divisor.
\item Let $R$ be a ring. For every integer $n<0$ and every $x\in R$, we define $nx$ to be the additive inverse of $(-n)x$. Prove the following for all $m,n\in\Z$ and all $a,b\in R$: $n(a+b)=na+nb$, $(n+m)a=na+ma$, and $(mn)a=m(na)$.
\item 
\begin{enumerate}
\item Let $p$ and $q$ be distinct prime integers. Show that the map $\F_q\to\F_q$ given by $x\mapsto px$ is bijective.
\item Let $n$ be a positive integer. Use part (a) to show that $\Z_n^{\times}$ is the set of elements $a\in\Z_n$ such that $\gcd(a,n)=1$.
\end{enumerate}
\item Let $F$ be a field and $x\in F\setminus\{1\}$. Prove the identity \[1+x+x^2+\cdots+x^{n-1}=(x^n-1)/(x-1)\] for every integer $n\ge 1$.
\item(\sel)\label{nilpotent_def} An element $x$ of a ring is called \col{nilpotent}\index{nilpotent element} if there exists an integer $n\ge 1$ such that $x^n=0$. Show that if $x$ is nilpotent then $1+x$ is a unit.
\item(\house)\label{boolean_def} A \col{Boolean ring}\index{Boolean ring} is a ring in which every element is idempotent. Show that every Boolean ring is commutative, and find an example of a Boolean ring among the rings already discussed.
\item Let $p$ be prime. Use the binomial theorem to prove the identity \[(x+y)^p=x^p+y^p\] for all $x,y\in\F_p$.
\end{exercises}

\chapter{Examples of Rings}\label{more_rings_chap}

In this chapter we focus on enlarging our catalog of known rings by studying several examples that may be less familiar than those discussed up to this point.

\subsection*{The trivial ring} Let $A=\{\bullet\}$ be a singleton set. There is only one binary operation on $A$, namely the map $(\bullet,\bullet)\mapsto\bullet$. Taking both $+$ and $\cdot$ to be this binary operation, the triple $(A,+,\cdot)$ satisfies all the ring axioms and is therefore a ring. This ring will be called the \col{trivial ring}\index{trivial ring} or the \col{zero ring}\index{zero ring}. Note that necessarily $0=1$ in this ring, as $A$ has only one element. In fact, any ring in which the zero element and the unity agree must be a trivial ring, as follows from the identities $x=x\cdot 1=x\cdot 0=0$. It will be useful in some definitions and theorems to exclude trivial rings from consideration; for instance, we may define a field as a nontrivial commutative ring in which every nonzero element is a unit.

\subsection*{Quadratic integer rings} Recall that a nonzero integer is called \col{squarefree}\index{squarefree} if it is not divisible by the square of any prime number. For example, 18 is not squarefree since it is divisible by $3^2$, but $6=2\cdot 3$ is squarefree. If $d\ne 1$ is any squarefree integer, we define
\[\Z[\sqrt d]=\{a+b\sqrt d:a,b\in\Z\},\] where $\sqrt d$ denotes a square root of $d$ in $\C$. Note that $\Z[\sqrt d]$ is contained in $\C$, and will be contained in $\R$ if $d>0$. The set $\Z[\sqrt d]$, equipped with the usual operations in $\C$, is a ring. (All the ring axioms follow from those in $\C$, so the only key point to make is that $\Z[\sqrt d]$ is closed under these operations.) Rings of the form $\Z[\sqrt d]$ are known as \col{quadratic integer rings}\index{quadratic integer rings}. In the special case of $d=-1$, we obtain the ring \[\Z[i]=\{a+bi:a,b\in\Z\},\] known as the \col{ring of Gaussian integers}\index{ring of Gaussian integers}. 

Quadratic integer rings are always domains but never fields; see Exercise \ref{quad_notfield}. However, each ring $\Z[\sqrt d]$ is contained in a slightly larger ring that \emph{is} a field, namely the set \[\Q(\sqrt d)=\{a+b\sqrt d: a,b\in\Q\}.\]

We can prove that $\Q(\sqrt d)$ is a field by using the common algebraic trick of rationalizing the denominator: if $a,b\in\Q$ are not both zero, then
\[\frac{1}{a+b\sqrt d}=\frac{1}{a+b\sqrt d}\cdot\frac{a-b\sqrt d}{a-b\sqrt d}=\frac{a-b\sqrt d}{a^2-b^2d}.\]
Letting $P=a/(a^2-b^2d)$ and $R=-b/(a^2-b^2d)$, the above equations imply that $P+R\sqrt d$ is a multiplicative inverse of $a+b\sqrt d$; hence every nonzero element of $\Q(\sqrt d)$ is a unit.

Fields of the form $\Q(\sqrt d)$ are called \col{quadratic number fields}\index{quadratic number field}. In the special case of $d=-1$ we obtain the field \[\Q(i)=\{a+bi:a,b\in\Q\}.\]

For any element $x=a+b\sqrt d\in\Q(\sqrt d)$, we define the \col{conjugate}\index{conjugate!in a quadratic number field} and the \col{norm}\index{norm!in a quadratic number field} of $x$, respectively, by \[\bar x=a-b\sqrt d\;\;\text{and}\;\;N(x)=x\bar x=a^2-b^2d.\] These concepts are quite useful for studying the properties of quadratic rings and fields; for example, the next theorem uses norms to describe the unit groups of quadratic rings. A proof of the theorem is left to the reader (Exercise \ref{quad_notfield}).

\begin{framed}
\begin{thm}\label{quad_int_unit}
Let $d\ne 1$ be a squarefree integer. 
\begin{enumerate}
\item[(a)] The norm map on $\Q(\sqrt d)$ is multiplicative: for all $x,y\in\Q(\sqrt d)$, \[N(xy)=N(x)N(y).\]
\item[(b)] An element $x\in\Z[\sqrt d]$ is a unit if and only if $N(x)=\pm 1$.
\end{enumerate}
\end{thm}
\end{framed}

\subsection*{Rings of functions}\index{ring of functions} If $R$ is any ring, we may utilize the structure of $R$ to create new rings consisting of functions whose codomain is $R$. The precise construction is explained below.

Letting $X$ be an arbitrary nonempty set, the set of all functions $X\to R$, denoted by $R^X$, can be made into a ring by using the operations of pointwise addition and multiplication: for functions $f,g\in R^X$, the functions $f+g$ and $fg$ are defined, respectively, by \[x\mapsto f(x)+g(x)\quad\text{and}\quad x\mapsto f(x)\cdot g(x).\] Note that $f(x)$ and $g(x)$ are added and multiplied using the operations of $R$. In this sense, $R^X$ inherits its ring structure from that of $R$. The zero element of $R^X$ is the constant map $x\mapsto 0$ for all $x\in X$, and the unity of $R^X$ is the constant map $x\mapsto 1$ for all $x\in X$.

As a variant of the above construction, we may restrict attention to functions with specified properties. One example, which will be referenced several times throughout the course, is the following: if $a<b$ are real numbers, let $\calC[a,b]$ denote the set of all continuous functions $[a,b]\to\R$, a type of function studied in one-variable calculus. The fact that $\calC[a,b]$ satisfies all the ring axioms follows easily from the known theorem that the sum and product of continuous functions are continuous. The ring $\calC[a,b]$ is not a domain and therefore not a field. One way to prove this is to construct two nonzero continuous functions such that each one is zero whenever the other is not; the product of the two functions will then be zero.

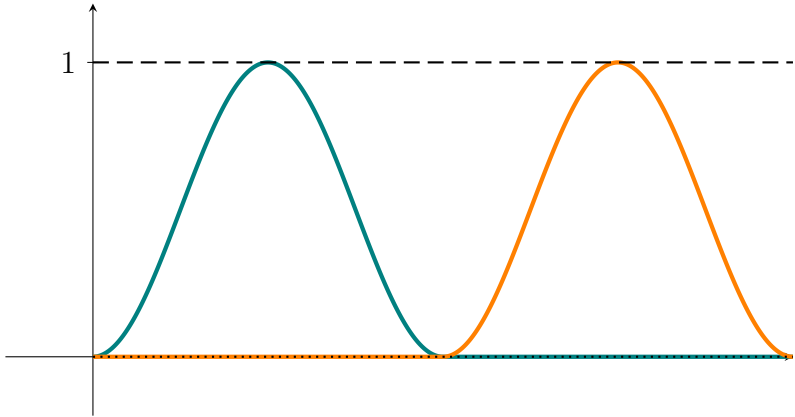
\begin{figure}[ht]
\centering
\begin{tikzpicture}
\begin{axis}[
    width=12cm,
    height=7cm,
    axis lines=middle,
    ymin=-0.2, ymax=1.2,
    xmin=-0.5, xmax=4,
    xtick=\empty,
    ytick={0,1},
    yticklabels={$0$,$1$},
    samples=200,
    tick label style={font=\normalsize, color=black},
    axis line style={black},
    tick style={black}
]

\addplot[domain=0:2, fone] {sin(deg(pi*x/2))^2};
\addplot[domain=2:4, fone] {0};

\addplot[domain=0:2, ftwo] {0};
\addplot[domain=2:4, ftwo] {sin(deg(pi*(x-2)/2))^2};

\addplot[domain=0:4, funity] {1};

\addplot[domain=0:4, fzero] {0};

\node[black] at (axis cs:4.05,1) {\small $1$};

\node[black] at (axis cs:4.05,0) {\small $0$};

\end{axis}
\end{tikzpicture}
\caption{The zero and unity elements of $\calC[a,b]$, together with two nonzero continuous functions 
$f_1$ and $f_2$ such that $f_1(x)f_2(x)=0$ for all $x$.}
\label{cont_zerodiv}
\end{figure}

Elements of $\calC[a,b]$ can be visualized via their graphs as continuous curves over the interval $[a,b]$. Figure \ref{cont_zerodiv} shows the zero element and unity of $\calC[a,b]$, together with a pair of zero divisors.

\subsection*{The ring of quaternions} The following example is perhaps the least intuitive among those we have seen thus far. The ring described here is called the \col{ring of Hamilton quaternions}\index{ring of quaternions}, and denoted by $\H$ in honor of \href{https://en.wikipedia.org/wiki/William_Rowan_Hamilton}{William Rowan Hamilton} (1805-1865), who first described it. 

Elements of $\H$, called \col{quaternions}\index{quaternions}, can be represented in the form \[a+bi+cj+dk,\] where $a,b,c,d\in\R$ and $i,j,k$ are specific elements of $\H$ satisfying the relations \begin{equation}\label{ham_eq}\tag{$\ast$}i^2=j^2=k^2=ijk=-1.\end{equation}

Note, in particular, that $\H$ contains $\R$ as the set of quaternions with $b=c=d=0$. 

Addition on $\H$ is quite natural and can be described as componentwise addition. For example, if $x,y\in\H$ are defined by \[x=2-3i+5j-k,\;\; y=7i-j+4k,\] then \[x+y=2+4i+4j+3k.\] 

\begin{figure}[ht]
\centering
\begin{tikzpicture}[scale=2, thick, >=Stealth]

  \draw[white!50!black] (0,0) circle(1);

  \node at (90:1.15) {\(i\)};
  \node at (330:1.15) {\(j\)};
  \node at (210:1.15) {\(k\)};

  \foreach \angle in {150, 270, 30} {
    \coordinate (A) at (\angle:1);
    \coordinate (B) at (\angle + 90:0.07);  
    \draw[->] (A) ++(B) -- (A);
  }

\end{tikzpicture}
\caption{Quaternion multiplication cycle.}
\label{quat_cycle}
\end{figure}
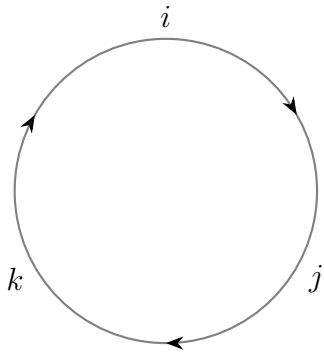

To multiply elements of $\H$, we first need to know the products of any two among $i,j,k$. These products can be deduced from the quaternion relations \eqref{ham_eq}. For example, we can obtain $ij$ as follows:
\[ijk=-1\implies ijk^2=-k\implies ij(-1)=-k\implies ij=k.\]
A simple way to remember these basic products is to place $i,j,k$ along a circle in a clockwise direction as in Figure \ref{quat_cycle}. When arranged this way, the product of two elements in a clockwise direction is the third element, and the product of two elements in a counterclockwise direction is the negative of the third element. For example, \[ki=j\;\;\text{and}\;\; ik=-j.\]

We can now compute any product in $\H$ by distributing and using the known products of $i,j,k$. For instance,
\[(2+3j)(5i-k)=10i-2k+15(-k)-3(i)=7i-17k.\]

In order to analyze the sets of units and zero divisors in $\H$, we will need two concepts analogous to those defined earlier for quadratic integer rings. If $q=a+bi+cj+dk$ is a quaternion, its \col{conjugate}\index{conjugate!of a quaternion} is defined to be the quaternion \[\bar q=a-bi-cj-dk.\] The product of any quaternion and its conjugate turns out to be a non-negative real number:
\[q\bar q=\bar q q=a^2+b^2+c^2+d^2.\]

The square root of the above real number is called the \col{norm}\index{norm!of a quaternion} of $q$ and denoted by $\|q\|$. Thus, \[\|q\|=\sqrt{a^2+b^2+c^2+d^2}\quad\text{if}\quad q=a+bi+cj+dk.\]

We can now prove that every nonzero element of $\H$ is a unit. Let \[q=a+bi+cj+dk\] be nonzero, so that $a,b,c,d$ are not all zero and therefore $\|q\|\ne 0$. Defining $p=\bar q/\|q\|^2$, we observe that \[pq=qp=q\bar q/\|q\|^2=1,\] so that $p$ is a multiplicative inverse of $q$ and thus $q\in\H^{\times}$. This proves the claim.

Note that it is not accurate to call $\H$ a field since it is not commutative; instead, $\H$ is an example of a \col{division ring}\index{division ring}, a ring in which every nonzero element is a unit. (Thus, fields are precisely the nontrivial commutative division rings.)

From a theoretical standpoint, it is important to prove that there does indeed exist a ring whose elements can be represented in the form $a+bi+cj+dk$, where $i,j,k$ satisfy \eqref{ham_eq}. We complete our discussion of quaternions by providing a rigorous definition of the ring $\H$ as a subset of $M_2(\C)$. 

We define $\H$ to be the set of matrices of the form \[\begin{pmatrix}z&w\\-\bar w&\bar z\end{pmatrix},\] where $z,w\in\C$ and $\bar z$ denotes the standard complex conjugate of $z$. This set can be shown to form a ring with the operations of matrix addition and multiplication; moreover, defining $i,j,k\in\H$ by
\[i=\begin{pmatrix}\sqrt{-1}&0\\0&-\sqrt{-1}\end{pmatrix},\;\;j=\begin{pmatrix}0&1\\-1&0\end{pmatrix},\;\;k=\begin{pmatrix}0&\sqrt{-1}\\\sqrt{-1}&0\end{pmatrix},\]
every element of $\H$ has a unique representation in the expected form $a+bi+cj+dk$, and indeed, $i^2=j^2=k^2=ijk=-1$.

\subsection*{Matrix rings}

The last example in this chapter is a straightforward generalization of the rings $M_n(\R)$ already introduced. If $R$ is an arbitrary ring and $n$ a positive integer, the set $M_n(R)$ of $n\times n$ matrices with entries in $R$ admits a natural ring structure with the usual addition and multiplication of matrices. Explicitly, if $A=(a_{ij})$ and $B=(b_{ij})$ are $n\times n$ matrices, then
\[A+B=(a_{ij}+b_{ij})\quad\text{and}\quad AB=(c_{ij}), \quad\text{where}\quad c_{ij}=\sum_{k=1}^na_{ik}b_{kj}.\] 

The zero element of $M_n(R)$ is the zero matrix, and the unity is the identity matrix. The rings $M_n(R)$ are typically not commutative even if the ring $R$ is commutative. If $R$ is a nontrivial ring and $n>1$, then $M_n(R)$ necessarily contains zero divisors; for instance, a matrix with exactly one nonzero entry is a zero divisor.

The remaining comments are non-essential to this course, but may be of interest to the reader familiar with the concepts from linear algebra being referenced. 

Standard concepts from linear algebra can be generalized to the context of matrix rings. The \col{trace}\index{trace of a square matrix} of a square matrix is defined as the sum of the diagonal entries. The trace function, $\text{tr}:M_n(R)\to R$, satisfies \[\text{tr}(A+B)=\text{tr}(A)+\text{tr}(B)\]
for all $A,B\in M_n(R)$.

If $R$ is commutative, the \col{determinant}\index{determinant of a square matrix} of a matrix in $M_n(R)$ is defined by the \href{https://en.wikipedia.org/wiki/Leibniz_formula_for_determinants}{Leibniz formula} and satisfies the expected property that \[\det(AB)=\det(A)\det(B)\]
for all $A,B\in M_n(R)$. Moreover, the \col{characteristic polynomial}\index{characteristic polynomial} of a matrix $A\in M_n(R)$ is defined as $\det(xI_n-A)$, where $I_n$ is the $n\times n$ identity matrix, and the classical \href{https://en.wikipedia.org/wiki/Cayley\%E2\%80\%93Hamilton_theorem}{Cayley-Hamilton theorem} applies.

 The following theorem describes the unit group of $M_n(R)$, denoted by $\GL_n(R)$, and called the \col{general linear group}\index{general linear group over a ring} of degree $n$ over $R$. This theorem generalizes the fact that a matrix with real entries is invertible if and only if its determinant is nonzero.

\begin{framed}
\begin{thm}\label{unitmatrix}
Let $R$ be a commutative ring. A matrix in $M_n(R)$ is a unit if and only if its determinant is a unit in $R$.
\end{thm}
\end{framed}
\begin{proof}
For the more substantial direction, suppose that $A\in M_n(R)$ is such that $\det(A)$ is a unit in $R$. To show that $A$ is a unit in $M_n(R)$, we use the \href{https://en.wikipedia.org/wiki/Adjugate_matrix}{adjugate matrix} of $A$, which satisfies the relations
\[A\cdot\adj(A)=\adj(A)\cdot A=\det(A)\cdot I.\] 
Letting $u=\det(A)^{-1}$, it follows that the matrix $u\cdot\adj(A)$ is the multiplicative inverse of $A$ in $M_n(R)$.
\end{proof}

Consider applying the above theorem to the matrix $A=\begin{pmatrix}17&29\\7&12\end{pmatrix}$. Since $\det(A)=1$ is a unit in $\Z$, then $A$ must have an inverse matrix whose entries are integers. The proof of Theorem \ref{unitmatrix} shows us how to find this inverse, yielding $A^{-1}=\begin{pmatrix}12&-29\\-7&17\end{pmatrix}$.

\section*{Exercises}

\begin{note}\;
\begin{itemize}
\item Solutions to starred exercises (\sel) are included in Chapter~\ref{review1chap}.
\item A house (\house) marks exercises suggested as sample homework.
\end{itemize}
\end{note}

\medskip

\begin{exercises}
\item \begin{enumerate}
\item Use Theorem \ref{quad_int_unit} to find all the units in the rings $\Z[i]$ and $\Z[\sqrt{-3}]$. 
\item(\house) If $d<0$ and $d\notin\{-1,-3\}$, show that $\Z[\sqrt d]^{\times}=\{\pm 1\}$.
\item(\sel)\label{chap3_sel1} Prove that the unit group of $\Z[\sqrt 2]$ is infinite.
\end{enumerate}
\item\label{quad_notfield} Prove Theorem \ref{quad_int_unit} and use it to show that no quadratic integer ring is a field.
\item(\house) In the ring $\calC[0,1]$, give an example of a nonzero function that is not a unit, and an example of a nonconstant function that is a unit. 
\item For $f\in\calC[a,b]$, show that $f$ is a zero divisor if and only if the set \[\{x\in [a,b]:f(x)=0\}\] contains an open interval. Use this fact to construct a nontrivial zero divisor in the ring $\calC[0,1]$.
\item(\sel)\label{chap3_sel2} Show that there exist infinitely many $x\in\H$ such that $x^2=-1$.
\item\begin{enumerate}
\item For all quaternions $x$ and $y$ show that $\overline{xy}=\bar y\bar x$.
\item For all $x,y\in\H$ show that $\|xy\|=\|x\|\cdot\|y\|$.
\item The integers 30 and 17 can be expressed as sums of four squares as follows: $30=1^2+2^2+3^2+4^2$ and $17=5^2+6^2+7^2+8^2$. Use properties of quaternions to find an expression for $510=30\cdot 17$ as a sum of four squares. For more on sums of four squares see \href{https://en.wikipedia.org/wiki/Lagrange\%27s\_four-square\_theorem}{Lagrange's four-square theorem}.
\end{enumerate}
\item Show that the determinant of the matrix $A=\begin{pmatrix}2&5\\8&6\end{pmatrix}$ in $M_2(\Z_9)$ is a unit. Use the proof of Theorem \ref{unitmatrix} to find the inverse of $A$.
\item Determine all elements of the group $\GL_2(\Z_4)$.
\item (\href{https://en.wikipedia.org/wiki/Cramer\%27s_rule}{Cramer's rule}) 
\begin{enumerate}
\item Let $R$ be a commutative ring, $A\in M_n(R)$ an $n\times n$ matrix, and $b\in R^n$ regarded as an $n\times 1$ matrix. If $\det(A)$ is a unit in $R$, show that there is a unique solution to the system of linear equations $Ax=b$, and it is given by 
\[x=\det(A)^{-1}\cdot(\adj A)\cdot b.\]
\item Show that the $i$th entry of the solution vector $x$ is given by $x_i=\frac{\det(A_i)}{\det(A)}$, where $A_i$ is the matrix obtained by replacing the $i$th column of $A$ with $b$.
\item Use the above rule to find $x,y,z\in\Z_{10}$ such that
\[\begin{cases}
2x-7y+4z&=3,\\
4x+9y+6z&=8,\\
6x+2y-z&=1.
\end{cases}
\]
\end{enumerate}
\end{exercises}

\chapter{Polynomial Rings}\label{poly_chap}

This chapter introduces a type of ring that plays a prominent role in ring theory. Our main goals here are to provide an informal description as well as a rigorous construction of these rings, and to explore some of their essential properties.

Let $A$ be a commutative ring. Informally, a \col{polynomial}\index{polynomial!in one variable} over $A$ can be defined as an expression of the form \[a_0+a_1x+a_2x^2+\cdots+a_nx^n,\]
where each $a_i$ is an element of $A$, and $x$ is a variable. We will discuss the formal definition later, but the key fact to know is that polynomials can be represented in the above form. The terms $a_i$ are the \col{coefficients}\index{coefficients!of a polynomial} of the polynomial, and $a_0$ is the \col{constant term}\index{constant term!of a polynomial}. If $a_n$ is nonzero, the integer $n$ is called the \col{degree}\index{degree!of a one-variable polynomial} of the polynomial and $a_n$ the \col{leading coefficient}\index{leading coefficient}. A \col{monic polynomial} is a polynomial whose leading coefficient is 1. A polynomial in which no positive exponent of $x$ appears is a \col{constant polynomial}\index{constant polynomial}. The constant polynomial with $a_0=0$ is the \col{zero polynomial}\index{zero polynomial}.

Polynomials are added and multiplied by distributing and combining like terms as in elementary algebra, e.g., \[(2-x+x^2)+(1+5x^2-4x^3)=3-x+6x^2-4x^3\] and
\[(2-x^2)(1+3x)=2+6x-x^2-3x^3.\]

The set of all polynomials over $A$ will be denoted by $A[x]$. This set, with the operations of addition and multiplication described above, forms a ring; a sketch of the proof is given at the end of the chapter. The choice of $x$ as the variable in our polynomials is entirely arbitrary, and any other symbol may be used to denote the variable. Elements of $A[x]$ can be denoted in the form $p(x)$ when the variable $x$ needs to be emphasized, but otherwise may be denoted without the variable, as $p$ rather than $p(x)$.

Every polynomial $p\in A[x]$ gives rise to a function $A\to A$ by \col{evaluation}\index{evaluation!of one-variable polynomial}: given any element $r\in A$, we define $p(r)$ to be the result of substituting $r$ for $x$ in the expression defining $p$ and evaluating the resulting expression using the ring operations in $A$. A \col{root}\index{root of a polynomial}, or \col{zero}\index{zero!of one-variable polynomial}, of $p$ is an element $r\in A$ such that $p(r)=0$.

\begin{ex}[Polynomial evaluation]\;
\begin{examples}
\item For any polynomial $p$, $p(0)$ is equal to the constant term of $p$.
\item Every element of $A$ is a root of the zero polynomial.
\item The complex number $i$ is a root of the polynomial $x^2+1\in\C[x]$.
\end{examples}
\end{ex}

Although every polynomial gives rise to a function, it is important to distinguish between polynomials and their associated functions, as different polynomials can correspond to the same function. For example, consider the polynomials $p,q\in\F_3[x]$ given by \[p=x^2+1,\quad q=x^3+x^2+2x+1.\] 

We compute $p(0)=q(0)=1$, $p(1)=q(1)=2$, and $p(2)=q(2)=2$, so the functions $\F_3\to\F_3$ corresponding to $p$ and $q$ are the same even though $p$ and $q$ are distinct polynomials. Thus, we will distinguish between the concepts of polynomial over $A$ and \col{polynomial function}\index{polynomial function!in one variable} on $A$, the latter term meaning any mapping $f:A\to A$ such that there exists a polynomial $p\in A[x]$ with $f(r)=p(r)$ for every $r\in A$. The above example shows that the function $\F_3\to\F_3$ defined by $0\mapsto 1$, $1\mapsto 2$, and $2\mapsto 0$ is a polynomial function.

Next, we prove a fundamental property of polynomials, namely the fact that one polynomial can always be divided by another leaving a small remainder (though the dividend may need to be scaled first---see details below). In order to state this result more precisely, it will be helpful to extend our notion of the degree of a polynomial to include the zero polynomial. Recall that the \href{https://en.wikipedia.org/wiki/Extended_real_number_line}{extended real number system}, often denoted by $\overline\R$, consists of $\R$ together with two additional elements, $\infty$ and $-\infty$, having the natural properties that $-\infty<r<\infty$ for every real number $r$ and also $r+\infty=\infty$ and $r-\infty=-\infty$ for every $r\in\R$; moreover, $\infty+\infty=\infty$ and $-\infty-\infty=-\infty$.

We define the degree of the polynomial 0 to be $-\infty$, thus obtaining a degree function $\deg: A[x]\to\overline\R$ whose image is the set $\N\cup\{-\infty\}$. Defining $\deg(0)$ to be $-\infty$ will prove to be especially useful when $A$ is an integral domain; see Theorem \ref{degthm}.

\begin{framed}
\begin{thm}[Division with remainder]\label{divrem}
Let $A$ be a commutative ring. Let $f,g\in A[x]$ be polynomials with $g$ nonzero, and let $b$ be the leading coefficient of $g$. There exist polynomials $q,r\in A[x]$ and an integer $m\ge 0$ such that \[b^mf=qg+r\quad\text{and}\quad\deg(r)<\deg(g).\] Moreover, we may assume that $m\le\max\{0,\deg(f)-\deg(g)+1\}$.
\end{thm}
\end{framed}
\begin{proof} We proceed by induction on $\deg(f)$. If $\deg(f)<\deg(g)$ we may take $q = 0$, $r = f$, and $m = 0$; this establishes a base case. For the inductive step, suppose that $\deg(f)\ge\deg(g)$. Let $f = ax^n + \cdots$ and $g = bx^d + \cdots$, where $a$ and $b$ are the leading coefficients of $f$ and $g$, respectively. Since $bf(x)$ and $ax^{n-d}g(x)$ both have leading term $abx^n$, the polynomial $f_1 = bf - ax^{n-d}g$ has degree less than $n$. By induction there exist $q_1, r_1 \in A[x]$ such that
  \[b^{m_1} f_1 = q_1 g + r_1 \quad\text{and}\quad \deg(r_1) < \deg(g),\] 
  where $m_1\le\max\{0,\deg f_1 - \deg g + 1\}$. It follows that
  \[b^{m_1 + 1} f = (q_1 + ab^{m_1}x^{n-d}) g + r_1.\] 
Setting $q=q_1 + ab^{m_1}x^{n-d}$, $r=r_1$, and $m=m_1+1$ gives the desired relation, $b^m f = qg + r$. Moreover, since $\deg f_1<\deg f$, we have $
m_1 \le \max\{0, \deg f - \deg g\}$ and therefore \[m = m_1 + 1 \le \max\{0, \deg f - \deg g + 1\}.\qedhere\]
\end{proof}

For practical applications it can be useful to explicitly compute polynomials $q$ and $r$ as in the above theorem. This can be done by setting $m = \max\{0, \deg f - \deg g + 1\}$ and applying \href{https://en.wikipedia.org/wiki/Polynomial_long_division}{polynomial long division} to divide $b^mf$ by $g$.

\begin{1ex} Suppose that $f,g\in\Z[x]$ are the polynomials
\[f(x) = 3x^3 + 2x^2 + x - 1,\quad g(x) = 2x^2 + x + 1.\] Scaling $f(x)$ by $2^2$ and applying long division, we find that
\[4f(x)=(6x+1)g(x)-3x-5.\]

In the notation of the theorem, we have $q=6x+1$ and $r=-3x-5$.
\end{1ex}

We now derive some important results from Theorem \ref{divrem}, the first one being the relation between roots of a polynomial and factors of degree one.

\begin{framed}
\begin{thm}[Factor Theorem]
Let $A$ be a commutative ring and $p(x)\in A[x]$. An element $a\in A$ is a root of $p(x)$ if and only if $x-a$ is a factor of $p(x)$, i.e., there exists $t(x)\in A[x]$ such that \[p(x)=(x-a)\cdot t(x).\]
\end{thm}
\end{framed}
\begin{proof}
We divide $p(x)$ by $x-a$ leaving a constant remainder: \[p(x)=t(x)(x-a)+r.\]

Evaluating at $a$, we have $p(a)=r$. Thus, $a$ is a root of $p(x)$ if and only if $r=0$, which is equivalent to $x-a$ being a factor of $p(x)$.
\end{proof}

Theorem \ref{divrem} takes an especially strong and simple form in the case where $A$ is a field. The following version of the theorem will be the most commonly used throughout the course. Note that, in contrast to Theorem \ref{divrem}, no scaling of the dividend is required here.

\begin{framed}
\begin{thm}[Division Algorithm]\label{div_algo}
Let $k$ be a field and let $f,g\in k[x]$ be polynomials with $g\ne 0$. There exists a unique pair of polynomials $q,r\in k[x]$ such that $f=qg+r$ and $\deg(r)<\deg(g)$.
\end{thm}
\end{framed}
\begin{proof}
Existence follows from Theorem \ref{divrem}. Uniqueness can be shown using Theorem \ref{degthm} below, and is left to the reader; see Exercise \ref{qr_unique}.
\end{proof}

The polynomials $q$ and $r$ in Theorem \ref{div_algo}\index{quotient (in division algorithm)} are called, respectively, the \col{quotient} and \col{remainder}\index{remainder (in division algorithm)} of $f$ divided by $g$. As before, $q$ and $r$ can be computed using long division. The division algorithm for polynomials is a close analogue of the division algorithm for integers, which lies at the foundation of modular arithmetic: if $a,b\in\Z$ with $b\ne 0$, there exists a unique pair of integers $q,r$ such that \[a=qb+r\;\;\text{and}\;\; 0\le r<|b|.\]

 (The integer $r$ is the reduction of $a$ modulo $b$.) 

For the most part, in this course we will work with polynomial rings $A[x]$ whose coefficient ring $A$ is a domain. This is done in order to avoid certain unexpected behaviors illustrated by the following examples:

\begin{itemize}
\item In general, we may expect that multiplying a polynomial of degree two and one of degree three, as in the product $(x^2-1)(4+x^3)$, should result in a polynomial of degree five. However, this may not always be the case. Consider, for example, the polynomials $2x^2+1$ and $4+3x^3$ in the ring $\Z_6[x]$. The product of these polynomials is equal to $2x^2+4+3x^3$, a polynomial of degree three, not five. What occurred here is that the leading coefficients of the polynomials multiplied to 0, so the degree of the product dropped unexpectedly. Theorem \ref{degthm} below shows that we can avoid this type of example by restricting attention to coefficient rings that are integral domains.
\item The familiar notion that a polynomial equation of degree $n$ has at most $n$ solutions may also fail to hold if we allow coefficients from a ring that is not an integral domain. For example, the polynomial $x^2+7\in\Z_8[x]$ has degree two and yet has four roots, namely $1,3,5$, and $7$. Theorem \ref{roots_number} below shows that this type of example can only occur for coefficient rings, such as $\Z_8$, that are not integral domains.
\end{itemize}

\begin{framed}
\begin{thm}\label{degthm}
If $A$ is an integral domain, then \[\deg(pq)=\deg(p)+\deg(q)\] for all $p,q\in A[x]$. 
\end{thm}
\end{framed}
\begin{proof}
If both $p$ and $q$ are the zero polynomial, the result follows from the fact that $-\infty=-\infty-\infty$. If $p=0$ and $q$ is nonzero (or vice versa), the theorem holds because $-\infty=-\infty+n$ for all $n\in\N$. If both $p$ and $q$ are nonzero, we may write $p(x)=ax^n+s(x)$ and $q(x)=bx^m+t(x)$, where $n=\deg(p)$, $m=\deg(q)$, $\deg(s)<n$, and $\deg(t)<m$. Then $p(x)q(x)=abx^{n+m}+h(x)$, where $\deg(h)<n+m$. Since $A$ is a domain and $a,b$ are nonzero, then $ab\ne 0$; hence $\deg(pq)=n+m$.
\end{proof}

\begin{framed}
\begin{thm}\label{roots_number}
The following hold if $A$ is an integral domain:
\begin{enumerate}
\item[(a)] The polynomial ring $A[x]$ is a domain and its unit group is equal to the unit group of $A$.
\item[(b)] If $p\in A[x]$ is a nonzero polynomial of degree $n$, then $p$ has at most $n$ distinct roots in $A$.
\end{enumerate}
\end{thm}
\end{framed}
\begin{proof}
Assuming $A[x]$ is not a domain, there exist nonzero polynomials $p,q\in A[x]$ such that $pq=0$. Then $\deg(p)+\deg(q)=\deg(0)=-\infty$, which is an immediate contradiction since the degrees of $p$ and $q$ are natural numbers. Hence $A[x]$ must be a domain.

Note that if $p\in A[x]$ is a unit and $q=p^{-1}$, then $p$ and $q$ are nonzero polynomials with $pq=1$, so $\deg(p)+\deg(q)=\deg(1)=0$. Since $\deg(p)$ and $\deg(q)$ must be nonnegative, this implies that $p$ and $q$ are both constant polynomials and thus elements of $A$. Moreover, since $pq=1$, $p$ is a unit in $A$. Conversely, every unit in $A$ is a unit in $A[x]$. This proves (a). Part (b) is proved by induction on $n$; see Exercise \ref{roots_deg}.
\end{proof}

We end this chapter by describing a formal construction of polynomial rings in order to establish a firm logical foundation for working with polynomials.

Let $A$ be an arbitrary commutative ring. Recall that $A^{\N}$ denotes the set of all functions $\N\to A$ (equivalently, the set of all sequences of elements of $A$). Define $P$ to be the subset of $A^{\N}$ consisting of sequences that are eventually zero, i.e., sequences $(a_0,a_1,a_2,\ldots)$ such that $a_i=0$ for all large enough indices $i$. We define addition in $P$ as pointwise addition of functions; representing sequences as infinite lists, this corresponds to componentwise addition. Multiplication in $P$ is defined as follows: for $f,g\in P$ we let $fg\in P$ be the function \[fg(n)=\sum_{i+j=n}f(i)g(j)\quad\text{for all }n\in\N,\] where the sum takes place over all ordered pairs $(i,j)\in\N^2$ such that $i+j=n$. Thus, if $f=(a_0,a_1,a_2,\ldots)$ and $g=(b_0,b_1,b_2,\ldots)$, then \[fg=(a_0b_0,\,a_0b_1+a_1b_0,\,a_0b_2+a_1b_1+a_2b_0,\ldots).\] 

This product is known as the \col{Cauchy product}\index{Cauchy product!of sequences} of sequences. All the ring axioms are satisfied, so $P$ is a ring.

It remains to show that elements of $P$ can be represented as polynomials. Elements of $A$ will henceforth be interpreted as elements of $P$ via the map $a\mapsto (a,0,0,\ldots)$ for all $a\in A$. Defining $x\in P$ to be the sequence $(0,1,0,0,\ldots)$, an induction argument shows that for all $n\in\N$, $x^n$ is the sequence with a 1 in the $(n+1)$-st entry and zeros elsewhere; from here it follows that a sequence $(a_0,a_1,\ldots, a_n,0,0,0,\ldots)$ in $P$ is equal to the element $a_0+a_1x+a_2x^2+\cdots+a_nx^n$, thus establishing the polynomial form of elements of $P$.  Moreover, if elements of $P$ are represented in this form, then addition and multiplication correspond to the process of distributing and combining like terms. In particular, multiplication satisfies the following: 
if $f=\sum_{i=0}^N{a_ix^i}$ and $g=\sum_{j=0}^M{b_jx^j}$, then \[fg=\sum_{k=0}^{M+N}c_kx^k,\quad\text{ where }\quad c_k=\sum_{i+j=k}a_ib_j.\]

The standard notation for the ring $P$ we have constructed is $A[x]$; this ring is called the \col{polynomial ring over $A$ in the variable $x$}. From now on we will write elements of $A[x]$ as polynomials rather than as sequences; in particular, elements of $A$ are interpreted as constant polynomials in $A[x]$.

\section*{Exercises}

\begin{note}\;
\begin{itemize}
\item Solutions to starred exercises (\sel) are included in Chapter~\ref{review1chap}.
\item A house (\house) marks exercises suggested as sample homework.
\end{itemize}
\end{note}

\medskip

\begin{exercises}
\item Let $f:\F_5\to\F_5$ be the function given by \[0\mapsto 1,\; 1\mapsto 2,\; 2\mapsto 4,\; 3\mapsto 1,\; 4\mapsto 0.\] Show that $f$ is a polynomial function; that is, find $p(x)\in \F_5[x]$ such that $f(a) = p(a)$ for all $a\in \F_5$.
\item Carry out division with remainder, in the sense of Theorem \ref{divrem}, applied to the polynomials $f=8x^5+2,\;g=3x^3+7\in\Z_{10}[x]$. Are the quotient and remainder unique in this case?
\item Let $A$ be an integral domain and $f\in A[x]$ a nonzero polynomial. If $r_1,\ldots, r_n$ are distinct roots of $f$ in $A$, show that there is a polynomial $g\in A[x]$ such that \[f(x)=(x-r_1)(x-r_2)\cdots(x-r_n)g(x).\]
\item(\sel)\label{chap4_sel} Let $p(x)=x^4-2x^3-3x^2+4x+2$. Show that $\pm\sqrt 2$ are roots of $p(x)$ in $\Z[\sqrt 2]$ and that $p(x)$ factors as $(x-r_1)(x-r_2)(x-r_3)(x-r_4)$ with each $r_i\in\Z[\sqrt 2]$. Find all the roots $r_i$ explicitly.
\item Let $A$ be a commutative ring. Show that \[\deg(f+g)\le\max\{\deg f, \deg g\}\] for all $f,g\in A[x]$, and that equality holds if $\deg f\ne\deg g$.
\item(\house)\label{qr_unique} Use Theorem \ref{degthm} to prove that the quotient and remainder in the Division Algorithm (Theorem \ref{div_algo}) are uniquely determined.
\item\label{roots_deg} Use induction to prove Theorem \ref{roots_number}(b).
\item(\house) Show by example that a polynomial over a commutative ring that is not a domain can have infinitely many roots.
\item\label{eisen_lem} This exercise is used in the proof of Eisenstein's Criterion (Theorem \ref{eisenstein}), an important result later in the course.

Let $A$ be an integral domain, and suppose that $s,t\in A[x]$ are polynomials such that \[ax^n=s(x)t(x),\] where $n\ge 2$ and $a\in A$ is nonzero. Show that both $s$ and $t$ have constant term 0.
\end{exercises}

\chapter{Power Series}\label{power_chap}

In this chapter and the next, we study two types of generalizations of polynomials in one variable: power series and polynomials in several variables. We begin with power series, which can be informally described as infinite polynomials, or sums of the form \[a_0+a_1x+a_2x^2+a_3x^3+\cdots\] that may continue indefinitely. One well-known power series is the Taylor series for the exponential function around $x=0$:
\[e^x=1+x+\frac{x^2}{2!}+\frac{x^3}{3!}+\frac{x^4}{4!}+\cdots.\]

In Calculus, identities involving power series (such as the above expansion of $e^x$) are interpreted in terms of convergence; the focus is on whether the infinite sum defines a meaningful function. In this course, we take a more general algebraic perspective. Power series will be defined and studied as algebraic objects, with particular emphasis on the fact that they form a ring, without requiring any notion of convergence. In specific cases, such as power series with complex coefficients, questions of convergence can be explored. We will show, in particular, that convergent power series also form a ring.

The formal construction of power series rings is quite straightforward now that we've seen the construction of polynomial rings. In the latter construction we began with the set $A^{\N}$ (the set of sequences in $A$), where $A$ is any commutative ring, and we defined addition and multiplication on a certain subset of $A^{\N}$, namely the set of sequences whose terms are eventually zero. To define power series rings we simply make the observation that these same operations extend to the entire set $A^{\N}$. Thus, the sum and product of sequences $a=(a_0,a_1,a_2,\ldots)$ and $b=(b_0,b_1,b_2,\ldots)$ are defined by
\begin{align*}
a+b&=(a_0+b_0,a_1+b_1,a_2+b_2,\ldots),\\
ab&=(a_0b_0,\,a_0b_1+a_1b_0,\,a_0b_2+a_1b_1+a_2b_0,\ldots).
\end{align*}

The set $A^\N$, equipped with these operations, satisfies all the ring axioms, and is called the \col{ring of power series over $A$}. Introducing the notation $\sum_{i=0}^\infty a_ix^i$ to represent the sequence $(a_0,a_1,\ldots)$, multiplication in this ring takes the form
\[\left(\sum_{i=0}^{\infty}a_ix^i\right)\cdot\left(\sum_{j=0}^{\infty}b_jx^j\right)=\sum_{k=0}^{\infty}c_kx^k,\quad\text{ where } c_k=\sum_{i+j=k}a_ib_j.\]

The ring of power series over $A$ is denoted by $A\pow x$. (Of course, we may use any other symbol for the variable $x$.) Note that this notation differs from $A[x]$, which represents a polynomial ring. Multiplication in $A\pow x$ is called the \col{Cauchy product}\index{Cauchy product!of power series} of power series. The terms $a_i$ in the power series $\sum a_ix^i$ are the \col{coefficients}\index{coefficients!of a power series} of the power series, and $a_0$ is the \col{constant term}\index{constant term!of a power series}. A power series  whose nonconstant terms are all 0 is a \col{constant power series}.

Polynomials can now be regarded as special cases of power series, namely as those power series whose coefficients are all zero after some point in the series. Hence, the ring $A\pow x$ contains $A[x]$.

The extra flexibility gained by extending polynomials to infinite series is reflected in special properties of the ring $A\pow x$ that are not shared with $A[x]$. In particular, we have seen that if $A$ is a domain, then the units of $A[x]$ must all be constant polynomials (Theorem \ref{roots_number}). In contrast, a polynomial of positive degree may happen to be a unit in $A\pow x$. One standard example is the fact that the multiplicative inverse of $1-x$ in $A\pow x$ is the geometric series
\[g=1+x+x^2+x^3+\cdots.\]
Indeed, the definition of the Cauchy product implies that $g\cdot(1-x)=1$.

Our first theorem concerning power series analyzes the units and zero divisors of the ring $A\pow x$. For this purpose, it is useful to introduce the \col{order of vanishing}\index{order of vanishing} of a power series, a certain element of $\overline\R$ associated to the power series. 

For a nonzero $f\in A\pow x$, we define $\ord(f)$ to be the least index of a nonzero coefficient of $f$. Thus, we have $\ord(f)=n$ if \[f=a_nx^n+a_{n+1}x^{n+1}+a_{n+2}x^{n+2}+\cdots\] with $a_n\ne 0$. Note that $f$ is then a multiple of $x^n$, and that $x^n$ is in fact the highest power of $x$ dividing $f$. For the zero power series we define $\ord(0)=\infty$, obtaining a map $\ord:A\pow x\to\overline\R$ whose image is the set $\N\cup\{\infty\}$. Using the notion of an \href{https://en.wikipedia.org/wiki/Infimum_and_supremum}{infimum}, a more succinct definition of the order of vanishing is  \[\ord f=\inf\{i:a_i\ne 0\},\] which applies to every power series $f=\sum a_ix^i$.

\begin{framed}
\begin{thm}\label{ordthm}
The following hold if $A$ is an integral domain:
\begin{enumerate}
\item[(a)] For all $p,q\in A\pow x$, $\ord(pq)=\ord(p)+\ord(q)$.
\item[(b)] The ring $A\pow x$ is a domain and its units are the power series whose constant term is a unit in $A$.
\end{enumerate}
\end{thm}
\end{framed}

\begin{proof}
Part (a) follows readily from the definitions, and the fact that $A\pow x$ is a domain follows from (a). Now suppose that $p=\sum a_ix^i$ is a unit in $A\pow x$. If $q=\sum b_jx^j$ is the multiplicative inverse of $p$, then $a_0b_0=1$, and in particular, $a_0$ is a unit in $A$. Conversely, suppose that $p=\sum a_ix^i$ is a power series such that $a_0\in A^{\times}$, and let $b_0=a_0^{-1}$. We may then recursively construct a power series $q=\sum b_j x^j$ that is a multiplicative inverse of $p$: assuming that $b_0,\ldots, b_{n-1}$ have been defined, we set
\begin{equation}\label{seriesinverse}
b_{n}=-a_0^{-1}\cdot\sum_{i=1}^{n}a_ib_{n-i}.\tag{$*$}
\end{equation}
The definition of $q$ ensures that $pq=1$, so that $p$ is a unit in $A\pow x$.
\end{proof}

The proof of Theorem \ref{ordthm} effectively provides an algorithm for inverting a power series. As an example of the method, take $f\in\R\pow x$ to be the Taylor series for $e^x$. This series is a unit in $\R\pow x$ since its constant term is 1, a unit in $\R$. Using the formula \eqref{seriesinverse}, we compute the first several terms of $f^{-1}$ and find that
\[f^{-1}=1 - x + \frac{x^2}{2} - \frac{x^3}{3!} + \frac{x^4}{4!} -\frac{x^5}{5!} + \cdots.\] 

We recognize this power series as the Taylor series for $e^{-x}$, which is consistent with the identities $e^xe^{-x}=e^0=1$.

To close this chapter, we restrict attention to power series with coefficients in $\C$, so that we may consider questions of convergence. Recall that every power series $f(z)\in\C\pow z$ has a \href{https://en.wikipedia.org/wiki/Radius_of_convergence}{radius of convergence} that is either a non-negative real number or $\infty$. We say that $f$ is \col{convergent}\index{convergent power series} if its radius of convergence is positive, so that there exists an open disk $D\subseteq\C$ centered at 0 such that $f(z)$ converges absolutely for all $z\in D$.

If $f$ and $g$ are convergent power series with radii of convergence $R_1$ and $R_2$, then $f$ and $g$ are both absolutely convergent on the open disk of radius $\min\{R_1,R_2\}$ centered at 0. It follows that $f+g$ and $fg$ also converge absolutely within this disk, the latter by a well-known \href{https://en.wikipedia.org/wiki/Cauchy_product#Convergence_and_Mertens'_theorem}{theorem of Mertens} concerning convergence of the Cauchy product of complex series. As a consequence, convergent power series form a ring, which we denote by $\C\conv z$. Note that the polynomial ring $\C[z]$ is contained in $\C\conv z$ since polynomials have infinite radius of convergence; hence, \[\C[z]\subset\C\conv z\subset\C\pow z.\]

Though the proof is an exercise in analysis rather than algebra, it is worth mentioning here that the multiplicative inverse of a convergent power series is necessarily convergent, so that the units of $\C\conv z$ are precisely the convergent power series with nonzero constant term. The inverse of any unit can be computed as above, using \eqref{seriesinverse}.

\newpage
\section*{Exercises}

\begin{note}\;
\begin{itemize}
\item Solutions to starred exercises (\sel) are included in Chapter~\ref{review1chap}.
\item A house (\house) marks exercises suggested as sample homework.
\end{itemize}
\end{note}

\medskip

\begin{exercises}
\item Let $f,g\in\Z_9\pow t$ be power series with initial terms \[f=2+7t+3t^2+5t^4+\cdots,\quad g=6+t^2+8t^4+\cdots.\] Compute the initial terms of the product series $fg$. Is $f$ a unit in $\Z_9\pow t$? Is $g$ a unit?
\item(\sel)\label{chap5_sel1} Let $f$ be the power series over $\F_7$ given by
\[f=\sum_{i=1}^{\infty}(i+1)i\,x^{i-1}=2+6x+5x^2+6x^3+\cdots.\]
Compute the first few terms of $f^{-1}$.
\item(\house) Find all invertible power series $f\in\Z\pow x$ such that $f=f^{-1}$.
\item(\sel)\label{chap5_sel2} Let $A$ be an integral domain. Show that
\[\ord(p+q)\ge\min\{\ord p,\ord q\}\] for all $p,q\in A\pow x$, and that equality holds if $\ord p\ne\ord q$. Find an example where equality fails when $\ord p=\ord q$.
\item(\house) Let $A$ be a nontrivial commutative ring that is not a domain. Show that $A\pow x$ contains infinitely many nonconstant zero divisors.
\item\label{pow_unit} Let $k$ be a field. Show that every nonzero element of $k\pow x$ can be represented uniquely in the form $x^nu(x)$, where $n\in\N$ and $u(x)\in k\pow x^{\times}$. 
\item In $\C\conv z$, let $f$ be the Taylor series of $\cos z$ centered at 0:
\[f(z)=1 - \frac{z^2}{2!} + \frac{z^4}{4!} - \frac{z^6}{6!} +\cdots.\] 

Compute the initial terms of $f^{-1}$ and use them to approximate $\sec(i/4)$.
\end{exercises}

\chapter{Polynomials in Several Variables}\label{multi_chap}

Let $A$ be a commutative ring and let $I$ be an arbitrary nonempty set. We will define a polynomial ring $A[X]$ having one variable for every element of $I$. In the case where $I$ is a singleton set, this construction yields the one-variable polynomial ring $A[x]$ defined in Chapter \ref{poly_chap}. For polynomial rings in two and three variables, we will use notation such as $A[x,y]$ and $A[x_1,x_2,x_3]$, respectively. Similarly, $A[x_1,x_2,x_3,\ldots]$ denotes a polynomial ring with a countably infinite number of variables.

Polynomials all have the same general form: a finite sum of terms, each of which is the product of an element of the base ring $A$ and finitely many non-negative powers of variables. For example, the following is an element of the ring $\Z[x,y,z]$: 
\[2x^3y+3128y^3z^5-219xy^2z^3+24.\]

The formal construction of polynomial rings is given at the end of the chapter; for now, let us focus on describing these rings, especially in the case of finitely many variables. If $A$ is a commutative ring and $n$ a positive integer, a \col{monomial}\index{monomial!in a polynomial ring} in $A[x_1,\ldots, x_n]$ is a product of the form \[x_1^{e_1}\cdots x_n^{e_n},\] where each $e_i\in\N$. A \col{polynomial}\index{polynomial!in several variables} is a finite sum of the form $\sum a_iM_i$, where each $M_i$ is a monomial and $a_i\in A$; in other words, a polynomial is a linear combination of monomials.

Every nonzero $f\in A[x_1,\ldots, x_n]$ has a unique representation as a linear combination of \underline{distinct} monomials $M_i$ with \underline{nonzero} coefficients $a_i$. The monomials in this representation will be called the \col{monomials of $f$}\index{monomial!of a polynomial}. For example, the polynomial $f=3xy+5x^2y-3xy+1$ can be simplified to $f=5x^2y+1$, so the monomials of $f$ are $x^2y$ and 1.

As in the one-variable case, every polynomial $f\in A[x_1,\ldots, x_n]$ gives rise to a function $A^n\to A$ by \col{evaluation}\index{evaluation!of multivariable polynomial}, i.e., by substituting the variables $x_1,\ldots, x_n$ with the entries in any given tuple $(a_1,\ldots, a_n)\in A^n$. For example, if $f(x,y,z)=2xy-y+7z^2x$, then \[f(1,2,0)=4-2+0=2.\] 

A \col{polynomial function in $n$ variables}\index{polynomial function!in several variables} over $A$ is a function $g:A^n\to A$ such that there exists $f\in A[x_1,\ldots, x_n]$ with $g(r)=f(r)$ for all $r\in A^n$. As in the one-variable case, it is important to distinguish between polynomials and their associated functions, as distinct polynomials can give rise to the same function. For instance, the function $\F_2^2\to\F_2$ given by
\[(0,0)\mapsto 1,\; (0,1)\mapsto 1,\; (1,0)\mapsto 0,\; (1,1)\mapsto 1\]
is a polynomial function, as it agrees with the function induced by the polynomial $f(x,y)=1+x+xy$. The same function corresponds to the polynomial $g(x,y)=1+x+xy+x^2+x$.

A \col{zero}\index{zero!of multivariable polynomial} of a polynomial $f\in A[x_1,\ldots, x_n]$ is any element $\alpha\in A^n$ such that $f(\alpha)=0$. (We reserve the term \emph{root} for a zero of a polynomial in one variable.) Unlike the case of one variable, a polynomial in several variables over an integral domain can have infinitely many zeros: consider, for instance, the polynomial $f(x,y)=x-y\in\Z[x,y]$.

For multivariate polynomials there are two natural types of degree. In all cases the degree of the zero polynomial is $-\infty$, so we focus on nonzero polynomials in the following definitions.

The first type of degree is the degree with respect to one of the variables. To define this, note that a polynomial $f(x_1,\ldots, x_n)$ can be expressed as a polynomial in any one variable, say $x_i$, whose coefficients are polynomials in the remaining $n-1$ variables. For example, if \[f(x,y)=y+2xy^4+3x^2y^3-5x^3-7,\] we may express $f$ as a polynomial in $x$ as \[f=-5x^3+(3y^3)x^2+(2y^4)x+(y-7)\] 
or as a polynomial in $y$ as 
\[f=(2x)y^4+(3x^2)y^3+y-(5x^3+7).\]

The \col{degree of $f$ with respect to $x_i$}\index{degree!with respect to a variable}, denoted by $\deg_{x_i}(f)$, is defined to be the degree of $f$ when expressed as a polynomial in $x_i$. In the above example, we have $\deg_x(f)=3$ and $\deg_y(f)=4$. In general, $\deg_{x_i}(f)$ is the largest exponent of $x_i$ that appears in the monomials of $f$.

The second type of degree, called the \col{total degree}\index{total degree}, or simply the \col{degree}\index{degree!of multivariable polynomial}, is defined as follows: we begin by setting the degree of a monomial $x_1^{e_1}\cdots x_n^{e_n}$ to be the sum $e_1+e_2+\cdots+e_n$. The degree of $f$, denoted by $\deg(f)$, is then defined to be the maximum of the degrees of the monomials of $f$. In the example above, $\deg(f)=5$. With this definition of degree, we have the following generalization of a theorem proved for one-variable polynomials in Chapter \ref{poly_chap}:

\begin{framed}
\begin{thm}\label{degthm_multi} If $A$ is an integral domain, then \[\deg(pq)=\deg(p)+\deg(q)\] for all $p,q\in A[x_1,\ldots, x_n]$. As a result, the ring $A[x_1,\ldots, x_n]$ is a domain and its unit group is $A^{\times}$.
\end{thm}
\end{framed}
\begin{proof}
Let $d=\deg(p)$ and $e=\deg(q)$. Let $p_d$ be the sum of all the monomials of $p$ having degree $d$, and similarly define $q_e$. Thus, we have $p=p_d+f$ and $q=q_e+g$, where $\deg(f)<d$ and $\deg(g)<e$. It follows that $pq=p_dq_e+h$, where $\deg(h)<d+e$. Moreover, every monomial in the product $p_dq_e$ has degree $d+e$, so the definition of degree implies that $\deg(pq)=d+e$. The second statement in the theorem follows from this property of the degree map (argue as in the proof of Theorem \ref{roots_number}(a)).
\end{proof}

The proof of Theorem \ref{degthm_multi} illustrates the usefulness of considering the sum of all the monomials of a polynomial having a fixed degree. If $f\in A[x_1,\ldots, x_n]$ is a nonzero polynomial of degree $d$, for every integer $0\le i\le d$ we will denote by $f_i$ the sum of all monomials of $f$ having degree $i$; if no such monomials exist, then $f_i=0$. Thus, \[f=f_0+f_1+\cdots+f_d,\quad d=\deg(f).\] The terms $f_i$ in this sum are called the \col{homogeneous components}\index{homogeneous components} of $f$. Note that if $f_i$ is nonzero, then $\deg(f_i)=i$. If all the monomials of $f$ have the same degree (equivalently, if $f=f_d$), then $f$ is called \col{homogeneous}\index{homogeneous polynomial}. For example, the polynomial \[f(x,y,z)=13x^2z-2xyz+y^3\] is homogeneous of degree three, and if \[g(x,y)=y+2xy^4+3x^2y^3-5x^3-7,\]  the nonzero homogeneous components of $g$ are \[g_0=-7,\; g_1=y,\;g_3=-5x^3,\;g_5=2xy^4+3x^2y^3.\]

We end this chapter by describing the formal construction of the general polynomial ring $A[X]=A[x_i:i\in I]$, which has one variable for every element of a nonempty index set $I$. As the particular details of the construction will not be needed later in this course, readers not interested in foundational considerations or infinite-variable cases are welcome to skip this final section. The key takeaway is that there is a logically sound way of working with polynomials in any number of variables. For example, we can meaningfully speak of a polynomial ring over $\Z$ having variables indexed by $\R$; one element of this ring is the polynomial \[2x_{\pi}^{20}\,x_{\sqrt 2}^7-31x_{-1}^3\,x_{2/3}^4+6.\]

Let $M$ be the set of all functions $I\to\N$ having finite support, i.e., functions $f$ such that $f(i)=0$ for all but finitely many $i\in I$. We equip the set $M$ with the binary operation of pointwise addition of functions; this operation is commutative and associative, and its identity is the function mapping every element of $I$ to 0.

Now let $A[X]$ be the set of all functions $M\to A$ having finite support. Addition in $A[X]$ is defined as pointwise addition of functions. For $f,g\in A[X]$, the product $fg\in A[X]$ is the map $M\to A$ given by
\[fg(s)=\sum_{k+m=s}f(k)g(m).\] Here, the sum takes place over all ordered pairs $(k,m)\in M^2$ such that $k+m=s$. Note that $f(k)$ and $g(m)$ are elements of the ring $A$, so the product $f(k)g(m)$ takes place in $A$, as does the sum of all such products.

With the above definitions, $A[X]$ is a commutative ring. It remains to show that this ring contains one variable $x_i$ for each $i\in I$, and that elements of $A[X]$ can be represented as polynomials in the variables $x_i$ with coefficients in $A$.

For each $a\in A$, let $\bar a\in A[X]$ denote the map  \[\bar a(s) = 
\begin{cases} 
a & \text{if } s = 0 \;\text{(the zero function)}, \\
0 & \text{otherwise.}
\end{cases}\]

Elements of $A$ will henceforth be interpreted as elements of $A[X]$ via the mapping $a\mapsto\bar a$. This will allow us to use elements of $A$ as coefficients in our polynomials. Next, we define the variables $x_i$.

For each $i\in I$, let $\delta_i\in M$ be the mapping $I\to\N$ given by \[\delta_i(j) = 
\begin{cases} 
1 & \text{if } j = i, \\
0 & \text{otherwise,}
\end{cases}\]

and let $x_i\in A[X]$ be the mapping $M\to A$ given by
 \[x_i(s) = 
\begin{cases} 
1 & \text{if } s = \delta_i, \\
0 & \text{otherwise.}
\end{cases}\]

We can now build polynomials in the variables $x_i$. To each $s\in M$ we associate the monomial
\[X^s=\prod_{i\in I}x_i^{s(i)},\]
which is well defined given that $s$ has finite support. Every element of $A[X]$ now has a unique representation in the form
\[\sum_{s\in M}a_sX^s,\] where $a_s\in A$ for each $s$ and only finitely many coefficients $a_s$ are nonzero. In terms of this polynomial representation, operations in $A[X]$ take a familiar form:
\begin{align*}
\left(\sum_{i\in M}a_iX^i\right)+\left(\sum_{i\in M}b_iX^i\right)&=\sum_{i\in M}(a_i+b_i)X^i,\\
\left(\sum_{i\in M}a_iX^i\right)\cdot\left(\sum_{j\in M}b_jX^j\right)&=\sum_{k\in M}\left(\sum_{i+j=k}a_ib_j\right)X^k.
\end{align*}

This completes the construction of $A[X]$.

\section*{Exercises}

\begin{note}\;
\begin{itemize}
\item Solutions to starred exercises (\sel) are included in Chapter~\ref{review1chap}.
\item A house (\house) marks exercises suggested as sample homework.
\end{itemize}
\end{note}

\medskip

\begin{exercises}
\item Show that if $A$ is an integral domain and $1\le i\le n$, then \[\deg_{x_i}(pq)=\deg_{x_i}(p)+\deg_{x_i}(q)\] for all $p,q\in A[x_1,\ldots, x_n]$.
\item(\house) Let $f(x, y) = 3x^5y + 2x^3y^3 - 7x^2y^2 + 4xy + 6x + 9y^2 + 1$. Identify the homogeneous components of $f$, the degree of each component, the total degree of $f$, and the degrees with respect to both variables, when $f$ is regarded (a) as a polynomial over $\Z_3$ and (b) as a polynomial over $\Z$.
\item(\house)
\begin{enumerate}
\item Find an example of a commutative ring $A$ and a nonzero polynomial $f\in A[x,y]$ such that every element of $A^2$ is a zero of $f$.
\item Let $k$ be an infinite field and $f\in k[x_1,\ldots, x_n]$ a nonzero polynomial. Show that $f(\alpha)\ne 0$ for some $\alpha\in k^n$.
\end{enumerate}
\item\label{finite_polynomial_fn}\begin{enumerate}
\item Let $F$ be a finite field and $n$ a positive integer. Show that every function $F^n\to F$ is a polynomial function by comparing the total number of functions with the number of polynomial functions.
\item Give an example of a function $\Q^2\to\Q$ that is not a polynomial function.
\end{enumerate}
\item(\sel)\label{chap6_sel} Let $A$ be a commutative ring and $f\in A[x_1,\ldots, x_n]$ a nonzero polynomial of degree $d$. Show that if $f$ is homogeneous, then \[f(\lambda x_1,\ldots,\lambda x_n)=\lambda^d f(x_1,\ldots, x_n)\] for every $\lambda\in A$. Conversely, show that if $A$ is an infinite integral domain and the above identity holds for every $\lambda\in A$, then $f$ is homogeneous.
\item Let $A$ be a commutative ring.
\begin{enumerate}
\item Show that a product of homogeneous polynomials in $A[x_1,\ldots, x_n]$ is homogeneous.
\item Assume that $A$ is a domain. Show that every factor of a homogeneous polynomial is homogeneous. In other words, if $F,G,H\in A[x_1,\ldots, x_n]$ are such that $F$ is homogeneous and $F=GH$, then both $G$ and $H$ are homogeneous.
\end{enumerate}
\end{exercises}

\chapter{The Field of Fractions of a Domain}\label{frac_chap}

Just as the integers can be extended to the field of rational numbers by adjoining fractions, any integral domain can be extended to a field by adjoining fractions of its elements. This brief chapter presents the construction in detail, generalizing the familiar construction of $\Q$ from $\Z$.

Let $D$ be a domain. We will show that there is a field $F$ whose elements can be represented in the form $a/b$ with $a,b\in D$ and $b\ne 0$, and where the usual rules of fraction arithmetic apply:
\begin{equation}\label{fracop}
\frac{a}{b}+\frac{c}{d}=\frac{ad+bc}{bd}\quad\text{and}\quad\frac{a}{b}\cdot\frac{c}{d}=\frac{ac}{bd}.\tag{$*$}
\end{equation}
Moreover, $F$ contains $D$ in the sense that every element $a\in D$ can be interpreted as the fraction $a/1\in F$. The field $F$ is known as the \col{field of fractions}\index{field of fractions of a domain} of $D$, or also the \col{quotient field}\index{quotient field of a domain} of $D$, and is denoted by $\Frac(D)$. Later in the course we will be able to show that $\Frac(D)$ is the smallest field containing $D$, in some precise sense---see Theorem \ref{frac_univ}. In the case $D=\Z$, the field $\Frac(D)$ is precisely $\Q$; in fact, this is the standard definition of $\Q$.

We now prove the existence of a field of fractions for every domain $D$. Let $S=D\times(D\setminus\{0\})$, and define an equivalence relation on the set $S$ by \[(a,b)\sim(c,d)\iff ad=bc.\] Denoting by $a/b$ the equivalence class of the pair $(a,b)$, we have \[a/b=c/d\iff ad=bc,\] as might be expected from prior experience with fractions.

Now let $F$ be the quotient set $S/{\sim}$, i.e., the set of all equivalence classes of elements of $S$. Thus, \[F=\{a/b:a,b\in D\text{ and } b\ne 0\}.\]

Note that $a/b=(xa)/(xb)$ for every nonzero $x\in D$, so the representation of an element of $F$ as a fraction is generally not unique. This reflects the same phenomenon observed in $\Q$, where, for instance, $2/3=4/6$.

Addition and multiplication on $F$ are defined by the identities \eqref{fracop}. These operations are well defined (see Exercise \ref{chap7_sel}) and make $F$ into a ring. The assumption that $D$ is a domain is needed to ensure that the denominator $bd$ in \eqref{fracop} is nonzero.

The zero element of $F$ is $0/1$ and the unity is $1/1$, and these elements are distinct, so $F$ is a nontrivial ring. Moreover, $F$ is commutative since $D$ is commutative, and every nonzero element of $F$ is a unit: any fraction $a/b$ with $a\ne 0$ has a multiplicative inverse in $F$, namely $b/a$. This shows that $F$ is a field.

Finally, to every element $a\in D$ we may associate the fraction $a/1\in F$, and the mapping $a\mapsto a/1$ is injective (see Exercise \ref{frac_map_inj}). Hence, the image of this map can be regarded as a copy of $D$ contained in $F$, and operations in this image set are equivalent to operations in $D$, in the sense that
\[\frac{x}{1}+\frac{y}{1}=\frac{x+y}{1}\quad\text{and}\quad\frac{x}{1}\cdot\frac{y}{1}=\frac{xy}{1}.\]

The above notion of forming a copy of one ring inside another ring is carefully developed in Chapter \ref{homo_chap}. In the terminology introduced there, $D$ is \emph{isomorphic} to a subring of $F$.

This concludes the construction of the field of fractions of $D$. The construction is important from a theoretical standpoint, as it establishes a solid logical foundation for working with fractions of ring elements. However, for practical purposes it suffices to know that elements of $\Frac(D)$ are represented as fractions $a/b$ with $b\ne 0$, that $a/b=c/d\iff ad=bc$, and that operations in $F$ are defined as in \eqref{fracop}.

We end this chapter with two standard examples of fields of fractions.

\begin{1ex}[Rational functions] If $k$ is a field, the ring $k[x_1,\ldots, x_n]$ is an integral domain (by Theorem \ref{degthm_multi}) and therefore has a field of fractions; this field is called the \col{rational function field in $n$ variables over $k$}\index{rational function field}, denoted by $k(x_1,\ldots, x_n)$. Thus,
\[k(x_1,\ldots, x_n)=\Frac k[x_1,\ldots, x_n].\] Elements of this field, called \col{rational functions}\index{rational function}, have the form \[\frac{p(x_1,\ldots, x_n)}{q(x_1,\ldots, x_n)},\] where $p$ and $q$ are polynomials in $n$ variables and $q\ne 0$. The following is a sample computation in the field $\F_5(x,y)$: \[\frac{2y+2}{x+3y}+\frac{4xy+3}{x}=\frac{4x^2y + 2xy^2 + 2xy + 4y}{x^2 + 3xy}.\]
\end{1ex}

\begin{1ex}[Laurent series] If $k$ is a field, the power series ring $k\pow x$ is an integral domain (Theorem \ref{ordthm}) and therefore has a field of fractions; this field is called \col{field of formal Laurent series} over $k$ and denoted by $k\lau x$. Thus,
\[k\lau x=\Frac k\pow x.\]
Note that since $k\pow x$ contains the polynomial ring $k[x]$, the field $k\lau x$ contains $k[x]$ and hence contains its field of fractions: $k(x)\subseteq k\lau x$.

By construction, elements of $k\lau x$, called \col{Laurent series}\index{Laurent series}, can be expressed as fractions whose numerator and denominator are both power series. However, the special properties of the power series ring lead to a simpler representation, namely as a single series in which finitely many negative exponents of the variable $x$ are allowed. Given $f=p/q\in k\lau x$, we may write $q(x)=x^n\cdot r(x)$, where $r(x)$ is a unit in $k\pow x$ and $n=\ord(q)$ is the order of vanishing of $q(x)$. Letting $s(x)=r(x)^{-1}\in k\pow x$, we have $f=x^{-n}ps$, and the product $ps$ is a power series, say
\[ps=a_0+a_1x+a_2x^2+\cdots.\]
We thus arrive at the standard representation of a Laurent series: \[f=\frac{a_0}{x^n}+\frac{a_1}{x^{n-1}}+\cdots+\frac{a_{n-1}}{x}+a_n+a_{n+1}x+a_{n+2}x^2+\cdots.\]

As a simple example, consider expanding the rational function $\frac{1}{x^2(1+3x)}$ as a Laurent series. By the geometric series formula, we have
\[\frac{1}{1+3x}=\sum_{i=0}^{\infty}(-3)^ix^i=1-3x+9x^2-27x^3+\cdots.\] Dividing both sides by $x^2$, we obtain a Laurent series: \[\frac{1}{x^2(1+3x)}=\frac{1}{x^2}-\frac{3}{x}+9-27x+\cdots.\]
\end{1ex}

\section*{Exercises}

\begin{note}\;
\begin{itemize}
\item Solutions to starred exercises (\sel) are included in Chapter~\ref{review1chap}.
\item A house (\house) marks exercises suggested as sample homework.
\end{itemize}
\end{note}

\medskip

\begin{exercises}
\item Which of the rings $\H$, $\R[x]$, $\calC[a,b]$ has a field of fractions?
\item(\house)\label{frac_map_inj} Let $D$ be a domain. Show that the canonical map $D\to\Frac(D)$ given by $a\mapsto a/1$ is injective.
\item(\sel)\label{chap7_sel} Let $D$ be a domain and $F=\Frac(D)$. Show that addition in $F$, defined by \eqref{fracop}, is well defined (independent of the choice of representatives). In other words, show that if $a/b=a'/b'$ and $c/d=c'/d'$, then \[\frac{ad+bc}{bd}=\frac{a'd'+b'c'}{b'd'}.\]
\item(\house) In the field $\F_5\lau t$, let \[f=\frac{1+4t+2t^2+\cdots}{t+3t^2+t^3+\cdots}.\] Find the first three terms of $f$ when expressed as a Laurent series, i.e., in the form $f=\displaystyle\sum_{i=-N}^{\infty}a_it^i$.
\item Let $k(x)$ be the field of rational functions in one variable over a field $k$. Let $v:k(x)\to\overline\R$ be the function defined by \[v(f)=\deg(q)-\deg(p)\]
if $f=p/q$ with $p,q\in k[x]$. Show that this function is well defined and satisfies \[v(fg)=v(f)+v(g)\quad\text{and}\quad v(f+g)\ge\min\{v(f),v(g)\}\] for all $f,g\in k(x)$. Functions with the above properties are called \col{valuations}, and play an important role in certain areas of ring theory. The next exercise includes another example of a valuation.
\item Let $k$ be a field and let $\ord: k\lau t\to\overline\R$ be the function defined by \[\ord(f)=\ord(a)-\ord(b)\] if $f=a/b$ with $a,b\in k\pow t$. Show that this function is well defined and extends the order of vanishing function $\ord: k\pow t\to\overline\R$. Moreover, show that \[\ord(fg)=\ord(f)+\ord(g)\] and \[\ord(f+g)\ge\min\{\ord f,\ord g\}\] for all $f,g\in k\lau t$.
\end{exercises}

\chapter{First Review on Rings}\label{review1chap}

This chapter provides a review of the material covered in Chapters \ref{intro_chap}--\ref{frac_chap}. The first two sections offer a concise summary of key definitions and theorems that the reader is encouraged to revisit before moving on to subsequent chapters. The next section presents detailed solutions to a selection of earlier exercises. The final section contains additional exercises that span a range of topics and levels of depth, and may be used for extra practice.

\section*{Core concepts}

The following is a list of the key concepts defined up to this point:

\begin{itemize}
\item The ring axioms. Commutative rings.
\item Units and zero divisors. The unit group of a ring.
\item Integral domains and fields.
\item Standard examples of rings:
\[\Z,\;\Z_n,\;\Z[\sqrt d],\;\Q(\sqrt d),\;\H,\;M_n(R),\;A[x],\;A\pow x,\;A[x_1,\ldots, x_n].\]
\item Norm functions on $\Q(\sqrt d)$ and on $\H$.
\item Polynomials vs. polynomial functions.
\item The order of vanishing of a power series.
\item Zeros of polynomials in several variables.
\item Degree and homogeneity of multivariable polynomials.
\item The field of fractions of an integral domain.
\item Rational functions and Laurent series over a field.
\end{itemize}

\section*{Essential theorems}

The following are the main theorems covered in Chapters 1--7. These results form the foundation for much of the material ahead.

\begin{itemize}
\item If $p$ is a prime, then $\Z_p$ is a field denoted by $\F_p$.
\item Units in various rings:
\begin{enumerate}
\item[(a)] Let $d\neq 1$ be a squarefree integer. An element of $\Z[\sqrt d]$ is a unit if and only if its norm is $\pm 1$.
\item[(b)] Let $R$ be a commutative ring. A matrix in $M_n(R)$ is invertible if and only if its determinant is a unit in $R$.
\item[(c)] Every nonzero element of $\H$ is a unit.
\item[(d)] If $A$ is an integral domain, the units of $A\pow x$ are the power series whose constant term is a unit in $A$.
\end{enumerate}
\item If $A$ is an integral domain and $p\in A[x]$ a nonzero polynomial of degree $n$, then $p$ has at most $n$ roots in $A$.
\item If $A$ is a domain, then $\deg(fg)=\deg(f)+\deg(g)$ for all $f,g\in A[x_1,\ldots,x_n]$.
\item The Factor Theorem: Let $A$ be a commutative ring and $p\in A[x]$. An element $r\in A$ is a root of $p(x)$ if and only if the polynomial $x-r$ is a factor of $p(x)$.
\item The Division Algorithm for polynomials: Let $k$ be a field and let $f,g\in k[x]$ be polynomials with $g\ne 0$. There exists a unique pair of polynomials $q,r\in k[x]$ such that $f=qg+r$ and $\deg(r)<\deg g$.
\end{itemize}

\section*{Selected solutions}

This section presents worked solutions to selected exercises from earlier chapters.

\begin{framed}
\noindent \textbf{Exercise \ref{chap1_sel}, Chapter \ref{intro_chap}:}\\ Let $(R,+,\cdot)$ be a ring and $x,y \in R$. Use the ring axioms to prove the following:
\begin{enumerate}
\item[(a)] $0 \cdot x = x \cdot 0 = 0$,
\item[(b)] $(-x)\cdot y = x \cdot (-y) = -(x \cdot y)$,
\item[(c)] $(-x)\cdot (-y) = x \cdot y$.
\end{enumerate}
\end{framed}
\begin{proof}[Solution]\;

\noindent \textbf{(a) $0\cdot x = 0$:} In $R$ we have $0 = 0 + 0$, so
\[
0 \cdot x = (0 + 0) \cdot x = 0\cdot x + 0\cdot x.
\]
Subtracting $0 \cdot x$ from both sides, we obtain $0 = 0\cdot x$. A similar argument proves $x\cdot 0=0$.

\medskip

\noindent \textbf{(b) $(-x)\cdot y = -(x \cdot y)$:} By distributivity,
\[
(-x)\cdot y + x\cdot y = (-x+x)\cdot y = 0 \cdot y = 0.
\]
Thus $(-x)\cdot y$ is the additive inverse of $x\cdot y$:
\[
(-x)y =  -(x\cdot y).
\]
The identity $x \cdot (-y) = -(x \cdot y)$ is proved similarly.

\medskip

\noindent \textbf{(c) $(-x) \cdot (-y) = x \cdot y$:} 
By part (b), we have

\[
(-x)\cdot(-y)=-\bigl(x\cdot(-y)\bigr)=-\bigl(-(x\cdot y)\bigr)=x\cdot y.
\qedhere
\]

\end{proof}

\begin{framed}
\noindent \textbf{Exercise \ref{chap2_sel}, Chapter \ref{rings_intro_chap}:}\\ If $a$ is an element of an integral domain such that $a^2=1$, show that $a=1$ or $a=-1$. Find a counterexample to this statement in some ring that is not a domain.
\end{framed}
\begin{proof}[Solution]
Suppose $a^2 = 1$ in an integral domain $D$. Then
\[
a^2 - 1 = 0 \quad \Rightarrow \quad (a - 1)(a + 1) = 0.
\]

Since $D$ is a domain, it has no zero divisors, so this implies
\[
a - 1 = 0 \quad \text{or} \quad a + 1 = 0,
\]
i.e., $a = 1$ or $a = -1$.

Now for a counterexample in a ring that is not a domain. In the ring $\Z_8$ we have $3^2 = 1$, and yet $3 \notin \{1, -1\}$. Thus, the result may fail when the ring in question is not a domain.
\end{proof}

\begin{framed}
\noindent\textbf{Exercise \ref{nilpotent_def}, Chapter \ref{rings_intro_chap}:}\\ An element $x$ of a ring is called \col{nilpotent} if there exists an integer $n\ge 1$ such that $x^n=0$. Show that if $x$ is nilpotent then $1+x$ is a unit.
\end{framed}
\begin{proof}[Solution]
Let $x$ be a nilpotent element of a ring $R$, so that $x^n = 0$ for some $n \ge 1$. Consider the element
\[
y = 1 - x + x^2 - x^3 + \cdots + (-1)^{n-1}x^{n-1}.
\]
We compute the product $(1+x)y$:
\begin{align*}
(1+x)y &= (1+x)\left(1 - x + x^2 - x^3 + \cdots + (-1)^{n-1}x^{n-1}\right) \\
&= 1 - x + x^2 - x^3 + \cdots + (-1)^{n-1}x^{n-1} \\
&\quad + x - x^2 + x^3 - \cdots + (-1)^{n-1}x^n.
\end{align*}
All terms cancel except the $1$ and $(-1)^{n-1}x^n$, which is 0 because $x^n=0$. So we obtain:
\[
(1+x)y = 1.
\]
Similarly, $y(1+x) = 1$, so $y$ is the multiplicative inverse of $1+x$. Hence, $1+x$ is a unit.
\end{proof}

\begin{framed}
\noindent\textbf{Exercise \ref{chap3_sel1}, Chapter \ref{more_rings_chap}:}\\ Prove that the unit group of $\Z[\sqrt 2]$ is infinite.
\end{framed}
\begin{proof}[Solution]
Let $u = 1 + \sqrt{2} \in \Z[\sqrt{2}]$. The norm of $u$ is
\[
N(u) = (1 + \sqrt{2})(1 - \sqrt{2}) = 1 - 2 = -1.
\]
Since the norm is a unit in $\Z$, $u$ is a unit in $\Z[\sqrt{2}]$, by Theorem \ref{quad_int_unit}. We claim that $u^n$ is a distinct unit for every $n\in\N$. Clearly $u^n$ is a unit for all $n$, since it is the product of units. It remains to show that $u^m \ne u^n$ if $m \ne n$. Suppose $u^m = u^n$ with $m < n$. Then $1 = u^{n - m}$. However, since $u=1+\sqrt 2>1$, all powers $u^s$ with $s\ge 1$ satisfy $u^s>1$. The identity $1 = u^{n - m}$ is thus impossible, so $u^m \ne u^n$ if $m \ne n$. It follows that the set $\{u^n:n\in\N\}$ is infinite, and therefore the unit group of $\Z[\sqrt 2]$ is infinite.
\end{proof}

\begin{framed}
\noindent\textbf{Exercise \ref{chap3_sel2}, Chapter \ref{more_rings_chap}:}\\ Show that there exist infinitely many $x\in\H$ such that $x^2=-1$.
\end{framed}
\begin{proof}[Solution]
Let $x = ai + bj$ for real numbers $a$ and $b$, and compute $x^2$:
\[
x^2 = (ai + bj)^2 = a^2i^2 + ab(ij + ji) + b^2j^2.
\]
Using the quaternion identities $i^2 = j^2 = -1$ and $ij = k = -ji$, we find:
\[
x^2 = -a^2 - b^2.
\]
Hence, if $(a,b)$ is any point on the unit circle in $\R^2$, then the quaternion $x=ai+bj$ satisfies $x^2=-1$. The infinitude of such elements $x$ follows immediately.
\end{proof}

\begin{framed}
\noindent\textbf{Exercise \ref{chap4_sel}, Chapter \ref{poly_chap}:}\\
Let $p(x)=x^4-2x^3-3x^2+4x+2$. Show that $\pm\sqrt 2$ are roots of $p(x)$ in $\Z[\sqrt 2]$ and that $p(x)$ factors as $(x-r_1)(x-r_2)(x-r_3)(x-r_4)$ with each $r_i\in\Z[\sqrt 2]$. Find all the roots $r_i$ explicitly.

\end{framed}
\begin{proof}[Solution]
Let $p(x) = x^4 - 2x^3 - 3x^2 + 4x + 2$. We work in the ring $\Z[\sqrt{2}]$. To check whether $\pm\sqrt{2}$ are roots, substitute each into $p(x)$:
\[
p(\sqrt{2}) = (\sqrt{2})^4 - 2(\sqrt{2})^3 - 3(\sqrt{2})^2 + 4\sqrt{2} + 2.
\]
Using $\sqrt{2}^2 = 2$ and $\sqrt{2}^3 = 2\sqrt{2}$, $\sqrt{2}^4 = 4$:
\[
p(\sqrt{2}) = 4 - 4\sqrt{2} - 6 + 4\sqrt{2} + 2 = 0.
\]
So $\sqrt{2}$ is a root. Similarly, we check that $-\sqrt 2$ is a root. Hence, the polynomial $(x - \sqrt{2})(x + \sqrt{2}) = x^2 - 2$ should be a factor of $p(x)$ in $\Z[\sqrt{2}][x]$.

Perform polynomial division to factor out $x^2 - 2$:
\[
p(x) = (x^2 - 2)(x^2 - 2x - 1).
\]

We now factor the second quadratic:
\[
x^2 - 2x - 1 = (x - r_1)(x - r_2),\quad\text{with } r = \frac{2 \pm \sqrt{(-2)^2 + 4}}{2} = \frac{2 \pm \sqrt{8}}{2} = 1 \pm \sqrt{2}.
\]
Therefore,
\[
p(x) = (x - \sqrt{2})(x + \sqrt{2})(x - (1 + \sqrt{2}))(x - (1 - \sqrt{2})).
\]

In particular, the four roots in $\Z[\sqrt{2}]$ are $\sqrt{2}, -\sqrt{2}, 1 + \sqrt{2}, 1 - \sqrt{2}$.
\end{proof}

\begin{framed}
\noindent\textbf{Exercise \ref{chap5_sel1}, Chapter \ref{power_chap}:}\\
Let $f$ be the power series over $\F_7$ given by
\[f=\sum_{i=1}^{\infty}(i+1)i\,x^{i-1}=2+6x+5x^2+6x^3+\cdots.\]
Compute the first few terms of $f^{-1}$.

\end{framed}
\begin{proof}[Solution]
We are given the power series
\[
f = \sum_{i=1}^\infty (i+1)i\, x^{i-1} = 2 + 6x + 5x^2 + 6x^3 + \cdots
\]
in $\F_7\pow x$. Since the constant term is 2, a unit in $\F_7$, $f$ is a unit in $\F_7\pow x$.

Let $f^{-1} = \sum_{n=0}^\infty b_nx^n$. To compute the first few $b_n$, we use the recursive formula from Chapter 5:
\[
b_0 = a_0^{-1}, \quad b_n = -a_0^{-1} \sum_{i=1}^n a_i b_{n-i},
\]
where $a_i$ are the coefficients of $f$.

We compute in $\F_7$:
\begin{align*}
a_0 &= 2 \quad\Rightarrow\quad b_0 = 2^{-1} = 4, \\
a_1 &= 6 \quad\Rightarrow\quad b_1 = -4 \cdot a_1 b_0 = 2, \\
a_2 &= 5 \quad\Rightarrow\quad b_2 = -4(a_1 b_1 + a_2 b_0)=5, \\
a_3 &= 6 \quad\Rightarrow\quad b_3 = -4(a_1 b_2 + a_2 b_1 + a_3 b_0)=3.
\end{align*}

Therefore, the first few terms of $f^{-1}$ are:
\[
f^{-1} = 4 + 2x + 5x^2 + 3x^3 + \cdots.\qedhere
\]
\end{proof}

\begin{framed}
\noindent\textbf{Exercise \ref{chap5_sel2}, Chapter \ref{power_chap}:}\\
Let $A$ be an integral domain. Show that
\[\ord(p+q)\ge\min\{\ord p,\ord q\}\] for all $p,q\in A\pow x$, and that equality holds if $\ord p\ne\ord q$. Find an example where equality fails when $\ord p=\ord q$.

\end{framed}
\begin{proof}[Solution]
Let $p = \sum a_i x^i$ and $q = \sum b_i x^i$ be elements of $A\pow x$. Define $r = p + q = \sum (a_i + b_i)x^i$. By definition,
\[
\begin{aligned}
\ord(p) &= \inf\{i : a_i \ne 0\},\\
\ord(q) &= \inf\{i : b_i \ne 0\},\\
\ord(r) &= \inf\{i : a_i + b_i \ne 0\}.
\end{aligned}
\]

Let $m = \ord(p)$ and $n = \ord(q)$, and suppose without loss of generality that $m \le n$. If $m < n$, then $a_m \ne 0$ and $b_m = 0$, so $a_m + b_m = a_m \ne 0$, and therefore $\ord(r) = m = \min\{\ord(p), \ord(q)\}$.

If $m > n$, the same argument shows $\ord(r) = n = \min\{\ord(p), \ord(q)\}$. If $m = n$, then it may happen that $a_m + b_m = 0$, in which case the order of $r$ is strictly greater than $m$.

In all cases, $\ord(r) \ge \min\{\ord(p), \ord(q)\}$, and equality holds when $\ord(p) \ne \ord(q)$.

For the required example, let $p = 3x + 2x^2 + \cdots$ and $q = -3x + 5x^3 + \cdots$ in $\Z\pow x$. Then $\ord(p) = \ord(q) = 1$, but $p + q = 2x^2 + 5x^3 + \cdots$ has $\ord(p + q) = 2 > 1$. Hence, equality can fail when $\ord(p) = \ord(q)$.
\end{proof}

\begin{framed}
\noindent\textbf{Exercise \ref{chap6_sel}, Chapter \ref{multi_chap}:}\\
Let $A$ be a commutative ring and $f\in A[x_1,\ldots, x_n]$ a nonzero polynomial of degree $d$. Show that if $f$ is homogeneous, then \[f(\lambda x_1,\ldots,\lambda x_n)=\lambda^d f(x_1,\ldots, x_n)\] for every $\lambda\in A$. Conversely, show that if $A$ is an infinite integral domain and the above identity holds for every $\lambda\in A$, then $f$ is homogeneous.
\end{framed}
\begin{proof}[Solution]
Suppose $f$ is homogeneous of degree $d$. Then every monomial of $f$ has total degree $d$. Let $m =  x_1^{e_1} \cdots x_n^{e_n}$ be such a monomial with $e_1 + \cdots + e_n = d$. For every $\lambda\in A$ we have
\[
\begin{aligned}
m(\lambda x_1, \ldots, \lambda x_n) 
&= (\lambda x_1)^{e_1} \cdots (\lambda x_n)^{e_n} \\
&= \lambda^{e_1 + \cdots + e_n} x_1^{e_1} \cdots x_n^{e_n} \\
&= \lambda^d m(x_1, \ldots, x_n).
\end{aligned}
\]
Writing $f$ as a linear combination of monomials and applying the above result shows that
\[
f(\lambda x_1, \ldots, \lambda x_n) = \lambda^d f(x_1, \ldots, x_n).
\]

Conversely, suppose that $A$ is an infinite integral domain and that the identity $f(\lambda x_1, \ldots, \lambda x_n) = \lambda^d f(x_1, \ldots, x_n)$ holds for all $\lambda \in A$. Write
\[
f = \sum_{i=0}^d f_i,
\]
where each $f_i$ is the homogeneous component of degree $i$. Then, for all $\lambda\in A$,
\[
f(\lambda x_1, \ldots, \lambda x_n) = \sum_{i=0}^d \lambda^i f_i(x_1, \ldots, x_n).
\]
It follows that
\[\lambda^df=\sum_{i=0}^d \lambda^i f_i\quad\text{for all}\;\lambda\in A.\]

Let $R=A[x_1, \ldots, x_n]$ and consider the polynomial $g\in R[X]$ given by
\[g(X)=X^d(f_d-f)+\sum_{i<d}X^i f_i.\]
Then we have $g(\lambda)=0$ for all $\lambda\in A$. In particular, $g$ has infinitely many roots in $R$, and since $R$ is a domain (by Theorem \ref{degthm_multi}), this implies that $g$ is the zero polynomial (by Theorem \ref{roots_number}). Hence, $f = f_d$ is homogeneous of degree $d$.
\end{proof}

\begin{framed}
\noindent\textbf{Exercise \ref{chap7_sel}, Chapter \ref{frac_chap}:}\\
Let $D$ be a domain and $F=\Frac D$. Show that addition in $F$ is well defined (independent of the choice of representatives). In other words, show that if $a/b=a'/b'$ and $c/d=c'/d'$, then \[\frac{ad+bc}{bd}=\frac{a'd'+b'c'}{b'd'}.\]
\end{framed}
\begin{proof}[Solution]
Suppose $a/b = a'/b'$ and $c/d = c'/d'$ in $F$. By definition of equality in $F$, we have
\[
ab' = a'b \quad \text{and} \quad cd' = c'd.
\]
We aim to show that
\[
\frac{ad + bc}{bd} = \frac{a'd' + b'c'}{b'd'},
\]
or equivalently, that \[(ad+bc)b'd'=(a'd'+b'c')bd.\] 

Computing the left-hand side of the above equation:
\[
(ad + bc)(b'd') = adb'd' + bcb'd' = a(b'd')d + b(c b'd') = a'd bd' + b c'd bd'.
\]

Computing the right-hand side:
\[
(a'd' + b'c')(bd) = a'd'bd + b'c'bd = a'd bd' + b c'd bd'.
\]

Since the two expressions obtained are equal, this proves the desired identity. Hence, addition is well defined.
\end{proof}

\section*{Additional practice}

\begin{exercises}
\item Let $C$ be the set of all continuous functions $\R\to\R$. Determine whether $(C,+,\circ)$ is a ring, where $\circ$ denotes composition of functions.
\item Let $R$ be a ring in which $x^3=x$ for every $x\in R$. Show that $R$ is commutative. Repeat the exercise if instead $x^4=x$ for all $x$.
\item Find all idempotent and all nilpotent elements in the rings $\Z_6$ and $\Z_8$, and find a nonzero nilpotent $2\times 2$ matrix over $\Z$.
\item Show that the set of nilpotent elements in a commutative ring is closed under addition.
\item Let $F$ be a field and $n$ a positive integer. Show that in the ring $M_n(F)$ the set of zero divisors is the complement of the set of units.
\item Let $S$ be a ring and $x,y\in S$. Show that if $1-xy\in S^{\times}$ then $1-yx\in S^{\times}$.
\item An \col{ordered ring}\index{ordered ring} is a pair $(R,\le)$, where $R$ is a ring and $\le$ is a \href{https://en.wikipedia.org/wiki/Partially_ordered_set}{partial order} on $R$ satisfying the following for all $a,b,c\in R$: if $a\ge 0$ and $b\ge 0$ then $ab\ge 0$, and if $a\le b$ then $a+c\le b+c$.
\begin{enumerate}
\item Show that the field $\C$ cannot be ordered, in the sense that there does not exist a partial order on $\C$ making $\C$ into an ordered ring.
\item Show that no Boolean ring can be ordered.
\end{enumerate}
\item Let $R$ be a commutative ring and $f\in R[x_1,\ldots, x_n]$ a polynomial  of degree $d$. The \col{homogenization}\index{homogenization of a polynomial} of $f$ with respect to a new variable, $z$ (or any other symbol), denoted by $f^h$, is the polynomial in $R[x_1,\ldots, x_n,z]$ obtained by replacing each monomial $m$ of $f$ with $m\cdot z^{d-\deg m}$.
\begin{enumerate}
\item[(a)] Prove that $f^h$ is homogeneous of degree $d$.
\item[(b)] If $R$ is a domain, show that \[f^h(x_1,\ldots,x_n,z)=z^df(x_1/z,\ldots,x_n/z).\]
\item[(c)] Compute $f^h$ if $f(x,y)=5x^3y-71xy+y^2+x^3+1$.
\end{enumerate}
\item Let $F$ be a finite field with $q$ elements. Suppose that $f,g\in F[x_1,\ldots, x_n]$ are polynomials of degree $<q$ in every variable, and that $f(a)=g(a)$ for all $a\in F^n$. Show that $f=g$.
\end{exercises}

\chapter{Subrings}\label{subrings_chap}

In this chapter and the next, we introduce two fundamental concepts in ring theory: subrings of a ring and ideals of a ring. Subrings and ideals are both defined as certain types of subsets of a ring, and as we will see, the only subset of a ring $R$ that is both a subring and an ideal is $R$ itself. Thus, ideals and subrings form essentially disjoint classes of subsets of $R$. By analyzing these types of subsets, we gain substantial information about $R$ and its ring-theoretic properties.

Let $R$ be a ring. A \col{subring}\index{subring} of $R$ is a subset $S\subseteq R$ that contains $0_R$ and $1_R$, is closed under addition and multiplication (meaning that these binary operations restrict to operations on $S$), and is a ring equipped with those operations (meaning that $(S,+,\cdot)$ satisfies all the ring axioms). Informally, $S$ is a ring contained inside the larger ring $R$, and operations in $S$ are carried out exactly as in $R$. 

Clearly, $R$ is a subring of itself. Any subring other than $R$ is a \col{proper subring}.

The following theorem provides a succinct criterion for checking whether a particular subset of a ring is a subring. We leave it to the reader to prove the theorem using the ring axioms.

\begin{framed}
\begin{thm}\label{subring_crit}
If $R$ is a ring and $S$ is a subset of $R$, then $S$ is a subring of $R$ if and only if it contains $0_R$ and $1_R$, and is closed under multiplication and subtraction (i.e., if $x,y\in S$ then $xy\in S$ and $x-y\in S$).
\end{thm}
\end{framed}

If $R$ is a field, a \col{subfield}\index{subfield} of $R$ is a subring closed under multiplicative inversion, i.e., a subring $S$ such that $x^{-1}\in S$ for every nonzero $x\in S$. In particular, $S$ is itself a field, and operations in $S$ agree with those in $R$. A \col{proper subfield} of $R$ is any subfield other than $R$.

\begin{ex}[Subrings and subfields]\; 
\begin{examples}
\item Each of these rings is a subring of the next: \[\Z,\,\Q,\,\R,\,\C,\,\C[z],\,\C\conv z,\,\C\pow z,\,\C\lau z.\]
Can you think of ways to continue this list indefinitely? Selecting the fields in the above list, we obtain examples of subfields: \[\Q\subset\R\subset\C\subset\C\lau z.\]
\item The field $\R$ is a subring of the quaternion ring $\H$.
\item The matrix ring $M_n(\Q)$ is a subring of $M_n(\C)$. More generally, if $S$ is a subring of $R$, then $M_n(S)$ is a subring of $M_n(R)$. 
\item Within the ring $\calC[a,b]$, the set of differentiable functions is a subring.
\item The quadratic integer ring $\Z[\sqrt d]$ is a subring of the quadratic field $\Q(\sqrt d)$, which is in turn a subfield of $\C$.
\item Every subring of an integral domain is itself a domain.
\item Within the quaternion ring $\H$, the set of elements $a+bi+cj+dk$ with $a,b,c,d\in\Z$ is a subring of $\H$ called the ring of \col{integer quaternions}. Similarly, the ring of \col{rational quaternions} is defined by the requirement that $a,b,c,d\in\Q$.
\end{examples}
\end{ex}

The following is a fundamental fact about subrings:

\begin{framed}
\begin{thm}\label{subring_int}
Let $R$ be a ring. Every intersection of subrings of $R$ is a subring of $R$. If $R$ is a field, every intersection of subfields is a subfield of $R$.
\end{thm}
\end{framed}

It follows from the above theorem that the collection of all subrings of $R$ is closed under arbitrary intersections. This allows us, in particular, to show that every ring $R$ has a unique smallest subring, called the \col{prime subring}\index{prime subring} of $R$, and defined to be the intersection of all subrings of $R$. This ring can be described more concretely as \[\{n\cdot 1_R: n\in\Z\}.\]

Thus, the prime subring of $R$ consists of all sums of the form \[1_R+1_R+\cdots+1_R\] and their additive inverses; see Exercise \ref{prime_sub}.

Similarly, if $F$ is a field, there is a unique smallest subfield of $F$, called the \col{prime subfield}\index{prime subfield} of $F$, and defined as the intersection of all subfields of $F$. Note that the prime subfield of $F$, being itself a subring of $F$, must contain the prime subring of $F$.
 
\begin{ex}[Prime subrings and subfields]\;
\begin{examples}
\item The prime subring of $\Z_n$ is $\Z_n$, and the prime subring of $\Z$ is $\Z$.
\item The prime subring of $\C$ is $\Z$, and the prime subfield of $\C$ is $\Q$.
\item For any ring $R$, the prime subring of $M_n(R)$ is the set of integer multiples of the identity matrix.
\item If $S$ is a subring of $R$, then $S$ and $R$ have the same prime subring. If $F$ is a subfield of $E$, then $F$ and $E$ have the same prime subfield.
\item Using the previous observation, we find that the prime subring of $\H$ is $\Z$ since $\R$ is a subring of $\H$ and its prime subring is $\Z$, and the prime subfield of $\C(z)$ is $\Q$ since $\C$ is a subfield of $\C(z)$.
\end{examples}
\end{ex}

Let $R$ be a ring and $R_0$ its prime subring. The \col{characteristic}\index{characteristic of a ring} of $R$, denoted by $\ch(R)$, is defined by

\[\ch(R)=
\begin{cases}
\#R_0, & \text{if}\; R_0\;\text{is finite};\\
0&\text{otherwise.}
\end{cases}\]

For example, $\ch(\Z_n)=n$ and $\ch(\Z)=0$. In later chapters, we will develop alternate descriptions of the characteristic and study its main properties. For the moment, it is sufficient to know the above definition.

\subsection*{Subrings generated by sets} The following is a common way of constructing subrings of a ring. Let $R$ be a ring and $S$ a subring of $R$. For any \underline{subset} $U\subseteq R$, let $S[U]$ denote the intersection of all subrings of $R$ that contain both $S$ and $U$. By Theorem \ref{subring_int}, $S[U]$ is in fact a subring of $R$, and is therefore the smallest subring of $R$ containing $S$ and $U$. The ring $S[U]$ is called the \col{subring generated by $S$ and $U$}. Note that if $U$ includes any elements of $S$, these elements may be removed from $U$ without changing the ring $S[U]$; put differently, $S[U]=S[U\setminus S]$. For this reason, we will often consider the ring $S[U]$ when $U\cap S=\emptyset$, in which case $S[U]$ is regarded as a ring obtained by adjoining to $S$ all the elements of $U$. Informally, the shorthand phrase `$S$ adjoin $U$' may be used to refer to $S[U]$.

 If $U=\{u_1,\ldots, u_n\}$ is a finite set, we denote $S[U]$ by $S[u_1,\ldots, u_n]$. Though this notation may seem to conflict with that of a polynomial ring over $S$, the two meanings are compatible, as the following theorem explains:

\begin{framed}
\begin{thm}\label{subring_adjoin}
Let $R$ be a ring, $S$ a subring of $R$, and $u_1,\ldots, u_n\in R$. The ring $S[u_1,\ldots, u_n]$ consists of all the elements of the form $p(u_1,\ldots, u_n)$, where $p\in S[x_1,\ldots, x_n]$ is an $n$-variable polynomial.
\end{thm}
\end{framed}
\begin{proof}
Any subring of $R$ containing $S$ as well as $u_1,\ldots, u_n$ will contain all powers of the elements $u_i$ and any product of such powers with elements of $S$; hence, it will contain the set $A$ defined by \[A=\{p(u_1,\ldots, u_n):p\in S[x_1,\ldots, x_n]\}.\] 

In particular, $S[u_1,\ldots, u_n]\supseteq A$. Now note that $A$ contains $S$ because any element $s\in S$ can be expressed as $p(u_1,\ldots, u_n)$, where $p$ is the constant polynomial $p(x_1,\ldots, x_n)=s$. Also, $A$ contains each $u_i$, as $u_i$ is obtained by evaluating the monomial $x_i$ at $(u_1,\ldots, u_n)$. It follows that $A$ contains the ring $S[u_1,\ldots, u_n]$, and is therefore equal to it.
\end{proof}

An entirely analogous process is used to construct subfields of any field $F$. Given a subfield $k$ of $F$ and a subset $U$ of $F$, the \col{subfield generated by $k$ and $U$}, denoted by $k(U)$, is the intersection of all the subfields of $F$ containing $k$ and $U$. If $U=\{u_1,\ldots, u_n\}$ is finite, then $k(U)$ consists of all elements of the form $p(u_1,\ldots, u_n)/q(u_1,\ldots, u_n)$, where $p,q\in k[x_1,\ldots, x_n]$ are polynomials such that $q(u_1,\ldots, u_n)\ne 0$.

\begin{ex}[Subrings generated by subsets]\;
\begin{examples}
\item In $\Q$, we may adjoin $1/3$ to $\Z$ to obtain the ring $\Z[1/3]$. Elements of this ring have the form $p(1/3)$, where $p\in\Z[x]$. Using this description, we deduce that
\[\Z[1/3]=\left\{\frac{n}{3^k}:n\in\Z,k\ge 0\right\}.\]
\item Let $A$ be any commutative ring. In the polynomial ring $A[x,y]$, consider the subring $A[x][y]$, i.e., the subring generated by $A[x]$ and $y$. Since this subring contains $A$, $x$, and $y$, it must be the entire ring $A[x,y]$. Thus, $A[x,y]=A[x][y]$. Similarly, $A[x,y]=A[y][x]$.
\item In $\C$, the subring generated by $\Z$ and $i$ is $\Z[i]$, the ring of Gaussian integers. To justify this, we begin with the description of $\Z[i]$ as elements of the form $p(i)$, where $p\in\Z[x]$. Since $i^2=-1$, the powers $i^n$, $n\in\N$, are all equal to $\pm 1$ or $\pm i$; thus, $p(i)$ can be expressed in the form $a+bi$ with $a,b\in\Z$, and therefore $\Z[i]=\{a+bi:a,b\in\Z\}$. More generally, if $d\ne 1$ is any squarefree integer, the subring of $\C$ generated by $\Z$ and $\sqrt d$ is the quadratic integer ring $\Z[\sqrt d]$. Similarly, the quadratic field $\Q(\sqrt d)$ is the subfield of $\C$ generated by $\Q$ and $\sqrt d$; see Exercise \ref{quad_subfield}.
\item In the quaternion ring $\H$, we may consider the subring $\Z[i,j]$, which consists of elements of the form $p(i,j)$, where $p\in\Z[x,y]$. Since $ij=k$, this ring contains $k$, and therefore contains every integer quaternion, i.e., every quaternion of the form $a+bi+cj+dk$ with $a,b,c,d\in\Z$. Conversely, the quaternion identities $i^2=j^2=-1$, $ij=k$, and $ji=-k$ imply that every element of $\Z[i,j]$ can be expressed in the above form; hence, $\Z[i,j]$ is equal to the ring of integer quaternions.
\end{examples}
\end{ex}

\subsection*{The lattice of subrings} For every ring $R$, let $\Sub(R)$ denote the set of all subrings of $R$, partially ordered by inclusion. This poset has an interesting structure in its own right: in the terminology of order theory, $\Sub(R)$ is a \href{https://en.wikipedia.org/wiki/Complete_lattice}{complete lattice}, meaning that every subset of $\Sub(R)$ has a join (or supremum) and a meet (or infimum) in $\Sub(R)$. The lattice $\Sub(R)$ is called the \col{lattice of subrings} of $R$.

The completeness of $\Sub(R)$ can be described explicitly using the construction of a ring generated by a subset. If $\mathcal{S}$ is any subset of $\Sub(R)$, i.e., any collection of subrings of $R$, the intersection of all rings in $\mathcal{S}$ is a subring of $R$ (by Theorem \ref{subring_int}), and is therefore the meet of $\mathcal{S}$ in $\Sub(R)$. The join of $\mathcal{S}$ is not the union of the rings in $\mathcal{S}$ (which may not be a subring), but instead is the intersection of all subrings containing this union. In the case of just two subrings $A,B\in\Sub(R)$, the meet is $A\cap B$ and the join is the ring $A[B]=B[A]$, in the notation of the earlier construction $S[U]$.

An analogous construction applies to fields. If $F$ is a field, the set of all subfields of $F$ forms a complete lattice called the \col{lattice of subfields} of $F$. The meet of two subfields $K$ and $L$ is $K\cap L$, and the join is the subfield $K(L)=L(K)$, the smallest subfield of $F$ containing both $K$ and $L$.

The figure below depicts a small portion of the subring lattice of $\C$.

\bigskip
\begin{figure}[h]
\centering
\begin{tikzpicture}[node distance=2cm, every node/.style={draw, rectangle, rounded corners, minimum size=1.2em, inner sep=4pt}]

\node (Z)         at (0,0)      {$\mathbb{Z}$};
\node (Zi)        at (-3,2)     {$\mathbb{Z}[i]$};
\node (Zsqrt2)    at (3,2)      {$\mathbb{Z}[\sqrt{2}]$};
\node (Zpi2)      at (5.5,2)    {$\mathbb{Z}[\pi^2]$};
\node (Zpi)       at (4,4)      {$\mathbb{Z}[\pi]$};
\node (ZiPi2)     at (-0.5,4)   {$\mathbb{Z}[i, \pi^2]$};
\node (Q)         at (0,2)      {$\mathbb{Q}$};
\node (Qi)        at (-2,4)     {$\mathbb{Q}(i)$};
\node (Qsqrt2)    at (2,4)      {$\mathbb{Q}(\sqrt{2})$};
\node (R)         at (0,6)      {$\mathbb{R}$};
\node (C)         at (0,8)      {$\mathbb{C}$};

\draw (Z) -- (Zi);
\draw (Z) -- (Zsqrt2);
\draw (Z) -- (Zpi2);
\draw (Z) -- (Q);
\draw (Zi) -- (Qi);
\draw (Zsqrt2) -- (Qsqrt2);
\draw (Zpi2) -- (Zpi);
\draw (Zpi2) -- (ZiPi2);
\draw (Zi) -- (ZiPi2);
\draw (Zpi) -- (R);
\draw (Q) -- (Qi);
\draw (Q) -- (Qsqrt2);
\draw (Qi) -- (C);
\draw (Qsqrt2) -- (R);
\draw (R) -- (C);

\end{tikzpicture}
\caption{
A portion of the lattice of subrings of \( \mathbb{C} \). Edges indicate inclusion.
}
\end{figure}
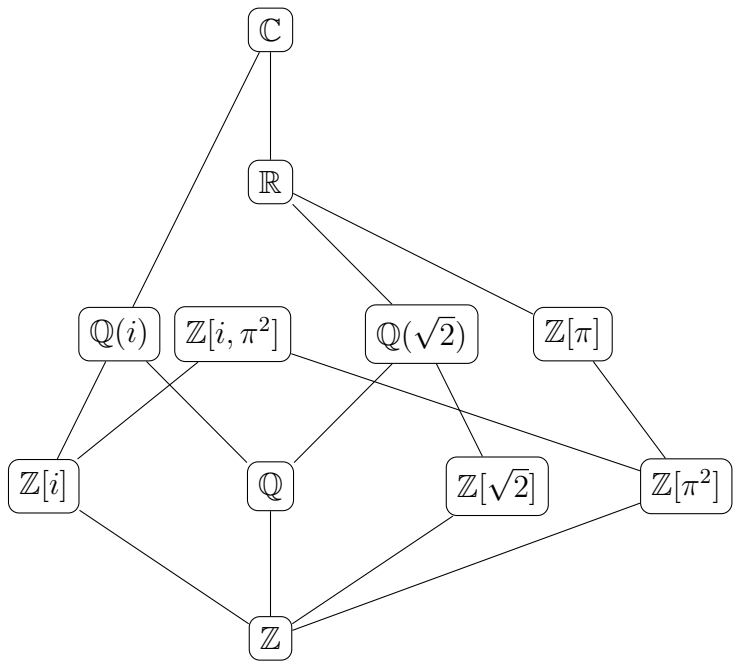

\newpage
\section*{Exercises}

\begin{note}\;
\begin{itemize}
\item Solutions to starred exercises (\sel) are included in Chapter~\ref{review2chap}.
\item A house (\house) marks exercises suggested as sample homework.
\end{itemize}
\end{note}

\medskip

\begin{exercises}
\item Let $R\subseteq M_2(\Z)$ be the set of matrices of the form \[\begin{pmatrix}a&a-b\\a-b&b\end{pmatrix}.\] Prove or disprove that $R$ is a subring of $M_2(\Z)$.
\item Determine which of the following subsets of $\calC[0,1]$ are subrings:
\begin{enumerate}
\item[(a)] The set of $f\in\calC[0,1]$ such that $0\le\int_0^1 f\le 1$.
\item[(b)] The set of functions in $\calC[0,1]$ having a finite number of zeros, together with the zero function.
\end{enumerate}
\item(\house)\label{int_poly_def} In the polynomial ring $\Q[x]$, let $R$ be the set of polynomials $f$ such that $f(n)\in\Z$ for all $n\in\Z$. Show that $R$ is a subring of $\Q[x]$, and give an example of an element of $R$ that is not in $\Z[x]$. The ring $R$ is known as the \col{ring of integer-valued polynomials}.
\item Let $S$ be a subring of a ring $R$. Prove or provide a counterexample: if $A$ is a subring of $S$, then $A$ is a subring of $R$.
\item Give an example to show that the union of two subrings of a ring is not necessarily a subring.
\item\label{prime_sub} Let $R$ be a ring. Show that the prime subring of $R$ is equal to the set $\{n\cdot 1_R: n\in\Z\}$.
\item Let $R=\Z[i]$. Show that the subring $\Z[2i]$ is a proper subring of $R$, and describe its elements.
\item(\house) Determine whether the subring $\Z[x^2,x^3]\subseteq\Z[x]$ is a proper subring, and describe its elements.
\item Show that the rings $\Z[x]$ and $\Z[\sqrt 5]$ have infinitely many subrings.
\item Prove that every subring of a field is an integral domain.
\item Prove that the prime subring of every ring is a commutative ring.
\item(\sel)\label{chap9_sel} Give an example (with proof) of a non-commutative ring whose prime subring is a field.
\item Prove that every nontrivial Boolean ring has characteristic 2.
\item Let $R$ be the ring of integer quaternions. Show that $R$ has exactly eight units. The group $R^{\times}$ is known as the \col{quaternion group}\index{quaternion group}, denoted $Q_8$.
\item Let $F$ be a field. In the field of Laurent series $F\lau x$, consider the subring $F[x,\frac{1}{x}]$ obtained by adjoining $x$ and $1/x$ to $F$. Show that $F[x,\frac{1}{x}]$ is the set of \col{Laurent polynomials}\index{Laurent polynomial}, i.e., Laurent series having only finitely many nonzero terms.
\item\label{quad_subfield} If $d\ne 1$ is a squarefree integer, show that the subfield of $\C$ generated by $\Q$ and $\sqrt d$ is the quadratic field $\Q(\sqrt d)$.
\item Find $\alpha\in\R$ such that $\Q(\alpha)=\Q(\sqrt 2,\sqrt 3)$.
\item(\house)\label{localization_exc} Let $R$ be an integral domain. A \col{multiplicative subset} of $R$ is a subset $S\subseteq R$ such that $1\in S$, $0\notin S$, and if $a\in S$ and $b\in S$, then $ab\in S$. Let $F=\Frac(R)$. The set 
\[R_S:=\left\{\frac{x}{y}\in F:y\in S\right\}\]

is called the \col{localization}\index{localization} of $R$ at $S$, and is also denoted by $S^{-1}R$.
\begin{enumerate}
\item Show that $R_S$ is a subring of $F$ containing $R$, and that every element of $S$ is a unit in $R_S$.
\item Show that every subring of $\Q$ is a localization of $\Z$ at some multiplicative subset.
\end{enumerate}
\end{exercises}

\chapter{Ideals}\label{ideals_chap}

Following our treatment of subrings in the previous chapter, we now discuss a different type of subset of a ring. Ideals lie at the core of ring theory, and will be used constantly for the remainder of the course.

An \col{ideal}\index{ideal} of a ring $R$ is a nonempty subset $I\subseteq R$ that is closed under addition and under multiplication by elements of $R$, meaning that
\[a\in R,\; x\in I\implies ax\in I\;\;\text{and}\;\;xa\in I.\]

Informally, we may view the above as a \emph{scaling property}, analogous to closure under scalar multiplication in a vector space.

Clearly, $R$ is an ideal of itself; any ideal other than $R$ is a \col{proper ideal}. Every ideal $I$ contains the element 0, as a consequence of the scaling property: for every $x\in I$, $0\cdot x\in I$. However, ideals---unlike subrings---usually do not contain 1. In fact, the only ideal that contains 1 is $R$ itself, again by the scaling property: if $1\in I$, then $a\cdot 1\in I$ for every $a\in R$, and thus $I=R$. By similar reasoning, $R$ is the only ideal that contains a unit.

One key reason for studying ideals is that every ideal $I$ of $R$ can be used to construct a simplified version of the ring $R$, denoted by $R/I$, whose properties may be easier to study than those of $R$ and yet provide nontrivial information about $R$. The details of this construction are given in the next chapter.

\begin{ex}[Ideals]\;
\begin{examples}
\item In any ring $R$, the subsets $R$ and $\{0\}$ are ideals; the latter is called the \col{trivial ideal} or the \col{zero ideal}\index{zero ideal}.
\item Every field $F$ has exactly two ideals, namely $F$ and $\{0\}$. Indeed, if $I$ is a nonzero ideal of $F$, let $x\in I$ be nonzero. Since $x$ is a unit, we must have $I=F$. 
\item In the polynomial ring $A[x]$, where $A$ is a commutative ring, the set of polynomials with constant term 0 is an ideal.
\item Let $n\in\Z$. In the ring of integer quaternions, the set of elements of the form $a+bi+cj+dk$ with $a,b,c,d$ all divisible by $n$ is an ideal.
\item For every $n\in\Z$, the set $n\Z=\{nk:k\in\Z\}$ is an ideal of $\Z$. We will show in Theorem \ref{Zideals} that every ideal of $\Z$ is equal to $n\Z$ for some integer $n$.
\item Let $I$ be an ideal of a commutative ring $A$, and let $I[x]$ denote the set of polynomials in $A[x]$ all of whose coefficients belong to $I$. Then $I[x]$ is an ideal of $A[x]$.
\item The ring $\calC[a,b]$ has uncountably many ideals: for $r\in [a,b]$, the set \[M_r=\{f\in\calC[a,b]:f(r)=0\}\] 

is an ideal. Graphically, elements of $M_r$ are continuous curves over $[a,b]$ that intersect the $x$-axis at $x=r$.
\end{examples}
\end{ex}

The following is a fundamental property of ideals, and is a straightforward consequence of the definition of an ideal.

\begin{framed}
\begin{thm}\label{ideal_int}
Let $R$ be a ring. Every intersection of ideals of $R$ is an ideal of $R$.
\end{thm}
\end{framed}

\subsection*{Ideals generated by sets} Theorem \ref{ideal_int} allows us to define the smallest ideal containing a given subset of a ring $R$. If $U\subseteq R$ is any subset, the \col{ideal generated by $U$}, denoted by $(U)$, is the intersection of all ideals of $R$ containing $U$. If $U=\{u_1,\ldots,u_n\}$ is a finite set, we may denote $(U)$ by $(u_1,\ldots, u_n)$. If $U$ is countably infinite, we may use notation such as $(u_1,u_2,u_3,\ldots)$.

An ideal $I$ is called \col{finitely generated}\index{finitely generated ideal} if it can be generated by a finite set; in other words, there exist elements $r_1,\ldots, r_n\in I$ such that $I=(r_1,\ldots, r_n)$. If $I$ can be generated by a single element, then $I$ is called a \col{principal ideal}\index{principal ideal}. Principal ideals can be described quite concretely: for any $r\in R$, we have \[(r)=\{arb:a,b\in R\}.\] If $R$ is commutative, the description is even simpler: \[(r)=\{ar:a\in R\}.\] More generally, if $R$ is commutative, the ideal generated by a subset $U$ has a fairly explicit description:

\begin{framed}
\begin{thm}\label{ideal_gen}
If $R$ is a commutative ring and $U$ is a nonempty subset of $R$, then the ideal $(U)$ is the set of all sums of the form \[a_1u_1+\cdots+a_nu_n\] with each $a_i\in R$ and $u_i\in U$. Thus, elements of $(U)$ are linear combinations of elements of $U$ with scalars in $R$.
\end{thm}
\end{framed}
\begin{proof}
Letting $J$ denote the set of all $R$-linear combinations of elements of $U$, we aim to prove that $J=(U)$. To begin, note that the set $J$ is nonempty and is closed under addition and scaling. Hence, $J$ is an ideal of $R$. Moreover, $J$ clearly contains $U$. As $(U)$ is defined to be the intersection of all ideals containing $U$, this implies that $J\supseteq (U)$.

To show that $J\subseteq(U)$, let $s=a_1u_1+\cdots+a_nu_n\in J$, where each $u_i\in U$ and each $a_i\in R$. If $I$ is any ideal of $R$ containing $U$, the scaling property implies that each product $a_iu_i$ belongs to $I$ since $u_i\in I$, and therefore $s\in I$. As this applies to every ideal $I$ containing $U$, $s$ belongs to the intersection of all ideals containing $U$, which is $(U)$ by definition. Hence, $J\subseteq (U)$ and thus $J=(U)$.
\end{proof}

As a consequence of Theorem \ref{ideal_gen} we obtain the following: if $R$ is commutative and $r_1,\ldots, r_n\in R$, then
\[(r_1,\ldots, r_n)=\{a_1r_1+\cdots+a_nr_n:a_i\in R\}.\]
In particular, for every $r\in R$, $(r)=\{ar:a\in R\}$, as observed earlier.

\begin{ex}[Ideals generated by sets]\;
\begin{examples}
\item In every ring $R$, $(0)=\{0\}$ and $(1)=R$. For every $u\in R^{\times}$, $(u)=R$.
\item The only ideals of $\H$ are $(0)$ and $(1)$ since every nonzero quaternion is a unit, as proved in Chapter \ref{more_rings_chap}.
\item In $\Z$, the principal ideal $(n)$ is equal to $n\Z$, the set of all integer multiples of $n$.
\item Let $R$ be a ring with ideal $I$. If $r_1,\ldots, r_n\in I$, then $(r_1,\ldots, r_n)\subseteq I$.
\item In $\Z$, $(6,15)=(3)$ because both 6 and 15 are multiples of 3, so $6,15\in(3)$, and $3=6(-2)+15$, so $3\in (6,15)$.
\item In the polynomial ring $A[x]$, where $A$ is a commutative ring, the principal ideal $(x)$ is the set of polynomials with constant term 0.
\item In the power series ring $A\pow x$, the ideal $(x)$ is equal to the set of power series $f$ with $\ord(f)>0$.
\end{examples} 
\end{ex}
 
A ring in which every ideal is principal is called a \col{principal ideal ring}\index{principal ideal ring (PIR)} (PIR). In the case where such a ring is an integral domain, it is called a \col{principal ideal domain}\index{principal ideal domain (PID)} (PID). As we will see in Theorems \ref{Zideals} and \ref{kxideals} below, examples of PID's include $\Z$, $k[x]$, and $k\pow x$, where $k$ is a field.

In the opposite extreme, there exist rings with plenty of ideals that are not finitely generated, much less principal. For example, let $R = \Z[x_1, x_2, x_3, \dots]$, a polynomial ring with infinitely many variables, and let $I$ be the ideal generated by all the variables:
\[I = (x_1, x_2, x_3, \dots).\]

We claim that $I$ is not finitely generated. Note, first of all, that every element $f\in I$ has constant term 0, as $f$ is a sum of terms of the form $x_ig_i$, by Theorem \ref{ideal_gen}, and each $x_ig_i$ has constant term 0. 

Suppose, by contradiction, that $I$ is generated by a finite set of elements $f_1, \dots, f_n\in I$. Since each $f_i$ is a polynomial, it involves only finitely many variables. Thus, there exists some index $N$ such that all of the generators $f_1, \dots, f_n$ lie in the subring $\mathbb{Z}[x_1, \dots, x_N]$. Now, since $x_{N+1}\in I$ by definition of $I$, and $I=(f_1,\ldots, f_n)$, there must exist $h_1,\ldots, h_n\in R$ such that
\[x_{N+1}=h_1f_1+\cdots+h_nf_n.\]

Since $x_{N+1}$ appears in the left-hand side of this equation, it must occur in at least one polynomial $h_i$; otherwise, it would not appear on the right-hand side. However, since $f_i$ has no constant term, any occurrence of $x_{N+1}$ in $h_if_i$ will appear multiplied by some of the variables $x_1,\ldots, x_N$. Hence, the monomial $x_{N+1}$ cannot appear in the sum $h_1f_1+\cdots+h_nf_n$. We have reached a contradiction, thus proving the claim that $I$ is not finitely generated.

\subsection*{Ideals in $\Z$, $k[x]$, and $k\pow x$} It will be important going forward to identify all the ideals of some basic examples of rings. We begin this discussion here by determining the ideals of $\Z$, $k[x]$, and $k\pow x$, where $k$ is a field. Further examples, including $\Z_n$ and $M_n(R)$, appear in the exercises. Interestingly, the argument used for $\Z$ will apply almost verbatim to the ring $k[x]$; this observation is the beginning of a long series of analogies between $\Z$ and polynomial rings in one variable over a field, particularly $\Q[x]$.

\begin{framed}
\begin{thm}\label{Zideals}
The ring $\Z$ is a PID, and if $k$ is a field, the polynomial ring $k[x]$ is a PID.
\end{thm}
\end{framed}
\begin{proof}
The zero ideal in $\Z$ is clearly principal, so we focus on proving that every nonzero ideal $I$ of $\Z$ is principal. The scaling property of ideals implies that $I$ must contain positive integers; letting $g$ be the smallest positive integer in $I$, we claim that $I=(g)$. Certainly $(g)\subseteq I$, so it suffices to prove the reverse inclusion. Let $a\in I$. By the Division Algorithm for integers, there exist $q,r\in\Z$ such that $a=qg+r$ and $0\le r<g$. Given that both $a$ and $g$ belong to $I$, and $I$ is an ideal, the integer $r=a-qg$ also belongs to $I$. The minimality of $g$ now implies that $r=0$, for otherwise $r$ would be a positive integer in $I$ which is less than $g$. Thus, $a=qg\in(g)$, proving that $I\subseteq(g)$. The proof for the ring $k[x]$ is nearly identical: using the Division Algorithm for polynomials (Theorem \ref{div_algo}), one can show that every nonzero ideal of $k[x]$ is generated by an element of least degree.
\end{proof}

The following theorem shows that the set of ideals of a power series ring over a field has a particularly simple structure.

\begin{framed}
\begin{thm}\label{kxideals}
Let $k$ be a field. The nonzero ideals of the power series ring $k\pow x$ are the principal ideals $(x^n)$ for $n\ge 0$. Hence, the ideals of $k\pow x$ lie in a descending chain \[(1)\supset (x)\supset (x^2)\supset(x^3)\supset\cdots\supset (0).\] 
\end{thm}
\end{framed}

\begin{proof}
Let $I$ be a nonzero ideal of $k\pow x$. Let $f(x)\in I$ have minimal order of vanishing among all nonzero elements of $I$. We claim that $I=(f(x))=(x^n)$, where $n=\ord(f)\ge 0$. Writing
\[f(x) = a_n x^n + a_{n+1} x^{n+1} + \cdots\]
with $a_n \ne 0$, we have $f(x) = x^n u(x)$, where $u(x)=a_n+a_{n+1}x+\cdots$ is a unit in $k\pow x$ since its constant term is nonzero (Theorem \ref{ordthm}). Now let $v(x)=u(x)^{-1}\in k\pow x$, and note that $x^n=f(x)v(x)\in I$ by the scaling property. Thus, $(x^n)\subseteq I$.

To show that $I\subseteq(x^n)$, let $g(x) \in I$ be nonzero and let $m=\ord(g)$, so that we may write $g(x) = x^m t(x)$ for some $t(x)\in k\pow x$. By minimality of $n$, we have $m \ge n$, hence $g(x)=x^n\cdot x^{m-n}t(x)\in(x^n)$. This proves that $I=(x^n)$, as claimed.

We have shown that the nonzero ideals of $k\pow x$ are the ideals $(x^n)$ for $n\ge 0$. A straightforward argument shows that these ideals form a strictly descending chain.
\end{proof}

\subsection*{The lattice of ideals} We conclude this chapter by discussing a lattice analogous to the lattice of subrings of a ring. For any ring $R$, let $\Id(R)$ denote the set of all ideals of $R$, partially ordered by inclusion. Every subset $\mathcal{I}\subseteq\Id(R)$ has a meet (infimum) in the lattice, namely the set $\cap_{A\in\mathcal{I}}A$, which is an ideal by Theorem \ref{ideal_int}. Moreover, $\mathcal{I}$ has a join (supremum) in $\Id(R)$, namely the ideal generated by the set $\cup_{A\in\mathcal{I}}A$. Thus, $\Id(R)$ is a complete lattice, and is called the \col{lattice of ideals}\index{lattice of ideals} of $R$.

The join described above is also called the \col{sum}\index{sum of ideals} of the ideals in $\mathcal{I}$ and denoted by $\sum_{A\in\mathcal{I}}A$, as it can be shown to consist of all finite sums $r_1+\cdots +r_n$ with each $r_i$ belonging to some ideal in $\mathcal{I}$ (not necessarily the same ideal for each $r_i$), See Exercise \ref{ideal_sum_exc}.

In the case of only two ideals, say $I$ and $J$, the sum is given by
\[I+J=\{a+b:a\in I,b\in J\}.\]

Analogously, there is a natural notion of \col{product}\index{product of ideals} of $I$ and $J$, denoted by $IJ$, though this product cannot simply be the set $\{ab:a\in I,\,b\in J\}$, as this set is not generally an ideal. Instead, $IJ$ is defined to be the ideal generated by this set. Explicitly: \[IJ:=\{a_1b_1+\cdots+a_nb_n: a_k\in I, b_k\in J\;\text{for all}\;k\}.\]
The defining properties of ideals imply that $IJ\subseteq I\cap J$, but equality does not hold in general, so $IJ$ is not the meet of $I$ and $J$ in the lattice $\Id(R)$.

\section*{Exercises}

\begin{note}\;
\begin{itemize}
\item Solutions to starred exercises (\sel) are included in Chapter~\ref{review2chap}.
\item A house (\house) marks exercises suggested as sample homework.
\end{itemize}
\end{note}

\medskip

\begin{exercises}
\item In the ring $R=M_2(\Z)$, let $V$ and $W$, respectively, be the sets of matrices of the form
\[\begin{pmatrix} a & 0 \\ b & 0 \end{pmatrix}\quad\text{and}\quad \begin{pmatrix} 2a & 2b \\ 2c & 2d \end{pmatrix}\]
with $a,b,c,d\in\Z$. Determine whether $V$ and $W$ are ideals of $R$.
\item In the ring $R=\calC[0,1]$, let $Z$ denote the set of functions having infinitely many zeros in $[0,1]$. Is $Z$ an ideal of $R$?
\item Find examples to show the following: if $I$ and $J$ are ideals of a ring, the sets $I\cup J$ and $\{ab:a\in I,\,b\in J\}$ are not necessarily ideals.
\item(\house)\label{asc_chain} Let $I_1\subseteq I_2\subseteq I_3\subseteq\cdots$ be an ascending chain of ideals of a ring. Show that the union $\cup_{n=1}^{\infty}I_n$ is an ideal.
\item Prove or provide a counterexample: if $S$ is a subring of a ring $R$, and $I$ is an ideal of $S$, then $I$ is an ideal of $R$.
\item Show that every ideal of a ring is an Abelian group with respect to addition.
\item If $R$ is a nontrivial commutative ring whose only ideals are $(0)$ and $(1)$, prove that $R$ is a field.
\item(\sel)\label{chap10_sel} Prove that in every commutative ring, the set of all nilpotent elements is an ideal. This ideal is called the \col{nilradical}\index{nilradical} of the ring. In the ring $M_2(\Z)$, show that the set of nilpotent elements is not an ideal.
\item Follow these steps to prove that $\Z_n$ is a PIR for every $n\ge 1$:
\begin{enumerate}
\item[(a)] Let $J$ be an ideal of $\Z_n$, and let $I$ be the set of integers whose reduction modulo $n$ belongs to $J$:\[I=\{a\in\Z:a \equiv b \;(\mathrm{mod}\ n)\;\text{for some}\;b\in J\}.\] Show that $I$ is an ideal of $\Z$.
\item[(b)] Letting $M$ be a generator of $I$, show that $J=(m)$, where $m$ is the reduction of $M$ modulo $n$, i.e., $m$ is the unique element of $\Z_n$ such that $M\equiv m\pmod n$.
\end{enumerate}
\item(\house) Let $I$ and $J$ be ideals of a ring. Show that the ideal generated by the set $\{ab:a\in I,\,b\in J\}$ is equal to the set of all finite sums $a_1b_1+\cdots+a_nb_n$ with each $a_k\in I$ and $b_k\in J$.
\item In $\H$, show that $(1+j)=(1)$. In the ring of integer quaternions, use the norm function to show that $(1+j)$ is a proper ideal.
\item\label{ideal_sum_exc} Let $R$ be a ring and $\mathcal{I}$ a subset of the lattice of ideals of $R$. Show that the join of the ideals in $\mathcal{I}$ is equal to the set of all finite sums $r_1+\cdots+r_n$ such that for each index $1\le s\le n$ there exists $A_s\in\mathcal{I}$ with $r_s\in A_s$.
\item\label{comax_exc} Let $I$ and $J$ be ideals of a ring $R$.
\begin{enumerate}
\item[(a)] Show that $IJ\subseteq I\cap J$ and give an example to show that equality need not hold.
\item[(b)] Show that if $R$ is commutative and $I+J=(1)$, then $IJ=I\cap J$.
\end{enumerate}
\item Prove the identity $(ab)=(a)(b)$ for principal ideals in a commutative ring.
\item\label{non_pid} Show that in the ring $\Z[x]$, the ideal $(2,x)$ is not principal. (Hence, the assumption in Theorem \ref{Zideals} that $k$ is a field cannot be weakened to assume only that $k$ is a domain.)
\item Let $A$ be a nontrivial commutative ring. Show that if $n>1$, then the polynomial ring $A[x_1,\ldots, x_n]$ is not a PIR.
\item(\house) Prove that the ring $\Z[\sqrt{-5}]$ is not a PID by showing that the ideal $(3,2+\sqrt{-5})$ is not principal.
\item Let $B$ be a Boolean ring. Prove that every finitely-generated ideal of $B$ is a principal ideal.
\item Let $R$ be a commutative ring. A \col{homogeneous ideal}\index{homogeneous ideal} of the polynomial ring $R[x_1,\ldots, x_n]$ is an ideal $I$ with the property that whenever $p\in I$, every homogeneous component of $p$ belongs to $I$. Show that an ideal is homogeneous if and only if it can be generated by a set of homogeneous polynomials. 
\item Let $I$ be an ideal of a ring $R$. In the matrix ring $M_n(R)$, let $M_n(I)$ be the set of matrices with entries in $I$. Show that $M_n(I)$ is an ideal of $M_n(R)$, and that the map $I\mapsto M_n(I)$ is a bijection between the set of ideals of $R$ and the set of ideals of $M_n(R)$. What are the ideals of $M_n(\C)$?
\item Give an example of a non-commutative PIR.
\end{exercises}

\chapter{Direct Products and Quotient Rings}\label{quo_chap}

In this chapter we study two standard constructions for creating new rings starting from known rings. The \emph{direct product} construction combines any collection of rings to form a larger ring, and the \emph{quotient ring} construction uses an ideal of a ring to construct a simplified version of that ring.

\section*{Direct products}

If $R$ and $S$ are rings, we may consider the Cartesian product $R\times S$, initially only a set, and use the operations of $R$ and $S$ to endow this set with a ring structure in a natural way: for pairs $(a,b)$ and $(c,d)$ in $R\times S$ we define
\[(a,b)+(c,d)=(a+c,b+d)\quad\text{and}\quad (a,b)\cdot(c,d)=(ac,bd).\]

Crucially, in this definition, the operations of $R$ are applied in the left-hand entry and the operations of $S$ in the right-hand entry. For example, in $\Z\times\F_5$ we have
\[(-3,2)+(7,3)=(4,0)\quad\text{and}\quad(-3,2)\cdot(7,3)=(-21,1).\]

With these componentwise operations, $R \times S$ acquires a natural ring structure called the \col{direct product}\index{direct product of rings} of $R$ and $S$. The zero element of $R\times S$ is $(0_R,0_S)$ and the unity is $(1_R,1_S)$.

The direct product construction can be generalized to any Cartesian product of rings $\prod_{i\in I}R_i$, once again defining operations componentwise. In particular, for any ring $R$ and positive integer $n$ we may form the ring $R^n$ consisting of ordered $n$-tuples of elements of $R$.

\begin{ex}[Direct products of rings]\;
\begin{examples}
\item In the ring $\H\times\Z_6\times\Q(\sqrt 2)$ we compute
\[(1+i+j,\;2,\;5-\sqrt 2)\cdot(k,\;3,\;\sqrt 2)=(i-j+k,\;0,\;5\sqrt 2-2).\]
\item Powers of elements of $R^n$ are computed componentwise: if $m\in\N$,
\[(r_1,\ldots, r_n)^m=(r_1^m,\ldots, r_n^m).\]
Similarly, integer multiples are computed componentwise: if $m\in\Z$,
\[m(r_1,\ldots, r_n)=(mr_1,\ldots, mr_n).\]
\item In the ring $\F_2^n$, every element is idempotent since this holds true in $\F_2$.  Thus, $\F_2^n$ is a Boolean ring.
\item Although $\R^2$ and $\C$ are both 2-dimensional vector spaces over $\R$, they are rather different as rings. For instance, $\R^2$  is not an integral domain (while $\C$ is even a field) because $(1,0)\cdot(0,1)=(0,0)$.
\item If $R$ and $S$ are nontrivial rings, then $R\times S$ is not a domain. Thus, direct products of rings are usually not domains.
\item The pair $(j,5)\in\H\times\Z_9$ is a unit with inverse $(j,5)^{-1}=(-j,2)$.
\item The unit group of a direct product is the Cartesian product of the individual unit groups: if $R=\prod_{i\in I}R_i$, then \[R^{\times}=\prod_{i\in I}R_i^{\times}.\]
\end{examples}
\end{ex}

The above examples address, in particular, the question of identifying the units and zero divisors in a direct product of rings. We turn now to the question of describing the subrings and ideals of a direct product, beginning with a simple observation:
\begin{framed}
\begin{thm}
Let $R=\prod_{i\in I}R_i$ be a direct product of rings.
\begin{enumerate}
\item[(a)] If $S_i$ is a subring of $R_i$ for each $i$, then the product $\prod_{i\in I}S_i$ is a subring of $R$.
\item[(b)] If $A_i$ is an ideal of $R_i$ for each $i$, then $\prod_{i\in I}A_i$ is an ideal of $R$.
\end{enumerate}
\end{thm}
\end{framed}

We are naturally led to ask whether the above theorem describes \emph{all} the subrings and ideals of a direct product of rings. The answer is different for subrings than for ideals.

There is no universal classification of all subrings of a direct product of rings. In particular, a subring of a direct product is not necessarily a direct product of subrings of the individual factors. For example, in the ring $\F_2\times\F_2$, the set $\{(0,0),(1,1)\}$ is a proper subring, but cannot be expressed as a product of subrings of $\F_2$, as $\F_2$ has no proper subrings.

Despite this complexity, there are certain special cases where a complete description of all subrings of a direct product can be achieved; some examples are included in the exercises.

Describing all the ideals of a direct product of rings is often more feasible than describing its subrings. The following theorem provides a complete description of the ideals in any direct product of finitely many rings.

\begin{framed}
\begin{thm}\label{prod_ideals}
Let $R$ and $S$ be rings. The ideals of $R\times S$ are the sets of the form $I\times J$, where $I$ is an ideal of $R$ and $J$ is an ideal of $S$. An analogous statement applies to any direct product of finitely many rings.
\end{thm}
\end{framed}
\begin{proof}
We show that if $A$ is an ideal of $R\times S$, then $A=I\times J$ for some ideals $I$ and $J$. Let $\pi_1$ and $\pi_2$ be the projection maps $R\times S\to R$ and $R\times S\to S$, respectively, and let $I=\pi_1(A)$ and $J=\pi_2(A)$. A straightforward argument shows that $I$ and $J$ are ideals.

Clearly, if $(x,y)\in A$, then $x\in I$ and $y\in J$, so $(x,y)\in I\times J$. Hence, $A\subseteq I\times J$. Now let $(x,y)\in I\times J$, so that there exist $r\in R$ and $s\in S$ such that $(x,s)\in A$ and $(r,y)\in A$. Since $A$ is an ideal of $R\times S$, we conclude from the identity \[(x,y)=(x,s)\cdot(1,0)+(r,y)\cdot(0,1)\]

that $(x,y)\in A$. Hence, $I\times J\subseteq A$.
\end{proof}

Theorem \ref{prod_ideals} does not extend to direct products of infinitely many rings. For example, let $R=\prod_{i=1}^\infty\Z$, the ring of all integer sequences, and let $J\subset R$ be the set of sequences having only finitely many nonzero terms. Then $J$ is a proper ideal of $R$, but cannot be expressed as a direct product of ideals of $\Z$: indeed, if $J=\prod_{i=1}^\infty n_i\Z$, then, since $J$ contains all the sequences having a 1 in one position and zeros elsewhere, it follows that $1\in n_i\Z$ for every $i$, and therefore $J=R$, a contradiction.

There is no general classification theorem describing all the ideals in an infinite direct product of rings. However, there are some instances where a complete description is possible. Exercise \ref{field_prod_ideal} discusses one notable case, namely any direct product of fields.

\section*{Quotient rings}

Recall that for any nonzero integer $n$, the equivalence relation of \emph{congruence modulo} $n$ in $\Z$ is defined by \[a\equiv b\ (\mathrm{mod}\ n)\iff n\;\;\text{divides}\;\; a-b.\] 

Moreover, this relation is compatible with addition and multiplication, in the sense that if $a\equiv b\pmod n$ and $c\equiv d\pmod n$, then \[a+c\equiv b+d\;\, (\mathrm{mod}\ n)\;\;\text{and}\;\; ac\equiv bd\;\, (\mathrm{mod}\ n).\] 

These facts lead to one common way of defining the ring $\Z_n$, namely, as the quotient set $\Z/{\equiv}$, i.e.,  the set of equivalence classes of integers under the relation of congruence modulo $n$.

In what follows, we will generalize the above notions to arbitrary rings: given any ring $R$ and any ideal $I$ in $R$, we will define the concept of \emph{congruence modulo} $I$ as well as a ring $R/I$ analogous to the ring $\Z_n$. Informally, we may describe the construction of $R/I$ as coloring the elements of $R$ with some number of colors, and then forming a new ring whose elements are the colors themselves. In this sense, $R/I$ is a simplification of the ring $R$, as $\Z_n$ is a simplification of $\Z$.

A \col{congruence}\index{congruence} on a ring $R$ is an equivalence relation $\equiv$ on $R$ that respects the ring operations: if $a\equiv b$ and $c\equiv d$, then $a+c\equiv b+d$ and $ac\equiv bd$. The quotient set of $R$ by the congruence will be denoted by $R/{\equiv}$, and the equivalence class of an element $a\in R$ by $[a]$. Explicitly,
\[[a]=\{b\in R:a\equiv b\}.\]
The set $R/{\equiv}$ can be made into a ring by defining
\[[x]+[y]=[x+y]\quad\text{ and }\quad [x]\cdot [y]=[xy]\]
for all $x,y\in R$. 

Note that the sum and product in the expressions $[x+y]$ and $[xy]$ are the operations of the ring $R$, so the ring $R/{\equiv}$ can be said to inherit its operations from $R$. The properties of a congruence ensure that the operations in $R/{\equiv}$ are well defined (independent of the representatives of the equivalence classes). Put differently, congruence classes can be added by choosing \emph{any} representatives of the classes, adding those representatives, and forming the corresponding congruence class. Likewise for multiplication. Clearly, the zero element of $R/{\equiv}$ is $[0]$ and the unity is $[1]$.

The following theorem shows that congruences on $R$ are in one-to-one correspondence with ideals of $R$; this is one way in which ideals arise naturally in ring theory.

\begin{framed}
\begin{thm}\label{cong_ideal}
Let $R$ be a ring. If $\equiv$ is a congruence on $R$, the equivalence class $[0]=\{r\in R:r\equiv 0\}$ is an ideal of $R$. Conversely, if $I$ is an ideal of $R$, then the relation defined by $a\equiv_I b\iff a-b\in I$ is a congruence on $R$. Moreover, the mappings \[\equiv\;\;\mapsto\; [0]\quad\text{and}\quad I\;\mapsto\;\;\equiv_I\] 

are inverse bijections between the set of ideals of $R$ and the set of congruences on $R$.
\end{thm}
\end{framed}
\begin{proof} Exercise \ref{ideal_to_congruence}.\end{proof}

The congruence $\equiv_I$ corresponding to an ideal $I$ is called \col{congruence modulo $I$}\index{congruence modulo an ideal}. Moving forward, we will write $a\equiv b\pmod I$ instead of $a\equiv_I b$. Thus,
\[a \equiv b \;(\mathrm{mod}\ I)\iff a-b\in I.\] 

In the case where $I=(x)$ is a principal ideal, we may ease notation by writing $a\equiv b\pmod x$ rather than $a\equiv b\pmod{(x)}$.

The ring $R/{\equiv_I}$ is denoted simply by $R/I$, and called \col{$R$ modulo $I$}. Elements of $R/I$ are the congruence classes $[x]$, where $x\in R$. Note that
\[[x]=\{y\in R:x-y\in I\}=\{x+r:r\in I\}.\]

For this reason, $[x]$ is commonly denoted by $x+I$, especially when the ideal $I$ needs to be emphasized in the notation. Elements of $R/I$ may thus be described as \emph{translates} of $I$ by elements of $R$, also called \emph{cosets} of $I$, borrowing a term from group theory. Using coset notation, the zero element of $R/I$ is $0+I=I$, and the unity is $1+I$. This completes the construction of the quotient ring $R/I$. 

To summarize the key points of the above construction: elements of $R/I$ have the form $x+I$ with $x\in R$, and operations are defined by \[(x+I)+(y+I)=(x+y)+I,\quad (x+I)(y+I)=xy+I.\] 

Moreover, we have \[x+I=y+I\iff x\equiv y\;(\mathrm{mod}\ I)\iff x-y\in I.\]

If the ideal $I$ is clear from context, we may simplify notation and write $[x]$ instead of $x+I$; in this case, operations are given by \[[x]+[y]=[x+y]\quad\text{ and }\quad [x]\cdot [y]=[xy].\]

\begin{ex}[Quotient rings]\;
\begin{examples}
\item For any nonzero integer $n$, the relation of congruence modulo the ideal $(n)=n\Z$ is precisely congruence modulo $n$ as defined at the beginning of this chapter. The quotient ring $\Z/(n)=\Z/n\Z$ consists of the equivalence classes $[0],[1],\ldots,[n-1]$, as every integer is congruent modulo $n$ to one of the integers $0,1,\ldots, n-1$. Operations on these classes correspond to operations in $\Z$ followed by reduction modulo $n$; for instance, $[1]+[n-1]=[n]=[0]$. Hence, the ring $\Z/(n)$ is simply an alternate description of $\Z_n$.
\item Elements of the quotient ring $\Z[x]/(x^2)$ can be uniquely represented in the form $[a+bx]$ with $a,b\in\Z$. Indeed, if \[f(x)=a_0+a_1x+a_2x^2+\cdots+a_nx^n\in\Z[x],\] then $f(x)\equiv a_0+a_1x\pmod{x^2}$, so $[f(x)]=[a_0+a_1x]$. Moreover, if $[a+bx]=[c+dx]$, then the polynomial $g(x)=a-c+(b-d)x$ is a multiple of $x^2$, which cannot occur unless $g(x)=0$ and therefore $a=c$ and $b=d$. Representing elements of $\Z[x]/(x^2)$ as above, multiplication takes the form \[[a+bx]\cdot [c+dx]=[ac+(ad+bc)x].\]

Consequently, the zero divisors in this quotient ring are the classes $[bx]$ with $b\in\Z$ nonzero. The units are the classes $[1+bx]$ and $[-1+bx]$ with $b\in\Z$ arbitrary.
\item We claim that the quotient ring $\Z[i]/(1+i)$ has exactly two elements. Note first that $[0]\ne[1]$ in this ring, as equality would imply that $1\equiv 0\pmod{(1+i)}$, hence $1\in(1+i)$ and therefore $1+i$ is a unit. However, the latter conclusion is false given that the norm of $1+i$ is $N(1+i)=2$. (Recall Theorem \ref{quad_int_unit}.) We now leave it to the reader to verify that a Gaussian integer $a+bi$ is congruent to either 0 or 1 modulo $(1+i)$, depending on whether $a+b$ is even or odd, respectively. Thus, $[0]$ and $[1]$ are the only two elements of the ring $\Z[i]/(1+i)$.
\item One standard use of quotient rings is to construct rings where a given polynomial equation has a solution. Let $R$ be the polynomial ring $\R[t]$ and let $I=(t^2+1)$. In the quotient $R/I$, let $\alpha=[t]$. Then \[\alpha^2=[t^2]=[-1],\] since $t^2\equiv-1\pmod I$. Hence, the equation $x^2=-1$ has a solution in $R/I$, namely, $x=\alpha$. We will see in the next chapter that $R/I$ is essentially the field $\C$, with $\alpha$ playing the role of the number $i$.
\end{examples}
\end{ex}

The following theorem is our main result concerning quotient rings. The theorem describes the subrings and ideals of a quotient $R/I$ in terms of the subrings and ideals of $R$.

For any subset $S\subseteq R$, let $S/I$ denote the set of translates of $I$ by elements of $S$:
\[S/I:=\{s+I:s\in S\}\subseteq R/I.\]

\begin{framed}
\begin{thm}\label{quo_sub_id}
Let $R$ be a ring and $I$ an ideal of $R$. The map $S\mapsto S/I$ induces inclusion-preserving bijections 
\[\{\text{Ideals of $R$ containing $I$}\}\longrightarrow\{\text{Ideals of $R/I$}\}\]
and
\[\{\text{Subrings of $R$ containing $I$}\}\longrightarrow\{\text{Subrings of $R/I$}\}.\]
In both cases, the inverse map is given by $A\mapsto\{r\in R:r+I\in A\}$.
\end{thm}
\end{framed}
\begin{proof} Exercise \ref{corresp_proof}.\end{proof}

As an application of the theorem, we describe all subrings and ideals of the ring $\Z_n$. Subrings of $\Z_n=\Z/(n)$ are in bijection with subrings of $\Z$ containing $(n)$. Since $\Z$ has only one subring (namely itself), the same must hold for $\Z_n$. Thus, $\Z_n$ has no proper subrings. Note that this same conclusion can easily be reached without using the theorem, so we are merely verifying that the theorem is consistent with previously established results.

The analysis of ideals of $\Z_n$ is more subtle. The ideals of $\Z/n\Z$ correspond to ideals of $\Z$ containing $n\Z$. Since $\Z$ is a PID (Theorem \ref{Zideals}), its ideals are of the form $d\Z$, and we may assume that $d$ is positive. In order for $d\Z$ to contain $n\Z$, we must have $n\in d\Z$, i.e., $n$ is a multiple of $d$, and thus $d$ is a divisor of $n$.

Theorem \ref{quo_sub_id} now implies that the ideals of $\Z/n\Z$ are precisely the sets of the form $d\Z/n\Z$, where $d$ is a positive divisor of $n$. Moreover, the set $d\Z/n\Z$ corresponds to the principal ideal $(d)$ in $\Z_n$. Hence, we conclude that the ideals of $\Z_n$ are the principal ideals $(d)=d\Z_n$, where $d$ is a positive divisor of $n$.

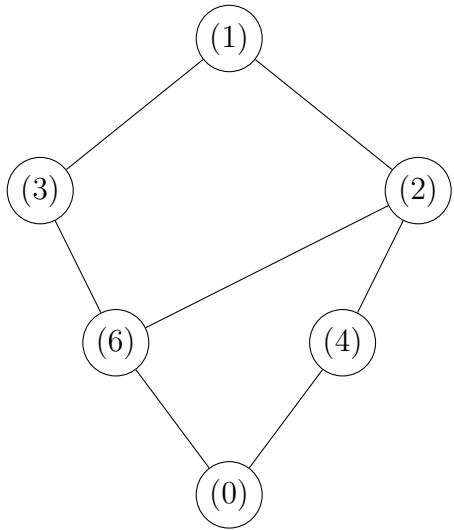
\begin{figure}[h]
\centering
\begin{tikzpicture}[
  every node/.style={draw, circle, inner sep=2pt},
  node distance=1.5cm and 1.8cm,
  ]

\node (0) at (0,0) {$(0)$};
\node (6) at (-1.5,2) {$(6)$};
\node (4) at (1.5,2) {$(4)$};
\node (3) at (-2.5,4) {$(3)$};
\node (2) at (2.5,4) {$(2)$};
\node (1) at (0,6) {$(1)$};

\draw (0) -- (6);
\draw (0) -- (4);
\draw (6) -- (3);
\draw (6) -- (2);
\draw (4) -- (2);
\draw (2) -- (1);
\draw (3) -- (1);

\end{tikzpicture}
\caption{The lattice of ideals of \( \mathbb{Z}_{12} \).}\label{Z12_lattice}
\end{figure}

Figure \ref{Z12_lattice} shows the full set of ideals of $\Z_{12}$. To explain the containment relations in the figure, note that if $d$ and $e$ are divisors of $n$, and $d$ divides $e$, then $(e)\subseteq (d)$.

\section*{Exercises}

\begin{note}\;
\begin{itemize}
\item Solutions to starred exercises (\sel) are included in Chapter~\ref{review2chap}.
\item A house (\house) marks exercises suggested as sample homework.
\end{itemize}
\end{note}

\medskip

\begin{exercises}
\item Find all the idempotent elements, nilpotent elements, units, and zero divisors in the ring $\F_3\times\F_5$.
\item Use a direct product to construct examples of the following:
\begin{enumerate}
\item A non-commutative ring having an infinite unit group and exactly ten ideals.
\item A ring having finite unit group and uncountably many ideals.
\end{enumerate}
\item Let $R=\prod_{i\in I} R_i$ be a direct product of rings. Let $S$ be the subset of $R$ consisting of tuples $(r_i)$ such that $r_i=0$ for all but finitely many $i\in I$. Show that $S$ is an ideal of $R$. This ideal is called the \col{direct sum}\index{direct sum of rings} of the rings $R_i$ and denoted by $\bigoplus_{i\in I} R_i$. (Of course, if the index set $I$ is finite, then $R=S$.)
\item(\house) Find all subrings and ideals of the ring $\Z_4\times\F_2$.
\item Show that a direct product of finitely many rings is a PIR if and only if every factor in the direct product is a PIR. Find a counterexample in the case of infinitely many rings.
\item\label{field_prod_ideal} Let $R=\prod_{i\in I}F_i$, where each $F_i$ is a field. In this exercise we determine all the ideals of $R$.
\begin{enumerate}
\item For any subset $S\subseteq I$, define \[I(S)=\{(r_i)\in R:r_i=0\,\text{ for all }\,i\in S\}.\] Show that $I(S)$ is an ideal of $R$.
\item For any ideal $J$ of $R$, define \[Z(J)=\{k\in I: r_k=0\,\text{ for all }\,(r_i)\in J\}.\] Show that $J\subseteq I(Z(J))$.
\item With $J$ as above, suppose that $(r_i)\in I(Z(J))$. Construct an element $(h_k)\in J$ such that $h_k\ne 0$ for all $k\notin Z(J)$. Use this to prove that $J=I(Z(J))$, and thus conclude that every ideal of $R$ is of the form $I(S)$ for some subset $S\subseteq I$.
\end{enumerate}
\item Determine which of the following relations on a ring are congruences:
\begin{enumerate}
\item On $\Q$, $a\sim b\iff a-b\in\Z$.
\item On $\Z$, $a\sim b\iff$ $a-b$ is even or $a=b=0$.
\item On $\Z\times\Z$, \[(a,b)\sim(c,d)\iff a\equiv c\;\,(\mathrm{mod}\ 3)\quad\text{and}\quad b\equiv d\;\,(\mathrm{mod}\ 5).\]
\end{enumerate}
\item Let $\equiv$ be a congruence on a commutative ring $R$. Show that if $a\equiv b$ in $R$, and $f\in R[x]$ is a polynomial, then $f(a)\equiv f(b)$.
\item Let $R=A\times B$ be a direct product of rings. Describe all congruences on $R$ in terms of congruences on $A$ and $B$.\item Let $R=\Z[x]$ and let $I$ be the ideal $(x^3+2x)$. Use polynomial long division to find $a,b,c\in\Z$ such that $x^5+x^6\equiv a+bx+cx^2\pmod I$.
\item Show that every quotient of a commutative ring is commutative.
\item(\house) Determine the cardinalities of the rings $\F_2[x]/(x^2+x+1)$ and $\F_2\pow x/(x^3)$.
\item\label{ideal_to_congruence} Prove Theorem \ref{cong_ideal}.
\item(\house) Use a quotient ring to construct examples of the following: 
\begin{enumerate}
\item For each positive integer $n$, a commutative ring having exactly $n$ ideals.
\item For each positive integer $n$, a non-commutative ring having exactly $n$ ideals.
\end{enumerate}
\item\label{corresp_proof} Prove the correspondence of ideals in Theorem \ref{quo_sub_id}.
\item Prove that every quotient ring of a PIR is a PIR. Give an example of a quotient of a PID that is not a PID.
\item(\sel)\label{chap11_sel} Prove that every quotient of the quadratic integer ring $\Z[\sqrt d]$ by a nonzero ideal is a finite ring.
\end{exercises}

\chapter{Isomorphisms and Homomorphisms}\label{homo_chap}

We introduce here a precise notion of structure-preserving map between rings. Such maps are used to establish a concrete link between two rings in such a way that the properties of one ring are reflected in those of the other.

A \col{homomorphism}\index{homomorphism} between rings $R$ and $S$ is a map $f:R\to S$ satisfying $f(0_R)=0_S$, $f(1_R)=1_S$, and
\[f(a+b)=f(a)+f(b),\quad f(ab)=f(a)f(b)\]
for all $a,b\in R$. An \col{isomorphism}\index{isomorphism} is a bijective homomorphism. Informally, a homomorphism represents a similarity of ring structures and an isomorphism represents an equivalence of structures.

The fundamental result of this chapter is the Isomorphism Theorem (Theorem \ref{iso_thm}), which associates to every homomorphism $f:R\to S$ an isomorphism that reveals the structure of the image of $f$. As applications of this theorem, we derive several results concerning the characteristic of a ring---most notably, that the prime subring of any ring is isomorphic either to $\Z$ or to $\Z_n$ for some $n$. We also introduce the \emph{Frobenius homomorphism}, a fundamental map for rings of prime characteristic (such as the fields $\F_p$), that plays a central role in understanding their algebraic properties. 

One classical consequence of the Frobenius homomorphism is a result from elementary number theory known as \href{https://en.wikipedia.org/wiki/Fermat\%27s_little_theorem}{Fermat's Little Theorem}, which states that for every prime $p$ and every integer $n$, $n^p \equiv n \pmod p$. A proof is given in Examples \ref{frob_hom_exs}.

\section*{Isomorphisms}

Let $R$ and $S$ be rings. We say that $R$ and $S$ are \col{isomorphic}, written $R\cong S$, if there exists an isomorphism $R\to S$. As the reader can readily verify, the relation of isomorphism is an equivalence relation on the class of all rings: the identity map $R\to R$ is an isomorphism, so $R\cong R$; the inverse of an isomorphism is an isomorphism, so if $R\cong S$ then $S\cong R$; and the composition of isomorphisms is an isomorphism, so if $R\cong S$ and $S\cong T$, then $R\cong T$.

When two rings are isomorphic, we may regard them as being essentially the same ring, merely described in different ways.

\begin{ex}[Ring isomorphisms]\;
\begin{examples}\label{iso_exs}
\item Let $R$ be a commutative ring and $n$ a positive integer. Let $D_n(R)$ be the ring of all $n\times n$ scalar matrices (diagonal matrices whose diagonal entries are all equal) with entries in $R$, a subring of the matrix ring $M_n(R)$. The natural map $R\to D_n(R)$ mapping $r\in R$ to $rI_n$ is an isomorphism. In particular, $D_n(R)\cong R$ for all $n$.
\item The quaternion ring $\H$ has a subring isomorphic to $\C$. To avoid confusion with notation, write quaternions in the form $a+bi+cj+dk$ and complex numbers in the form $a+bI$. The set $C=\{a+bi:a,b\in\R\}$ is a subring of $\H$ isomorphic to $\C$ via the map $a+bi\mapsto a+bI$. Thus $\H$ may be regarded as an extension of $\C$, an enlargement of $\C$ into a broader ring.
\item The quotient ring $\R[x]/(x^2+1)$ is isomorphic to $\C$. Indeed, the map $a+bi\;\mapsto\; [a+bx]$ is an isomorphism $\C\to\R[x]/(x^2+1)$. Surjectivity follows from the Division Algorithm for polynomials: every $p(x)\in\R[x]$ can be expressed in the form \[p(x)=q(x)(x^2+1)+r(x),\]

where $\deg r<2$. Hence $r(x)=a+bx$ for some $a,b\in\R$, and in the quotient ring $\R[x]/(x^2+1)$ we have $[p(x)]=[a+bx]$.
\item Isomorphisms map units to units: if $u\in R$ is a unit and $f:R\to S$ is an isomorphism, then $f(u)$ is a unit and $f(u)^{-1}=f(u^{-1})$. In particular, if $R$ and $S$ are fields, then $f(a/b)=f(a)/f(b)$ for all $a,b\in R$ with $b\ne 0$. Similarly, isomorphisms map zero divisors to zero divisors.
\end{examples}
\end{ex}

\subsection*{Proving non-isomorphism} Isomorphic rings share all set-theoretic and ring-theoretic properties. For instance, if the rings $R$ and $S$ are isomorphic, then
\begin{align*}
&R \text{ has cardinality } 6 \iff S \text{ has cardinality } 6;\\
&R \text{ has exactly five units } \iff S \text{ has exactly five units};\\
&R \text{ is commutative} \iff S \text{ is commutative};\\
&R \text{ is a domain} \iff S \text{ is a domain}.\\
\end{align*}

\vspace{-5mm}
\noindent This observation suggests a general strategy for proving that two rings are \emph{not} isomorphic: find a property possessed by one ring but not by the other. For example, $\Z_4$ and $\Z_5$ are not isomorphic because they have different cardinalities; $\Q$ and $\Z$ have the same cardinality but are not isomorphic because only one of them is a field.

A slightly more subtle strategy to distinguish between two ring structures is to use polynomial equations. For instance, the rings $\Q$ and $\Q(i)$ have the same cardinality, and both are fields, yet they are not isomorphic because the equation $x^2+1=0$ has a solution in only one of them. Indeed, if there existed an isomorphism $f:\Q(i)\to\Q$, then applying $f$ to both sides of the equation $i^2+1=0$ would yield $f(i)^2+1=0$, so $f(i)$ would be a rational number whose square is $-1$; a contradiction.

As a final example, we show that \[\Z_4\not\cong\Z_2\times\Z_2.\]
Both rings are commutative rings and have cardinality four; however, $\Z_4$ has a nonzero nilpotent element, namely $2$, whereas $\Z_2\times\Z_2$ has no such element. Hence, the two rings cannot be isomorphic, since any isomorphism would map nonzero nilpotent elements to nonzero nilpotent elements.

\subsection*{Automorphisms} An important special case of an isomorphism arises when the domain and codomain are the same ring. An \col{automorphism}\index{automorphism} of a ring $R$ is an isomorphism $R\to R$. The set of all automorphisms of $R$, equipped with the operation of composition, forms a group called the \col{automorphism group}\index{automorphism group} of $R$, denoted by $\Aut(R)$. The identity element of $\Aut(R)$ is the identity map $R\to R$, which is clearly an automorphism. Elements of $\Aut(R)$ may be regarded as \emph{symmetries} of the ring $R$.

\begin{ex}[Ring automorphisms]\;
\begin{examples}
\item Complex conjugation, i.e., the map $a+bi\mapsto a-bi$, is an element of $\Aut(\C)$, as follows from the identities \[\overline{z+w}=\bar z+\bar w,\quad\overline{zw}=\bar z\bar w\quad\text{for all}\;\;z,w\in\C.\] Note that $\bar{\bar z}=z$ for all $z\in\C$, so conjugation is an example of an \col{involution}\index{involution}, i.e., a map whose composition with itself is the identity map (equivalently, an automorphism that is equal to its inverse).
\item The group $\Aut(\Z)$ has only one element, namely the identity map. Any automorphism $f$ must map 0 to 0 and 1 to 1, and by induction, $f(n)=n$ for all $n\in\N$. Since $f(-n)=-f(n)=-n$, it follows that $f$ is the identity map on all of $\Z$. By similar reasoning, $\Aut(\Q)$ has only one element, as does $\Aut(\Z_n)$ for every $n\ge 1$.
\item In contrast to the examples of $\Z$, $\Z_n$, and $\Q$, the rational function field $\Q(x)$ has infinitely many automorphisms. For example, the map $f(x)\mapsto f(ax)$ is an automorphism for every nonzero $a\in\Q$; its inverse is the map $g(x)\mapsto g(x/a)$.
\item Recall that $\GL_n(\R)$ denotes the unit group of $M_n(\R)$, i.e., the set of $n\times n$ invertible matrices with real entries. For every matrix $U\in\GL_n(\R)$, the conjugation map $A\mapsto UAU^{-1}$ is an automorphism of $M_n(\R)$.
\end{examples}
\end{ex}

\section*{Homomorphisms}

Let $R$ and $S$ be rings. Recall that a map $f:R\to S$ is a homomorphism if it satisfies $f(0)=0$, $f(1)=1$,  and \[f(a+b)=f(a)+f(b),\quad f(ab)=f(a)f(b)\]
for all $a,b\in R$. Note that homomorphisms are not assumed to be bijective, injective, or surjective. Having discussed isomorphisms in the previous section, we now focus on more general homomorphisms. 

We begin with two minor observations: the identity \[f(a-b)=f(a)-f(b)\quad\text{for all}\;\;a,b\in R\] follows immediately from the definition, and the condition $f(0)=0$ may be dropped from the definition, as it follows from $f(a+b)=f(a)+f(b)$.

\begin{ex}[Homomorphisms]\label{homo_exs}\;
\begin{examples}
\item If $A$ is a commutative ring, there is a surjective homomorphism $A\pow x\to A$ mapping every power series $a_0+a_1x+a_2x^2+\cdots$ to its constant term, $a_0$.
\item If $R$ and $S$ are any rings, the projection maps $R\times S\to R$ and $R\times S\to S$ given by $(a,b)\mapsto a$ and $(a,b)\mapsto b$, respectively, are surjective homomorphisms.
\item Let $R$ be a ring and $I$ an ideal of $R$. The natural map $R\to R/I$ mapping $a\in R$ to its equivalence class $[a]=a+I$ is a surjective homomorphism called \col{reduction modulo $I$}\index{reduction modulo an ideal}.
\item Let $f:R\to S$ be a homomorphism. For every $n\ge 1$, there is an induced homomorphism $M_n(R)\to M_n(S)$ given by $(a_{ij})\mapsto(f(a_{ij}))$.
\item\label{unique_homZ} If $R$ is any ring, there exists a unique homomorphism $\Z\to R$, namely the map $f_R(n)=n\cdot 1_R$. The fact that $f_R$ is a homomorphism follows from Theorem \ref{identi}(b). Uniqueness follows from the observation that any homomorphism $\Z\to R$ maps $1$ to $1_R$, and this implies $n\mapsto n\cdot 1_R$ for all $n\in\Z$.
\end{examples}
\end{ex}

The \col{kernel}\index{kernel of a homomorphism} of a homomorphism $f:R\to S$ is the set \[\ker(f):=\{r\in R:f(r)=0_S\}=f^{-1}(0_S),\] and the  \col{image}\index{image of a homomorphism} of $f$ is the image of $f$ as a function, i.e., \[\im(f):=f(R)=\{s\in S:s=f(r)\text{ for some }r\in R\}.\]

The proof of the following theorem is left to the reader (Exercise \ref{ker_prop}).

\begin{framed}
\begin{thm}\label{im_ker}
Let $f:R\to S$ be a ring homomorphism. 
\begin{enumerate}
\item[(a)] The kernel of $f$ is an ideal of $R$ (the domain).
\item[(b)] The image of $f$ is a subring of $S$ (the codomain).
\item[(c)] The map $f$ is injective if and only if $\ker(f)=(0)$.
\end{enumerate}
\end{thm}
\end{framed}

A couple of observations regarding Theorem \ref{im_ker}:

\begin{itemize}
\item The theorem provides a new description of ideals: the ideals of a ring $R$ are precisely the kernels of homomorphisms whose domain is $R$. Indeed, every such kernel is an ideal, and conversely, every ideal $I$ is the kernel of a homomorphism, namely the reduction map $R\to R/I$.
\item An injective homomorphism between rings is called an \col{embedding}\index{embedding}. According to the theorem, a homomorphism is an embedding if and only if its kernel is trivial. Informally, an embedding is a way of placing one ring inside another while preserving its structure. Examples include the canonical map $a\mapsto a/1$ from a domain $D$ to its field of fractions, $\Frac(D)$, and the map $R\to R^2$, where $R$ is any ring, given by $x\mapsto(x,x)$; this map is called the \emph{diagonal embedding} of $R$ in $R^2$.
\end{itemize}

The following theorem is arguably the most consequential fact concerning isomorphisms and homomorphisms.
\begin{framed}
\begin{thm}[Isomorphism Theorem]\label{iso_thm}\index{Isomorphism Theorem}
If $f:R\to S$ is a ring homomorphism, then the map $[a]\mapsto f(a)$ is a well-defined isomorphism $R/\ker(f)\longrightarrow\im(f)$. In particular, \[R/\ker(f)\;\cong\;\im(f).\]
\end{thm}
\end{framed}
\begin{proof}
We begin by showing that the map $F:[a]\mapsto f(a)$ is well defined. If $[a]=[b]$ in $R/\ker(f)$, then $a-b\in\ker(f)$, so $f(a-b)=0$; hence $f(a)=f(b)$, as required. Next, we note that $F$ is a homomorphism because $f$ is. That $F$ is surjective follows immediately from the definitions; moreover, $F$ is injective because $F([a])=0$ implies $a\in\ker(f)$, hence $[a]=0$ in $R/\ker(f)$. This completes the proof that $F$ is an isomorphism.
\end{proof}

If $R$ is any ring, a \col{homomorphic image}\index{homomorphic image} of $R$ is the image of $R$ under any ring homomorphism. With this terminology, the Isomorphism Theorem states that every homomorphic image of $R$ is isomorphic to a quotient of $R$.

\begin{note} Theorem \ref{iso_thm} is sometimes called the \emph{First Isomorphism Theorem}, as there are also a \emph{Second} and \emph{Third} Isomorphism Theorems; see Exercises \ref{2nd_isothm} and \ref{third_iso_thm}. In this course we apply the First Isomorphism Theorem frequently, whereas the second and third will not be needed for developing the theory.
\end{note}

\begin{ex}[Using the Isomorphism Theorem]\label{iso_thm_exs}\;
\begin{examples}
\item Let $R$ and $S$ be rings. Since the projection map $R\times S\to R$ given by $(a,b)\mapsto a$ is a surjective homomorphism, $R$ is a homomorphic image of $R\times S$ and is therefore isomorphic to a quotient of $R\times S$. Explicitly, $R\cong (R\times S)/((0)\times S)$, as the ideal $(0)\times S$ is the kernel of the projection map. Similarly, $S$ is isomorphic to a quotient of $R\times S$.
\item For any ring $R$, $R/(0)\cong R$. This follows from the Isomorphism Theorem applied to the identity map $R\to R$.
\item\label{cab_max_idls} In the ring $\calC[a,b]$, we have previously defined ideals \[M_r=\{f\in\calC[a,b]:f(r)=0\}\] for every $r\in[a,b]$. For every such $r$, let $e_r:\calC[a,b]\to\R$ be the map $f\mapsto f(r)$. Then $e_r$ is a surjective homomorphism whose kernel is $M_r$, so $\calC[a,b]/M_r\cong\R$. We will use this observation later in the course to show that $M_r$ is a \emph{maximal ideal}, meaning that there is no proper ideal larger than $M_r$.
\item\label{prod_quo_id} Let $A$ be a ring with ideal $I$, and $B$ a ring with ideal $J$. We claim that \[(A\times B)/(I\times J)\cong(A/I)\times(B/J).\] To prove this, let $f:A\times B\to (A/I)\times(B/J)$ be defined by $f(a,b)=(a+I,b+J)$. Then $f$ is a surjective homomorphism with kernel equal to $I\times J$. The claim now follows from the Isomorphism Theorem.
\end{examples}
\end{ex}

\subsection*{Prime subring and characteristic} We end this chapter with some applications of the Isomorphism Theorem to questions concerning the prime subring and the characteristic of a ring. Recall that if $R$ is a ring with prime subring $R_0$, we defined $\ch(R)=\#R_0$ if $R_0$ is finite, and otherwise $\ch(R)=0$.

\begin{framed}
\begin{thm}\label{prime_sub_thm}\;
\begin{enumerate}
\item[(a)] The prime subring of any ring $R$ is isomorphic to either $\Z$ or $\Z_n$ for some $n>0$. In the latter case, we have $\ch(R)=n$, and $n$ is the least positive integer such that $n\cdot 1_R=0_R$.
\item[(b)] The characteristic of an integral domain is either 0 or a prime number.
\end{enumerate}
\end{thm}
\end{framed}
\begin{proof}
Let $R$ be a ring with prime subring $R_0$. Let $f_R:\Z\to R$ be the unique homomorphism $\Z\to R$, which is given by $f_R(n)=n\cdot 1_R$ (see Examples \ref{homo_exs}). Recall from Exercise \ref{prime_sub}, Chapter \ref{subrings_chap} that $R_0=\{n\cdot 1_R:n\in\Z\}$, so that $R_0$ is equal to the image of $f_R$. Applying the Isomorphism Theorem to $f_R$, we conclude that $R_0\cong\Z/\ker(f_R)$. Recalling that the ideals of $\Z$ are the principal ideals $(n)$ with $n\ge 0$, we conclude that $R_0$ is isomorphic to $\Z$ if $\ker(f_R)=(0)$, and to $\Z_n$ if $\ker(f_R)=(n)$, $n>0$. In the latter case, it follows immediately that $n=\#R_0=\ch(R)$. Moreover, since $n$ is the least positive integer in $(n)=\ker(f_R)$, we have $n=\min\{m>0: m\cdot 1_R=0_R\}$. This proves (a).

We prove (b) by contradiction. Suppose that $R$ is a domain whose characteristic is a composite number $n$. Then $n=ab$ with $1<a,b<n$, so $0_R=f_R(n)=f_R(a)f_R(b)$. By the minimality of $n$ proved in part (a), $f_R(a)$ and $f_R(b)$ are nonzero. Hence, the equation $0_R=f_R(a)f_R(b)$ expresses $0_R$ as a product of nonzero elements of $R$, contradicting the assumption that $R$ is a domain.
\end{proof}

The final theorem of this chapter is a fundamental fact concerning rings with prime characteristic. The map referenced in the theorem is called the \col{Frobenius homomorphism}\index{Frobenius homomorphism} in honor of \href{https://en.wikipedia.org/wiki/Ferdinand_Georg_Frobenius}{Ferdinand Georg Frobenius} (1849-1917).

\begin{framed}
\begin{thm}\label{frob_thm}
If $R$ is a commutative ring with prime characteristic $p$, the map $x\mapsto x^p$ is a homomorphism $R\to R$. In particular, the identity $(a+b)^p=a^p+b^p$ holds for all $a,b\in R$.
\end{thm}
\end{framed}
\begin{proof}
The essential detail to prove is the identity $(a+b)^p=a^p+b^p$. By the binomial theorem, we have
\[(a+b)^p=\sum_{i=0}^pa^ib^{p-i}{p\choose i}.\]

For $0<i<p$, the binomial coefficient ${p\choose i}$ is an integer multiple of $p$, hence ${p\choose i}\cdot 1_R=0_R$ and the binomial identity simplifies to $(a+b)^p=a^p+b^p$.
\end{proof}

Colloquially, the identity $(a+b)^p=a^p+b^p$ is sometimes called the ``Freshman's Dream," as it would appear to oversimplify the process of calculating $(a+b)^p$.

\begin{ex}[Frobenius homomorphism]\label{frob_hom_exs}\;
\begin{examples}
\item With $R$ as in Theorem \ref{frob_thm}, let $\sigma$ be its Frobenius homomorphism. Since the composition of homomorphisms is a homomorphism, the map $\sigma^n:=\sigma\circ\sigma\circ\sigma\cdots\circ\sigma$ is a homomorphism for every $n>0$. Note that, by definition of $\sigma$, we have $\sigma^n(x)=x^{p^n}$ for all $x\in R$. We thus obtain the following generalization of the identity $(a+b)^p=a^p+b^p$:
\[(a+b)^{p^n}=a^{p^n}+b^{p^n}\quad\text{for all}\;a,b\in R,\, n\in\N.\]
\item Let $k$ be a finite field and $p=\ch(k)$. Let $\sigma:k\to k$ be the Frobenius homomorphism. Since $k$ is a domain, $\sigma$ has trivial kernel (if $x^p=0$, then $x=0$) and is therefore injective. Moreover, since $\sigma$ is a map from a finite set to itself, it must also be surjective. Hence, $\sigma$ is an automorphism of $k$. This element $\sigma\in\Aut(k)$ is naturally called the \emph{Frobenius automorphism} of $k$. For an arbitrary ring of characteristic $p$, the Frobenius map is not generally an automorphism; see Exercise \ref{frob_inj}.
\item (\href{https://en.wikipedia.org/wiki/Fermat\%27s_little_theorem}{Fermat's little theorem}) Let $p$ be prime. We will prove by induction on $n$ that $n^p\equiv n\pmod{p}$ for every positive integer $n$. The base case $n=1$ is clear. For the inductive step, we use the fact that $\F_p$ has characteristic $p$ and thus, by Theorem \ref{frob_thm},
\[(a+1)^p=a^p+1\;\;\text{for all}\;\; a\in\F_p.\]
Assuming $n^p\equiv n\pmod{p}$, it follows from the above identity that \[(n+1)^p\equiv n^p+1\equiv n+1\;\text{(mod $p$)},\] completing the induction.
\end{examples}
\end{ex}

\section*{Exercises}

\begin{note}\;
\begin{itemize}
\item Solutions to starred exercises (\sel) are included in Chapter~\ref{review2chap}.
\item A house (\house) marks exercises suggested as sample homework.
\end{itemize}
\end{note}

\medskip

\begin{exercises}
\item Show that the fields $\Q(\sqrt 2)$ and $\Q(\sqrt 3)$ are not isomorphic.
\item Show that the quotient ring $\Z[i]/(3)$ is not isomorphic to $\F_3\times\F_3$.
\item Prove that the inverse map of an isomorphism is an isomorphism.
\item\label{frac_iso} Let $A$, $B$ be isomorphic integral domains. Show that $\Frac(A)\cong\Frac(B)$.
\item Let $R$ be a ring. Suppose that $f:R\to X$ is a bijection, where $X$ is a set. We will define a ring structure on $X$ by using the ring structure of $R$; this is a process known as  \col{transporting}\index{transporting ring structure} the ring structure of $R$. Explicitly, operations in $X$ are defined as follows: given $x,y\in X$, we may write $x=f(a)$ and $y=f(b)$ for uniquely determined elements $a,b\in R$; we then define
\[x+y=f(a+b)\quad\text{and}\quad xy=f(ab).\]

With these operations, show that $X$ is a ring, and that moreover, $X$ is isomorphic to $R$. This construction is described in the following figure.
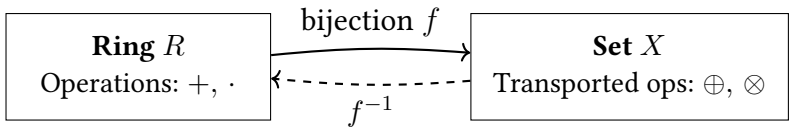
\begin{figure}[h!]
\centering
\begin{tikzpicture}[node distance=6.5cm, auto, font=\small]

  \node (R) [draw, rectangle, minimum width=3.5cm, minimum height=1.4cm, align=center] 
    {\textbf{Ring $R$}\\ Operations: $+,\, \cdot$};

  \node (X) [right of=R, draw, rectangle, minimum width=4.2cm, minimum height=1.4cm, align=center] 
    {\textbf{Set $X$} \\ Transported ops: $\oplus,\, \otimes$};

  \draw[->, thick] (R) to[bend left=5] node[above, font=\normalsize] {bijection $f$} (X);
  \draw[->, thick, dashed] (X) to[bend left=5] node[below] {$f^{-1}$} (R);

  \node at ($(R)!0.5!(X)+(0,-2.4)$) [align=left] {
    Define operations in $X$ by pullback: \\
    \quad $x \oplus y := f(f^{-1}(x) + f^{-1}(y))$ \\
    \quad $x \otimes y := f(f^{-1}(x) \cdot f^{-1}(y))$
  };
\end{tikzpicture}
\caption{Transporting a ring structure from a ring $R$ to a set $X$.}
\label{fig:transport_structure}
\end{figure}
\item Let $f:R\to S$ be a ring isomorphism. Show that $f$ induces a bijection between the set of units in $R$ and the set of units in $S$, and similarly for the sets of zero divisors.
\item Prove that isomorphic rings have the same characteristic.
\item Show that $\Z_8$ is not isomorphic to $\Z_2\times\Z_4$.
\item Show that $M_2(\R)$ has a subring isomorphic to $\C$.
\item Suppose that \( R \cong R_1 \times \cdots \times R_n \). Show that \( R \) contains nontrivial idempotents \( e_1, \ldots, e_n \) such that
\[
1 = e_1 + \cdots + e_n,\quad e_ie_j = 0\; \text{ for }\; i \ne j.
\]
\item Determine all the automorphisms of the quadratic field $\Q(\sqrt d)$.
\item Which of the following maps are homomorphisms $M_2(\R)\to\R$? 
\begin{enumerate}
\item The map $\begin{pmatrix}a&b\\c&d\end{pmatrix}\mapsto a$;
\item The trace map, $\begin{pmatrix}a&b\\c&d\end{pmatrix}\mapsto a+d$;
\item The determinant map, $\begin{pmatrix}a&b\\c&d\end{pmatrix}\mapsto ad-bc$.
\end{enumerate}
\item\label{ker_prop} Prove Theorem \ref{im_ker}.
\item Does there exist a commutative ring that is not an integral domain and yet has prime characteristic?
\item Let $A$ be a commutative ring. Show that for every $n\ge 1$, there is an embedding $A\to M_n(A)$ via scalar matrices, $a\mapsto aI_n$.
\item\label{complex_autos} The purpose of this exercise is to show that $\Aut(\R)$ consists of only the identity map. Let $\sigma\in\Aut(\R)$.
\begin{enumerate}
\item Show that $\sigma(x)=x$ for every $x\in\Q$.
\item Use the fact that positive real numbers are squares to show that $\sigma$ is a strictly increasing function.
\item Let $x\in\R$. Show that $\sigma(x)=x$ by considering a strictly increasing sequence $(s_n)$ of rational numbers converging to $x$, and strictly decreasing sequence $(t_n)$ of rational numbers converging to $x$.
\end{enumerate}
\item Let $f:R\to S$ be a ring homomorphism. Determine which of the following are necessarily true: (a) The preimage of a subring of $S$ under $f$ is a subring of $R$; (b) The image of an ideal is an ideal; (c) The preimage of an ideal is an ideal.
\item(\house) Let $k$ be a field and $R$ a nontrivial ring. Show that every homomorphism $k\to R$ is an embedding. In particular, every homomorphism between fields is an embedding.
\item Let $\Z_{(p)}$ be the subring of $\Q$ consisting of fractions whose denominator is not divisible by $p$. Find a surjective homomorphism $\Z_{(p)}\to\F_p$ whose kernel is $(p)$. Conclude that $\Z_{(p)}/(p)\cong\F_p$.
\item Let $f:R\to S$ be a ring homomorphism. The \col{graph}\index{graph of a homomorphism} of $f$ is the set $\Gamma_f=\{(r,f(r)):r\in R\}\subseteq R\times S$. Show that $\Gamma_f$ is a subring of $R\times S$ isomorphic to $R$. This construction provides a new example of a type of subring of a direct product.
\item(\house) Let $A$ be a ring and $I$ an ideal of $A$. Let $I[x]$ be the ideal of $A[x]$ consisting of polynomials all of whose coefficients belong to $I$. Show that \[A[x]/I[x]\cong (A/I)[x].\]
\item Recall that a PIR is a ring in which every ideal is principal. Show that every homomorphic image of a PIR is a PIR.
\item What is the characteristic of the direct product $\Z_n\times\Z_m$?
\item(\sel)\label{2nd_isothm} (Second Isomorphism Theorem) Let $R$ be a ring and $S$ a subring of $R$. Let $I$ be an ideal of $R$. Show that $S\cap I$ is an ideal of $S$, that $S+I$ is a subring of $R$, and that $I$ is an ideal of $S+I$. Moreover, \[S/(S\cap I)\cong (S+I)/I.\]
 Here,  $S+I=\{s + i : s \in S,\, i \in I\}$.
\item(\house)\label{third_iso_thm} (Third Isomorphism Theorem) Let $R$ be a ring and let $I \subseteq J \subseteq R$ be ideals. Show that \[(R/I)\big/\!(J/I)\;\cong\; R/J.\]
\item Find all homomorphisms $\Z\times\Z\to\Z$. For each homomorphism, describe its kernel and image.
\item Find a familiar ring that is isomorphic to the quotient $\Z[i]/(2+5i)$.
\item(\sel)\label{chap12_sel} Let $X$ be a set and $R(X)$ the ring of sets $(P(X),\Delta,\cap)$ in which addition and multiplication of subsets of $X$ are, respectively, the symmetric difference and intersection of sets. Show that $R(X)$ is isomorphic to the ring of functions $\F_2^X$.
\item(\house) Let $R$ be a commutative ring with prime characteristic. If $a\in R$ is nilpotent, show that $(1+a)^n=1$ for some integer $n\ge 1$.
\item Prove that in every commutative ring of characteristic 2, the set of idempotent elements is a subring.
\item Let $F$ be a field of characteristic 2 with more than two elements. Show that there exist $x,y\in F$ such that $(x+y)^3\ne x^3+y^3$.
\item Determine the characteristic of the ring $\Z[i]/(2+i)$.
\item Let $T$ be the subset of $M_2(\Z)$ consisting of matrices of the form $\begin{pmatrix} a&2b\\b&a\end{pmatrix}$. Show that $T$ is a subring of $M_2(\Z)$ isomorphic to $\Z[\sqrt 2]$.
\item For every positive integer $n$, give an example of a non-commutative ring of characteristic $n$.
\item\label{frob_inj} Let $p$ be a prime number.
\begin{enumerate}
\item Let $R$ be a commutative ring of characteristic $p$. Show that the Frobenius homomorphism $R\to R$ is injective if and only if $R$ has no nonzero nilpotent elements.
\item Let $R=\F_p[x]$, the ring of polynomials over $\F_p$. Show that the Frobenius map of $R$ is injective but not surjective.
\item A field with prime characteristic is called \col{perfect}\index{perfect field} if its Frobenius map is surjective (and therefore an automorphism). Show that every finite field is perfect, and that the field $\F_p(x)$ of rational functions over $\F_p$ is not.
\end{enumerate}
\end{exercises}

\chapter{Universal Properties}\label{univ_chap}

A \col{universal property}\index{universal property} of a ring is a formal description of how the ring interacts with other rings via homomorphisms. Importantly, a universal property determines the ring uniquely up to isomorphism. This chapter is devoted to studying the universal properties satisfied by three fundamental constructions: the field of fractions of a domain, the quotient of a ring by an ideal, and the univariate polynomial ring over a commutative ring. Universal properties arise throughout mathematics; their systematic treatment belongs to \href{https://en.wikipedia.org/wiki/Category_theory}{Category Theory}.

Throughout this chapter, we write $\Hom(A,B)$ for the set of all ring homomorphisms from a ring $A$ to a ring $B$. Thus,
\[\Hom(A,B):=\{\text{Homomorphisms}\;A\to B\}.\]

\subsection*{Fields of fractions}

As a first example of a universal property, consider the field of fractions of a domain. We have seen in Chapter \ref{frac_chap} that every domain $D$ can be embedded in its field of fractions, $F=\Frac(D)$. The following theorem shows that there is a precise sense in which every field that contains $D$ must also contain $F$, so that $F$ is the smallest field containing $D$.

\newpage

\begin{framed}
\begin{thm}[Universal property of the field of fractions]\label{frac_univ}
 Let $D$ be an integral domain and $K$ a field. For every embedding $f:D\to K$ there exists a unique embedding $\tilde f:\Frac(D)\to K$ such that $f(x)=\tilde f(x/1)$ for all $x\in D$. Explicitly, $\tilde f$ is given by \[\tilde f(a/b)=f(a)/f(b)\] for all $a,b\in D$ with $b\ne 0$.
\end{thm}
\end{framed}

\begin{proof}
Let $F=\Frac(D)$ and define $\tilde f:F\to K$ as in the theorem. This map is well defined, as if $a/b=c/d$ in $F$, then $ad=bc$ in $D$; applying $f$ to this equation we obtain $f(a)f(d)=f(b)f(c)$. Since $b$ and $d$ are nonzero and $f$ is an embedding, $f(b)$ and $f(d)$ are nonzero elements of $K$, hence we may divide by these elements to obtain $f(a)/f(b)=f(c)/f(d)$ and therefore $\tilde f(a/b)=\tilde f(c/d)$.

We leave it to the reader to verify that $\tilde f$ is a homomorphism and injective. Now, for $x\in D$ we have $\tilde f(x/1)=f(x)/f(1)=f(x)$. This proves the existence of $\tilde f$ with the required properties. Finally, if $h:F\to K$ is any embedding such that $f(x)=h(x/1)$ for all $x\in D$, then the fact that $h$ is a homomorphism implies that $h(a/b)=h(a/1)h(1/b)=f(a)/f(b)=\tilde f(a/b)$ for all $a,b\in D$ with $b\ne 0$. Hence $h=\tilde f$, proving uniqueness of $\tilde f$.
\end{proof}

\begin{note}
If we identify $D$ with the subring $\{x/1:x\in D\}\subseteq\Frac(D)$, then Theorem \ref{frac_univ} states that every embedding of $D$ into a field $K$ extends uniquely to an embedding of $\Frac(D)$ into $K$.
\end{note}

\begin{1ex}[The field of fractions of a field] For any field $k$ we have $k\cong\Frac(k)$, as is to be expected given that $k$ is the smallest field containing $k$. For a formal proof of this claim, apply the universal property to the identity map $f:k\to k$. This embedding extends uniquely to an embedding $\tilde f:\Frac(k)\to k$ satisfying $\tilde f(x/1)=f(x)=x$ for all $x\in k$. It follows that $\tilde f$ is surjective and therefore an isomorphism.
\end{1ex}

As an application of Theorem \ref{frac_univ}, we derive the following result describing the prime subfield of any field.

\begin{framed}
\begin{thm} Let $k$ be a field.
\begin{enumerate}
\item[(a)] If $k$ has prime subring $R$ and prime subfield $F$, then $F\cong\Frac(R)$.
\item[(b)] The prime subfield of $k$ is isomorphic to either $\Q$ or $\F_p$ for some prime $p$.
\end{enumerate}
\end{thm}
\end{framed}

\begin{proof}
To prove (a), we begin by noting that $R$ is a domain since it is a subring of a field, so $R$ has a well-defined field of fractions. Now let $i:R\to k$ be the inclusion map $x\mapsto x$, clearly an embedding. By the universal property of $\Frac(R)$, there is an extension $h:\Frac(R)\to k$ such that $h(a/1)=a$ for all $a\in R$.  Letting $E$ be the image of $h$, we claim that $E=F$, from which the desired result follows. We have $E\supseteq F$ since $E$ is a subfield of $k$ and therefore contains the prime subfield. Moreover, we have $R\subseteq F$ since $F$ is a subring of $k$ and $R$ is the prime subring of $k$. To show that $E\subseteq F$, let $x\in E$ and write $x=h(a/b)$, where $a,b\in R$ and $b\ne 0$. Then $x=h((a/1)(1/b))=h(a/1)h(b/1)^{-1}=ab^{-1}\in F$ since $F$ is a field and $a,b\in F$. This proves the claim that $E=F$ and thus $\Frac(R)\cong F$.

Now (b) follows from (a) together with Theorem \ref{prime_sub_thm}: the prime subring $R$ is isomorphic to either $\Z$ or $\Z_p=\F_p$ for some prime $p$, and therefore $\Frac(R)$ is isomorphic to $\Frac(\Z)\cong\Q$ or $\Frac(\F_p)\cong\F_p$. Note that we have used the fact that isomorphic domains have isomorphic fields of fractions, which is the content of Exercise \ref{frac_iso}, Chapter \ref{homo_chap}. Alternatively, this fact can be derived from the universal property---see Exercise \ref{frac_univ_iso}.
\end{proof}

\subsection*{Quotient rings} Our next theorem is the main tool needed to construct homomorphisms whose domain is a quotient ring, such as $\Z_n$.

\begin{framed}
\begin{thm}[Universal property of quotient rings]\label{quo_univ}
Let $R$ and $S$ be rings, and $I$ an ideal of $R$. If $f:R\to S$ is a homomorphism such that $I\subseteq\ker(f)$, then there exists a unique homomorphism $\bar f:R/I\to S$ such that $\bar f([x])=f(x)$ for all $x\in R$.
\end{thm}
\vspace{-2mm}
\end{framed}

\begin{proof}
The main claim to prove is that $\bar f$ is well defined. The fact that $\bar f$ is a homomorphism, and the uniqueness of $\bar f$, follow readily from the definitions. Assuming $[a]=[b]$ in $R/I$, we have $a-b\in I$ and therefore $a-b\in\ker(f)$. From here we deduce that $f(a)=f(b)$, hence $\bar f([a])=\bar f([b])$, showing that $\bar f$ is well defined.
\end{proof}

\begin{note} In the context of Theorem \ref{quo_univ}, let $\pi:R\to R/I$ be reduction modulo $I$. Then $f=\bar f\circ\pi$, so these maps form the \href{https://en.wikipedia.org/wiki/Commutative_diagram}{commutative diagram} shown below.

\[
\begin{tikzcd}
R \arrow[r, "f"] \arrow[d, "\pi"'] & S \\
R/I \arrow[ur, dashed, "\bar f"'] &
\end{tikzcd}
\]

The map $f$ is said to \emph{descend to the quotient}, as its domain can be changed from $R$ to $R/I$, a vertical descent in the above diagram.
\end{note}

The following theorem is equivalent to the universal property of quotient rings, and may itself be referred to as the universal property. It provides a concise formulation of the property as a correspondence between two sets of homomorphisms: specifying a homomorphism $R/I\to S$ is equivalent to specifying a homomorphism $R\to S$ that maps $I$ to 0.

\begin{framed}
\begin{thm}[Universal property of quotient rings]\label{quo_ring_bij}
Let $R$ and $S$ be rings, and $I$ an ideal of $R$. The map $f\mapsto\bar f$ (in the notation of Theorem \ref{quo_univ}) is a bijection
\[\{f\in\Hom(R,S):f(I)=(0)\}\;\xlongrightarrow{\sim}\;
\Hom(R/I,S).\]
\end{thm}
\end{framed}

\begin{proof}
Let $\Phi$ be the map $f\mapsto\bar f$. We claim that the map $\Psi(g)=g\circ\pi$ is the inverse of $\Phi$. Let $X$ be the domain of $\Phi$ and $Y=\Hom(R/I,S)$. Note, first of all, that $\Psi$ is a well-defined map from $Y$ to $X$. Indeed, if $g:R/I\to S$ is a homomorphism, then $g\circ\pi$ is a homomorphism $R\to S$, and for every $x\in I$ we have $g(\pi(x))=g([0])=0_S$, so $g\circ\pi$ maps $I$ to 0. Now, for $f\in X$, the defining property of $\bar f$ implies that $\Psi(\Phi(f))=\bar f\circ\pi=f$. Finally, for $g\in Y$ we claim that $\Phi(\Psi(g))=g$. By definition of $\Phi$ and $\Psi$, the following holds for every $x\in R$: the value of $\Phi(\Psi(g))$ at $[x]$ is $\overline{g\circ\pi}[x]=g\circ\pi(x)=g([x])$, proving the claim. This shows that $\Phi$ and $\Psi$ are inverses, and in particular, that $\Phi$ is a bijection.
\end{proof}

\begin{1ex}[Maps between quotient rings] Letting $m$ and $n$ be positive integers, we will determine all homomorphisms $\Z_m\to\Z_n$. By the universal property of quotient rings, the set $\Hom(\Z_m,\Z_n)$ is in bijection with the set of homomorphisms $f:\Z\to\Z_n$ such that $f(m)=[0]$. Now, there is exactly one homomorphism $f:\Z\to\Z_n$, namely the map $s\mapsto s\cdot [1]$ (see \ref{unique_homZ} in Examples \ref{homo_exs}). Hence $\Hom(\Z_m,\Z_n)$ is empty if $f(m)\ne [0]$, and otherwise is a singleton set. Since $f(m)=m\cdot[1]=[m]$, we have $f(m)=[0]$ in $\Z_n$ precisely when $n$ divides $m$. Therefore, \[
\Hom(\Z_m,\Z_n) =
\begin{cases}
\{\,x+(m)\mapsto x+(n)\,\}, & n\mid m,\\[4pt]
\emptyset, & n\nmid m.
\end{cases}
\] 

Concretely, when $n\mid m$, the map $\Z_m\to\Z_n$ is simply reduction modulo $n$. For instance, the image of $[8]$ under the homomorphism $\Z_{12}\to\Z_6$ is $[2]$, since $8\equiv 2\pmod 6$.
\end{1ex}

\subsection*{Polynomial rings} The final theorem of this chapter is the essential tool for constructing homomorphisms whose domain is a polynomial ring. Put simply, the theorem states that specifying a homomorphism $A[x]\to B$ is equivalent to specifying a homomorphism $A\to B$ as well as an element of $B$, namely the image of $x$.

\begin{framed}
\begin{thm}[Universal property of polynomial rings]\label{poly_univ}
Let $A$ and $B$ be commutative rings. For every $f\in\Hom(A,B)$ and every $b\in B$ there exists a unique  $\tilde f\in\Hom(A[x],B)$ extending $f$ and mapping $x$ to $b$. Explicitly, $\tilde f$ is given by
\[\tilde f\left(\sum a_ix^i\right)=\sum f(a_i)b^i.\]
\end{thm}
\end{framed}

\begin{proof}
Let $f:A\to B$ be a homomorphism, let $b\in B$, and define $\tilde f$ as in the theorem. Proving that $\tilde f$ is a homomorphism is Exercise \ref{poly_hom_thm}. The definition of $\tilde f$ implies that $\tilde f(x)=f(1)\cdot b=b$ and that $\tilde f$ extends $f$, as $\tilde f(a)=f(a)\cdot b^0=f(a)$ for every $a\in A$. Hence, $\tilde f$ is a homomorphism $A[x]\to B$ extending $f$ and mapping $x$ to $b$. Uniqueness of $\tilde f$ follows immediately from the defining properties of a homomorphism; any extension of $f$ mapping $x$ to $b$ will necessarily map a polynomial $\sum a_ix^i$ to $\sum f(a_i)b^i$.
\end{proof}

The following is an equivalent formulation of the universal property proved above. The proof is left to the reader as Exercise \ref{prove_poly_univ2}.

\begin{framed}
\begin{thm}[Universal property of polynomial rings]\label{poly_univ2}
Let $A$ and $B$ be commutative rings. The map $h\mapsto (h|_{A},h(x))$ is a bijection
\[\Hom(A[x],B)\;\xlongrightarrow{\sim}\;\Hom(A,B)\times B.\]
Here, $h|_A$ denotes the restriction of $h$ to $A$.
\end{thm}
\end{framed}

\begin{ex}[Homomorphisms on polynomial rings]\;
\begin{examples}
\item The universal property implies that every evaluation map $A[x]\to A$ is a homomorphism. For $a\in A$, let $e_a:A[x]\to A$ be defined by $e_a(\sum a_ix^i)=\sum a_ia^i$. Applying the universal property to the identity map $f:A\to A$ and the element $a\in A$, we deduce that there is a unique homomorphism $\tilde f:A[x]\to A$ given by $\sum a_ix^i\mapsto\sum a_ia^i$. Note that $\tilde f$ is precisely the map $e_a$, so we conclude that $e_a$ is a homomorphism. This argument is merely a formal verification of an idea that is intuitively very simple: polynomial expressions such as $p(x)^3+5q(x)$ evaluate naturally to $p(a)^3+5q(a)$.

\item Homomorphisms between coefficient rings extend to homomorphisms between polynomial rings. Given any homomorphism $f:A\to B$ between commutative rings, there is a corresponding homomorphism $A[x]\to B[x]$ given by $\sum a_ix^i\mapsto\sum f(a_i)x^i$. To derive this fact, we apply the universal property to the composition of the homomorphisms $A\xrightarrow{f} B\hookrightarrow B[x]$ and the element $x\in B[x]$. As a specific instance of this construction, the reduction map $\Z\to\Z_6$ gives rise to a homomorphism $\Z[x]\to\Z_6[x]$; the image of $8x^{10}-2x+20$ under this map is $2x^{10}+4x+2$.
\end{examples}
\end{ex}

\section*{Exercises}

\begin{note}\;
\begin{itemize}
\item Solutions to starred exercises (\sel) are included in Chapter~\ref{review2chap}.
\item A house (\house) marks exercises suggested as sample homework.
\end{itemize}
\end{note}

\medskip

\begin{exercises}
\item Let $k$ be a field. Show that $\ch(k)=0$ if and only if $k$ has a subfield isomorphic to $\Q$, and show that $\ch(k)=p$, where $p$ is a prime number, if and only if $k$ has a subfield isomorphic to $\F_p$.
\item\label{frac_univ_iso} Use the universal property of the field of fractions to show that if $A$ and $B$ are isomorphic domains, then $\Frac(A)\cong\Frac(B)$.
\item Use the universal property of polynomial rings to show that if the commutative rings $R$ and $S$ are isomorphic, then the polynomial rings $R[x]$ and $S[x]$ are isomorphic.
\item Let $d\ne 1$ be a squarefree integer. Show that $\Q(\sqrt d)$ is isomorphic to the field of fractions of $\Z[\sqrt d]$.
\item(\house) Use evaluation maps to prove the following: if $A$ is a commutative ring and $a\in A$, then \[A[x]/(x-a)\cong A.\]
\item(\house) Prove the following universal property for direct products of rings: for rings $R$, $S$, and $T$, there is a bijection \[\Hom(T,R\times S)\,\xlongrightarrow{\sim}\,\Hom(T,R)\times \Hom(T,S).\]
\item(\sel)\label{chap13_sel} Let $R$ be a ring and $n\ge 1$. Show that there exists a homomorphism $\Z_n\to R$ if and only if $\ch(R)$ divides $n$.
\item List the image of each element of $\Z_{10}$ under the homomorphism $\Z_{10}\to\Z_5$.
\item\label{prove_poly_univ2} Prove Theorem \ref{poly_univ2}.

\item\label{poly_hom_thm} Show that the map $\tilde f$ from Theorem \ref{poly_univ}, needed for proving the universal property of polynomial rings, is a homomorphism.
\end{exercises}

\chapter{Prime and Maximal Ideals}\label{prime_chap}

Prime and maximal ideals serve as the basic building blocks of ring theory, much as prime numbers do in the integers. They capture, in algebraic form, the idea of “irreducible pieces” from which more complicated structures can be assembled. In particular, we will prove in this chapter that the maximal ideals of $\Z$ correspond exactly to the prime numbers.

Every maximal ideal of a ring turns out to be a prime ideal, so the collection of prime ideals always includes the maximal ones. In $\Z$ the two notions essentially coincide, but in general rings there are many more prime ideals than maximal ones, and studying this larger set reveals much about the ring’s internal structure.

The set of all prime ideals of a ring $R$, called the \col{prime spectrum}\index{prime spectrum} and denoted by \(\Spec(R)\), plays a central role far beyond elementary ring theory. It serves as a unifying object across algebraic geometry, number theory, order theory, and even mathematical logic. Laying the groundwork for understanding \(\Spec(R)\) is therefore one of the main goals of this chapter. See Exercise \ref{spec_def} for some of its basic properties.

\section*{Basic concepts}

Let $R$ be a commutative ring. A \col{maximal ideal}\index{maximal ideal} of $R$ is a proper ideal that is not contained in any other proper ideal. Thus, an ideal $M$ of $R$ is maximal if and only if $M\ne R$ and the only ideals $I$ with $M\subseteq I\subseteq R$ are $M$ and $R$.  A \col{prime ideal}\index{prime ideal} of $R$ is a proper ideal $P$ with the property that for all $a,b\in R$, $ab\in P$ implies $a\in P$ or $b\in P$. Thus (by induction), whenever a finite product $a_1a_2\cdots a_n$ lies in $P$, at least one of the factors $a_i$ must belong to $P$.

\begin{ex}[Prime and maximal ideals]\label{pr_max_id_exs}\;
\begin{examples}
\item In $\Z$, every principal ideal $(p)$ generated by a prime number $p$ is a prime ideal, because if $p$ divides a product $ab$, then $p$ must divide one of the factors. (The latter claim is known as \href{https://en.wikipedia.org/wiki/Euclid\%27s_lemma}{Euclid's Lemma}, and will be clear to the reader willing to assume the fundamental theorem of arithmetic. Otherwise, it follows from the \emph{B\'ezout Identity}, Theorem \ref{bezout_id}.)
\item For every commutative ring $A$, $A$ is an integral domain if and only if $(0)$ is a prime ideal of $A$. This follows immediately from the definition of a domain as a nonzero commutative ring in which, for all $a,b\in A$, $ab=0\implies a=0$ or $b=0$.
\item In any field, $(0)$ is a maximal ideal because the only other ideal is $(1)$.
\item Let $k$ be a field. In the power series ring $k\pow x$, the unique maximal ideal is $(x)$, as we have seen in Theorem \ref{kxideals} that the ideals of $k\pow x$ lie in a chain $(x)\supset(x^2)\supset(x^3)\supset\cdots\supset(0)$. The prime ideals of this ring are $(0)$ and $(x)$. Indeed, $(0)$ is prime because $k\pow x$ is a domain, and $(x)$ is readily shown to be prime. Moreover, an ideal $I=(x^n)$ with $n\ge 2$ is not prime because $x\cdot x^{n-1}\in I$ but $x\notin I$ and $x^{n-1}\notin I$. 
\end{examples}
\end{ex}

The following theorem provides a criterion for determining whether an ideal $I$ is prime or maximal by analyzing the quotient ring $R/I$.
\begin{framed}
\begin{thm}\label{quo_prime}
For every commutative ring $R$ and ideal $I$ of $R$:
\begin{enumerate}
\item[(a)] $I$ is prime if and only if the quotient ring $R/I$ is a domain;
\item[(b)] $I$ is maximal if and only if $R/I$ is a field. 
\end{enumerate}
In particular, every maximal ideal is prime.
\end{thm}
\end{framed}

\begin{proof}
We prove one direction of (a) and (b); the converse directions are proved similarly and left to the reader. If $R/I$ is a domain, then $I$ must be a prime ideal: assuming $ab\in I$, in $R/I$ we have $[a][b]=0$, hence $[a]=0$ or $[b]=0$, and thus $a\in I$ or $b\in I$. Now, if $R/I$ is a field, we will show that $I$ is maximal. Suppose that $J$ is an ideal properly containing $I$, and let $x\in J\setminus I$. Then $[x]\ne 0$ in $R/I$, so there exists $[y]\in R/I$ such that $[x][y]=1$ and hence $xy=1+i$ for some $i\in I$. It follows that $1=xy-i\in J$ since both $x$ and $i$ lie in $J$, and thus $J=R$. This proves that $I$ is maximal.
\end{proof}

\begin{ex}[Criteria for primality and maximality]\;
\begin{examples}
\item In $\Z$, every prime ideal $(p)$ is maximal, since $\Z/(p)=\F_p$ is a field. Thus, the prime ideals of $\Z$ are $(0)$ and $(p)$ for primes $p$, and the maximal ones are the ideals $(p)$.
\item In $\Z[x]$, the ideal $(x)$ is prime but not maximal, as the quotient \[\Z[x]/(x)\cong\Z\] is a domain but not a field. The non-maximality of $(x)$ can also be proved directly by noting that the ideal $(2,x)$ is a proper ideal strictly containing $(x)$.
\item Recall from \ref{cab_max_idls} in Examples \ref{iso_thm_exs} that the ring $\calC[a,b]$ contains an ideal $M_r$ for every $r\in[a,b]$, namely the kernel of evaluation at $r$, and that \[\calC[a,b]/M_r\cong\R\] as a consequence of the Isomorphism Theorem. Since $\R$ is a field, Theorem \ref{quo_prime} implies that $M_r$ is a maximal ideal. This shows that $\calC[a,b]$ is a ring with uncountably many maximal ideals.
\item In the polynomial ring $k[x]$, where $k$ is a field, every principal ideal $(x-a)$ is maximal. Indeed, the Isomorphism Theorem applied to the evaluation map $e_a:k[x]\to k$ yields an isomorphism \[k[x]/(x-a)\cong k.\] Since $k$ is a field, Theorem \ref{quo_prime} implies that $(x-a)$ is a maximal ideal. Taking $k=\R$, we obtain a second example of a ring with uncountably many maximal ideals, namely $\R[x]$.
\item Let $k$ be a field and $a_1,a_2,\ldots, a_n\in k$. Generalizing the previous example, we claim that the ideal $M=(x_1-a_1,x_2-a_2,\ldots, x_n-a_n)$ is a maximal ideal of the polynomial ring $k[x_1,x_2,\ldots, x_n]$. This observation lies at the foundation of modern algebraic geometry, as it establishes a clear link between points $(a_1,\ldots, a_n)\in k^n$ and maximal ideals in the ring $k[x_1,\ldots,x_n]$; put loosely, it provides an algebraic way of describing geometry. See Exercise \ref{alggeo} for a guided proof of the claim. 
\end{examples}
\end{ex}

\section*{Existence of maximal ideals}

Having defined the concepts of maximal/prime ideals, we now consider the question of whether such ideals always exist in any given commutative ring. The answer is yes, with one trivial exception: the zero ring has no proper ideal and therefore no maximal or prime ideal.

\begin{framed}
\begin{thm}[Using the axiom of choice]\label{exists_max}
Every nonzero commutative ring has at least one maximal (and therefore prime) ideal; in fact, every proper ideal is contained in a maximal ideal.
\end{thm}
\end{framed}
\begin{proof}[Proof sketch] Let $R$ be a nonzero commutative ring and let $I$ be a proper ideal of $R$.  The set of proper ideals containing $I$ forms a poset under inclusion. Every chain in this poset has an upper bound (namely, the union of the chain), so by \href{https://en.wikipedia.org/wiki/Zorn\%27s_lemma}{Zorn's Lemma}, a maximal element exists.
\end{proof}

In many concrete examples of rings, there is no need to rely on the axiom of choice to prove the existence of maximal/prime ideals or even to identify all such ideals of a ring. For instance, we have already determined (in \ref{pr_max_id_exs}) all the prime and maximal ideals of the rings $\Z$ and $k\pow x$, where $k$ is a field. The toolkit we are developing will enable us to do the same for other familiar rings, such as $\Z_n$ and $\Q[x]$.

There is no restriction on the number of maximal ideals that a ring could have; for example, this number could be uncountable, as in the case of $\R[x]$, or countably infinite, as in the case of $\Z$, or it could be exactly one, as in the case of $k\pow x$. See Examples \ref{prod_id_exs} for a construction of rings having $n$ maximal ideals, for any positive integer $n$.
 
A \col{local ring}\index{local ring} is a commutative ring having exactly one maximal ideal. Rings of this type have the key property that their unit groups are easily determined: if $R$ is a local ring whose unique maximal ideal is $M$, then \[R^{\times}=R\setminus M.\]
This property greatly simplifies questions about invertibility and ideal containment. To prove the identity, note first that no unit can belong to $M$ because $M$ is a proper ideal, so $R^{\times}\subseteq R\setminus M$. For the opposite inclusion, apply Theorem \ref{exists_max}: every non-unit generates a proper ideal and is thus contained in some maximal ideal, which must be $M$.

\begin{ex}[Local rings]\;
\begin{examples}\label{local_ex}
\item If $k$ is a field, then both $k$ and $k\pow x$ are local rings. The maximal ideal of $k$ is $(0)$, and that of $k\pow x$ is $(x)$.
\item The polynomial ring $k[x]$ is not local, as it has at least two maximal ideals: $(x)$ and $(x-1)$. Similarly, $\Z$ is not a local ring.
\item For every $n\ge 1$, the ring $\Z_{p^n}$ is local. By Theorem \ref{quo_sub_id}, the ideals of $\Z/(p^n)$ have the form $(d)/(p^n)$, where $d$ is a divisor of $p^n$, i.e., $d=p^k$ for some $0\le k\le n$. Hence, the ideals of $\Z_{p^n}$ form a chain \[(p)/(p^n)\supset(p^2)/(p^n)\supset\cdots\supset(p^n)/(p^n).\] Clearly, there is a single maximal ideal, namely, $(p)/(p^n)$.
\end{examples}
\end{ex}

The following theorem provides a criterion for checking whether a given ring is local. See Exercise \ref{Zplocal} for a proof of the theorem as well as an additional example of a local ring.

\begin{framed}
\begin{thm}\label{local_crit}
Let $A$ be a commutative ring.
\begin{enumerate}
\item[(a)] Suppose that $I$ is an ideal of $A$ such that every element of $A\setminus I$ is a unit. Then $A$ is a local ring and $I$ is its unique maximal ideal.
\item[(b)] $A$ is a local ring if and only if the set $A\setminus A^{\times}$ is an ideal of $A$.
\end{enumerate}
\end{thm}
\end{framed}

Using the above theorem, the conclusions of the examples in \ref{local_ex} can be easily re-derived. For instance, a field is a local ring because $I=(0)$ is an ideal meeting the requirements in part (a) of the theorem.

\section*{Primes in a direct product or quotient ring}

We now turn to the question of how to determine the prime ideals in two fundamental constructions: a direct product of rings and a quotient of a ring by an ideal. The conclusions of our analysis will allow us, for instance, to identify the prime and maximal ideals of $\Z_n$ and of $\Z^n$ for every $n\ge 1$.

Recall that the ideals of a quotient ring $R/I$ can be uniquely represented in the form $J/I$, where $J$ is an ideal of $R$ containing $I$ (Theorem \ref{quo_sub_id}). The following theorem shows that this representation preserves the property of being prime or maximal, thus providing an explicit description of the prime and maximal ideals of $R/I$ in terms of the corresponding ideals of $R$.

\begin{framed}
\begin{thm}\label{quo_max}
Let $R$ be a commutative ring and $I$ an ideal of $R$. The maximal ideals of $R/I$ are the ideals of the form $M/I$, where $M$ is a maximal ideal of $R$ containing $I$. The prime ideals of $R/I$ are of the form $P/I$, where $P$ is a prime ideal of $R$ containing $I$.
\end{thm}
\end{framed}
\begin{proof}
Let $J$ be an ideal containing $I$. The reduction homomorphism $R\to R/J$ maps $I$ to 0, hence descends to a homomorphism $R/I\to R/J$ given by $a+I\mapsto a+J$ (by the universal property of quotient rings, \ref{quo_univ}). The kernel of this homomorphism is precisely $J/I$, hence the Isomorphism Theorem yields 
\[(R/I)/(J/I)\cong R/J.\] 
Applying Theorem \ref{quo_prime}, we have $J/I$ is maximal $\iff$ $(R/I)/(J/I)$ is a field $\iff$ $R/J$ is a field $\iff$ $J$ is maximal. Similarly, $J/I$ is prime $\iff$ $R/J$ is a domain $\iff$ $J$ is prime. 
\end{proof}
\newpage
\begin{ex}[Prime ideals in a quotient ring]\;
\begin{examples}
\item Theorem \ref{quo_max} enables us to determine the prime and maximal ideals of $\Z_n$ for any $n\ge 1$. Recall that ideals of $\Z_n$ have the form $(d)/(n)$, where $d$ is a divisor of $n$. The maximal ideals are those having $(d)$ a maximal ideal of $\Z$, hence $d$ being prime. Thus, the maximal ideals of $\Z_n$ are the ideals $(p)/(n)$, where $p$ is a prime dividing $n$. These ideals are also precisely the prime ideals of $\Z_n$, as the only non-maximal prime ideal of $\Z$ is $(0)$, which does not contain $(n)$.
\item Using the previous example we can re-derive the fact that $\Z_{p^n}$ is a local ring. Its maximal ideals have the form $(q)/(p^n)$, where $q$ is prime and divides $p^n$. It follows that $q=p$, showing that $(p)/(p^n)$ is the unique maximal ideal of $\Z_{p^n}$.
\end{examples}
\end{ex}

We conclude this chapter by describing the prime and maximal ideals of a direct product of finitely many rings. For this purpose, we introduce the following notation. Let $R_1,R_2,\ldots,R_n$ be rings and 
\[\pi_k:R_1\times\cdots\times R_n\to R_k\] the projection map for each index $1\le k\le n$. For any ideal $I\subseteq R_k$, set
\[I^{(k)}:=\pi_k^{-1}(I)=R_1\times\cdots\times I\times\cdots\times R_n.\]
For instance, in a product $A\times B$ we have $I^{(1)}=I\times B$ and $I^{(2)}=A\times I$.

\begin{framed}
\begin{thm}\label{prod_max}
Let $R_1,R_2,\ldots, R_n$ be commutative rings.
\begin{enumerate}
\item[(a)] The prime ideals of $\prod_{i=1}^n R_i$ are the ideals of the form $P^{(k)}$, where $1\le k\le n$ and $P$ is a prime ideal of $R_k$.
\item[(b)] The maximal ideals of $\prod_{i=1}^n R_i$ are the ideals of the form $M^{(k)}$, where $1\le k\le n$ and $M$ is a maximal ideal of $R_k$.
\end{enumerate}
\end{thm}
\end{framed}
\begin{proof}
Let $R=R_1\times R_2\times\cdots\times R_n$ and let $I=I_1\times I_2\times\cdots\times I_n$ be an ideal of $R$. (Recall from \ref{prod_ideals} that every ideal of $R$ has this form.) By \ref{prod_quo_id} in Examples \ref{iso_thm_exs}, the quotient $R/I$ decomposes as
\[R/I\cong(R_1/I_1)\times(R_2/I_2)\times\cdots\times(R_n/I_n).\]

A finite direct product of rings is an integral domain if and only if exactly one factor is a domain and all the others are zero rings. Hence $R/I$ is a domain precisely when there exists an index $k$ such that:
$R_k/I_k$ is a domain and $R_j/I_j$ is a zero ring for all $j\ne k$. Equivalently, $I_j = R_j$ for $j\ne k$ and $I_k$ is a prime ideal of $R_k$. Thus $I = P^{(k)}$ for some prime ideal $P \subseteq R_k$. This proves (a). The proof of (b) is entirely analogous, noting that a finite direct product of rings is a field if and only if exactly one factor is a field and all the others are zero rings.
\end{proof}

\begin{ex}[Prime ideals in a direct product]\label{prod_id_exs}\;
\begin{examples}
\item The maximal ideals of $\Z^2=\Z\times\Z$ are the ideals $(p)\times\Z$ and $\Z\times(p)$ for all primes $p$. The prime ideals of $\Z^2$ are all the maximal ideals as well as $(0)\times\Z$ and $\Z\times(0)$.
\item We can construct a ring having exactly $n$ maximal ideals for any positive integer $n$. Fixing a prime $p$, let $R=\F_p^n$. The maximal ideals of $R$ are the ideals $(0)^{(k)}$ for $1\le k\le n$, as $(0)$ is the unique maximal ideal of $\F_p$.
\item In an infinite direct product of rings, prime and maximal ideals need not be of the simple form $P^{(k)}$. Describing these more general ideals requires additional set-theoretic machinery; readers interested in this topic may consult resources on \emph{ultrafilters} and \emph{ultraproducts}.
\end{examples}
\end{ex}

\section*{Exercises}

\begin{note}\;
\begin{itemize}
\item Solutions to starred exercises (\sel) are included in Chapter~\ref{review2chap}.
\item A house (\house) marks exercises suggested as sample homework.
\end{itemize}
\end{note}

\medskip

\begin{exercises}
\item Show that in $\Z[i]$, the ideal $(3)$ is maximal, while $(5)$ is not.
\item Use Theorem \ref{prod_max} to describe all the prime ideals of $\Z^3$ and all the maximal ideals of $\Z\times\Q$.
\item Draw a diagram of the ideal lattice of $\Z_{42}$. Label all prime and maximal ideals.
\item Prove, directly from the definitions of prime and maximal ideals, that every maximal ideal is prime.
\item Show by example that the intersection of two prime ideals need not be prime.
\item Let $f:A\to B$ be a ring homomorphism. Show that the preimage of a prime ideal is prime; in other words, if $P$ is a prime ideal of $B$, then $f^{-1}(P)=\{a\in A:f(a)\in P\}$ is a prime ideal of $A$. Give an example to show that the same is not generally true for maximal ideals.
\item In the ring $\Z^\N$, let $S$ be the set of sequences with only finitely many nonzero terms. Show that $S$ is an ideal. Is it prime? Is it maximal?
\item(\house) Prove the following:
\begin{enumerate}
\item In a finite commutative ring, every prime ideal is maximal.
\item In a Boolean ring, every prime ideal is maximal.
\item In a PID, every nonzero prime ideal is maximal.
\end{enumerate}
\item Let $F$ be a field and $R=F[x_1,\ldots,x_n]$. Use Theorem \ref{quo_prime} to show that in the following chain of ideals of $R$, every ideal is prime, but only the largest is a maximal ideal:\[(x_1)\subset(x_1,x_2)\subset(x_1,x_2,x_3)\subset\cdots\subset (x_1,x_2,\ldots,x_n).\]
\item Let $A$ be a commutative ring. Show that the union of all the maximal ideals of $A$ is equal to the complement of the unit group, $A\setminus A^{\times}$. (In the case where $A$ is a local ring with maximal ideal $M$, this proves that $A^{\times}=A\setminus M$.)
\item In $\calC[0,1]$, let $I$ be the set of all functions $f$ such that \[f(1/3)=f(1/2)=0.\] Show that $I$ is a proper ideal that is neither prime nor maximal.
\item Classify, up to isomorphism, all local rings of cardinality $\le 8$.
\item(\sel)\label{chap14_sel} Let $P$ be a prime ideal of a commutative ring $A$. Recall that $P[x]$ denotes the set of polynomials whose coefficients all lie in $P$. Is $P[x]$ a prime ideal of $A[x]$?
\item In a local ring, show that the only idempotents are 0 and 1.
\item(\house) Let $R$ be a commutative ring and $I$ an ideal of $R$. The \col{radical}\index{radical of an ideal} of $I$ is the set \[\sqrt{I}=\{x\in R:x^n\in I\;\text{for some}\;n\in\N\}.\]

Show that $\sqrt I$ is an ideal, and is equal to the intersection of all prime ideals of $R$ containing $I$.
\item(\house)\label{spec_def} The \col{prime spectrum}\index{prime spectrum} of a commutative ring $R$, denoted $\Spec(R)$, is the set of prime ideals of $R$.
\begin{enumerate}
\item Show that every ring homomorphism $f:A\to B$ induces a map $f^{\ast}:\Spec(B)\to\Spec(A)$ given by $P\mapsto f^{-1}(P)$.
\item Let $I$ be an ideal of the ring $A$, and let $\pi:A\to A/I$ be the reduction modulo $I$. Show that the map $\pi^{\ast}:\Spec(A/I)\to\Spec(A)$ is injective.
\item For commutative rings $A$ and $B$, show that $\Spec(A\times B)$ is in bijection with the disjoint union $\Spec(A) \sqcup \Spec(B)$.
\end{enumerate}
\item\label{Zplocal} Prove Theorem \ref{local_crit} and use it to show that the ring $\Z_{(p)}$, the set of rational numbers whose denominator is not divisible by $p$, is a local ring whose unique maximal ideal is $(p)$.
\item Determine whether the ring of all infinitely differentiable functions $\R\to\R$ is a local ring.
\item Let $R$ be an integral domain and $P$ a prime ideal of $R$. Show that the set $R\setminus P$ is a multiplicative subset as defined in Exercise \ref{localization_exc}, Chapter \ref{subrings_chap}, and that the localization $R_S$ is a local ring.
\item Let $A$ be a commutative ring.
\begin{enumerate}
\item Show that if $M$ is a maximal ideal of the power series ring $A\pow x$, then $M\cap A$ is a maximal ideal of $A$. 
\item Show that for every prime ideal $P \subset A$, there exists a prime ideal $Q \subset A\pow x$ such that $Q \cap A = P$.
\end{enumerate}
\item\label{alggeo} Let $k$ be a field and $a_1,a_2,\ldots, a_n\in k$. Let \[M=(x_1-a_1,x_2-a_2,\ldots, x_n-a_n),\] an ideal in $k[x_1,x_2,\ldots, x_n]$. Denote by $e:k[x_1,x_2,\ldots, x_n]\to k$ the evaluation map $x_i\mapsto a_i$ for all $i$.
\begin{enumerate}
\item Noting that $x_i\equiv a_i\pmod M$ for every $i$, show that
\[f(x_1,\ldots,x_n)\equiv f(a_1,\ldots,a_n)\pmod M\] for every polynomial $f\in k[x_1,x_2,\ldots, x_n]$.
\item Show that $M=\ker(e)$ and that $M$ is a maximal ideal.
\end{enumerate}
\end{exercises}

\chapter{Second Review on Rings}\label{review2chap}

This chapter provides a review of the material covered in Chapters \ref{subrings_chap}--\ref{prime_chap}. The first two sections provide concise summaries of the core definitions, constructions, and theorems introduced in this part of the text. The next section contains detailed solutions to a selection of exercises from earlier chapters. The final section presents additional exercises that draw on ideas from across the chapters under review.

\section*{Core concepts}

\begin{itemize}
\item Subring of a ring, subfield of a field.
\item Subring generation: the construction $R[U]$.
\item Prime subring of a ring; isomorphism with either $\Z$ or $\Z_n$.
\item Prime subfield of a field; isomorphism with either $\Q$ or $\F_p$.
\item Ideals, finitely-generated ideals, principal ideals.
\item PIR and PID.
\item Sum and product of ideals.
\item The lattice of subrings of a ring; the lattice of ideals.
\item Direct product of rings.
\item Quotient rings: elements of $R/I$ have the form $[a]=a+I$ for $a\in R$.
\item Congruence modulo an ideal $I$: $a\equiv b\pmod I\iff a-b\in I$.
\item Homomorphisms, isomorphisms, automorphisms, embeddings.
\item The Frobenius homomorphism on a ring of prime characteristic.
\item Prime and maximal ideals. Local rings.
\end{itemize}

\section*{Essential theorems}

\begin{itemize}
\item Let $R$ be a ring. Every intersection of subrings of $R$ is a subring of $R$. Every intersection of ideals of $R$ is an ideal of $R$.
\item The ring $\Z$ is a PID. If $k$ is a field, the polynomial ring $k[x]$ is a PID, and the power series ring $k\pow x$ is a PID in which the nonzero ideals are $(x^n)$, $n\ge 0$.
\item Let $R$ and $S$ be rings. The ideals of $R\times S$ are the sets of the form $I\times J$, where $I$ is an ideal of $R$ and $J$ is an ideal of $S$. An analogous statement applies to direct products of finitely many rings.
\item For every ring $R$ and ideal $I$ of $R$, there are bijections
\[\{\text{Ideals of $R/I$}\}\longrightarrow\{\text{Ideals of $R$ containing $I$}\},\]
\[\{\text{Subrings of $R/I$}\}\longrightarrow\{\text{Subrings of $R$ containing $I$}\}.\]
\item The kernel of a ring homomorphism $f:R\to S$ is an ideal of $R$, and the image of $f$ is a subring of $S$. The map $f$ is injective if and only if $\ker(f)=(0)$.
\item (Isomorphism Theorem) For every homomorphism $f:R\to S$ there is a corresponding isomorphism \[R/\ker(f)\;\cong\;\im(f).\]
\item (Universal property of the field of fractions) Let $D$ be a domain and $K$ a field. Every embedding $D\to K$ extends uniquely to an embedding $\Frac(D)\to K$.
\item (Universal property of quotient rings) Let $R$ be a ring and $I$ an ideal of $R$. Every ring homomorphism $R\to S$ mapping $I$ to $0$ descends to a homomorphism $R/I\to S$.
\item (Universal property of polynomial rings) Let $A,B$ be commutative rings. The data of a homomorphism $A[x]\to B$ is equivalent to the data of a homomorphism $A\to B$ and an element $b\in B$.
\item For every commutative ring $R$ and ideal $I$ of $R$, $I$ is prime if and only if $R/I$ is a domain, and $I$ is maximal if and only if $R/I$ is a field.
\end{itemize}

\section*{Selected solutions}

\begin{framed}
\noindent \textbf{Exercise \ref{chap9_sel}, Chapter \ref{subrings_chap}:}\\
Give an example of a non-commutative ring whose prime subring is a field.
\end{framed}
\begin{proof}[Solution]
Let $R = M_n(\F_p)$ for any $n \ge 2$ and any prime $p$, clearly a non-commutative ring. The prime subring $R_0$ of $R$ is the set of integer multiples of the identity matrix, and hence is equal to the set of scalar matrices in $R$. The map $g:\F_p\to R_0$ mapping $x\mapsto xI_n$ is thus surjective; moreover, a simple argument shows that $g$ is a homomorphism and that $\ker(g)=(0)$. Thus, $g$ is an isomorphism. This proves that $R_0$ is isomorphic to $\F_p$ and is therefore a field.
\end{proof}

\begin{framed}
\noindent \textbf{Exercise \ref{chap10_sel}, Chapter \ref{ideals_chap}:}\\
Prove that in every commutative ring, the set of all nilpotent elements is an ideal. (This ideal is called the \col{nilradical} of the ring.) In the ring $M_2(\Z)$, show that the set of nilpotent elements is not an ideal. 
\end{framed}
\begin{proof}[Solution]
Let $A$ be a commutative ring. Define
\[
\text{Nil}(A) = \{x \in A : x^n = 0 \text{ for some } n \ge 1\}.
\]
We will show that $\text{Nil}(A)$ is an ideal. Clearly $0 \in \text{Nil}(A)$. Now, if $x, y \in \text{Nil}(A)$, then $x^m = 0$ and $y^n = 0$ for some $m,n \ge 1$; we claim that \[(x+y)^{m+n} = 0,\]

 so that $x+y \in \text{Nil}(A)$ and thus $\text{Nil}(A)$ is closed under addition. Using the binomial theorem,
\[
(x+y)^{m+n} = \sum_{k=0}^{m+n} \binom{m+n}{k} x^k y^{m+n-k}.
\]

Each term $x^k y^{m+n-k}$ vanishes because either $k \ge m$ implies $x^k = 0$ or $m+n-k \ge n$ implies $y^{m+n-k} = 0$. Thus $(x+y)^{m+n} = 0$, as claimed.

Next, we show that $\text{Nil}(A)$ is closed under scaling. If $x \in \text{Nil}(A)$ and $a \in A$, then $x^n = 0$ for some $n \ge 1$, and since $A$ is commutative,
\[
(ax)^n = a^n x^n = a^n \cdot 0 = 0,
\]
so $ax \in \text{Nil}(A)$. This completes the proof that $\text{Nil}(A)$ is an ideal.

Now consider the ring $M_2(\Z)$. We will show that the set of nilpotent matrices is not closed under addition, so it does not form an ideal. Let
\[
A = \begin{pmatrix} 0 & 1 \\ 0 & 0 \end{pmatrix}, \quad B = \begin{pmatrix} 0 & 0 \\ 1 & 0 \end{pmatrix}.
\]
Then $A^2 = B^2 = 0$, so $A$ and $B$ are nilpotent. We compute
\[
A + B = \begin{pmatrix} 0 & 1 \\ 1 & 0 \end{pmatrix}, \quad (A+B)^2 = \begin{pmatrix} 1 & 0 \\ 0 & 1 \end{pmatrix} = I.
\]
It follows that $A + B$ is not nilpotent, as every power of $A+B$ is equal to either $A+B$ or $I$. This proves that the sum of nilpotent matrices is not necessarily nilpotent.
\end{proof}

\begin{framed}
\noindent \textbf{Exercise \ref{chap11_sel}, Chapter \ref{quo_chap}:}\\
Prove that every quotient of the quadratic integer ring $\Z[\sqrt d]$ by a nonzero ideal is a finite ring.

\end{framed}
\begin{proof}[Solution]
Let $R = \Z[\sqrt{d}]$ and let $I$ be a nonzero ideal of $R$. We aim to show that the quotient ring $R/I$ is finite.

Recall that the norm map $N: R \to \Z$ is given by
\[N(a + b\sqrt{d}) = a^2 - db^2.\]

Fix a nonzero element $z = a_0 + b_0\sqrt{d} \in I$ and let $n=|N(z)|$, a positive integer. Note that $n\in I$ since $n=\pm z(a_0-b_0\sqrt d)$ and $z\in I$. This implies that if two integers are congruent modulo $n$, then they are congruent modulo $I$. Now define a map $\Z_n\times\Z_n\to R/I$ by \[(a+(n),b+(n))\mapsto (a+b\sqrt d)+I.\] 

This map is well defined and surjective, so $\#(R/I)\le n^2$.
\end{proof}

\begin{framed}
\noindent\textbf{Exercise \ref{2nd_isothm}, Chapter \ref{homo_chap}}\\
Let $R$ be a ring and $S$ a subring of $R$. Let $I$ be an ideal of $R$. Show that $S\cap I$ is an ideal of $S$, that $S+I$ is a subring of $R$, and that $I$ is an ideal of $S+I$. Moreover, \[S/(S\cap I)\cong (S+I)/I.\]
 Here,  $S+I=\{s + i : s \in S,\, i \in I\}$.
\end{framed}
\begin{proof}[Solution]
We verify each claim in turn.

First, \(S \cap I\) is an ideal of \(S\): both $S$ and $I$ contain 0 and are closed under addition, so the same holds true for $S\cap I$. Now, if \(s \in S\) and \(x \in S \cap I\), then \(sx \in S\) since \(S\) is a subring, and \(sx \in I\) since \(x \in I\) and \(I\) is an ideal of \(R\). So \(sx \in S \cap I\). This proves that $S\cap I$ is closed under scaling by elements of $S$, and therefore an ideal of $S$.

Next, $S + I$ is a subring of \(R\): it contains \(0 = 0 + 0\) and $1=1+0$, and is closed under addition and additive inverses because both $S$ and $I$ are. If \(s_1 + i_1, s_2 + i_2 \in S + I\), then their product is
\[
(s_1 + i_1)(s_2 + i_2) = s_1s_2 + s_1i_2 + i_1s_2 + i_1i_2.
\]
Now \(s_1s_2 \in S\), \(s_1i_2, i_1s_2, i_1i_2 \in I\) since \(I\) is an ideal and \(S \subseteq R\). Therefore the product lies in \(S + I\), showing that \(S + I\) is a subring.

Next, we show that \(I\) is an ideal of \(S + I\), the key detail being closure under scaling. For any \(r \in S + I\), say \(r = s + i\), and \(x \in I\), we have
\[
r x = (s + i)x = sx + ix \in I,
\]
since \(sx \in I\) (as \(s \in S \subseteq R\), \(x \in I\)) and \(ix \in I\). Thus \(I\) absorbs multiplication from \(S + I\).

Finally, we prove the isomorphism \(S/(S \cap I) \cong (S + I)/I\). Define
\[
\varphi : S \to (S + I)/I, \quad \varphi(s) = s + I.
\]
The definitions easily imply that $\varphi$ is a surjective ring homomorphism whose kernel is $S\cap I$. By the Isomorphism Theorem, we conclude that
\[
S/(S \cap I) \cong (S + I)/I.\qedhere
\]
\end{proof}

\begin{framed}
\noindent \textbf{Exercise \ref{chap12_sel}, Chapter \ref{homo_chap}:}\\
Let $X$ be a set and $R(X)$ the ring of sets $(P(X),\Delta,\cap)$ in which addition and multiplication of subsets of $X$ are, respectively, the symmetric difference and intersection of sets. Show that $R(X)$ is isomorphic to the ring of functions $\F_2^X$.
\end{framed}
\begin{proof}[Solution]
Define a map \(\varphi : R(X) \to \F_2^X\) as follows: for each subset \(A \subseteq X\), the function \(\varphi(A)\) is given by
\[
\varphi(A)(x) = \begin{cases}
1 & \text{if } x \in A, \\
0 & \text{if } x \notin A.
\end{cases}
\]

In other words, $\varphi(A)$ is the \href{https://en.wikipedia.org/wiki/Indicator_function}{characteristic function} of $A$.

The map $\varphi$ is surjective: given a function \(f \in\F_2^X\), we have $\varphi(A)=f$, where \(A=\{x \in X : f(x) = 1\}\).

Next, we check that \(\varphi\) is a ring homomorphism. Recall that $X$ is the unity of the ring $R(X)$. By definition, $\varphi(X)$ is the constant function $x\mapsto 1$ for all $x\in X$, which is the unity in $\F_2^X$. Similarly, $\emptyset$ is the additive identity of $R(X)$, and $\varphi(\emptyset)$ is the constant function $x\mapsto 0$, which is the zero element of $\F_2^X$. This shows that $\varphi(0)=0$ and $\varphi(1)=1$.

Now, for subsets \(A, B \subseteq X\), we have
\[
\varphi(A \Delta B)(x) = \begin{cases}
1 & \text{if } x \in A \text{ or } x \in B, \text{ but not both}, \\
0 & \text{otherwise},
\end{cases}
\]

which agrees with $\varphi(A)(x)+\varphi(B)(x)$. Thus, $\varphi(A \Delta B)=\varphi(A)+\varphi(B)$. Similarly,
\[
\varphi(A \cap B)(x) = \varphi(A)(x) \cdot \varphi(B)(x),
\]

so $\varphi(A\cap B)=\varphi(A)\varphi(B)$. This shows that $\varphi$ is a homomorphism. Finally, we have $\ker(\varphi)=\{\emptyset\}$, which is the zero ideal of $R(X)$.

Therefore, \(\varphi\) is a ring isomorphism \(R(X)\to\F_2^X\).
\end{proof}

\begin{framed}
\noindent \textbf{Exercise \ref{chap13_sel}, Chapter \ref{univ_chap}:}\\
Let $R$ be a ring and $n\ge 1$. Show that there exists a homomorphism $\Z_n\to R$ if and only if $\ch(R)$ divides $n$.
\end{framed}
\begin{proof}[Solution]
Since \(\mathbb{Z}_n=\mathbb{Z}/(n)\), a homomorphism \(\varphi : \mathbb{Z}_n \to R\) corresponds to a homomorphism \(\mathbb{Z} \to R\) that sends \(n\) to \(0\), by the universal property of quotient rings. Let \(f_R:\mathbb{Z} \to R\) be the unique homomorphism, which is given by \(k \mapsto k \cdot 1_R\). As shown in the proof of Theorem \ref{prime_sub_thm}, the kernel of $f_R$ is the principal ideal generated by $\ch(R)$. Hence, we have $f_R(n)=0$ if and only if $n$ is a multiple of $\ch(R)$.
\end{proof}

\begin{framed}
\noindent \textbf{Exercise \ref{chap14_sel}, Chapter \ref{prime_chap}:}\\
Let $P$ be a prime ideal of a commutative ring $A$. Recall that $P[x]$ denotes the set of polynomials whose coefficients all lie in $P$. Is $P[x]$ a prime ideal of $A[x]$?
\end{framed}
\begin{proof}[Solution]
We claim that \( P[x] \) is a prime ideal of \( A[x] \). The reduction map \( A \to A/P \) induces a  homomorphism  \(f:A[x] \to (A/P)[x]\). Note that $f$ is surjective and that its kernel is \( P[x] \). By the Isomorphism Theorem, $A[x]/P[x]$ is isomorphic to $(A/P)[x]$. Since \( A/P \) is an integral domain (because \( P \) is prime), the ring \( (A/P)[x] \) is also an integral domain. Therefore, $A[x]/P[x]$ is a domain, and \( P[x] \) is a prime ideal of \( A[x] \), as claimed.
\end{proof}

\section*{Additional Practice}

\begin{exercises}
 \item Describe all subrings of $\Z[x]$ of the form \( \Z[f(x)] \) where \( f(x) \in \Z[x] \) has degree 2. How do these subrings relate to one another?
\item Let $x_1,\ldots, x_n$ be elements of a commutative ring. Show that \[(x_1,\ldots, x_n)=(x_1)+\cdots+(x_n).\]
\item Determine all ideals and all subrings of $\Q\times\R$.
\item Find all idempotent elements of the ring $\Z[x]/(x^2-1)$.
\item Let $T$ be the set of matrices of the form $\begin{pmatrix}a&0\\b&c\end{pmatrix}$ with $a,b,c\in\Z$. Show that $T$ is a subring of $M_2(\Z)$ and find all its ideals.
\item In the ring $R=M_2(\R)$, let $S$ and $T$, respectively, be the sets of matrices of the form \[\begin{pmatrix}a&b\\0&a\end{pmatrix}\quad\text{and}\quad\begin{pmatrix}0&b\\0&0\end{pmatrix}.\] Are $S$ and $T$ subrings of $R$? Is $T$ an ideal of $S$? Is $T$ an ideal of $R$? Is $S$ an ideal of $R$?
\item Give an example of a non-commutative ring having a commutative quotient.
\item The \col{center}\index{center of a ring} of a ring $A$, denoted by $Z(A)$, is the set of all elements of $A$ that commute with every element of $A$: \[Z(A)=\{a\in A:ax=xa\;\text{ for every }\;x\in A\}.\]
\begin{enumerate}
\item Show that $Z(A)$ is a subring of $A$ and that it is commutative.
\item Show that the center of a division ring is a field. (Recall that a division ring is a nontrivial ring in which every nonzero element is a unit. 
\item What is the center of $\H$? Of $M_n(\Z)$?
\end{enumerate}
\item Recall that in $\calC[a,b]$ there is a maximal ideal $M_r$ for every $r\in[a,b]$.
\begin{enumerate}
\item Show that $M_r\ne (x-r)$.
\item Show that $M_r$ is not finitely generated.
\end{enumerate}
\item Let $A$ be a commutative ring. Suppose that $P$, $I$, $J$ are ideals of $A$ with $P$ prime. Prove the following:
\begin{enumerate}
\item If $IJ\subseteq P$, then $I\subseteq P$ or $J\subseteq P$.
\item If $P=IJ$, then $P=I$ or $P=J$.
\end{enumerate}
\item An \col{irreducible ideal}\index{irreducible ideal} in a ring is an ideal $H$ such that whenever $H=I\cap J$, where $I$ and $J$ are ideals, necessarily $H=I$ or $H=J$. Prove that prime ideals are irreducible, and give an example of an irreducible ideal that is not prime.
\item Let $A$ be an integral domain and $a,b$ integers with $\gcd(a,b)=1$. Show that the ideal $(x^a-y^b)$ is a prime ideal of $A[x,y]$.
\item Let $A$ be a commutative ring with prime ideals $P_1,\ldots, P_n$, and let $I$ be an ideal of $A$. Show that if $I\subseteq\cup_{i=1}^n P_i$, then $I\subseteq P_i$ for some $i$.
\item Show that for every commutative ring $A$, the set of all nilpotent elements is equal to the intersection of all prime ideals of $A$.
\item Let $R$ be a finite ring. Show that there exist integers $a>b$ such that $x^a=x^b$ for every $x\in R$. (Use the ring $R^n$, where $n=\#R$.)
\item Show that the rings $\Z[x]$ and $\Q[x]$ are not isomorphic.
\item \begin{enumerate}
\item Let $R$ be a ring and $h:R\to R$ a homomorphism. Show that the set of elements fixed by $h$, i.e., the set $\{a\in R:h(a)=a\}$ is a subring of $R$.
\item Let $R=\F_p[x]$, where $p$ is prime, and let $h:R\to R$ be the Frobenius homomorphism, the map $f(x)\mapsto f(x)^p$. Show that the subring of $R$ fixed by $h$ is precisely the subring $\F_p[x^p]$.
\end{enumerate}
\item Show that, up to isomorphism, there is only one Boolean field.
\item For every prime $p$, we define the rings $R_1=\F_p[x]/(x^2-2)$ and $R_2=\F_p[x]/(x^2-3)$. In the cases $p=2,5$, and $11$, determine whether $R_1\cong R_2$.
\item A \col{simple ring}\index{simple ring} is a ring having exactly two ideals, namely $(0)$ and $(1)$. Show that for every field $k$, the matrix ring $M_n(k)$ is simple.
\item Let $R$ be a commutative ring with ideal $I$. Recall that $M_n(I)$ denotes the set of matrices with entries in $I$, and that $M_n(I)$ is an ideal of $M_n(R)$. Use the Isomorphism Theorem to show that \[M_n(R)/M_n(I)\cong M_n(R/I).\]
\item Let $F$ be a field. Determine whether the rings \[F[x,y]/(y^2-x)\quad\text{and}\quad F[x,y]/(x^2-y)\] are isomorphic.
\item Let $R=\F_2[x,y]$ and $I=(x^2,xy,y^2)$. Show that the quotient ring $R/I$ has cardinality 8.
\item Let $R=\Z[x]$ and $I=(x^2-2x)$.
\begin{enumerate}
\item Find a nontrivial zero divisor in $R/I$.
\item Find a proper ideal of $R$ strictly containing $I$.
\item Determine whether there exist integers $a,b$ such that \[x^5\equiv a+bx\;\,\text{(mod $I$)}.\]
\end{enumerate}
\item A \col{von Neumann regular ring}\index{von Neumann regular ring} is a ring $R$ with the following property: for every $a\in R$ there exists $x\in R$ such that $a=axa$. We will write VNR to denote this type of ring. Prove the following:
\begin{enumerate}
\item Every field is a VNR. A domain is a VNR if and only if it is a field.
\item Any direct product of VNR is a VNR.
\item Every VNR is \emph{reduced}, i.e., has no nonzero nilpotent elements.
\item If $R$ is commutative VNR, show that every finitely-generated ideal of $R$ is principal and generated by an idempotent element.
\end{enumerate}
\item Let $d$ be a squarefree integer with $d\not\equiv 1\pmod 4$, and let $R=\Z[\sqrt d]$.
\begin{enumerate}
\item For every prime $p\in\Z$, show that $R/(p)\cong\F_p[x]/(x^2-d)$.
\item If $p\mid d$, show that $(p)=P^2$ for some maximal ideal $P$ of $R$.
\item If $d$ is not a square modulo $p$, then $(p)$ is a maximal ideal of $R$.
\item If $d$ is a square modulo $p$, then $(p)=PQ$, where $P$ and $Q$ are distinct maximal ideals of $R$.
\end{enumerate}
\end{exercises}

\chapter{The Remainder Theorem}\label{crt_chap}

This chapter introduces a powerful structural theorem that expresses certain quotient rings as direct products of simpler ones. Such representations make it easier to analyze the quotient ring—its units, zero divisors, ideals, and prime ideals, among others. Two classical consequences are developed here: the \href{https://en.wikipedia.org/wiki/Chinese_remainder_theorem}{Chinese Remainder Theorem}, fundamental in number theory and modular arithmetic, and the \href{https://en.wikipedia.org/wiki/Lagrange_polynomial}{Lagrange Interpolation Formula}, an important tool in numerical analysis, coding theory, and cryptography.

The following terminology is needed in order to state the main theorem. Let $I$ and $J$ be ideals of a commutative ring $R$. We say that $I$ and $J$ are \col{comaximal}\index{comaximal ideals} if $I+J=R$, or equivalently, if $R$ is the smallest ideal containing both $I$ and $J$. A finite list of ideals $I_1,\ldots,I_n$ is \col{pairwise comaximal}\index{pairwise comaximal ideals} if every pair of distinct ideals in the list is comaximal.

\begin{ex}[Comaximal ideals]\;
\begin{examples}
\item If $a$ and $b$ are \emph{coprime integers}---that is, integers having no prime factors in common---then the ideals $(a)$ and $(b)$ of $\Z$ are comaximal. Indeed, if $(m)$ is an ideal of $\Z$ containing both $a$ and $b$, then $a$ and $b$ are divisible by $m$. It follows that $m=\pm 1$ and thus $(m)=\Z$, for otherwise any prime factor of $m$ would divide both $a$ and $b$, contradicting coprimality. More generally, if $a_1,\ldots, a_n$ are pairwise coprime integers, then the principal ideals they generate are pairwise comaximal.
\item In the ring $\Z[i]$ of Gaussian integers, the ideals $(3+2i)$ and $(4+i)$ are comaximal. Indeed, the sum of these ideals is the full ring $\Z[i]$, as the following equation shows that 1 can be expressed as a linear combination of $3+2i$ and $4+i$:
\[(1+2i)(3+2i)+(-2i)(4+i)=1.\]
The method used to find such an identity is known as the \emph{extended Euclidean algorithm}, discussed in the next chapter.
\item In any ring $R$, suppose that $M$ and $N$ are distinct maximal ideals. Then these ideals are comaximal, as the ideal $M+N$ is strictly larger than $M$ and hence equal to $R$ by maximality of $M$. More generally, any finite list of distinct maximal ideals is pairwise comaximal.
\end{examples}
\end{ex}

Recall from Chapter \ref{ideals_chap} that if $I$ and $J$ are ideals of $R$, then the product ideal $IJ$ is the ideal generated by the set of elements of the form $ij$, where $i\in I$ and $j\in J$. It follows from this definition that $IJ\subseteq I\cap J$. If, moreover, $I$ and $J$ are comaximal, then Exercise \ref{comax_exc} in the same chapter shows that $IJ=I\cap J$. We now generalize this result to any finite collection of comaximal ideals.
\begin{framed}
\begin{thm}\label{comax_prod} If $I_1,\ldots,I_n$ are pairwise comaximal ideals of a commutative ring, where $n\ge 2$, then \[I_1\cap I_2\cap\cdots\cap I_n=I_1\cdot I_2\cdots I_n.\]\end{thm}
\end{framed}
\begin{proof}
It suffices to prove that the intersection is contained in the product, as the reverse inclusion holds in general. We proceed by induction on $n$, with the base case $n=2$ being the exercise cited above. Assume $n\ge 3$ and that the result is true for $n-1$ ideals. Let $I_1,\ldots, I_n$ be pairwise comaximal ideals. Setting $J=\bigcap_{k=2}^nI_k$, we claim that $I_1$ and $J$ are comaximal. For $2\le k\le n$ we have $I_1+I_k=R$, so there exist elements $a_k\in I_1$ and $b_k\in I_k$ such that $a_k+b_k=1$. It follows that
\[\prod_{k=2}^n(a_k+b_k)=1.\]
Expanding this product, let $a$ be the sum of all terms involving at least one $a_k$, and let $b=b_2\cdots b_n$ be the unique term with no $a_k$. The properties of ideals then imply that $a\in I_1$ and $b\in J$. Since $a+b=1$, this proves the claim that $I_1$ and $J$ are comaximal. Now, since $I_2,\ldots, I_n$ are pairwise comaximal, we have $J=I_2\cdots I_n$ by inductive hypothesis. Moreover, since $I_1$ and $J$ are comaximal, then $I_1J=I_1\cap J$ by the base case $n=2$. Therefore,
\[I_1\cap I_2\cap\cdots\cap I_n=I_1\cap J=I_1J=I_1\cdot I_2\cdots I_n,\]
completing the induction.
\end{proof}

The following theorem captures an essential feature of comaximal ideals: the existence of elements that are congruent to 1 modulo one ideal and to 0 modulo all the others. 

\begin{framed}
\begin{thm}\label{rt_idempo}
Let $R$ be a commutative ring and let $I_1,\ldots,I_n$ be pairwise comaximal ideals of $R$, with $n\ge 2$. There exist elements $e_1,\ldots, e_n\in R$ satisfying, for each $1\le k\le n$,
\[e_k\equiv 1\;{\text{(mod $I_k$)}}\quad\text{and}\quad e_k\equiv 0\;{\text{(mod $I_j$)}}\quad\text{for}\;\;j\ne k.\]
\end{thm}
\end{framed}
\begin{proof}
Fix $1\le k\le n$. For each $j\ne k$ we have $I_k+I_j=R$, so there exist $s_j\in I_k$ and $t_j\in I_j$ such that $s_j+t_j=1$. Let $e_k:=\prod_{j\ne k}t_j$ and consider the product \[1=\prod_{j\ne k}(s_j+t_j).\] Expanding this product, let $x$ be the sum of all terms containing at least one $s_j$. Then $x\in I_k$ since each $s_j$ belongs to $I_k$ and $I_k$ is an ideal. Moreover, the above equation implies that $1=x+e_k$. Reducing modulo $I_k$, we obtain $e_k\equiv 1\pmod{I_k}$. Now, the definition of $e_k$ implies that $e_k\in I_j$ for each $j\ne k$, as $e_k$ is a multiple of $t_j$ and $t_j\in I_j$. Hence $e_k\equiv 0\pmod{I_j}$. Thus the elements $e_1,\dots,e_n$ satisfy all the required congruences.
\end{proof}

We may now state the main result of this chapter.

\begin{framed}
\begin{thm}[Remainder Theorem (\textnormal{isomorphism form})]\label{RT_thm} Let $R$ be a commutative ring and $I_1,\ldots,I_n$ pairwise comaximal ideals of $R$. There is an isomorphism
\[R/(I_1\cdots I_n)\;\longrightarrow\;(R/I_1)\times(R/I_2)\times\cdots\times(R/I_n)\] given by $[r]\mapsto(r+I_1,r+I_2,\ldots,r+I_n)$. In particular,
\[R/(I_1\cdots I_n)\;\cong\;(R/I_1)\times(R/I_2)\times\cdots\times(R/I_n).\]
\end{thm}
\end{framed}
\begin{proof}
The map $f:R\to(R/I_1)\times(R/I_2)\times\cdots\times(R/I_n)$ defined by
\[f(r)=(r+I_1,r+I_2,\ldots,r+I_n)\]
is a well-defined homomorphism whose kernel is
\[\ker(f) = \{ r \in R : r \in I_k \text{ for all } k\} = I_1 \cap \cdots \cap I_n.\]
Applying Theorem \ref{comax_prod} gives $\ker(f)=I_1 \cdots I_n$.

We claim that $f$ is surjective. Let $a=(a_1+I_1,\ldots,a_n+I_n)$ be an arbitrary element of the codomain of $f$. With elements $e_1,\ldots, e_n$ as in Theorem \ref{rt_idempo}, let \[x = e_1 a_1 + e_2 a_2 + \cdots + e_n a_n.\] The properties of the $e_k$'s imply that $f(x)=a$, proving the claim. Finally, the Isomorphism Theorem (\ref{iso_thm}) implies that $f$ induces the desired isomorphism between $R/(I_1\cdots I_n)$ and $(R/I_1)\times\cdots\times(R/I_n)$.
\end{proof}

An alternate formulation of the Remainder Theorem is in terms of congruences. The following theorem uses the proofs of theorems \ref{rt_idempo} and \ref{RT_thm} to construct a solution to a given system of congruences.

\newpage

\begin{framed}
\begin{thm}[Remainder Theorem (\textnormal{congruence form})]\label{rt_cong} Let $R$ be a commutative ring and $I_1,\ldots,I_n$ pairwise comaximal ideals of $R$. For any $b_1,b_2,\ldots,b_n\in R$, there exists $x\in R$ such that \[x\equiv b_k\;\text{(mod $I_k$)}\quad\text{for}\quad 1\le k\le n.\] Moreover, such an element $x$ is unique modulo $I_1\cdots I_n$. One explicit solution is given by $x = \sum_{k=1}^ne_k b_k$, where $e_k$ is defined as follows: for each $j\ne k$, choose $t_j\in I_j$ such that $t_j\equiv 1\pmod{I_k}$. Then $e_k=\prod_{j\ne k}t_j$.
\end{thm}
\end{framed}

\begin{1ex}[Solving a system of congruences] Consider the problem of finding $x\in\Z$ satisfying \[
\begin{cases}
\begin{aligned}
x &\equiv 3\;\text{(mod 4)}, \\
x &\equiv 8\;\text{(mod 13)}.
\end{aligned}
\end{cases}
\]

Since $\gcd(4,13)=1$, the ideals $(4)$ and $(13)$ are comaximal, so there exists $x\in\Z$ solving the above system, and $x$ is unique modulo 52. An explicit solution is given by $3e_1+8e_2$ for some integers $e_1$ and $e_2$. To construct $e_1$, we seek $t_2\in(13)$ such that $t_2\equiv 1\pmod{4}$, so we take $t_2=13$ and then $e_1=t_2$. To construct $e_2$, we seek $t_1\in(4)$ such that $t_1\equiv 1\pmod{13}$, so we take $t_1=40$ and then $e_2=t_1$. Our solution to the system is thus $x=3(13)+8(40)=359$. Reducing $x$ modulo 52  gives $x\equiv 47\pmod{52}$. The full solution set is then described by $x=47+52n$, $n\in\Z$.

\end{1ex}

\subsection*{Lagrange Interpolation} 

\begin{framed}
\begin{thm}[Lagrange Interpolation]\label{lagrange} Let $F$ be a field and let $(a_1,\ldots, a_n)$ and $(b_1,\ldots, b_n)$ be elements of $F^n$ such that the $a_k$'s are pairwise distinct. There exists a unique polynomial $p(x)\in F[x]$ with $\deg(p)<n$ such that $p(a_k)=b_k$ for $1\le k\le n$. Explicitly,
\[p(x)=\sum_{k=1}^nb_ke_k(x),\quad\text{where}\quad e_k(x)=\prod_{j\ne k}\frac{x-a_j}{a_k-a_j}.\]
\end{thm}
\vspace{-4mm}
\end{framed}

\begin{proof} To prove the uniqueness of $p(x)$, suppose that $g(x)$ is also a polynomial of degree less than $n$ satisfying $g(a_k)=b_k$ for every $k$. Then the polynomial $p-g\in F[x]$ has degree less than $n$ and has each $a_k$ as a root, a total of $n$ roots. This implies that $p-g$ is the zero polynomial, by Theorem \ref{roots_number}, proving the uniqueness of $p(x)$. 

Existence of $p(x)$ follows from the Remainder Theorem applied in the ring $F[x]$. Every principal ideal $(x-a_k)\subset F[x]$ is maximal, so the ideals \[(x-a_1),\ldots, (x-a_n)\] are pairwise comaximal. By Theorem \ref{rt_cong}, there exists $p(x)\in F[x]$ such that, for $1\le k\le n$,
\begin{equation}\label{lag_poly}\tag{$\ast$}
p(x)\equiv b_k\;\text{(mod $(x-a_k)$)}.
\end{equation}
Moreover, one such $p(x)$ is given by \[p(x)=\sum_{k=1}^ne_k(x)b_k,\] where $e_k$ is constructed as follows: for $j\ne k$, let $t_j(x)=\frac{1}{a_k-a_j}(x-a_j)$, which lies in $(x-a_j)$. The relation
\[t_j(x)=1+\frac{x-a_k}{a_k-a_j}\] implies that $t_j\equiv 1\pmod{x-a_k}$. The recipe in Theorem \ref{rt_cong} now leads to defining \[e_k(x)=\prod_{j\ne k}t_j=\prod_{j\ne k}\frac{x-a_j}{a_k-a_j}.\]

For each $1\le k\le n$, the above polynomial $p(x)$ satisfies \eqref{lag_poly}, so we may write $p(x)=b_k+(x-a_k)q(x)$ for some polynomial $q(x)$. Evaluating at $x=a_k$ gives $p(a_k)=b_k$, showing that $p(x)$ meets the requirements of the theorem.
\end{proof}

\begin{1ex}[Functions on finite fields]
As an application of Lagrange Interpolation, we prove that every function on a finite field is a polynomial function. This result also follows from Exercise \ref{finite_polynomial_fn}, Chapter \ref{multi_chap}. Let $k$ be a finite field and $f:k\to k$ any function. Let $a_1,\ldots, a_n$ be an enumeration of the elements of $k$. Let $b_i=f(a_i)$ for all $i$. Then Theorem \ref{lagrange} yields the existence of a polynomial $p(x)$ satisfying $p(a_i)=b_i=f(a_i)$ for all $i$. Hence $f$ coincides with the polynomial function induced by $p(x)$.
\end{1ex}

\subsection*{Chinese Remainder Theorem}

Specializing to the case $R=\Z$ in theorems \ref{RT_thm} and \ref{rt_cong}, we obtain a result classically known in number theory as the Chinese Remainder Theorem (CRT)\index{Chinese Remainder Theorem}, an essential tool for solving systems of congruences in $\Z$.

\begin{framed}
\begin{thm}[CRT]\label{CRT}
Let $m_1,\ldots, m_n$ be pairwise coprime integers. For any integers $a_1,\ldots, a_n$, the system of congruences \[x\equiv a_i\;\text{(mod $m_i$)}\quad\text{for}\quad i=1,\ldots, n\]
has a solution $x\in\Z$, and such a solution is unique modulo $m:=m_1\cdots m_n$. In particular, there exists a ring isomorphism
\[\Z_m\cong\Z_{m_1}\times\cdots\times\Z_{m_n}.\]
\end{thm}
\end{framed}

\begin{1ex}[Counting units modulo $n$] Since $30=2\cdot 3\cdot 5$, CRT gives \[\Z_{30} \cong \Z_2 \times \Z_3 \times \Z_5.\]

Using this description of $\Z_{30}$, we can easily determine, for instance, the number of units in $\Z_{30}$. Units in $\Z_2 \times \Z_3 \times \Z_5$ have the form $(a,b,c)$, with each entry being a unit in the appropriate ring. Since each ring in this direct product is a field, the number of units is equal to $p-1$, where $p$ is the cardinality of the field. Hence, the number of units in $\Z_{30}$ is $(2-1)(3-1)(5-1)=8$. See Exercise \ref{eulerphi} concerning the \emph{Euler phi function}, which generalizes the observations in this example.
\end{1ex}

\newpage
\section*{Exercises}

\begin{note}\;
\begin{itemize}
\item Solutions to starred exercises (\sel) are included in Chapter~\ref{review3chap}.
\item A house (\house) marks exercises suggested as sample homework.
\end{itemize}
\end{note}

\medskip

\begin{exercises}
\item Let $R = \F_3[x]$, and consider the ideals $I = (x)$ and $J = (x+1)$. Show that $I$ and $J$ are comaximal in $R$, and use the Remainder Theorem to prove that there is an isomorphism
\[
\F_3[x]/(x(x+1)) \cong \F_3\times\F_3.
\]

Is the ring \( \F_3[x]/(x(x+1)) \) a field?
\item Show that the system of congruences
\[
\begin{cases}
x \equiv 2\;\text{(mod 4)}, \\
x \equiv 3\;\text{(mod 6)}
\end{cases}
\]
has no solution $x\in\Z$, and relate this to $(4)$ and $(6)$ not being comaximal.
\item Find all solutions $x\in\Z$ to the system of congruences
 \[
\begin{cases}
\begin{aligned}
x &\equiv 1\;\text{(mod 5)}, \\
x &\equiv 2\;\text{(mod 6)},\\
x &\equiv 3\;\text{(mod 7)}.
\end{aligned}
\end{cases}
\]
\item Find the unique cubic polynomial $p\in\F_7[x]$ satisfying
\[p(0)=1,\;p(1)=3,\;p(2)=2,\;p(3)=6.\]
\item(\house) Let \( R = \Z_{105} \). Use CRT to show that $R\cong \Z_3\times\Z_5\times\Z_7$, and use this isomorphism to determine:
\begin{enumerate}
    \item The number of units in $R$;
    \item The number of zero divisors in $R$;
    \item The number of maximal ideals of $R$.
\end{enumerate}
\item(\sel)\label{chap16_sel} Let $p_1,\ldots, p_s$ be distinct odd prime integers, and $n=p_1\cdots p_s$. Determine the number of solutions of the equation $x^2=x$ in the ring $\Z_n$.
\item(\house)\label{eulerphi} The \col{Euler $\phi$ function}\index{Euler $\phi$ function} is defined for $n\ge 1$ by
\[\phi(n):=\#\Z_n^{\times}=\#\{a\in\{1,2,\ldots,n\}: \gcd(a,n)=1\}.\]
\begin{enumerate}
\item Suppose $n$ is factored as $n=p_1^{e_1}\cdots p_s^{e_s}$, where the $p_i$'s are distinct primes. Use CRT to show that $\phi(n)=\phi(p_1^{e_1})\cdots\phi(p_s^{e_s})$.
\item Let $p$ be prime and $e\ge 1$. Show that $\phi(p^e)=p^e-p^{e-1}$.
\item Calculate $\phi(n)$ for $10\le n\le 20$ by factoring $n$.
\end{enumerate}
\end{exercises}

\chapter{Euclidean and GCD Domains}\label{euc_chap}

In this chapter and the next, we generalize several concepts that are well known in the context of $\Z$: divisibility, greatest common divisor (gcd), division algorithm, prime number, and factorization as a product of primes.  One novelty in this generalization is the concept of an \emph{irreducible element} of an integral domain, which becomes especially relevant for rings with a more intricate structure than $\Z$. The present chapter develops terminology for divisibility in integral domains, and studies \emph{Euclidean domains}, a class of rings in which a form of division with remainder applies, as well as \emph{GCD domains}, rings in which every pair of elements has a well-defined gcd. The next chapter defines irreducibility and studies \emph{UFDs}, rings whose elements admit unique factorizations into irreducibles.

\section*{Divisibility in integral domains}
Let $D$ be an integral domain. We define a relation of divisibility on $D$ as follows: for elements $x,y\in D$, we say that $x$ \col{divides}\index{divisor of a ring element} $y$, written $x\mid y$, if $y=xa$ for some $a\in D$. In that case, $y$ is a \col{multiple}\index{multiple of a ring element} of $x$, and $x$ is a \col{divisor} of $y$. The relation of divisibility is reflexive and transitive: every element divides itself, and if $a\mid b$ and $b\mid c$, then $a\mid c$. The units of $D$ are precisely the divisors of 1, and every divisor of a unit is a unit.

A \col{greatest common divisor}\index{greatest common divisor} (gcd) of elements $a,b\in D$ is an element $d$ satisfying $d\mid a$, $d\mid b$, and if $x\mid a$ and $x\mid b$, then $x\mid d$. Similarly, a \col{least common multiple}\index{least common multiple} (lcm) of $a$ and $b$ satisfies $a\mid m$, $b\mid m$, and if $a\mid x$ and $b\mid x$, then $m\mid x$. We say that $a$ and $b$ are \col{coprime}\index{coprime  elements} if 1 is a gcd of $a$ and $b$; equivalently, the only common divisors of $a$ and $b$ are the units of $D$.

The gcd and lcm can be defined not only for pairs of elements of $D$, but in fact for any nonempty subset of $D$. The definitions are analogous to those above: for instance, a gcd of a subset $S\subseteq D$ is an element $g\in D$ dividing every element of $S$, and such that if $x\in D$ satisfies $x\mid s$ for all $s\in S$, then $x\mid g$.

Gcds and lcms do not necessarily exist in a given domain---see \ref{no_gcd} in Examples \ref{gcd_ass} below, which exhibits a pair of elements that do not have a gcd. Moreover, even when gcds do exist, they are typically not unique; for instance, in $\Z$, both 2 and $-2$ satisfy the definition of a gcd of 4 and 6. This failure of uniqueness is clarified by the following concept: elements $a,b\in D$ are called \col{associates}\index{associates}, written \[a\sim b,\] if there is a unit $u\in D^{\times}$ such that $a=ub$. Equivalently, $a$ and $b$ are associates if they generate the same ideal: $(a)=(b)$. For example, in $\Z$, 2 and $-2$ are associates. The relation $\sim$ is easily seen to be an equivalence relation on $D$.

The following theorem shows that when a gcd or lcm exists, it is unique up to associates.

\begin{framed}
\begin{thm}\label{gcd_ass_thm}
Let $D$ be a domain and let $a,b\in D$. Suppose that $d$ is a gcd of $a$ and $b$ in $D$. Then the set of all gcds of $a$ and $b$ agrees with the set of associates of $d$. The same holds true for lcm.
\end{thm}
\end{framed}
\begin{proof}
Exercise \ref{gcd_ass_exc}.
\end{proof}

In some domains it is common to ensure uniqueness of gcds by imposing additional requirements in the definition. In $\Z$, this is achieved by choosing the unique \emph{positive} gcd, and in $k[x]$, where $k$ is a field, by choosing the unique \emph{monic} gcd. We will follow this convention for the rest of the course. In domains where gcds are not uniquely determined, we will write $\gcd(a,b)$ to represent \emph{any} gcd of the elements $a$ and $b$.

\begin{ex}[Associates and gcds]\label{gcd_ass}\;
\begin{examples}
\item For all $a\in D$ and $u\in D^{\times}$, $\gcd(a,0)\sim a$ \;and\; $\gcd(a,u)\sim 1$.
\item In $\Z[i]$, the elements $1-i$ and $1+i$ are associates, as $1+i=i(1-i)$ and $i$ is a unit.
\item In the power series ring $\Q\pow x$, let \[a=5x^2-x^4+3x^5+\cdots,\quad b=17x^2+x^4+5x^{12}+\cdots.\] Then $\gcd(a,b)\sim x^2$, as both $a$ and $b$ are associates of $x^2$.
\item\label{x_gcd} In $\Q[x]$, $\gcd(x^2+x,x^3+x)=x$. Clearly $x$ is a common factor, so $x$ divides this $\gcd$. Moreover, any polynomial dividing both $x^2+x$ and $x^3+x$ must also divide the linear combination
\[ (1-x)(x^2 + x) + (x^3+x)=2x.\]

This shows that $x\mid\gcd(x^2+x,x^3+x)\mid 2x$, which implies the claim. The above linear combination was obtained using the \emph{extended Euclidean algorithm}, a method explained later in this chapter.
\item\label{no_gcd} In the ring $\Z[\sqrt{-5}]$, the elements $a=6$ and $b=2+2\sqrt{-5}$ have no gcd. To prove this, we begin by finding the full set of common divisors. Suppose that $d=x+y\sqrt{-5}$ is a common divisor of $a$ and $b$. Then $N(d)=x^2+5y^2$ is a positive integer dividing both $N(a)=36$ and $N(b)=24$, hence $N(d)\in\{1,2,3,4,6,12\}$. However, the integers $2,3,12$ are not of the form $x^2+5y^2$, so $N(d)\in\{1,4,6\}$. Solving \[x^2+5y^2=n\]

for $n\in\{1,4,6\}$, we conclude that $d\in\{\pm 1,\pm 2,\pm 1\pm\sqrt{-5}\}$. Clearly, $\pm 1$ and $\pm 2$ are common divisors of $a$ and $b$. The element $1-\sqrt{-5}$ is not, as \[\frac{2+2\sqrt{-5}}{1-\sqrt{-5}}=-\frac{4}{3}+\frac{2}{3}\sqrt{-5}\notin\Z[\sqrt{-5}].\]

Similar calculations show that the full set of common divisors of $a$ and $b$ is \[C=\{\pm 1,\pm 2,\pm(1+\sqrt{-5})\}.\]

Now, a $\gcd$ of $a$ and $b$ would have to be an element of $C$ that is divisible by every element of $C$. However, simple calculations show that no such element exists.
\end{examples}
\end{ex}

\section*{GCD domains}

Having already proved the uniqueness of gcds up to associates, we now turn to the question of existence. A \col{GCD domain}\index{GCD domain} is an integral domain in which every pair of elements has a gcd. The following theorem shows that this class of rings includes all PIDs; in particular, $\Z$ is a GCD domain, as is the polynomial ring $k[x]$ for every field $k$.

\begin{framed}
\begin{thm}[Bézout Identity]\label{bezout_id} Let $D$ be a PID. For all elements $a_1,\ldots, a_n\in D$ there exist $x_1,\ldots, x_n\in D$ such that
\[\gcd(a_1,\ldots,a_n)=a_1x_1+\cdots+a_nx_n.\]

\noindent(In particular, $D$ is a GCD domain.)
\end{thm}
\end{framed}
\begin{proof}
We prove the case $n=2$, the general case following by induction. Let $a,b\in D$. Since every ideal of $D$ is principal, there exists $d\in D$ such that $(d)=(a,b)$. We have \[d=ax+by\]
for some $x,y\in D$ because $d\in (a,b)$; thus, it suffices to show that $d\sim\gcd(a,b)$. Since both $a$ and $b$ belong to $(d)$, then $d\mid a$ and $d\mid b$, hence $d$ is a common divisor of $a$ and $b$. Now, if $g$ is any common divisor, the identity $d=ax+by$ implies that $g\mid d$ since $g\mid a$ and $g\mid b$. This proves that $d$ is a gcd of $a$ and $b$.
\end{proof}

\begin{ex}[GCD Domains]\;
\begin{examples}
\item Let $R=\Z$ or $k[x]$, where $k$ is a field. Then $R$ is a PID by Theorem \ref{Zideals} and hence a GCD domain by Theorem \ref{bezout_id}. By choosing the positive gcd when $R=\Z$, and the monic gcd when $R=k[x]$, we conclude that every pair of elements of $R$ has exactly one gcd in $R$. Moreover, the gcd is a linear combination of the two elements. For instance, the integers 7 and 10 are coprime, and $1=7a+10b$ for $a=3$ and $b=-2$.
\item There are many GCD domains that are not PIDs. Two examples are $\Z[x]$ and $\Q[x,y]$, though proving this requires material from later chapters.
\item The B\'ezout Identity holds in any PID, but not in every GCD domain. In the ring $\Z[x]$ we have $\gcd(2,x)\sim 1$, but 1 is not a linear combination of 2 and $x$. To prove this, we use the homomorphism $\Z[x]\to\F_2[x]$ induced by reduction modulo 2. Suppose there exist $a,b\in\Z[x]$ such that \[2a(x)+xb(x)=1.\] Reduction yields a factorization $x\bar b(x)=1$ in $\F_2[x]$. It follows that $x$ is a unit in $\F_2[x]$, a clear contradiction since units have degree 0.
\item There is no need for a separate concept of \emph{LCM domain}. In a GCD domain, every pair of elements has an lcm---see Exercise \ref{lcm_domain}.
\end{examples}
\end{ex}

To illustrate the utility of gcds, we prove the Rational Root Theorem below, a useful aid in finding roots of polynomials. We begin by proving the following generalization of \href{https://en.wikipedia.org/wiki/Euclid\%27s_lemma}{Euclid's Lemma}, a classical result in number theory.

\begin{framed}
\begin{thm}[Euclid's Lemma]\label{euclidlem} Let $a,b,c$ be elements of a GCD domain. If $a$ and $b$ are coprime and $a\mid bc$, then $a\mid c$.
\end{thm}
\end{framed}
\begin{proof} The statement is trivial if $c=0$, so assume $c\ne 0$. Let $g\sim\gcd(ac,bc)$. We claim that $g$ is associate to $c$. Granting this, we complete the proof as follows: since $a\mid ac$ and $a\mid bc$, then $a\mid \gcd(ac,bc)$, hence $a\mid c$, as required.

To prove the claim, note first that $c$ divides both $ac$ and $bc$, hence divides their gcd, which is $g$. Thus $c\mid g$. Write $g=ct$ for some $t$. We will show that $t$ is a unit by showing that $t\mid a$ and $t\mid b$. 

Write $ac=gu$ for some $u$. Cancellation of $c$ yields:
\[ac=gu=ctu\implies a=tu\implies t\mid a.\]
Similarly, writing $bc=gv$ for some $v$, we show that $t\mid b$. It follows that $t\mid\gcd(a,b)\sim 1$, so $t$ is a unit and $g=ct$ is associate to $c$, as claimed.
\end{proof}
 
\begin{framed}
\begin{thm}[Rational Root Theorem]\label{rat_root}
Let \( D \) be a GCD domain with field of fractions \( F = \Frac(D) \). Let
\[
p(x) = c_0 + c_1x + \cdots + c_nx^n \in D[x]
\]
be a polynomial with \( c_0 \ne 0 \) and \( c_n \ne 0 \). If \( r = a/b \in F \) is a root of \( p(x) \) with \( a, b \in D \) and \( \gcd(a,b)\sim 1\), then $a \mid c_0$ and $b \mid c_n$ in $D$.
\end{thm}
\end{framed}
\begin{proof} We leave it to the reader to show that $\gcd(a,b^n)\sim 1$ for all $n\ge 1$; see Exercise \ref{euclid_lemlem}. Since $r=a/b$ is a root of $p(x)$, the following holds in $D$:
\[0=b^np(a/b)=c_0b^n+c_1ab^{n-1}+\cdots+c_{n-1}a^{n-1}b+c_na^n.\]
Reducing modulo $a$ yields $0\equiv c_0b^n\pmod a$, so $a\mid c_0b^n$. By Euclid's Lemma, this implies $a\mid c_0$. Similarly, reduction modulo $b$ yields $0\equiv c_na^n\pmod b$, hence $b\mid c_na^n$ and thus $b\mid c_n$.
\end{proof}

\begin{ex}[Rational Root Theorem]\;
\begin{examples}
\item We find all rational roots of $p(x):=3x^3-5x^2+5x-2$. For any root $a/b$ in lowest terms, we have $b\mid 3$ and $a\mid 2$. This yields eight options for the root, namely $\pm1,\pm 2,\pm1/3,\pm 2/3$. Evaluating $p(x)$ at these eight values shows that only $2/3$ is a root and hence the unique rational root.
\item We show that the polynomial $f(X):=X^3+(t^2+1)X-t\in\Q[t][X]$ has no root in $\Q(t)$. Note that $f$ has coefficients in $\Q[t]$, which is a GCD domain. Clearly $X=0$ is not a root of $f$. Suppose that $a(t)/b(t)$ is a nonzero root of $f(X)$, where $a,b\in\Q[t]$ are coprime polynomials. Then $a(t)\mid t$ and $b(t)\mid 1$. It follows that $b$ is a nonzero constant and $a$ is either a nonzero constant or is associate to $t$. This shows that the roots of $f(X)$ in $\Q(t)$ are either nonzero constants or have the form $ct$ for some nonzero $c\in\Q$. However, evaluating $f(c)$ for any nonzero $c\in\Q$ yields a quadratic polynomial in $\Q[t]$, hence $c$ is not a root of $f$. Evaluating $f(ct)$ yields a cubic polynomial, hence $ct$ is not a root. Therefore, $f$ has no root in $\Q(t)$.
\end{examples}
\end{ex}

We end our discussion of GCD domains with the following property, which will be used in subsequent sections. A simple application of this property is the parametrization of \emph{Pythagorean triples}, i.e., triples of integers such that $a^2+b^2=c^2$; see \ref{pyth_trip} in Examples \ref{ufd_exs11}.

\begin{framed}
\begin{thm}[$n$th power property]\label{nth_pow}
Let $D$ be a GCD domain. Suppose that $a_1,\ldots, a_k$ are pairwise coprime elements whose product is the $n$th power of an element of $D$:
\[a_1\cdots a_k=x^n,\quad x\in D.\] Then each $a_i$ is associate to an $n$th power in $D$: 
\[a_i=u_iy_i^n,\quad u_i\in D^{\times}, \;y_i\in D.\]
\end{thm}
\end{framed}
\begin{proof}
We prove the case $k=2$; the general case follows by a simple induction. Suppose that $ab = x^n$ with $\gcd(a,b)\sim 1$. We will show that $a$ is associate to an $n$th power; by symmetry, the same will hold for $b$.

Let $g\sim\gcd(b,x)$ and write $b=gc$, $x=gd$ with $\gcd(c,d)\sim 1$. Note that $\gcd(a,g)\sim 1$ since $g\mid b$ and $\gcd(a,b)\sim 1$. Substituting $b=gc$ and $x=gd$ into $ab=x^n$ gives
\[
ac = g^{n-1} d^n.
\]
Since $d^n$ is coprime to $c$ and divides $ac$, Euclid’s Lemma implies $d^n\mid a$.  
Write $a = u d^n$. Substituting this into the above displayed equation yields
\[
u c = g^{n-1},
\]
hence $u\mid g^{n-1}$. Because $u\mid a$ and $\gcd(a,g)\sim 1$, we have $\gcd(u,g^{n-1})\sim 1$, so the relation $u\mid g^{n-1}$ implies that $u$ is a unit. Thus $a = u d^n$ is associate to an $n$th power, as claimed.
\end{proof}

\section*{Euclidean domains}

We now introduce a class of rings in which a division algorithm applies. Naturally, this includes $\Z$ and $k[x]$ for any field $k$, as division algorithms for these rings were discussed previously.

A \col{Euclidean function}\index{Euclidean function} on a domain $D$ is a function $f:D\setminus\{0\}\to\N$ with the following property: for all $a,b\in D$ with $b\ne 0$, there exist $q,r\in D$ such that \[a=qb+r\quad\text{with}\;\;r=0\;\;\text{or}\;\; f(r)<f(b).\]

A domain admitting a Euclidean function is called a \col{Euclidean domain}\index{Euclidean domain}. If $f$ is a Euclidean function on $D$, we may also refer to the pair $(D,f)$ as a Euclidean domain.

\begin{ex}[Euclidean domains]\;
\begin{examples}
\item The usual absolute value is a Euclidean function on $\Z$. (This follows from the division algorithm for $\Z$.) Under a certain notion of equivalence (see Exercise \ref{euc_equiv}), the absolute value is in fact the \emph{only} Euclidean function on $\Z$.
\item If $k$ is a field, the degree function on the polynomial ring $k[x]$ is a Euclidean function, by Theorem \ref{div_algo}.
\item Every field is a Euclidean domain in a rather trivial way: if $k$ is a field, \emph{every} function $f:k\setminus\{0\}\to\N$ is Euclidean. Indeed, given $a,b\in k$ with $b\ne 0$, we have $a=(a/b)b+0$. 
\item The norm function on $\Z[i]$, given by $N(a+bi)=a^2+b^2$, is Euclidean; see Exercise \ref{GaEuc}. It is far from true, however, that the norm function on \emph{every} quadratic integer ring $\Z[\sqrt d]$ is Euclidean, the most immediate reason being that $N(a+b\sqrt d)=a^2-b^2d$ need not take values in $\N$. But more significantly, even the function $|a^2-b^2d|$ is generally not Euclidean; the example of $\Z[\sqrt{14}]$ is discussed in Exercise \ref{quadEuc}. See also Exercise \ref{non_ufd_quad}, which shows that $\Z[\sqrt d]$ is not Euclidean if $d\equiv 1\pmod 4$.
\item In general, the elements $q$ and $r$ in the definition of a Euclidean function are not unique. For example, in $\Z$ we have the following valid Euclidean divisions of 10 by 4: \[10=2\cdot 4+2\quad\text{and}\quad 10=3\cdot 4+(-2).\]
Thus, additional restrictions need to be placed in order for quotients and remainders to be uniquely determined. In the case of $\Z$, this is done by requiring the remainder $r$ to satisfy $0 \le r < |b|$; this condition uniquely determines both $q$ and $r$, with $q$ being the integer part (floor) of $a/b$, and $r = a - qb$.
\end{examples}
\end{ex}

The following theorem is our main result concerning Euclidean domains.

\begin{framed}
\begin{thm}\label{eucpid}
Every Euclidean domain is a PID.
\end{thm}
\end{framed}
\begin{proof}
Let $(D,f)$ be a Euclidean domain and let $I$ be an ideal of $D$. If $I=(0)$, clearly $I$ is a principal ideal, so suppose $I$ is nonzero. Then the set $f(I\setminus\{0\})$ is a nonempty subset of $\N$ and therefore has a least element; let $x\in I\setminus\{0\}$ satisfy $f(x)\le f(y)$ for all nonzero $y\in I$. We claim that $I=(x)$, making $I$ a principal ideal. To prove the claim, we show that $I\subseteq (x)$. Let $a\in I$, and divide $a$ by $x$: there exist $q,r\in D$ such that $a=qx+r$ with $r=0$ or $f(r)<f(x)$. We must have $r=0$, since otherwise the inequality $f(r)<f(x)$ contradicts the minimality of $x$. Thus $a=qx\in (x)$, proving that $I\subseteq (x)$.
\end{proof}

\begin{ex}[Euclidean domains]\label{euc_exs}\;
\begin{examples}
\item If $k$ is a field, the polynomial ring $k[x]$ is a Euclidean domain, by Theorem \ref{div_algo}, therefore a PID and a GCD domain. By convention, the gcd of two polynomials is taken to be the unique monic gcd. By the B\'ezout Identity, the gcd of two polynomials can be expressed as a linear combination of them. An effective method for doing this in practice is detailed in the next section.
\item Using Theorem \ref{eucpid}, we can give non-examples of Euclidean domains. For instance, the ring $\Z[x]$ cannot be Euclidean because it is not a PID (see Exercise \ref{non_pid}, Chapter \ref{ideals_chap}). A much more difficult proposition is finding a PID that is not Euclidean, as this requires ruling out the existence of \emph{any} Euclidean function on some ring. The earliest published example of a non-Euclidean PID is the ring $\Z[(1+\sqrt{19})/2]$, due to Motzkin \cite{Motzkin}.
\end{examples}
\end{ex}

\subsection*{The Euclidean Algorithm} Let $D$ be a Euclidean domain. Since $D$ is a PID by Theorem \ref{eucpid}, the B\'ezout Identity (\ref{bezout_id}) implies that $D$ is a GCD domain, so every pair of elements of $D$ has a gcd; moreover, any gcd of elements $a$ and $b$ can be expressed in the form $ax+by$ for some $x,y\in D$. We will now describe the \col{Euclidean algorithm}\index{Euclidean algorithm}, a method for computing $\gcd(a,b)$, and the \emph{extended} Euclidean algorithm, which computes $\gcd(a,b)$ as well as elements $x,y$ such that $\gcd(a,b)=ax+by$.

The basis for both algorithms is the observation that if $a=qb+r$, then $\gcd(a,b)\sim\gcd(b,r)$, a fact that follows directly from the definition of a gcd. Using this, we can describe $\gcd(a,b)$ as the end result of a process detailed below:

\begin{framed}
\begin{thm}[Euclidean Algorithm]\label{euc_algo}
Let $D$ be a Euclidean domain and $a,b\in D$ with $b\ne 0$. Define a sequence \[a,b,r_1,r_2,r_3,\ldots\] by setting \[r_{-1}=a,\;\; r_0=b,\] and for $i\ge 0$, defining $r_{i+1}$ as follows: if $r_i=0$, then $r_{i+1}=0$; otherwise, $r_{i+1}$ is a remainder obtained upon dividing $r_{i-1}$ by $r_i$.

\noindent Then the following hold: 
\begin{enumerate}
\item[(a)] The terms $r_i$ are eventually 0.
\item[(b)] Letting $n=\inf\{i\ge 1:r_i=0\}$, we have $\gcd(a,b)=r_{n-1}$.

Thus, $\gcd(a,b)$ is the last nonzero term of the sequence $r_1,r_2,r_3,\ldots$.
\end{enumerate}
\end{thm}
\end{framed}
\begin{proof}
Exercise \ref{euc_algo_proof}.
\end{proof}

\begin{ex}[Euclidean Algorithm]\;
\begin{examples}
\item We compute $\gcd(252,198)$. The terms of the sequence $r_1,r_2,r_3,\ldots$ from Theorem \ref{euc_algo} are the integers boxed in the following calculations:
\begin{align*}
252 &= 1 \cdot 198 + \boxed{54} \\
198 &= 3 \cdot 54 + \boxed{36} \\
54 &= 1 \cdot 36 + \boxed{18} \\
36 &= 2 \cdot 18 + \boxed{0}
\end{align*}

Since 18 is the last nonzero remainder in the sequence, we conclude that \[\gcd(252,198)=18.\] 

\item Continuing with the previous example, we now illustrate the \col{extended Euclidean algorithm}\index{extended Euclidean algorithm} by finding integers $x$ and $y$ such that \[252x+198y=18.\]

In the notation of Theorem \ref{euc_algo}, the technique consists of expressing every term of the sequence $r_1,r_2,r_3,\ldots$ as a linear combination of $a$ and $b$. After $n-1$ steps, this expresses $\gcd(a,b)=r_{n-1}$ as such a linear combination.

In the example of $\gcd(252,198)$, let us set $a=252$ and $b=198$. The above equations with boxed integers lead to the following expressions of $r_1=54$, $r_2=36$, and $r_3=18$ as linear combinations of $a$ and $b$:
\begin{align*}
54 &=a-b, \\
36 &= b-3\cdot 54=b-3\cdot (a-b)=4b-3a,\\
18 &= 54-36=(a-b)-(4b-3a)=4a-5b.\\
\end{align*}
In particular, we obtain the desired relation:
\[18=\gcd(252,198)=252\cdot 4-198\cdot 5.\]

\item One common use of the extended Euclidean algorithm is for computing modular inverses. As a simple example, we can compute the inverse of 3 in $\F_{19}$ by expressing $1=\gcd(3,19)$ as a linear combination of 3 and 19: \[1=19-3\cdot 6.\]

Reducing modulo 19, this shows that \[1=3\cdot 13\;\;\text{in}\;\;\F_{19},\] hence $3^{-1}=13$ in $\F_{19}$.
\end{examples}
\end{ex}

\newpage

\section*{Exercises}

\begin{note}\;
\begin{itemize}
\item Solutions to starred exercises (\sel) are included in Chapter~\ref{review3chap}.
\item A house (\house) marks exercises suggested as sample homework.
\end{itemize}
\end{note}

\medskip

\begin{exercises}
\item Determine the values of $n\ge 1$ for which $\Z_n$ is a Euclidean domain.
\item Find all associates of $2+3i$ in $\Z[i]$ and in $\C$.
\item Let $k$ be a field. Show that the order of vanishing is a Euclidean function on the power series ring $k\pow x$.
\item Find a gcd of $a=1+3i$ and $b=1-3i$ in $\Z[i]$, and express $\gcd(a,b)$ as a linear combination of $a$ and $b$.
\item Find the (monic) $\gcd$ of $a=3x^2 + 2x + 1$ and $b=x^2 + 1$ in $\Q[x]$, and express $\gcd(a,b)$ as a linear combination of $a$ and $b$.
\item\label{gcd_ass_exc} Prove Theorem \ref{gcd_ass_thm}.
\item(\house)\label{lcm_domain} Show that every pair of elements $a,b$ of a GCD domain has an lcm, and that if $a$ and $b$ are not both 0, then \[\lcm(a,b)\sim \frac{ab}{\gcd(a,b)}.\]
\item\label{finite_gcd} Show that every nonempty finite subset of a GCD domain has a gcd in the domain.
\item\label{scale_gcd} Let $D$ be a GCD domain and $s,a_1,\dots,a_n\in D$.
Show that \[\gcd(sa_1,\dots,sa_n)\sim s\cdot\gcd(a_1,\dots,a_n).\]
\item Let $R$ be a GCD domain and let $a,b\in R$. Show that $a$ and $b$ are coprime if and only if 1 is a linear combination of $a$ and $b$.
\item Let $F=\Frac(R)$, where $R$ is a GCD domain. Show that every element of $F$ can be expressed in lowest terms, i.e., in the form $a/b$ with $a,b\in R$ and $\gcd(a,b)\sim 1$. Moreover, show that such a representation is unique up to multiplication of $a$ and $b$ by the same unit in $R$.
\item Prove that if an integer is the sum of two rational squares, then it is a sum of two integer squares. For example, \[13=(1/5)^2+(18/5)^2=2^2+3^2.\]
\item Use the Rational Root Theorem to find all roots of the given polynomial $p(x)$ in the field of fractions of the given domain $D$.
\begin{enumerate}
\item $p(x)=x^3-(1+\sqrt 2)x^2+(\sqrt 2-1)x-2$,\;\; $D=\Z[\sqrt 2]$
\item $p(x)=(1+i)x^3-(2+i)x^2+(3-i)x-(1+i)$,\;\; $D=\Z[i]$.
\end{enumerate}
\item Suppose $r$ is a real number such that $r+1/r$ is an odd integer. Show that $r$ must be irrational.
\item Let $a,b$ be elements of a GCD domain. Show that $(a)\cap(b)=(\lcm(a,b))$.
\item(\house)\label{euclid_lemlem} Let $a,b$ be coprime elements of a GCD domain. Show that $a$ and $b^n$ are coprime for every positive integer $n$.
\item Let $(D,f)$ be a Euclidean domain. If $a,b\in D$ are associates, show that $f(a)=f(b)$.
\item\label{euc_algo_proof} Prove Theorem \ref{euc_algo}.
\item Show that the mapping $x\mapsto x^2$ is a Euclidean function on $\Z$.
\item\label{euc_equiv} Let $D$ be a domain. We define a relation on the set of all Euclidean functions on $D$ by $f\approx g$ if there exists a strictly increasing function $\phi:\N\to\N$ such that $f(x)=\phi(g(x))$ for all $x\in D\setminus\{0\}$. 
\begin{enumerate}
\item Show that $\approx$ is an equivalence relation.
\item Show that any two Euclidean functions on $\Z$ are equivalent.
\end{enumerate}
\item(\house)\label{quadEuc} Show that the mapping $a+b\sqrt{14}\mapsto |a^2-14b^2|$ is not a Euclidean function on $\Z[\sqrt{14}]$. (Attempt a Euclidean division of 1 by $\sqrt{14}$.)
\item\label{GaEuc} Recall that the norm on $\Z[i]$ is given by $N(a+bi)=a^2+b^2$.
\begin{enumerate}
\item Show that $N$ is a Euclidean function. (Given $\alpha,\beta\in\Z[i]$, consider $x=\alpha/\beta \in \Q(i)$ and round the real and imaginary parts of $x$ to the nearest integers.)
\item Find $q$ and $r$ in $\Z[i]$ such that $3-4i=(2+5i)q+r$ and $N(r)<N(2+5i)$. 
\end{enumerate}
\item Let $A$ be a domain and $R$ a subring of $A$ that is a PID. Show that if $x,y\in R$ have a gcd in $R$, then they have a gcd in $A$.
\item\label{bez_dom} A \col{B\'ezout domain}\index{B\'ezout domain} is a domain in which every finitely-generated ideal is principal. Show that every PID is a B\'ezout domain, and every B\'ezout domain is a GCD domain.
\item(\house)\label{non_ufd_quad} Let $A$ be an integral domain and $F=\Frac(A)$. We say that $A$ is \col{integrally closed}\index{integrally closed} if the following holds: if $p(x)\in A[x]$ is a monic polynomial and $r\in F$ is a root of $p(x)$, then $r\in A$.
\begin{enumerate}
\item Show that every GCD domain is integrally closed.
\item If $d\ne 1$ is a squarefree integer with $d\equiv 1\pmod 4$, show that the ring $\Z[\sqrt d]$ is not integrally closed. Conclude that $\Z[\sqrt d]$ is not a Euclidean domain.
\end{enumerate}
\end{exercises}

\chapter{Unique Factorization Domains}\label{ufd_chap}

This chapter introduces the concept of a UFD, an integral domain in which an analogue of the Fundamental Theorem of Arithmetic holds: every nonzero non-unit element can be expressed as a product of prime factors in an essentially unique way. The theory of UFDs provides a unified framework for studying  factorization in rings such as $\Z$, the Gaussian integers $\Z[i]$, and the polynomial ring $k[x_1,\ldots, x_n]$ over any field $k$.

A central theme will be the relationship between irreducible and prime elements. In general domains these notions need not coincide, but in a UFD they do, and this equivalence lies at the heart of unique factorization.

The theory has striking applications to Diophantine equations. Unique factorization was historically used to settle several special cases of \href{https://en.wikipedia.org/wiki/Fermat\%27s_Last_Theorem}{Fermat's Last Theorem} (the claim that $x^n+y^n=z^n$ has no nontrivial integer solutions for $n>2$). The idea goes back to \href{https://en.wikipedia.org/wiki/Gabriel_Lam\%C3\%A9}{Gabriel Lam\'e} (1795--1870), who gave an argument that would prove Fermat’s theorem for prime exponents $n$ under the assumption that certain rings admit unique factorization. Although this assumption turns out to be false for most primes, it does hold for the exponents $n\in\{3,5,7,11,13,17,19\}$, so Lam\'e’s method is valid in these cases.

\section*{Primes and irreducibles}
 
Let $D$ be a domain. An \col{irreducible}\index{irreducible element} element of $D$ is a nonzero non-unit element $\pi\in D$ satisfying \[\pi=ab\;\;\text{implies}\;\;a\in D^{\times}\;\;\text{or}\;\;b\in D^{\times}.\]
Equivalently, every divisor of $\pi$ is a unit or an associate of $\pi$.
 
 A \col{reducible}\index{reducible element} element of $D$ is a nonzero non-unit element $f\in D$ that is expressible as a product of two non-units.
 
 A \col{prime}\index{prime element} element of $D$ is a nonzero non-unit element $p\in D$ satisfying \[p\mid ab\;\;\text{implies}\;\;p\mid a\;\;\text{or}\;\; p\mid b.\]
 
 Note that 0 and $1$ are not classified as irreducible, reducible, or prime.
 
 The following theorem clarifies the relation between prime and irreducible elements, in particular showing that every prime is irreducible, and that the two concepts agree in GCD domains.

\begin{framed}
\begin{thm}\label{prime_irred}
Let $D$ be an integral domain.
\begin{enumerate}
\item[(a)] Every associate of an irreducible (respectively, prime) element of $D$ is irreducible (respectively, prime).
\item[(b)] Every prime element of $D$ is irreducible.
\item[(c)] If $D$ is a GCD domain, then every irreducible element is prime.
\end{enumerate}
\end{thm}
\end{framed}
\begin{proof}
Part (a) follows immediately from the definitions. For (b), suppose that $p\in D$ is a prime element and that $p=ab$ with $a,b\in D$. Since $p\mid ab$, we may assume that $p\mid a$, say $a=pc$. Then $p=pcb$. By the cancellation property of domains, this implies that $cb=1$, so $b\in D^{\times}$. This proves that $p$ is irreducible.

Suppose now that $D$ is a GCD domain and let $f\in D$ be irreducible. To show that $f$ is prime, suppose that $f\mid xy$ for some $x,y\in D$. Since $f$ is irreducible, every divisor of $f$ is a unit or an associate of $f$; in particular, $\gcd(f,x)$ is either a unit or is associate to $f$. In the latter case, clearly $f\mid x$. In the former case, Euclid's Lemma (Theorem \ref{euclidlem}) implies $f\mid y$.
\end{proof}

\begin{ex}[Prime and irreducible elements]\label{prime_irr_ex}\;
\begin{examples}
\item In $\Z$, the irreducible elements are precisely the prime integers (including negative primes). This follows from Theorem \ref{prime_irred}.
\item A field has no reducible or irreducible elements, as every nonzero element is a unit. All nonzero elements are associates.
\item Let $F$ be a field. A polynomial $f\in F[x]$ is reducible if and only if $f=gh$ for some nonconstant polynomials $g$ and $h$. The associates of $f$ are the polynomials $cf(x)$ for $c\in F\setminus\{0\}$. Since $F[x]$ is a PID, the prime polynomials are precisely the irreducible polynomials.
\item In $\Z[x]$, a constant polynomial can be irreducible: if $p$ is a prime integer, then $p$ is irreducible in $\Z[x]$. Contrast this with the previous example, in which a constant polynomial is either 0 or a unit, hence is neither reducible nor irreducible.
\item Irreducible elements of a domain need not be prime. For example, $2\in\Z[\sqrt{-5}]$ is not prime because it divides \[(1+\sqrt{-5})(1-\sqrt{-5})=6\] but does not divide either factor (for instance, if \[2(x+y\sqrt{-5})=(1+\sqrt{-5})\] with $x,y\in\Z$, then $2x=2y=1$, a contradiction.) Thus 2 is not prime. However, 2 is irreducible: it is a nonzero non-unit, and if $2=ab$ with $a,b$ non-units, then $4=N(2)=N(a)N(b)$ forces $N(a)=N(b)=2$. But no element of $\Z[\sqrt{-5}]$ has norm 2, as the equation $x^2+5y^2=2$ has no integer solutions.
\end{examples}
\end{ex}

The next two theorems provide simple criteria for proving irreducibility of elements in two types of rings: quadratic integer rings and polynomial rings over fields.

\begin{framed}
\begin{thm}[Prime Norm Criterion]\label{prime_norm}
Let $R=\Z[\sqrt d]$ be a quadratic integer ring, and $x\in R$. If $N(x)$ is a prime integer, then $x$ is irreducible.
\end{thm}
\end{framed}

\begin{proof} Recall from Theorem \ref{quad_int_unit} that an element of $R$ is a unit if and only if its norm is $\pm 1$, a unit in $\Z$. Let $x\in R$ have prime norm. By Theorem \ref{prime_irred}, $N(x)$ is irreducible in $\Z$. It follows that $x$ cannot be zero or a unit, for otherwise $N(x)$ would be zero or $\pm 1$. Now assume, by contradiction, that $x$ is reducible, say $x=ab$ with $a,b$ non-units. Then $N(x)=N(a)N(b)$ is a product of nonzero non-unit integers, contradicting irreducibility of $N(x)$. Hence $x$ must be irreducible.
\end{proof}

\begin{note} With more work, it is possible to prove a statement stronger than the above theorem: if $N(x)$ is prime in $\Z$, then $x$ is prime in $R$; see Exercise \ref{prime_quad}.
\end{note}

\begin{framed}
\begin{thm}[Low-degree Test]\label{lowdeg}
If $F$ is a field and $p(x)\in F[x]$ is a polynomial of degree two or three, then $p(x)$ is reducible if and only if $p(x)$ has a root in $F$.
\end{thm}
\end{framed}
\begin{proof}
If $p(x)$ has a root $a\in F$, then $p(x)$ is divisible by $x-a$ and is thus reducible. Conversely, if $p(x)$ is reducible, a degree argument implies that $p(x)$ has a factor of degree one in $F[x]$ and therefore has a root in $F$.
\end{proof} 

\section*{UFDs}

A \col{unique factorization domain}\index{unique factorization domain} (UFD) is a domain $D$ satisfying:
\begin{itemize}
\item[(a)] Every nonzero non-unit element $f\in D$ can be factored as \[f=p_1\cdots p_n,\] where each $p_i$ is irreducible.
\item[(b)] With $f$ as above, if \[f=q_1\cdots q_m\] is any other factorization into irreducibles, then there is a bijection \[\sigma:\{1,\ldots,n\}\to\{1,\ldots, m\}\] such that $p_i$ is associate to $q_{\sigma(i)}$ for $1\le i\le n$. (In particular, $m=n$.)
\end{itemize}

In short, elements of a UFD have a unique factorization into irreducible elements, up to reordering of factors and multiplication by units.

\begin{ex}[Factorization into irreducibles]\;
\begin{examples}
\item In $\Z$, the integer $n=-294$ can be factored as $2\cdot 3\cdot 7\cdot(-7)$. Every other factorization of $n$ as a product of primes is obtained by writing the above factors in a different order or with a different sign, for instance, $n=(-2)\cdot 7\cdot 3\cdot 7$.
\item In $\Z[i]$, the element $n=13$ factors as $n=(2+3i)(2-3i)$. Each factor in this product is irreducible since its norm is $2^2+3^2=13$, a prime integer. (Here, we are using Theorem \ref{prime_norm}.) Other factorizations of $n$ into irreducibles are obtained by scaling these factors by a unit (namely, $\pm 1$ or $\pm i$) and/or reordering the factors; for instance, $n=-(3-2i)(2i+3)$. As follows from Theorem \ref{pid_is_ufd} below, $\Z[i]$ is a UFD, hence every factorization of $n$ into irreducibles is obtained in this way.
\item The domain $\Z[\sqrt{-6}]$ is not a UFD, as $n=10$ is a nonzero non-unit which factors in two ways not obtainable from each via reordering and scaling by units: \[10=2\cdot 5=(2+\sqrt{-6})(2-\sqrt{-6}).\] Irreducibility of all factors involved can be proved by assuming a factorization into non-units and using norms to reach a contradiction. Further, 2 and 5 are not associate to $2\pm\sqrt{-6}$ because they have different norms. Hence, we have an example of non-uniqueness of factorization into irreducibles.
\end{examples}
\end{ex}

Theorem \ref{gcd_form} below establishes the existence of gcds and lcms, a key property of UFDs.  To motivate the type of factorization assumed in the theorem, consider the following example in $\Z$: the integers $a=-294$ and $b=120$ can be factored as \[a=-2\cdot 3\cdot 7^2,\quad b=2^3\cdot 3\cdot 5,\]

so they can be expressed as products of powers of a common set of primes: \[a=-2\cdot 3\cdot 5^0\cdot 7^2,\quad b=2^3\cdot 3\cdot 5\cdot 7^0.\]

\newpage

\begin{framed}
\begin{thm}\label{gcd_form} Let $D$ be a UFD. Every pair of elements of $D$ has a gcd and an lcm in $D$. Moreover, if $a=up_1^{e_1}p_2^{e_2}\cdots p_n^{e_n}$ and $b=vp_1^{f_1}p_2^{f_2}\cdots p_n^{f_n}$, where $u,v\in D^{\times}$, the $p_i$'s are pairwise non-associate irreducibles in $D$, and $e_i,f_i\ge 0$ for each $i$, then
\[\gcd(a,b)\sim\prod_{i=1}^np_i^{\min(e_i,f_i)}\quad\text{and}\quad\lcm(a,b)\sim\prod_{i=1}^np_i^{\max(e_i,f_i)}.\]
(In particular, $D$ is a GCD domain.)
\end{thm}
\end{framed}
\begin{proof} Exercise \ref{gcd_form_pf}.\end{proof}

The following theorem provides a criterion for showing that a particular domain is a UFD. Applying the criterion, we will be able to prove that every PID is a UFD, thus providing a large collection of examples of UFDs.

\begin{framed}
\begin{thm}\label{ufd_crit}
A domain $D$ is a UFD if and only if every irreducible element of $D$ is prime, and every ascending chain of principal ideals in $D$ is eventually constant.
\end{thm}
\end{framed}
\begin{proof}
We prove the ``if" direction of the theorem, leaving the other direction as an exercise. We begin by noting that the chain condition in the theorem implies that every nonempty set of principal ideals in $D$ has a maximal element under inclusion, i.e., an ideal $(f)$ for which there is no strictly larger ideal in the set.

Define $S$ to be the set of all principal ideals $(x)\subseteq D$ such that $x$ is a nonzero non-unit that cannot be expressed as a product of irreducibles. We prove by contradiction that $S=\emptyset$, thus proving that every nonzero non-unit is a product of irreducibles.

Assuming $S$ is nonempty, let $(f)$ be a maximal element. Since $f$ is not irreducible (otherwise it is trivially a product of irreducibles), there exist nonzero non-units $a,b\in D$ such that $f=ab$. At least one of $a$ or $b$ must fail to be a product of irreducibles; otherwise $f$ would be such a product. Without loss of generality, assume $(a)\in S$. Now, since $f=ab\in (a)$ and $b$ is not a unit, the ideal $(f)$ is properly contained in $(a)$, contradicting the maximality of $(f)$. This proves that $S=\emptyset$, as required.

It remains to prove uniqueness of factorizations into irreducibles. We prove the following statement by induction on $n\ge 1$: if $p_1,\ldots, p_n$ are irreducibles and \[p_1\cdots p_n=q_1\cdots q_m,\] where each $q_i$ is irreducible, then there is a bijection \[\sigma:\{1,\ldots,n\}\to\{1,\ldots, m\}\] such that $p_i$ is associate to $q_{\sigma(i)}$ for each $i$.

The base case $n=1$ follows immediately from the definition of an irreducible element. Assuming now that $n\ge 2$ and that the conclusion holds true for $n-1$ irreducibles, suppose that $p_1\cdots p_n=q_1\cdots q_m$, with all elements involved being irreducible. Note that $m>1$, as otherwise $q_1$ is reducible. Since $p_n\mid q_1\cdots q_m$ and $p_n$ is prime (by hypothesis), then $p_n\mid q_j$ for some $j$. Since $q_j$ is irreducible, this implies that $p_n$ is associate to $q_j$, say $up_n=q_j$ with $u\in D^{\times}$. By cancellation, we obtain
\[p_1\cdots p_{n-1}=(uq_1)q_2\cdots\hat q_j\cdots q_m,\] where the hat over $q_j$ indicates that $q_j$ is omitted from the product.

By induction, there is a bijection $s:\{1,\ldots, n-1\}\to\{1,\ldots, m\}\setminus\{j\}$ such that $p_i$ is associate to $q_{s(i)}$ for $1\le i<n$. Defining $\sigma(k)=s(k)$ for $k<n$ and $\sigma(n)=j$, we obtain a bijection $\sigma:\{1,\ldots,n\}\to\{1,\ldots, m\}$ such that $p_i$ is associate to $q_{\sigma(i)}$ for $1\le i\le n$.
\end{proof}

Applying Theorem \ref{ufd_crit} we now prove the main result of this chapter.

\begin{framed}
\begin{thm}\label{pid_is_ufd}
Every PID (and every Euclidean domain) is a UFD.
\end{thm}
\end{framed}
\begin{proof}
Let $D$ be a PID. We will apply Theorem \ref{ufd_crit} to show that $D$ is a UFD. By Theorem \ref{prime_irred}, every irreducible element of $D$ is prime (recall that PIDs are GCD domains, by Theorem \ref{bezout_id}). Now let \[(a_1)\subseteq(a_2)\subseteq(a_3)\subseteq\cdots\]
be a chain of principal ideals in $D$. Let $I=\cup_{i=1}^{\infty}(a_i)$. Then $I$ is an ideal of $D$ (see Exercise \ref{asc_chain}, Chapter \ref{ideals_chap}), hence $I=(f)$ for some $f\in D$. Since $f\in I$, there exists $N$ such that $f\in (a_N)$. It follows that $(a_i)=(a_N)$ for all $i\ge N$, as $f\in(a_i)$ and $(a_i)\subseteq (f)$. By Theorem \ref{ufd_crit}, this proves that $D$ is a UFD.
\end{proof}

\begin{ex}[UFDs]\;
\begin{examples}\label{ufd_exs11}
\item Since $\Z$ is a Euclidean domain, it is a UFD. This implies the Fundamental Theorem of Arithmetic: every integer $n\notin\{0,\pm 1\}$ has a unique factorization of the form \[n=\pm p_1^{e_1}\cdots p_r^{e_r},\] where the $p_i$'s are primes with $0<p_1<p_2<\cdots< p_r$, and $e_i\ge 1$ for every $i$.
\item If $k$ is a field, the polynomial ring $k[x]$ is a Euclidean domain and therefore a UFD. Hence, every nonconstant polynomial $f(x)\in k[x]$ has a factorization of the form \[f(x)=c\cdot p_1(x)^{e_1}\cdots p_r(x)^{e_r},\] where the $p_i$'s are distinct \emph{monic} irreducible polynomials; moreover, such a factorization is unique up to reordering of the factors $p_i(x)^{e_i}$. For purposes of computing the $\gcd$ of two polynomials, the factorization into irreducibles should be used (as in Theorem \ref{gcd_form}) only if it is already available. If the factorizations of $f(x)$ and $g(x)$ are not both known, then the Euclidean algorithm is arguably the best tool to compute $\gcd(f,g)$, and to express it as a linear combination of $f$ and $g$ if necessary.
\item\label{pyth_trip} A \col{Pythagorean triple}\index{Pythagorean triple} is a triple of positive integers $a,b,c$ satisfying $a^2+b^2=c^2$. For instance, $(3,4,5)$ and $(5,12,13)$ are such triples. We will parametrize \emph{primitive} Pythagorean triples, meaning triples that are pairwise coprime, by using the fact that $\Z[i]$ is a UFD, hence a GCD domain, and therefore has the $n$th power property (Theorem \ref{nth_pow}).

Let $(a,b,c)$ be a primitive Pythagorean triple. In $\Z[i]$, we have the factorization
\[c^2=a^2+b^2=(a+bi)(a-bi).\] We claim that $a+bi$ and $a-bi$ are coprime in $\Z[i]$. The proof uses basic properties of the Gaussian prime $1+i$; see Exercise \ref{conj_copr}.

The $n$th power property now implies that $a+bi$ is associate to a square in $\Z[i]$, say $a+bi=u(m+ni)^2$, $u\in\{\pm 1,\pm i\}$. Conjugation yields $a-bi=\bar u(m-ni)^2$. Hence,
\begin{align*}a+bi&=u(m^2-n^2)+2iu(mn),\\
a-bi&=\bar u(m^2-n^2)-2i\bar u(mn).
\end{align*}
Consideration of the real and imaginary parts above rules out the possibilities $u=\pm i$. Hence $u=\pm 1$.

Setting $u=1$ leads to the parametrization \[a=m^2-n^2,\;\;b=2mn,\;\;c=m^2+n^2,\quad m,n\in\Z.\] Setting $u=-1$ leads to the same parametrization upon applying the transformation $(m,n)\mapsto(n,-m)$. 

The above description allows us to easily generate Pythagorean triples: taking $(m,n)=(3,2)$ we obtain the known triple $(5,12,13)$, and taking $(m,n)=(12,5)$ we arrive at the new triple $(119,120,169)$.

\end{examples}
\end{ex}

\subsection*{Classes of domains}

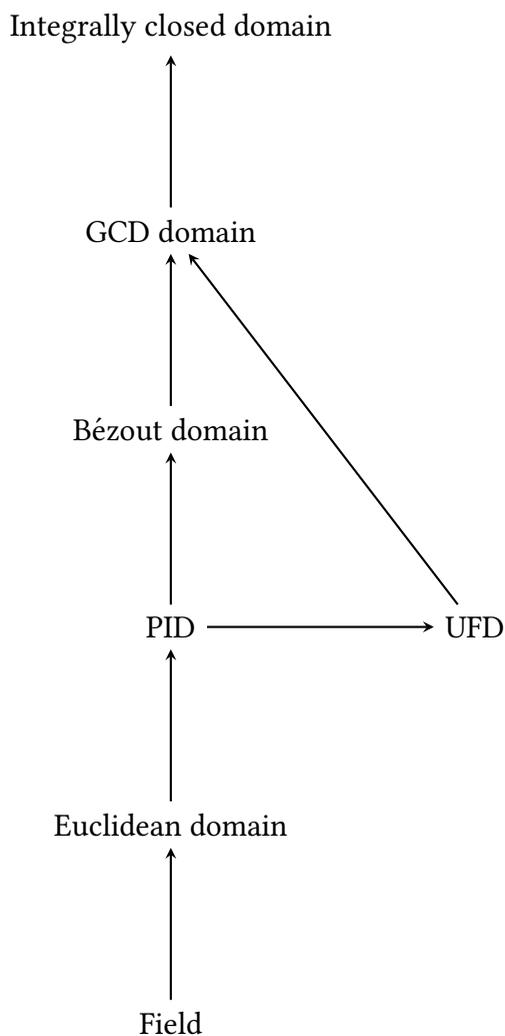
\begin{figure}[h!]
\centering
\begin{tikzpicture}[node distance=2cm, >=stealth, thick]

  \node (intcl) {Integrally closed domain};
  \node (gcd)   [below=of intcl] {GCD domain};
  \node (bez)   [below=of gcd]   {B\'ezout domain};
  \node (pid)   [below=of bez]   {PID};
  \node (euc)   [below=of pid]   {Euclidean domain};
  \node (field) [below=of euc]   {Field};

  \node (ufd) [right=3cm of pid] {UFD};

  \draw[->] (field) -- (euc);
  \draw[->] (euc)   -- (pid);
  \draw[->] (pid)   -- (bez);
  \draw[->] (bez)   -- (gcd);
  \draw[->] (gcd)   -- (intcl);

  \draw[->] (pid) -- (ufd);
  \draw[->] (ufd) -- (gcd);

\end{tikzpicture}

\bigskip
\caption{Relations among important classes of integral domains.}
\label{fig:domain_inclusions}
\end{figure}

We conclude this chapter with some remarks on various classes of domains that have appeared throughout this book. Figure \ref{fig:domain_inclusions} summarizes their containment relations. For the definitions of B\'ezout domains and integrally closed domains, see Exercises \ref{bez_dom} and \ref{non_ufd_quad} in Chapter \ref{euc_chap}.

Although the figure suggests a tidy hierarchy, it represents only a very small portion of the universe of rings. There are countless integral domains that lie outside all of the classes shown, and even within each class there is a rich supply of examples exhibiting rather different behavior. Nevertheless, the study of these classes has already yielded several powerful tools---Euclid’s Lemma, the existence of gcds, the $n$th power property, and unique factorization---and will continue to play a central role in the chapters ahead.

\newpage

\section*{Exercises}

\begin{note}\;
\begin{itemize}
\item Solutions to starred exercises (\sel) are included in Chapter~\ref{review3chap}.
\item A house (\house) marks exercises suggested as sample homework.
\end{itemize}
\end{note}

\medskip

\begin{exercises}
\item Show that in the ring $\Z[\sqrt{-3}]$, the elements 5 and $2+3\sqrt{-3}$ are irreducible, while 7 is reducible.
\item Find all monic irreducible polynomials in $\F_3[x]$ having degree $\le 3$.
\item Find all irreducible polynomials of degrees 3 and 4 in $\F_2[x]$.
\item(\sel)\label{chap17_sel1} Find all units of the ring $R=\Z[\sqrt{-5}]$. Is 3 a prime element of $R$? Is $1+\sqrt{-5}$ an irreducible element of $R$?
\item(\house) Let $R$ be a PID and $x\in R$ a nonzero element. Show that $x$ is irreducible $\iff$ $(x)$ is a maximal ideal of $R$ $\iff$ $(x)$ is a prime ideal of $R$.
\item Factor the polynomials $f=x^5+x^3+6x^2+6$ and $g=x^4 + x^3 + x^2 + x + 6$ into irreducibles in $\F_7[x]$. Compute $\gcd(f,g)$ in two ways; (a) by applying Theorem \ref{gcd_form}, and (b) by using the Euclidean algorithm.
\item Let $p,q$ be irreducible elements of a UFD, and let $e,f$ be positive integers. Show that $p^e$ is associate to $q^f$ if and only if $p$ and $q$ are associates.
\item Show that the polynomial $x^2-\sqrt 2$ is irreducible in $\Z[\sqrt 2][x]$.
\item Prove or provide a counterexample: every subring of a UFD is a UFD. Repeat the exercise for PIDs and Euclidean domains.
\item(\sel)\label{prime_id_el} Let $D$ be a domain and $x\in D$ a nonzero non-unit element. Show that $x$ is prime if and only if the ideal $(x)$ is a prime ideal.
\item\label{reverse_poly} Let $k$ be a field and $f\in k[x]$ a polynomial of degree $n$. The polynomial $g(x)=x^nf(1/x)$ is called the \col{reverse}\index{reverse of a polynomial} of $f$.
\begin{enumerate}
\item Describe the coefficients of $g$ in terms of those of $f$.
\item Show that $g$ is irreducible if and only if $f$ is irreducible.
\end{enumerate}
\item(\house) Let $R$ be a UFD. Show that every nonzero prime ideal of $R$ contains a prime element of $R$.
\item Show that for every field $F$ there exist infinitely many irreducible polynomials in $F[x]$. (If $F$ is finite, imitate Euclid's proof of the infinitude of prime numbers.)
\item(\house) Find all integer solutions of the equation $x^4+y^4=z^2$ by using a factorization in $\Z[i]$.
\item Prove the ``only if" direction of Theorem \ref{ufd_crit}.
\item Let $D$ be a UFD and $F=\Frac(D)$. Show that every nonzero element $r\in F$ can be written in the form $r=u\cdot p_1^{e_1}\cdots p_n^{e_n}$, where each $p_i$ is an irreducible element of $D$, and each $e_i$ is a nonzero integer. Moreover, this representation is unique up to reordering the factors and multiplication by units in $D$.
\item\label{gcd_form_pf} Prove Theorem \ref{gcd_form}.
\item Let $D$ be a UFD with field of fractions $F$. Show that $\Frac(D[x])$ is isomorphic to the function field $F(x)$.
\item(\house) Let $D$ be a UFD. A nonzero element $a\in D$ is called \col{squarefree}\index{squarefree} if there is no irreducible element $p\in D$ such that $p^2\mid a$. In particular, units are squarefree.
\begin{enumerate}
\item Let $a\in D$ be squarefree and let $b\in D$ satisfy $b^2\mid a$. Show that $b \in D^\times$.
\item Show that every nonzero element $a \in D$ is associate to a product of a square and a squarefree element. That is, show there exist $u \in D^\times$, a squarefree $s \in D$, and $q \in D$ such that $a = u \cdot q^2 s$. Moreover, show that $s$ is uniquely determined up to associates. The element $s$ is called the \col{squarefree part}\index{squarefree part} of $a$.
    \item Compute the squarefree part of $180\in\Z$ and of \[f(x) = x^4 - x^3 + x^2-2x+1\in\Q[x].\]
\end{enumerate}
\item\label{conj_copr} Let $(a,b,c)$ be pairwise coprime positive integers such that $c^2=a^2+b^2$. Prove the following:
\begin{enumerate}
\item $c\equiv 1\pmod 2$ and $a\not\equiv b\pmod 2$.
\item The element $p=1+i$ of $\Z[i]$ is prime.
\item The Gaussian integers $p$, $a+bi$, and $a-bi$ are pairwise coprime.
\end{enumerate}
\item Let $R$ be a commutative ring and $a,b\in R$. Show that the element $[a]\in R/(b)$ is prime if and only if $[b]\in R/(a)$ is prime.
\item(\sel)\label{prime_quad} Let $R=\Z[\sqrt d]$, where $d\ne 1$ is a squarefree integer, and suppose that $x\in R$ has prime norm $p=N(x)$.
\begin{enumerate}
\item Show that $R/(x)$ is isomorphic to $\F_p$.
\item Show that $x$ is a prime element of $R$.
\end{enumerate}
\end{exercises}

\chapter{Factorization in Polynomial Rings}\label{fac_chap}

In the previous chapter we developed the abstract theory of UFDs and saw that several familiar rings, such as $\Z$, $\Z[i]$, and $k[x]$ for a field $k$, admit unique factorizations. In this final chapter we prove that every polynomial ring over a UFD is itself a UFD. Thus, for instance, the ring $\Z[x]$ is a UFD, as is $k[x_1,\dots,x_n]$ for every field $k$. This important fact underlies the study of \emph{varieties} in algebraic geometry, where irreducible polynomials define irreducible geometric objects, the basic building blocks of the subject.

The key tool in this chapter is \emph{Gauss's Lemma}, which separates a polynomial into its \emph{content} and its \emph{primitive part}. This allows us to reduce questions about factorization in $R[x]$, where $R$ is a UFD, to questions in the Euclidean domain $\Frac(R)[x]$, a ring whose factorization theory is well understood. Once this structural result is in place, we develop two practical irreducibility tests---reduction modulo a prime ideal and Eisenstein's criterion---which often allow us to prove that a polynomial is irreducible.

These ideas form a natural culmination of our study of rings and factorization. They also serve as a gateway to further topics such as \href{https://en.wikipedia.org/wiki/Galois_theory}{Galois Theory}, where irreducible polynomials play a central role.

\section*{Gauss's Lemma}

Let $R$ be a UFD and $f\in R[x]$ a nonzero polynomial. The \col{content}\index{content of a polynomial} of $f$, denoted by $c(f)$, is any $\gcd$ of the coefficients of $f$. (This is well defined by Exercise \ref{finite_gcd} in Chapter \ref{euc_chap}: every nonempty finite subset of a GCD domain has a gcd.) We say that $f\in R[x]$ is \col{primitive}\index{primitive polynomial} if $c(f)\sim 1$; equivalently, no prime element of $R$ divides every coefficient of $f$. By factoring out the content, $f(x)$ can be expressed in the form \[f(x)=c\cdot\hat f(x),\]

where $c=c(f)$ and $\hat f(x)\in R[x]$ is a primitive polynomial called the \col{primitive part}\index{primitive part} of $f$. For example, the content and primitive part of \[f(x)=6x^7-38x^5+4x+2\in\Z[x]\] 

are
\[c=2,\;\;\hat f(x)=3x^7-19x^5+2x+1.\]

Note that, since gcds are unique up to multiplication by units in $R$, the same is true of $c(f)$ and $\hat f(x)$.

Theorem \ref{gauss_lemma} below, the multiplicativity of the content of polynomials, is a core result of this chapter and will be used to show that $R[x]$ is a UFD whenever $R$ is a UFD. The following notation and terminology will be used in the proof of the theorem. Recall that if $A$ and $B$ are commutative rings, any homomorphism $h:A\to B$ extends to a homomorphism $A[x]\to B[x]$ whose action on a given polynomial is to apply $h$ to every coefficient. In particular, if $I$ is any ideal of $A$, the quotient map $A\to A/I$ extends to a homomorphism $A[x]\to(A/I)[x]$ whose action on a polynomial $f(x)$ is to reduce its coefficients modulo $I$. The resulting polynomial is often denoted by $f(x)\!\mod I$, but we will write $\tilde f(x)$ instead when $I$ is clear from context. For example, the reduction of \[f(x)=12x^4-7x^3+29\in\Z[x]\]

 modulo 5 is \[\tilde f(x)=2x^4+3x^3+4\in\F_5[x].\]
 
\begin{framed}\index{Gauss's Lemma}
\begin{thm}[Gauss's Lemma]\label{gauss_lemma}
If $R$ is a UFD, then \[c(fg)\sim c(f)c(g)\] for all nonzero $f,g\in R[x]$. In particular, the product of primitive polynomials is primitive.
\end{thm}
\end{framed}
\begin{proof}
We begin by showing that the product of primitive polynomials is primitive. Suppose that $f,g\in R[x]$ are primitive polynomials and let $h=fg$. We will show that $h$ is primitive by showing that $c(h)$ is not divisible by any prime element of $R$. 

Suppose, by contradiction, that $p\in R$ is prime and $p\mid c(h)$. Note that the ideal $(p)$ of $R$ is a prime ideal because $p$ is a prime element of $R$ (see Exercise \ref{prime_id_el}, Chapter \ref{euc_chap}: every principal ideal generated by a prime element is prime). Now consider reductions modulo $(p)$: reducing coefficients modulo $(p)$, the relation $h=fg$ in $R[x]$ yields $\tilde h=\tilde f\tilde g$ in $S[x]$, where $S:=R/(p)$. The ring $S$ is a domain because $(p)$ is a prime ideal, hence $S[x]$ is a domain.

Noting that $\tilde h=0$ because $p$ divides every coefficient of $h(x)$, the equation $\tilde h=\tilde f\cdot\tilde g$ becomes $0=\tilde f\cdot\tilde g$. Since $S[x]$ is a domain, this implies that either $\tilde f=0$ or $\tilde g=0$, say $\tilde f=0$. By definition of $\tilde f$, this implies that every coefficient of $f$ is a multiple of $p$ in $R$, hence $p\mid c(f)$, a contradiction since $c(f)\sim 1$. We conclude that $c(h)$ is not divisible by any prime, as claimed. 

Now, for arbitrary nonzero $f,g\in R[x]$, write $f=c(f)\hat f$ and $g=c(g)\hat g$. Then $fg=c(f)c(g)\hat f\hat g$, and the product $\hat f\hat g$ is primitive. Letting $a_1,\ldots, a_n$ be the coefficients of $\hat f\hat g$ and $s=c(f)c(g)$, the previous factorization shows that the coefficients of $fg$ are $sa_1,\ldots, sa_n$. Hence \[c(fg)=\gcd(sa_1,\ldots,sa_n)\sim s\cdot\gcd(a_1,\ldots, a_n)= s\cdot c(\hat f\hat g)\sim s.\] Note that we have used the result of Exercise \ref{scale_gcd}, Chapter \ref{euc_chap} to factor $s$ out of the $\gcd$. This proves that $c(fg)\sim c(f)c(g)$.
\end{proof}

\section*{Unique factorization of polynomials}

Our next goal is to use Gauss's Lemma to show that the polynomial ring $R[x]$ is a UFD whenever $R$ is a UFD. A key step in this direction is Theorem \ref{hat_def} below, which relates irreducible elements of $R[x]$ to those of $\Frac(R)[x]$.

Consider $f(x):=3x^2+6\in\Z[x]$ as an example of the distinction between irreducibility in $R[x]$ and in $F[x]$. The factorization \[f(x)=3(x^2+2)\] expresses $f(x)$ as a product of two non-units in $\Z[x]$, hence $f(x)$ is reducible in $\Z[x]$. However, 3 is a unit in $\Q$, so $f(x)$ is associate to $x^2+2$ in $\Q[x]$, hence is irreducible in $\Q[x]$ by the Low-degree Test (Theorem \ref{lowdeg}). Thus, $f(x)$ is reducible in $\Z[x]$ but irreducible in $\Q[x]$.

Let $R$ be a UFD, $F=\Frac(R)$, and $f\in F[x]$ a nonzero polynomial. A \col{primitive associate}\index{primitive associate} of $f$ is a primitive polynomial $\hat f\in R[x]$ that is associate to $f$ in $F[x]$, i.e.,
\[f(x)=c\hat f(x)\quad\text{for some}\;\;c\in F^{\times}.\]

\begin{1ex}[Primitive associate] Consider $R=\Z$ and \[
f(x) = \frac{6}{5}x^3 - \frac{3}{10}x^2 + \frac{9}{25}x \in \Q[x].
\] To compute a primitive associate $\hat f(x)\in\Z[x]$, we clear denominators by multiplying $f(x)$ by $50=\lcm(5,10,25)$:
\[
50f(x) = 60x^3 - 15x^2 + 18x.
\]
Factoring out the content on the right-hand side:
\[
50f(x) = 3 \cdot (20x^3 - 5x^2 + 6x).
\]

A primitive associate of $f$ is therefore
\[
\hat f(x) = 20x^3 - 5x^2 + 6x.
\]
\end{1ex}

\begin{framed}
\begin{thm}[Properties of primitive associates]\label{hat_def} Let $R$ be a UFD and $F=\Frac(R)$. 
\begin{enumerate}
\item[(a)] Every nonzero $f\in F[x]$ has a primitive associate $\hat f\in R[x]$, unique up to multiplication by units in $R$.
\item[(b)] If $f\in R[x]$, then the primitive part of $f$ is a primitive associate. (Hence the notation $\hat f$ is unambiguous.)
\item[(c)] For all nonzero $f,g\in F[x]$,\;\; $\hat{\hat f}\sim\hat f\;\;\text{and}\;\;\widehat{fg}\sim\hat f\hat g$.
\item[(d)] For all nonconstant $f\in F[x]$,\[f\,\;\text{is irreducible in}\, F[x]\iff \hat f \;\text{is irreducible in}\;\, R[x].\]
\end{enumerate}
\end{thm}
\end{framed}
\begin{proof}
Obtaining a primitive associate requires only clearing denominators in the coefficients of $f(x)$ and then factoring out the content of the resulting polynomial in $R[x]$, as in the example preceding this theorem. To prove uniqueness of $\hat f(x)$ up to associates in $R$, suppose that $p,q\in R[x]$ are primitive polynomials associate to $f(x)$ in $F[x]$. Then $p(x)=t\cdot q(x)$ for some $t\in F^{\times}$. Clearing denominators, we obtain \[r\cdot p(x)=s\cdot q(x)\;\;\text{with}\;\; r,s\in R.\] 

Considering the content of both sides of this equation, we obtain $r\sim s$, hence $p(x)$ and $q(x)$ differ by a unit in $R$. This completes the proof of (a).

Part (b) is immediate from the definitions, as is the relation $\hat{\hat f}\sim\hat f$. The equivalence $\widehat{fg}\sim\hat f\hat g$ follows from Gauss's Lemma: the polynomial $\hat f\hat g$ is primitive and clearly associate to $fg$ in $F[x]$. This proves (c).

To prove (d) we show that $f$ is reducible in $F[x]$ if and only if $\hat f$ is reducible in $R[x]$. Assume first that $f=qh$, where $q,h\in F[x]$ are nonconstant. Then (c) implies $\hat f\sim\hat q\hat h$, so $\hat f$ is a product of non-units in $R[x]$, as $\hat q$ and $\hat h$ are nonconstant. Therefore, $\hat f$ is reducible in $R[x]$.

Conversely, assume that $\hat f$ is reducible in $R[x]$, say $\hat f(x)=a(x)b(x)$, where $a,b\in R[x]$ are non-units. Since $\hat f$ is primitive, neither $a(x)$ nor $b(x)$ can be constant. Writing $f(x)=c\cdot\hat f(x)$ for some $c\in F^{\times}$,  the identity \[f(x)=c\cdot\hat f(x)=c\cdot a(x)b(x)\] shows that $f(x)$ can be expressed as a product of nonconstant polynomials in $F[x]$, so $f(x)$ is reducible in $F[x]$.
\end{proof}

We can now prove the central theorem of this chapter.

\begin{framed}
\begin{thm}\label{poly_is_ufd}
Let $R$ be a UFD. For every $n\ge 1$, the polynomial ring $R[x_1,\ldots, x_n]$ is a UFD.
\end{thm}
\end{framed}
\begin{proof}
We argue by induction on $n$, the base case $n=1$ being the substantial step. For the inductive step, assume that \( R[x_1, \dots, x_n] \) is a UFD. Then \[R[x_1, \dots, x_{n+1}]\cong R[x_1, \dots, x_n][x_{n+1}]\] is a polynomial ring in one variable over a UFD, hence again a UFD by the case $n=1$.

It remains to prove that $R[x]$ is a UFD. The overall strategy is to use primitive associates, via Theorem \ref{hat_def}, to relate the questions of existence and uniqueness of factorizations in $R[x]$ to the same questions over $F[x]$, where $F=\Frac(R)$.

Letting $f(x)\in R[x]$ be a nonzero non-unit polynomial, we aim to show that $f(x)$ can be expressed as a product of irreducible elements of $R[x]$, and that such a representation is unique up to reordering and multiplication by units in $R$. The uniqueness statement is left as Exercise \ref{ufd_uniq_fac}; we provide the details of existence here.

 Write $f(x)=c\cdot \hat f(x)$, where $c\in R$ is the content of $f(x)$. If $f(x)$ is a constant polynomial, then $f(x)=c$ is a non-unit in $R$, hence has a factorization of the form \[f=\pi_1\cdots\pi_s,\]
 
 where each $\pi_i$ is irreducible in $R$ and therefore in $R[x]$.

Now suppose that $f(x)$ is nonconstant, so that $\hat f(x)$ has positive degree. Since $F[x]$ is a UFD, there exists a factorization \[\hat f=p_1\cdots p_n,\] 

where each $p_i$ is irreducible in $F[x]$. Applying Theorem \ref{hat_def}, we have \[\hat f\sim\hat p_1\cdots\hat p_n,\] with each $\hat p_i$ being primitive and irreducible in $R[x]$. Thus there exists a unit $u\in R^\times$ such that
\[
f(x)=c\,\hat f(x) = cu\,\hat p_1(x)\cdots \hat p_n(x).
\]

If $c\in R^{\times}$, absorbing $cu$ into one of the irreducible factors we obtain a factorization of $f$ as a product of irreducibles in $R[x]$.

If $c\notin R^{\times}$, then $c=\pi_1\cdots\pi_s$, where each $\pi_j$ is irreducible in $R$ and in $R[x]$, and the identity \[f(x)=u\pi_1\cdots\pi_s\cdot\hat p_1(x)\cdots\hat p_n(x)\]

expresses $f(x)$ as a product of irreducibles in $R[x]$.
\end{proof}

\begin{note} The statement of Theorem \ref{poly_is_ufd} remains true even if the polynomial ring has infinitely many variables; see Exercise \ref{poly_ufd_general}.
\end{note}

\begin{ex}[Polynomial UFDs]\;
\begin{examples}
\item Let $k$ be a field. Since fields are trivially UFDs, Theorem \ref{poly_is_ufd} implies that the ring $k[x]$ is a UFD. This had already been shown indirectly by proving that $k[x]$ is Euclidean, hence a PID, hence a UFD.
\item There is a crucial distinction between one-variable and multivariable polynomial rings. If $k$ is a field and $R=k[x_1,\ldots, x_n]$, then $R$ is a Euclidean domain when $n=1$, but is not even a PID (hence not Euclidean) when $n>1$. However, $R$ is always a UFD, by Theorem \ref{poly_is_ufd}.
\item Since $\Z$ is a UFD, so is $\Z[x]$. This shows that $\Z[x]$ is an example of a UFD that is not a PID, already mentioned in Chapter \ref{euc_chap}.
\end{examples}
\end{ex}

\section*{Irreducibility criteria}

Let $R$ be a UFD. The fact that $R[x]$ is a UFD raises several natural questions, among them: 
\begin{itemize}
\item In practice, how can we decide whether a given polynomial is irreducible?
\item How can we express a given polynomial as a product of irreducible polynomials?
\end{itemize} 

We will focus here on developing irreducibility criteria for polynomials. Regarding the second question above, we encourage the interested reader to explore \href{https://en.wikipedia.org/wiki/Berlekamp\%27s_algorithm}{Berlekamp's algorithm} and the \href{https://en.wikipedia.org/wiki/Cantor\%E2\%80\%93Zassenhaus_algorithm}{Cantor–Zassenhaus algorithm} for polynomials over finite fields, and the \href{https://en.wikipedia.org/wiki/Lenstra\%E2\%80\%93Lenstra\%E2\%80\%93Lov\%C3\%A1sz_lattice_basis_reduction_algorithm}{LLL algorithm} for polynomials over $\Q$. See Chapters 14--16 in von zur Gathen \& Gerhard \cite{Gathen-Gerhard} for a careful treatment of these topics.

The following theorem introduces a general technique for proving irreducibility of \emph{primitive} polynomials, namely, to reduce the coefficients modulo a prime ideal.

\begin{framed}
\begin{thm}[Reduction mod $P$ Criterion]\label{modP_crit}
Let $R$ be a UFD and $P$ a prime ideal of $R$. Suppose that $f(x)\in R[x]$ is a \underline{primitive} polynomial such that the reduction $\tilde f(x)\in(R/P)[x]$ is irreducible and has the same degree as $f(x)$. Then $f(x)$ is irreducible in both $R[x]$ and $\Frac(R)[x]$.
\end{thm}
\end{framed}
\begin{proof}
We show that $f(x)$ is irreducible in $R[x]$, which implies that it is irreducible in $\Frac(R)[x]$, by Theorem \ref{hat_def}.

Note that $\tilde f$ being irreducible implies that $f$ is nonconstant, as otherwise $f$ would be a primitive constant polynomial, hence a unit in $R[x]$, and then $\tilde f$ would be a unit.

Suppose, by contradiction, that $f$ is reducible. Then $f=ab$ with $a,b$ being non-units in $R[x]$. Since $f$ is primitive, neither $a$ nor $b$ can be constant polynomials. Reducing coefficients modulo $P$, the relation $f=ab$ becomes $\tilde f=\tilde a\tilde b$, which implies that one of $\tilde a,\tilde b$ must be a unit since $\tilde f$ is irreducible. Assuming $\tilde a$ is a unit in $(R/P)[x]$, we have $\deg(\tilde a)=0$ and therefore
\[\deg(\tilde f)=\deg(\tilde b)\le\deg(b)<\deg(f),\]
contradicting the hypothesis that $f$ and $\tilde f$ have the same degree. Hence $f$ must be irreducible.
\end{proof}

\begin{1ex}[Reduction mod $P$ Criterion] Let
\[
f(x) = 18x^3 + 27x^2 + 15x + 14 \in \Z[x].
\]
We will show that \( f(x) \) is irreducible in \( \Q[x] \). Since $f(x)$ is a cubic polynomial, we could, in principle, apply the Rational Root Theorem to show that $f(x)$ has no rational root and is therefore irreducible. However, this approach is tedious due to the number of divisors of the coefficients 18 and 14 of $f(x)$. Thus, we take a different approach.

First, observe that \( f(x) \) is primitive, so we may apply Theorem \ref{modP_crit}. To this effect, we choose a prime $p$ such that the reduction of $f(x)$ modulo $p$ remains a cubic polynomial. The smallest $p$ for which this holds is $p=5$. Reducing the coefficients of \( f(x) \) modulo \(5\):
\[
\tilde f(x)=3x^3+2x^2+ 4 \in \F_5[x].
\]
By calculating $\tilde f(a)$ for $a\in\F_5$ we conclude that $\tilde f(x)$ has no root in $\F_5$ and, being a cubic polynomial, must be irreducible. Hence, the criterion applies and allows us to conclude that $f(x)$ is irreducible in both $\Z[x]$ and $\Q[x]$.
\end{1ex}

The next and final theorem is a powerful irreducibility criterion for polynomials over a UFD. Though the theorem is stated for primitive polynomials, it can readily be applied to arbitrary polynomials by first factoring out the content. Historically, the criterion was first discovered by Sch\"onemann \cite{Cox}.

\begin{framed}\index{Eisenstein's Criterion}
\begin{thm}[Eisenstein's Criterion]\label{eisenstein}
Let $R$ be a UFD and \[f(x)=a_0+a_1x+\cdots+a_nx^n\in R[x]\] a nonconstant \underline{primitive} polynomial. Suppose that there exists a prime element $p$ of $R$ such that \[p\mid a_i\;\;\text{for}\;\; i<n,\;\;\text{and}\;\; p^2\nmid a_0.\] Then $f(x)$ is irreducible in $R[x]$ and in $\Frac(R)[x]$.
\end{thm}
\end{framed}

\begin{proof} We recall the result of Exercise \ref{eisen_lem}, Chapter \ref{poly_chap}: if $A$ is an integral domain and $s,t\in A[x]$ are polynomials such that \[ax^n=s(x)t(x),\]

where $n\ge 2$ and $a\in A$ is nonzero, then $s$ and $t$ have constant term 0.

To prove the theorem it suffices to show that $f(x)$ is irreducible in $R[x]$, by Theorem \ref{hat_def}. Suppose, by contradiction, that $f(x)=q(x)h(x)$, where $q$ and $h$ are non-units in $R[x]$. Since $f$ is primitive, both $q$ and $h$ must have positive degree, hence $n\ge 2$. Reducing the equation $f(x)=q(x)h(x)$ modulo $p$, we obtain \[\tilde a_nx^n=\tilde q(x)\cdot\tilde h(x),\]
since $p$ divides every coefficient $a_i$ with $i<n$. Note that $\tilde a_n\ne 0$ because $p\nmid a_n$. Indeed, if $p\mid a_n$, then $p$ divides every coefficient of $f$, hence $p$ divides $c(f)\sim1$, a contradiction.

Since $R/(p)$ is a domain, the exercise cited above implies that both $\tilde q(x)$ and $\tilde h(x)$ must have constant term 0. It follows that $q(x)$ and $h(x)$ have constant term divisible by $p$, and thus $p^2\mid a_0$, contradicting the hypotheses of the theorem. Hence, $f(x)$ must be irreducible in $R[x]$.
\end{proof}

\begin{note} Eisenstein's Criterion is often stated for arbitrary (not necessarily primitive) polynomials, and using prime ideals instead of prime elements. See Exercise \ref{eisen_gen} for a generalized version of the criterion incorporating both aspects.
\end{note}

We say that a polynomial $f(x)$ is \col{$p$-Eisenstein} if $f(x)$ and $p$ satisfy the hypotheses of Theorem \ref{eisenstein}. Showing that $f(x)$ is $p$-Eisenstein for some prime $p$ immediately proves that $f(x)$ is irreducible. More generally, if $f(x)\in R[x]$ is not necessarily primitive, then showing that the primitive part $\hat f(x)$ is $p$-Eisenstein implies that $f(x)$ is irreducible in $\Frac(R)[x]$, by Theorem \ref{hat_def}.

\begin{ex}[Eisenstein's Criterion]\;
\begin{examples}
\item Consider the polynomial
\[
f(x) = \frac{1}{2}x^4 + \frac{3}{4}x^3 + \frac{6}{5}x^2 + \frac{9}{10}x + \frac{12}{25} \in \Q[x].
\]
We will prove that \( f(x) \) is irreducible in \( \Q[x] \) by computing its primitive associate and applying the criterion. A straightforward calculation shows that
\[
\hat f(x)=100f(x) = 50x^4 + 75x^3 + 120x^2 + 90x + 48 \in \Z[x].
\]

Observe that $\hat f(x)$ is 3-Eisenstein and thus irreducible in \( \Z[x] \). It follows that \( f(x) \) is irreducible in \( \Q[x] \).
\item The following observation is a useful supplement to Eisenstein's Criterion. Suppose $f(x)\in R[x]$ is a primitive polynomial we suspect is irreducible but cannot find a prime $p$ such that $f(x)$ is $p$-Eisenstein. Then there may exist $a\in R$ such that the polynomial $f(x+a)$ is $p$-Eisenstein. In that case, $f(x+a)$ is proved to be irreducible, and thus $f(x)$ is irreducible. (For the last point, note that any factorization $f(x)=g(x)h(x)$ would yield a factorization $f(x+a)=g(x+a)h(x+a)$.) Consider \[f(x) = x^4+6x^3+16x^2+26x-9\in \Z[x].\] This polynomial is not Eisenstein at any prime; however, we consider the translate
\[f(x+1)= x^4+10x^3+40x^2+80x+40.\]
This polynomial is 5-Eisenstein and thus irreducible. Hence \( f(x) \) is irreducible in \(\Z[x]\) and in $\Q[x]$.
\item The idea of using translates as in the previous example may not always succeed in producing an Eisenstein polynomial. For example, for the polynomial $f(x)=x^4+x+1\in\Z[x]$, no translate $f(x+a)$ with $a\in\Z$ is Eisenstein at any prime---we leave the proof of this claim to the reader. However, $f(x)$ can be proved to be irreducible by reducing modulo 2; see Exercise \ref{chap18_sel2}.
\end{examples}
\end{ex}

\newpage

\section*{Exercises}

\begin{note}\;
\begin{itemize}
\item Solutions to starred exercises (\sel) are included in Chapter~\ref{review3chap}.
\item A house (\house) marks exercises suggested as sample homework.
\end{itemize}
\end{note}

\medskip

\begin{exercises}
\item(\house) Compute the content and primitive part of the following polynomials with coefficients in a UFD:
\begin{align*}
p(x) &= 6x^3 - 12x^2 + 18x\in\Z[x],\\
q(x) &= (2t^2 + 4t)x^2 + (6t)x + 10\in\Q[t][x],\\
r(x) &= (4 + 2t^2)x^2 + (6t^3 + 2)x + 2t\in\Q\pow t[x],\\
s(x)&= (3 + 6i)x^2 + (9 - 3i)x + 6\in\Z[i][x].
\end{align*}
\item Let $k$ be a field and $f(x)\in k[x]$ nonzero. Determine the content of $f(x)$. (Note that fields are trivially UFDs, so the content in question is well defined.)
\item(\sel)\label{chap18_sel1} Let $R$ be a PID and $f\in R[x]$ a nonzero polynomial. Suppose that $I$ is an ideal of $R$ containing every coefficient of $f(x)$. Show that $c(f)\in I$.
\item Regarding the rational number $12/35$ as a constant polynomial in $\Q[x]$, calculate its primitive associate in $\Z[x]$.
\item Let $F=\Q(t)$ be a function field in one variable over $\Q$.
\begin{enumerate}
\item Show that the rings $A=\Q[1/t]$ and $B=\Q[t]$ (subrings of $F$) are isomorphic, and are both UFDs with field of fractions $F$.
\item Let \[
f(x) = \frac{1}{t}x^4 + \frac{2}{t^2}x^3 + \frac{3}{t^3}x^2 + \frac{4}{t^4}x + \frac{5}{t^5}\in F[x].
\] Compute the primitive associate of $f$ in $A[x]$ and in $B[x]$.
\item Determine whether $f(x)$ is irreducible in $F[x]$.
\end{enumerate}
\item(\house) Factor the polynomials $f(x)=x^8-1$ and $g(x)=x^6-1$ as products of irreducibles in the UFDs $\Z[x]$, $\Q[x]$, $\F_2[x]$, and $\F_3[x]$.
\item\label{ufd_uniq_fac} Let $R$ be a UFD. Letting $f(x)\in R[x]$ be a nonzero non-unit polynomial, we showed in the proof of Theorem \ref{poly_is_ufd} that $f(x)$ can be expressed as a product of irreducible elements of $R[x]$. Show that such a factorization is unique up to reordering and multiplication by units in $R$.
\item Show that the polynomial
\[
f(x) = x^5 + 10x^4 + 25x^3 + 50x^2 + 100x + 35 \in \Z[x]
\]
is irreducible in $\Q[x]$ by applying Eisenstein’s Criterion at the prime $p = 5$.
\item Let \( f(x) = x^4 + 4x^3 + 6x^2 + 4x + 3 \in \Z[x] \).
\begin{enumerate}
\item Show that \( f(x) \) is not Eisenstein at any prime.
\item Find a translated polynomial \( f(x+a) \) that is $2$-Eisenstein.
\item Show that \( f(x) \) is irreducible in \( \Q[x] \).
\end{enumerate}
\item Let $R = \Q[t]$, $F = \Frac(R) = \Q(t)$, and consider the polynomial
\[
f(x) = x^3 + (t+1)x + (t+1) \in R[x].
\]
\begin{enumerate}
\item Show that $f(x)$ is primitive and that $(t)$ is a prime ideal of $R$.
\item Compute the reduction of $f(x)$ modulo $(t)$ and show that
\[
\tilde f(x)\in (R/(t))[x] \cong \Q[x]
\]
is irreducible in $\Q[x]$.
\item Show that $f(x)$ is irreducible in $R[x]$ and in $F[x] = \Q(t)[x]$.
\end{enumerate}
\item(\house) Let $f(x) = x^4 + 5x^3 + 10x^2 + 10x + 5 \in \Z[x]$. Use both the Eisenstein Criterion and the Reduction mod $P$ Criterion to show that $f$ is irreducible.
\item(\sel)\label{chap18_sel2} This exercise shows that the Reduction mod \( P \) Criterion can sometimes succeed where Eisenstein’s Criterion fails.

Let \( f(x) = x^4 + x + 1 \in \Z[x] \).
\begin{enumerate}
\item Show that \( f(x) \) is not Eisenstein at any prime.
\item Show that the reduction \( \tilde f(x) \in \F_2[x] \) is irreducible.
\item Conclude that \( f(x) \) is irreducible in \( \Z[x] \) and \( \Q[x] \).
\end{enumerate}
\item Define $f(x)\in\Z[i][x]$ by \[f(x) = x^5 + (1+i)x^4 + 2x^3 -10x + (3+3i).\]
\begin{enumerate}
\item Show that \( 1+i \) is a prime element of \( \Z[i] \).
\item Apply Eisenstein's Criterion to show that \( f \) is irreducible in \( \Q(i)[x] \).
\end{enumerate}
\item Let $k$ be a field. Show that the polynomial $x^n-t$ is irreducible in $k\pow t[x]$.
\item(\house) Let $f(x) = 2x^4 + 6x^3 + 4x^2 + 2x + 1 \in \Z[x]$.
\begin{enumerate}
\item Show that $f(x)$ is primitive but not Eisenstein at any prime. Moreover, the reductions of $f(x)$ modulo $2$, $3$, and $5$ are reducible or drop in degree.
\item Let $g(x)$ be the reverse of $f(x)$ as defined in Exercise \ref{reverse_poly}, Chapter \ref{ufd_chap}. Show that \( g(x) \), and hence \( f(x) \), is irreducible in \( \Q[x] \).
\end{enumerate}
\item Prove that the polynomial $1+x+x^2+\cdots+x^n\in\Z[x]$ is irreducible if and only if $n+1$ is prime.
\item Show that the polynomial $x^3+nx+2\in\Z[x]$ is irreducible for all $n\notin\{1,-3,-5\}$.
\item Show that there exist infinitely many integers $a$ such that the polynomial \[x^7+15x^2-30x+a\] is irreducible in $\Q[x]$.
\item\label{poly_ufd_general}
\begin{enumerate}
\item[(a)] Do some independent reading on directed unions of rings.
\item[(b)] Suppose that a ring \( R \) is a directed union of UFD's $R_i$, and that every irreducible in each $R_i$ remains irreducible in $R$. Show that if every nonzero non-unit in $R$ lies in some $R_i$, then $R$ is a UFD.
\item[(c)] Let $I$ be an arbitrary nonempty set, and let $D$ be a UFD. Show that the polynomial ring \( D[x_i : i \in I] \) is a UFD.
\end{enumerate}
\item\label{eisen_gen} The following is a generalization of the Eisenstein Criterion as stated in Theorem \ref{eisenstein}.
\begin{enumerate}
\item Let $f(x)\in R[x]$, where $R$ is a UFD, be a nonzero polynomial with leading coefficient $a_n$ and constant term $a_0$. Suppose that there exists a prime ideal $P$ of $R$  containing every coefficient of $f$ except for $a_n$. Moreover, suppose that $a_0\notin P^2$. Show that $\hat f(x)$ is irreducible in $R[x]$ and $f(x)$ is irreducible in $F[x]$.
\item Use the above criterion to show that the polynomial \[z^3 + xz^2 + y^2 z + xy^2 \in \Z[x, y,z]\] is irreducible.
\end{enumerate}
\end{exercises}

\chapter{Third Review on Rings}\label{review3chap}

This third and final review of the course covers Chapters \ref{crt_chap}–\ref{fac_chap}. The overall structure of the review matches that of the previous two.

\section*{Core concepts}

\begin{itemize}
\item Comaximal ideals.
\item Euclidean function, Euclidean domain.
\item Greatest common divisor, least common multiple.
\item Coprime elements, associate elements.
\item Extended Euclidean algorithm.
\item Irreducible and prime elements.
\item Unique factorization domain (UFD).
\item GCD domain.
\item Squarefree element of a UFD.
\item Content and primitive part of a polynomial.
\item $p$-Eisenstein polynomial.
\end{itemize}

\section*{Essential theorems}

\begin{itemize}
\item The Remainder Theorem: if $R$ is a commutative ring and $I_1,\ldots,I_n$ are pairwise comaximal ideals of $R$, then
\[R/(I_1\cdots I_n)\;\cong\;(R/I_1)\times(R/I_2)\times\cdots\times(R/I_n).\]
\item Every Euclidean domain is a PID.
\item Every PID is a UFD.
\item Every UFD is a GCD domain.
\item The Rational Root Theorem: let \( D \) be a GCD domain with field of fractions \( F = \Frac(D) \). Let
\[
p(x) = c_0 + c_1x + \cdots + c_nx^n \in D[x]
\]
be a polynomial with \( c_0 \ne 0 \) and \( c_n \ne 0 \). If \( r = a/b \in F \) is a root of \( p(x) \) with \( a, b \in D \) and \( \gcd(a,b)\sim 1\), then $a \mid c_0$ and $b \mid c_n$ in $D$.
\item If $R$ is a UFD, then $R[x]$ is a UFD.
\item The B\'ezout Identity: if $R$ is a PID and $a,b\in R$, then $\gcd(a,b)=ax+by$ for some $x,y\in R$.
\item Every prime element of a domain is irreducible. In a GCD domain, the two concepts coincide.
\item Gauss's Lemma: if $R$ is a UFD, then $c(fg)\sim c(f)c(g)$ for all nonzero $f,g\in R[x]$.
\end{itemize}

\subsection*{Irreducibility criteria}

\begin{itemize}
\item Low-degree Test: if $F$ is a field and $p(x)\in F[x]$ is a polynomial of degree two or three, then $p(x)$ is reducible if and only if $p(x)$ has a root in $F$.
\item Reduction mod $P$ Criterion: Let $R$ be a UFD and $P$ a prime ideal of $R$. Suppose that $f(x)\in R[x]$ is a \underline{primitive} polynomial such that the reduction $\tilde f(x)\in(R/P)[x]$ is irreducible and has the same degree as $f(x)$. Then $f(x)$ is irreducible in $R[x]$ and in $\Frac(R)[x]$.
\item Eisenstein's Criterion: Let $R$ be a UFD and \[f(x)=a_0+a_1x+\cdots+a_nx^n\in R[x]\] a nonconstant \underline{primitive} polynomial. Suppose that there exists a prime element $p$ of $R$ such that \[p\mid a_i\;\;\text{for}\;\; i<n,\;\;\text{and}\;\; p^2\nmid a_0.\] Then $f(x)$ is irreducible in $R[x]$ and in $\Frac(R)[x]$.
\end{itemize}

\section*{Selected solutions}

\begin{framed}
\noindent \textbf{Exercise \ref{chap16_sel}, Chapter \ref{crt_chap}:}\\
Let $p_1,\ldots, p_s$ be distinct odd prime integers, and $n=p_1\cdots p_s$. Determine the number of solutions of the equation $x^2=x$ in the ring $\Z_n$.
\end{framed}
\begin{proof}[Solution]
Since the $p_i$ are pairwise coprime, CRT gives an isomorphism:
\[
\Z_n \cong \Z_{p_1} \times \cdots \times \Z_{p_s}.
\]
Under this isomorphism, an element $x \in \Z_n$ corresponds to a tuple $(x_1, \ldots, x_s)$ with $x_i \in \Z_{p_i}$. The condition $x^2 = x$ in $\Z_n$ translates to
\[
(x_1^2, \ldots, x_s^2) = (x_1, \ldots, x_s),
\]
which holds if and only if $x_i^2 = x_i$ in $\Z_{p_i}$ for each $i$.

Now, for each odd prime $p_i$, the ring $\Z_{p_i}$ is a field, and in any field the only solutions to the equation $x^2 = x$ are $x = 0$ and $x = 1$. Therefore, each coordinate $x_i$ has exactly 2 possible values satisfying $x_i^2 = x_i$.

It follows that the required number of solutions is $2^s$.
\end{proof}

\begin{framed}
\noindent \textbf{Exercise \ref{chap17_sel1}, Chapter \ref{euc_chap}:}\\
Find all units of the ring $R=\Z[\sqrt{-5}]$. Is 3 a prime element of $R$? Is $1+\sqrt{-5}$ an irreducible element of $R$?
\end{framed}
\begin{proof}[Solution]
We first determine the units of $R$. An element \( u = a + b\sqrt{-5} \in R \) is a unit if and only if \[a^2 + 5b^2=N(u)=1,\]

 by Theorem \ref{quad_int_unit}. The only solution to the above equation is $b=0$, $a=\pm 1$, so we conclude that \(R^{\times}=\{\pm1\}\).

Next, we determine whether \( 3 \in R \) is a prime element. We show it is not. Note that
\[
(1+\sqrt{-5})(1-\sqrt{-5})= 6,
\]
and \( 3 \mid 6 \), so \( 3 \mid (1+\sqrt{-5})(1-\sqrt{-5}) \). However, we claim that \( 3 \nmid (1+\sqrt{-5}) \) and \( 3 \nmid (1-\sqrt{-5}) \), so 3 divides a product without dividing either factor, meaning 3 is not prime.

To prove this, suppose for contradiction that \( 3 \mid(1+\sqrt{-5}) \), i.e., \[1+\sqrt{-5} = 3(a + b\sqrt{-5})\]

for some \( a, b \in \Z \). Comparing real and imaginary parts, this implies that $3a=1$, a contradiction since $a\in\Z$. Therefore, \( 3\in R \) is not prime.

Finally, we show that \( 1 + \sqrt{-5} \) is irreducible in \( R \). Consider its norm:
\[
N(1 + \sqrt{-5}) = 1^2 + 5 = 6.
\]

(Alas, the norm is not prime, so we cannot apply Theorem \ref{prime_norm}.)
Suppose \( 1 + \sqrt{-5} = ab \) for non-units \( a, b \in R \). Then
\[
N(a)N(b) = 6.
\]
Since \( a \) and \( b \) are non-units, necessarily \( N(a), N(b) > 1 \), so one of $N(a)$ or $N(b)$ must equal 2; we may assume that $N(a)=2$. However, there is no element in \( R \) of norm 2, as the norm equation \( x^2 + 5y^2 = 2 \) has no integer solutions. This contradiction shows that \( 1+\sqrt{-5} \) is irreducible.
\end{proof}

\begin{framed}
\noindent \textbf{Exercise \ref{prime_id_el}, Chapter \ref{ufd_chap}:}\\
Let $D$ be a domain and $x\in D$ a nonzero non-unit element. Show that $x$ is prime if and only if the ideal $(x)$ is a prime ideal.
\end{framed}
\begin{proof}[Solution]
Assuming first that \( x \in D \) is a prime element, we show that the principal ideal \( (x) \) is prime. Let \( a, b \in D \) such that \( ab \in (x) \). We must show that $a\in (x)$ or $b\in (x)$. To see this, note that \( ab = x d \) for some \( d \in D \), so \( x \mid ab \). Since \( x \) is prime, it follows that \( x \mid a \) or \( x \mid b \), i.e., \( a \in (x) \) or \( b \in (x) \). Thus, \( (x) \) is a prime ideal.

Conversely, suppose that \( (x) \) is a prime ideal. Then \( x \) is nonzero and not a unit (since the ideal generated by a unit is the entire ring, which is not a proper ideal). We must show that \( x \) is a prime element, i.e., that \( x \mid ab \) implies \( x \mid a \) or \( x \mid b \).

So suppose \( x \mid ab \). Then \( ab \in (x) \), and since \( (x) \) is prime, we have \( a \in (x) \) or \( b \in (x) \), i.e., \( x \mid a \) or \( x \mid b \). Hence \( x \) is a prime element, as required.
\end{proof}

\begin{framed}
\noindent \textbf{Exercise \ref{chap18_sel1}, Chapter \ref{fac_chap}:}\\
 Let $R$ be a PID and $f\in R[x]$ a nonzero polynomial. Suppose that $I$ is an ideal of $R$ containing every coefficient of $f(x)$. Show that $c(f)\in I$.
\end{framed}
\begin{proof}[Solution] Let \( f(x) = a_0 + a_1x + \cdots + a_nx^n \in R[x] \). By assumption, we have \( a_i \in I \) for all $i$, so the ideal $A=(a_0, a_1, \ldots, a_n)$ is contained in $I$. By definition of content of a polynomial:
\[
c(f) = \gcd(a_0, a_1, \ldots, a_n).
\]
Since $R$ is a PID, the B\'ezout Identity implies that $c(f)\in A$, and therefore $c(f)\in I$.
\end{proof}

\begin{framed}
\noindent \textbf{Exercise \ref{chap18_sel2}, Chapter \ref{fac_chap}:}\\
This exercise shows that the Reduction mod \( P \) Criterion can sometimes succeed where Eisenstein’s Criterion fails.

 Let \( f(x) = x^4 + x + 1 \in \Z[x] \).
\begin{enumerate}
\item[(a)] Show that \( f(x) \) is not Eisenstein at any prime.
\item[(b)] Show that the reduction \( \tilde f(x) \in \F_2[x] \) is irreducible.
\item[(c)] Conclude that \( f(x) \) is irreducible in \( \Z[x] \) and \( \Q[x] \).
\end{enumerate}
\end{framed}
\begin{proof}[Solution] Note, first of all, that $f(x)$ is primitive. If $p$ is any prime, then $f(x)$ is not $p$-Eisenstein, as that would require $p$ dividing the non-leading coefficients of $f(x)$, and in particular $p\mid 1$, a contradiction. This proves (a).

Reducing \( f(x) \) modulo 2:
\[
\tilde f(x) = x^4 + x + 1 \in \F_2[x].
\]

To show that \( \tilde f(x) \) is irreducible in \( \F_2[x] \), we first rule out any linear factors by computing $\tilde f(0) = \tilde f(1)=1$, so $\tilde f$ has no roots in $\F_2$.

Next, check whether \( \tilde f(x) \) has an irreducible quadratic factor. In \( \F_2[x] \), the quadratic polynomials are:
\[
x^2,\quad x^2 + 1,\quad x^2 + x,\quad x^2 + x + 1.
\]
The only irreducible quadratic is $g(x)=x^2 + x + 1$. We claim that  $g(x)$ does not divide \( \tilde f(x) \). One way to prove this is using congruences: let $I$ be the principal ideal $(g)$ in $\F_2[x]$. Since $x^2\equiv x+1\pmod I$, then \[x^4\equiv x^2+2x+1\equiv x^2+1\equiv x+1+1\equiv x \;\text{(mod $I$)},\] hence \[\tilde f(x)=x^4+x+1\equiv 2x+1\equiv 1\not\equiv 0\;\text{(mod $I$)},\]

showing that $g(x)$ does not divide $\tilde f(x)$. This proves (b).

Part (c) follows immediately from the Reduction mod $P$ Criterion: we have verified that $f(x)\in\Z[x]$ is a primitive polynomial whose reduction modulo the prime ideal $(2)$ is irreducible and has the same degree as $f(x)$.
\end{proof}

\section*{Additional Practice}

\begin{exercises}
\item Let $R=\Z[2i]$, a subring of $\Z[i]$.
\begin{enumerate}
\item Show that in $R$, the element $2i$ is irreducible but not prime. Conclude that $R$ is not a UFD.
\item In $R$, find two factorizations of 4 that are not related via reordering and multiplication by units.
\end{enumerate}
\item Let $f\in\Z[x]$ be a monic polynomial. Show that if $f$ has a rational root, then it has an integer root. Using this, show that the real number $\sqrt 2+\sqrt[3] 2$ is irrational. 
\item Let $n,k>1$ be integers. Show $\sqrt[n]{k}$ is either an integer or is irrational.
\item Show that the prime ideals of $\Z[x]$ are $(0)$, $(f(x))$ for irreducible polynomials $f(x)$, and maximal ideals. Moreover, show that every maximal ideal has the form $(p,g(x))$, where $g\in\Z[x]$ is a polynomial whose reduction modulo $p$ is irreducible in $\F_p[x]$.
\item Let $R$ be an integral domain. Show that $R$ is a UFD if and only if every nonzero prime ideal of $R$ contains a prime element.
\item Let $p,q\in\Z[x]$. Show that $p$ and $q$ are coprime in $\Q[x]$ if and only if the ideal $(p,q)$ in $\Z[x]$ contains an integer.
\item For every positive integer $n$, construct a UFD with exactly $n$ irreducible elements up to associates.
\item Show that if $k$ is a field and $f(x)\in k[x]$ irreducible, then the quotient ring $k[x]/(f)$ is a field. If $f(x)$ is any squarefree polynomial, show that the ring $k[x]/(f)$ is isomorphic to a direct product of fields.
\item Let \( \omega = e^{2\pi i/3} = \frac{-1 + \sqrt{-3}}{2} \), and let
\[
\Z[\omega] = \{ a + b\omega \mid a, b \in \Z \}
\]
denote the ring of \col{Eisenstein integers}. The goal of this exercise is to show that \( \Z[\omega] \) is a Euclidean domain.

\begin{enumerate}
\item Show that \( \Z[\omega] \subseteq \C \) is a subring of $\C$.  
\item Define a norm function \( N: \Z[\omega] \to \N \) by
\[
N(a + b\omega) = a^2 - ab + b^2.
\]
Show that the norm is multiplicative: for all \( z, w \in \Z[\omega] \),
\[
N(zw) = N(z)N(w).
\]
\item Let \( \alpha, \beta \in \Z[\omega] \) with \( \beta \ne 0 \). Show that there exist \( q, r \in \Z[\omega] \) such that
\[
\alpha = q\beta + r \quad \text{with} \quad N(r) < N(\beta).
\]
(Approximate $\alpha/\beta\in\C$ by a nearest Eisenstein integer \( q \).)
\end{enumerate}
\item Let $R$ be a commutative ring. The \col{Krull dimension} of $R$, denoted $\dim(R)$, is the supremum (in $\overline\R$) of the set of natural numbers $n$ such that there exists a strictly increasing chain
\[P_0\subset P_1\subset\cdots\subset P_n,\]
where each $P_i$ is a prime ideal of $R$. (Note that $n$ is the number of links in the chain, rather than the number of ideals.)

Prove the following:
\begin{enumerate}
\item Every field has dimension 0.
\item Every PID has dimension 1.
\item The ring $\Z[x]$ has dimension 2.
\end{enumerate}
\end{exercises}

\appendix
\chapter{Thematic Quizzes}

This appendix gathers a collection of short, focused quizzes designed to test mastery of the core ideas developed throughout the course. Each section corresponds to a major topic, and the problems emphasize conceptual understanding as well as essential computational skills. The quizzes are suitable for self-study, review sessions, or in-class assessment.

\section{Rings and Subrings}

\begin{enumerate}
  \item \textbf{Definition and Examples of Rings.}  
  Give an example of a ring that is not a domain, and explain why it fails to be one.

  \item \textbf{Zero Divisors.}  
  Give an example of a zero divisor in \( \Z_6 \) and prove that your example satisfies the definition.

  \item \textbf{Units and Invertibility.}  
  List all units in \( \Z_8 \) and verify that each is invertible.

  \item \textbf{Polynomial Rings.}  
  Is \( \mathbb{Z}[x] \) an integral domain? Is it a field?

  \item \textbf{Formal Power Series.}  
  Let \( R \) be a commutative ring. Prove that the units in the power series ring \( R\llbracket x \rrbracket \) are precisely the power series whose constant term is a unit in \( R \).

    \item \textbf{Multivariable Polynomials.}  
    Let \( f(x,y) = 3x^2y + 2xy^2 - y \in \mathbb{Z}[x,y] \).  Determine the total degree of \( f \), the degree in \( x \), and the degree in \( y \). Is $f$ a homogeneous polynomial?
  
  \item \textbf{Subrings.}  
  Explicitly describe all elements of the subring of \( \mathbb{Q} \) generated by \( \frac{1}{2} \).

  \item \textbf{Prime Subring.}  
  Define the prime subring of a ring. Choose two rings and describe their prime subrings.

\end{enumerate}

\section{Ideals}

\begin{enumerate}
  \item \textbf{Ideal Containment and Equality.}  
  Let \( I = (6) \), \( J = (2, 3) \), and \( K = (12) \) in \( \mathbb{Z} \).  
  \begin{enumerate}
    \item Determine which of the containments \( I \subseteq J \), \( J \subseteq I \), \( K \subseteq I \) hold.  
    \item Are any of the ideals equal?
  \end{enumerate}

  \item \textbf{Ideal Generation.}  
  In the ring \( \mathbb{Z}[x] \), consider the ideal \( (2, x) \). Is this ideal principal? Explain why or why not.

\item \textbf{Ideals in \( \Z\pow x \).}  
Let \( A = \Z\pow x\), the ring of formal power series.  
\begin{enumerate}
  \item Show that the set of power series whose constant term is divisible by 3 forms an ideal in \( A \).  
  \item Is this ideal principal? If so, give a generator.
\end{enumerate}

\item \textbf{Ideal Membership.}  
Let \( f(x) = 6x + 4 \in \mathbb{Z}[x] \), and let \( I = (2, x) \subseteq \mathbb{Z}[x] \). Show that \( f(x) \in I \), and determine whether \( g(x) = x^2 + 1 \) belongs to \( I \).

  \item \textbf{Lattice of Ideals.}  
  Let \( R = \mathbb{Z} \), and let \( I = (4) \), \( J = (6) \).  
  \begin{enumerate}
    \item Compute \( I + J \) and \( I \cap J \).  
    \item Verify that \( I + J \) is the smallest ideal containing both \( I \) and \( J \).
  \end{enumerate}
\end{enumerate}

\section{Quotients and Direct Products}

\begin{enumerate}
  \item \textbf{Direct Product Structure.}  
  Let \( R = \mathbb{Z} \times \mathbb{Z} \).  
  \begin{enumerate}
    \item Is \( R \) an integral domain? Explain.  
    \item List all the units in \( R \).
  \end{enumerate}

  \item \textbf{Zero Divisors in a Product Ring.}  
  Let \( A = \Z_6\times\Z_{10} \).  
  Find a zero divisor in \( A \) and justify your answer.

  \item \textbf{Ideal in a Product Ring.}  
  Let \( R = \mathbb{Z} \times \mathbb{Z} \), and let  
  \[
    I = \{ (a,0) \mid a \in \mathbb{Z} \}.
  \]  
  \begin{enumerate}
    \item Show that \( I \) is an ideal of \( R \).  
    \item Is \( I \) a principal ideal? Explain.
  \end{enumerate}

    \item \textbf{Quotient Ring Construction.}  
    Let \( R = \mathbb{Z}[x] \), and let \( I = (x^2 + 1) \).  
    \begin{enumerate}
     \item Are \( x + I \) and \( -x + I \) equal in \( R/I \)? Justify your answer.  
      \item Compute \( (x + I)^2 \) and \( (2x + 3 + I)(x + 1 + I) \) in \( R/I \).  
    \end{enumerate}

    \item \textbf{Zero Divisors in a Quotient Ring.}  
    Let \( R = \mathbb{Z}[x] \), and let \( I = (x^2) \subseteq R \). Consider the quotient ring \( R/I \).  
    \begin{enumerate}
      \item Show that every element of \( R/I \) can be written in the form \( a + bx + I \), where \( a, b \in \mathbb{Z} \).  
      \item Prove that \( x + I \) is a zero divisor in \( R/I \).
    \end{enumerate}
\end{enumerate}

\section{Homomorphisms and Isomorphisms}

\begin{enumerate}
    \item \textbf{Homomorphisms.}  
    Let \( f: \mathbb{Z} \to M_2(\mathbb{Z}) \) be defined by
    \[
    f(n) = \begin{bmatrix} n & 0 \\ 0 & n \end{bmatrix}.
    \]
    
    Show that \( f \) is a ring homomorphism.

    \item \textbf{Kernel and Image.}  
    Let \( \varphi: \mathbb{Z}[x] \to \mathbb{Z} \) be defined by \( \varphi(f(x)) = f(1) \), i.e., evaluation at \( x = 1 \).  
 Describe the kernel and image of \( \varphi \). Is \( \varphi \) injective? Is it surjective? Justify your answers.

\item \textbf{Isomorphism Theorem.}  
Let \( \varphi: \mathbb{Z}[x] \to \mathbb{Z} \) be the evaluation map defined by \( \varphi(f(x)) = f(1) \).
\begin{enumerate}
  \item Show that \( \varphi \) is a surjective ring homomorphism.  
  \item Determine \( \ker(\varphi) \).  
  \item Use the isomorphism theorem to describe the ring \( \mathbb{Z}[x]/\ker(\varphi) \).
\end{enumerate}

  \item \textbf{Substitution Homomorphisms.}  
  Let 
  \[
    \varphi : \mathbb{Z}[x,y] \longrightarrow \mathbb{Z}[x]
  \]
  be the ring homomorphism defined by  
  \[
    \varphi\big(f(x,y)\big) = f\big(x,\,x+1\big).
  \]
  \begin{enumerate}
    \item Prove that \( \varphi \) is a ring homomorphism.
    \item Compute \( \varphi(x) \) and \( \varphi(y) \).
    \item Find two nonzero elements of \( \ker(\varphi) \).
  \end{enumerate}

\end{enumerate}

\section{Prime and Maximal Ideals}

\begin{enumerate}
  \item \textbf{Prime vs. Maximal Ideals.}  
  Let \( R = \mathbb{Z} \).  
  \begin{enumerate}
    \item Which of the following ideals are prime: \( (0) \), \( (4) \), \( (7) \)?  
    \item Which are maximal? Justify your answers.
  \end{enumerate}

  \item \textbf{Characterizing Ideals via Quotients.}  
  Let \( R = \mathbb{Z}[x] \), and let \( I = (x^2 + 1) \).  
  \begin{enumerate}
    \item Is \( I \) a maximal ideal?  
    \item Is \( I \) a prime ideal?  
    \item Justify your answers using properties of the quotient ring \( R/I \).
  \end{enumerate}
  
    \item \textbf{Prime and Maximal Ideals in Product Rings.}  
  Let \( R = \mathbb{Z} \times \mathbb{Z} \).  
  \begin{enumerate}
    \item Describe all prime ideals of the form \( P \times \mathbb{Z} \).  
    \item Is \( (2) \times \mathbb{Z} \) a maximal ideal of \( R \)? Why or why not?
  \end{enumerate}

  \item \textbf{Local Rings.}  
  Let \( R = \mathbb{Z}_{(5)} = \{ \frac{a}{b} \in \mathbb{Q} : 5\nmid b \} \).  
  \begin{enumerate}
    \item Show that the set \( M = \{ \frac{a}{b} \in R : 5\mid a \} \) is a maximal ideal.  
    \item Explain why \( R \) is a local ring.
  \end{enumerate}

  \item \textbf{CRT.}  
  Describe all solutions $x\in\Z$ of the system of congruences
  \[x \equiv 2\;\text{(mod 5)},\quad x \equiv 1\;\text{(mod 7)}.\]
  
\end{enumerate}

\section{Euclidean Domains and UFDs}

\begin{enumerate}
  \item \textbf{Divisibility in Euclidean Domains.}  
  Let \( D = \mathbb{Z} \).  
  \begin{enumerate}
    \item Use the Euclidean algorithm to compute \( \gcd(84, 30) \).  
    \item Express the gcd as a linear combination of 84 and 30.
  \end{enumerate}

\item \textbf{Examples of Euclidean Domains and UFDs.}  
Provide one example for each of the following. In each case, justify your answer briefly.
\begin{enumerate}
  \item A ring that is a Euclidean domain but not a field.  
  \item A ring that is a UFD but not a Euclidean domain.  
  \item A ring that is an integral domain but not a UFD.
\end{enumerate}

  \item \textbf{Irreducibility in \( \mathbb{Z}[x] \).}  
  Let \( f(x) = x^3 + 2x^2 + 4x + 8 \in \mathbb{Z}[x] \).  
  \begin{enumerate}
    \item Use the Rational Root Theorem to determine whether \( f \) has any rational roots.  
    \item Use Eisenstein’s Criterion (after an appropriate substitution, if needed) to test whether \( f \) is irreducible in \( \mathbb{Q}[x] \).
  \end{enumerate}

  \item \textbf{Content and Primitive Part.}  
  Let \( f(x) = 6x^2 + 9x + 3 \in \mathbb{Z}[x] \).  
  \begin{enumerate}
    \item Compute the content and the primitive part of \( f \).  
    \item Factor the primitive part into irreducibles in \( \mathbb{Z}[x] \).
  \end{enumerate}

\item \textbf{Irreducibility via Reduction mod \( P \).}  
Let \[f(x) = x^4 + 2x^2 + 9 \in \mathbb{Z}[x].\]  
\begin{enumerate}
  \item Reduce \( f(x) \) modulo 5 and factor it in \( \mathbb{F}_5[x] \).  
  \item Determine whether $f(x)$ is irreducible in \( \mathbb{Q}[x] \).  
\end{enumerate}

\end{enumerate}

\backmatter

\addcontentsline{toc}{chapter}{Index}
\printindex

\end{document}